\documentclass[a4paper,oneside]{article}

\usepackage{amsmath,amsthm}

\usepackage[pdftitle={Tables of Calabi--Yau equations},
pdfauthor={Gert Almkvist, Christian van Enckevort, Duco van Straten,
  and Wadim Zudilin}]{hyperref}
\newcommand{\hlink}[1]{\href{#1}{\texttt{#1}}}
\newcommand{\mailto}[1]{\href{mailto:#1}{\texttt{#1}}}
\newcommand{\dburl}{http://enriques.mathematik.uni-mainz.de/CYequations/}
\newcommand{\eqnno}[2][4]{\href{\dburl/summary.php?org=cytab\&order=#1\&eqnno=#2}{#2}}
\newcommand{\eqnnol}[3][4]{\href{\dburl/summary.php?org=cytab\&order=#1\&eqnno=#3}{#2}}
\newcommand{\eqnnoltmp}[2]{{#1}}

\theoremstyle{remark}
\newtheorem{cond}{Condition}
\newtheorem{step}{Step}
\newtheorem{remarka}{Remark}
\newtheorem*{remark}{Remark}

\newtheorem*{comments}{Comments}

\usepackage{calc}
\usepackage{array}
\usepackage{tabularx}

\newcolumntype{C}{>{$}c<{$}}
\newcolumntype{L}{>{$}l<{$}}
\newcolumntype{R}{>{$}r<{$}}

\usepackage{myref}
\usepackage{mysymb}

\ifx\pdfoutput\undefined\else
\advance\evensidemargin by -28mm
\advance\oddsidemargin by -28mm
\advance\topmargin by -32mm
\fi

\newcommand{\wh}{\widehat}
\renewcommand{\d}{\mathrm{d}}
\newcommand{\ph}{\vphantom{\bigg|_1^1}}
\newcommand{\dsize}{\displaystyle}

\begin{document}
\title{Tables of Calabi--Yau equations}
\author{Gert Almkvist, Christian van Enckevort,\\
Duco van Straten, and Wadim Zudilin}
\date{9 October 2010}
\maketitle

\begin{abstract}
  The main part of this paper is a big table 
  containing what we believe to be a complete list of all fourth order
  equations of Calabi--Yau type known so far. In the text preceding
  the tables we explain what a differential equation of Calabi--Yau
  type is and we briefly discuss how we found these equations. We also
  describe an electronic version of this list.
\end{abstract}

\tableofcontents

\section{Differential equations of Calabi--Yau type}
The differential equations we will be investigating are all of the
following form
\begin{equation}
\label{eq:Dnk}
  D_{n,k}:=\sum_{i=0}^k z^i P_i(\theta).
\end{equation}
Here $\theta:=z\frac{\partial}{\partial z}$ and the $P_i(\theta)$ are
polynomials with integral coefficients and order $n$. We will mostly
restrict ourselves to equations of order 4, but the number of terms
$k+1$ will vary. Independently of the number of terms we can rewrite the
differential equation $D_{n,k} y=0$ in the following way
\[
  y^{(n)}+a_{n-1}(z)y^{(n-1)}+ \dots +a_2(z)y'' +a_1(z)y'+a_0(z)y=0,
\]
where the $a_i(z)$ are rational functions of $z$. Recall that at any
point $z=a$ (including $z=\infty$ via the transformation $w=z^{-1}$)
we can define the so called indicial equation which governs the
existence of solutions of the form $z^\lambda \times (\text{power
  series in $z$})$. The set of solutions (with multiplicities) to the
indicial equation at a certain point will be called the spectrum at
that point.

In \cite{AZ} the authors define Calabi--Yau equations to be
differential equations of the form $D_{4,k} y = 0$ satisfying certain
extra conditions. Let us briefly list the conditions satisfied by the
equations from \pref{app:CYlist}. For an explanation of these
conditions we refer to \cite{AZ,ES}. Together they can be considered
to define the notion of a Calabi--Yau equation.

\begin{cond}
\label{en:condMUM}
  The singular point $z=0$ is a point of maximal
  unipotent monodromy, i.e., the indicial equation at $z=0$ should
  have $0$ as its only solution.
\end{cond}

\begin{cond}
\label{en:condAZ}
  The coefficients $a_i(z)$ satisfy the following
  equation
\[
  a_1 = \frac{1}{2}a_2a_3 - \frac{1}{8}a_3^3 + a_2' -
  \frac{3}{4}a_3a_3' - \frac{1}{2}a_3''.
\]
\end{cond}

\begin{cond}
\label{en:condspecinf}
  The solutions $\lambda_1 \le \lambda_2 \le
  \lambda_3 \le \lambda_4$ of the indicial equation at $z=\infty$ are
  positive rational numbers satisfying $\lambda_1 + \lambda_4 =
  \lambda_2 + \lambda_3=s$ for some $s \in \Q$ (symmetry). We also
  suppose that the eigenvalues $e^{2\pi i \lambda}$ of the monodromy
  around $z=\infty$ are the zeros (counted with multiplicity) of a
  product of cyclotomic polynomials, which can be interpreted as the
  characteristic polynomial of the monodromy around $z=\infty$.
\end{cond}

\begin{cond}
  The power series solution near $z=0$ has integral coefficients.
\end{cond}

\begin{cond}
  The genus zero instanton numbers computed by the standard recipe
  are integral (up to multiplication by an overall positive integer
  corresponding to the degree, which is not apparent from the
  differential equation alone).
\end{cond}

\section{Construction of Calabi--Yau equations}
Some of the equations we collected in this paper have a geometric
origin. These account for equations \#1--28 except \#9 in our list
(see introduction to the main table for the precise references). Such
geometric equations served as the starting point of our list. They are
automatically of Calabi--Yau type because of their geometric origin. We
observed that the expressions for the differential operators are very
similar. This encouraged us to start looking for more equations that
share the properties listed above. To find such equations we used various
constructions:
\begin{itemize}
\item Translation to other singular points (not $z=0$) with spectrum
  $\{a,a,a,a\}$
\item Manipulation of binomial expressions for the coefficients of
  the power series solution
\item Hadamard product
\item `Pullbacks' from fifth order equations
\item Computer searches
\end{itemize}
The first method of obtaining new equations is rather trivial.  Note
that we defined a Calabi--Yau equation to have a point of maximal
unipotent monodromy at $z=0$. However, it can have other singular
points where the spectrum consists of a single value $a$ with
multiplicity $4$. These points can be found by computing the spectrum
at all singular points. By translating to such a point (i.e., by
substitution $z=z_0+w$ if $z_0$ is the MUM-point or $z=1/w$ if the
MUM-point is at infinity) and simultaneously replacing $y$ by
$w^{-a}\tilde{y}$ we get a new equation for $\tilde{y}$, that has
maximal unipotent monodromy at $w=0$.  Of course in some sense it is
still the same equation. However, it can happen that one arrives at
the same equation at a different point of maximal unipotent monodromy
in entirely different ways. Then it may not always be clear that two
equations are just translates of each other (in the sense described
above).  An example of a Picard--Fuchs equation with two points of
maximal unipotent monodromy was discussed in \cite{Ro,Tj}. In that
case both points have a geometric interpretation.  By including
translates it is easier to check if we include a certain equation.

The next three constructions are discussed in more detail in
\cite{AZ}. We will restrict ourselves to a few brief remarks.  In many
cases explicit formulas are known for the coefficients of the power
series solution around $z=0$. They satisfy a recursion relation coming
from the differential equation. As can be seen from the list, in many
cases the coefficients can be expressed in terms of sums or double
sums of binomial expressions.  There are many ways of transforming
such expressions to other expressions which still satisfy a recursion
relation. Translating this recursion relation back to a differential
equation we obtain a new differential equation which in many cases
turns out to be of Calabi--Yau type again.

The Hadamard product of two equations with power series solutions
around $z=0$ given by $\sum_{n=0}^\infty A_n z^n$ and
$\sum_{n=0}^\infty B_n z^n$ is the differential equation that has
$\sum_{n=0} A_n B_n z^n$ as its power series solution. As explained
in \cite{AZ} the Hadamard product can be used to construct higher
order equations from lower order ones. The lower order ones come from
\cite{Za,Go} and computer searches.

As discussed in \cite{AZ} we can construct a fifth order equation from
a fourth order one. This construction also goes the other way: given a
suitable fifth order equation the corresponding fourth order equation
can be constructed. The fifth order equations are found by other
methods similar to the ones discussed above for fourth order equations.

\subsection{Computer searches}
A look at the list of differential equations shows that certain types
occur very often. If we restrict to equations with $k=1$, then the
list reduces to the 14 hypergeometric cases. These are all of the form
\begin{equation}
\label{eq:Dhg}
  D_{4,1}^{\text{hg}} = \theta^4 -c z (v \theta + u) (x \theta + w)
    (x \theta + x-w) (v \theta + v-u),
\end{equation}
for integers $c$, $u$, $v$, $w$, and $x$.  The term without $z$ has to
be $a \theta^4$ for some integer $a \ne 0$ because we want the
equation to have maximal unipotent monodromy at $z=0$. For some reason
only $a=1$ occurs. The last term determines the indicial equation at
$z=\infty$, as can be seen by transforming \pref{eq:Dnk} to
$w=z^{-1}$.
\begin{equation}
\label{eq:Dnkw}
  D_{n,k}^w:=\sum_{i=0}^k w^i P_{k-i}(-\theta_w),
\end{equation}
where $\theta_w:=w \frac{\partial}{\partial w}$. Therefore the indicial
equation at $z=\infty$ is
\[
  P_k(-\lambda)=0.
\]
Hence Condition~\ref{en:condspecinf} tells us among other things that
$P_k(-\lambda)$ has rational roots. Another part of this condition is
the symmetry $\lambda_1 + \lambda_4 = \lambda_2 + \lambda_3 =s$ for
some $s \in \Q$. Take $s=1$ and write $\lambda_1=\frac{u}{v}$ and
$\lambda_2=\frac{w}{x}$, then we find the expression from
\pref{eq:Dhg}.

However, we still have not fully exploited Condition
\ref{en:condspecinf}. The characteristic polynomial has to be a
product of cyclotomic polynomials. As it has fixed degree $4$, there
are only finitely many possibilities (see \sref{tab:spectra}). The
$\lambda_i$ are defined by the property that $e^{2\pi i
  \lambda_1},\dots, e^{2\pi i\lambda_4}$ is the full set of zeros of
the characteristic polynomial. This determines the $\lambda_i$ up to
integer shifts. Because they also have to be positive and we have the
symmetry condition $\lambda_1 + \lambda_4 = \lambda_2 + \lambda_3 =s$,
there are only finitely many choices for the $\lambda_i$'s for a fixed
characteristic polynomial and fixed $s$. As mentioned above, for the
hypergeometric equations $s$ turns out to be $1$. We listed the 14
possible spectra in \sref{tab:spectra}. These correspond exactly to
the first 14 equations in our list. Given the spectrum at $z=\infty$
the only thing that is not fixed is the constant $c$. This constant
can easily be found by requiring the coefficients of the power series
solution around $z=0$ to be integral.

\begin{table}
\begin{center}
\setlength{\extrarowheight}{2pt}
\begin{tabular}{|C|l|>{\raggedright\arraybackslash}m{6cm}|}
\hline
\multicolumn{1}{|l|}{char.\ polyn.} & spectra ($s=1$) & spectra ($s=2$)
\\ \hline
\phi_1\phi_2\phi_3 & --- & --- \\ \hline
\phi_1\phi_2\phi_6 & --- & --- \\ \hline
\phi_1\phi_2\phi_6 & --- & --- \\ \hline
\phi_1\phi_2^3     & --- & --- \\ \hline
\phi_1^2\phi_2^2   & --- & $\{\frac12,1,1,\frac32\}^*$
\rule[-5pt]{0pt}{0pt} \\ \hline
\phi_1^2\phi_3     & --- & $\{\frac23,1,1,\frac43\}^*$,
$\{\frac13,1,1,\frac53\}^*$ \rule[-5pt]{0pt}{0pt} \\ \hline
\phi_1^2\phi_4 & --- & $\{\frac34,1,1,\frac54\}$,
$\{\frac14,1,1,\frac74\}^*$ \rule[-5pt]{0pt}{0pt} \\ \hline
\phi_1^2\phi_6 & --- & $\{\frac56,1,1,\frac76\}^*$,
$\{\frac16,1,1,\frac{11}{6}\}^*$ \rule[-5pt]{0pt}{0pt} \\ \hline
\phi_1^3\phi_2 & --- & --- \\ \hline
\phi_1^4       & --- & $\{1,1,1,1\}^*$ \\ \hline
\phi_2^2\phi_3 & $\{\frac13,\frac12,\frac12,\frac23\}$ &
$\{\frac12,\frac23,\frac43,\frac32\}$,
$\{\frac13,\frac12,\frac32,\frac53\}$ \rule[-5pt]{0pt}{0pt} \\  \hline
\phi_2^2\phi_4 & $\{\frac14,\frac12,\frac12,\frac34\}$ &
$\{\frac12,\frac34,\frac54,\frac32\}$,
$\{\frac14,\frac12,\frac32,\frac74\}$ \rule[-5pt]{0pt}{0pt} \\ \hline
\phi_2^2\phi_6 & $\{\frac16,\frac12,\frac12,\frac56\}$ &
$\{\frac16,\frac12,\frac32,\frac{11}{6}\}$,
$\{\frac12,\frac56,\frac76,\frac32\}$, \rule[-5pt]{0pt}{0pt} \\ \hline
\phi_2^4 & $\{\frac12,\frac12,\frac12,\frac12\}$ &
$\{\frac12,\frac12,\frac32,\frac32\}$ \rule[-5pt]{0pt}{0pt} \\ \hline
\phi_3\phi_4 & $\{\frac14,\frac13,\frac23,\frac34\}$ &
\rule[-5pt]{0pt}{0pt}$\{\frac23,\frac34,\frac54,\frac43\}$,
$\{\frac14,\frac13,\frac53,\frac74\}$,
$\{\frac13,\frac34,\frac54,\frac53\}$,
$\{\frac14,\frac23,\frac43,\frac74\}$ \rule[-5pt]{0pt}{0pt} \\ \hline
\phi_3\phi_6 & $\{\frac16,\frac13,\frac23,\frac56\}$ &
\rule[-5pt]{0pt}{0pt}$\{\frac23,\frac56,\frac76,\frac43\}$,
$\{\frac16,\frac13,\frac53,\frac{11}{6}\}$,
$\{\frac13,\frac56,\frac76,\frac53\}$,
$\{\frac16,\frac23,\frac43,\frac{11}{6}\}$ \rule[-5pt]{0pt}{0pt} \\ \hline
\phi_3^2 & $\{\frac13,\frac13,\frac23,\frac23\}$ &
$\{\frac13,\frac23,\frac43,\frac53\}$, $\{\frac23,\frac23,\frac43,\frac43\}$,
$\{\frac13,\frac13,\frac53,\frac53\}$ \rule[-5pt]{0pt}{0pt} \\ \hline
\phi_4\phi_6 & $\{\frac16,\frac14,\frac34,\frac56\}$ &
\rule[-5pt]{0pt}{0pt}$\{\frac16,\frac34,\frac54,\frac{11}{6}\}$,
$\{\frac34,\frac56,\frac76,\frac54\}$,
$\{\frac14,\frac56,\frac76,\frac74\}$,
$\{\frac16,\frac14,\frac74,\frac{11}{6}\}$ \rule[-5pt]{0pt}{0pt} \\ \hline
\phi_4^2 & $\{\frac14,\frac14,\frac34,\frac34\}$ &
$\{\frac14,\frac34,\frac54,\frac74\}$,
$\{\frac34,\frac34,\frac54,\frac54\}$,
$\{\frac14,\frac14,\frac74,\frac74\}$ \rule[-5pt]{0pt}{0pt} \\ \hline
\phi_5 & $\{\frac15,\frac25,\frac35,\frac45\}$ &
\rule[-5pt]{0pt}{0pt}$\{\frac25,\frac45,\frac65,\frac85\}$,
$\{\frac35,\frac45,\frac65,\frac75\}$,
$\{\frac15,\frac35,\frac75,\frac95\}$,
$\{\frac15,\frac25,\frac85,\frac95\}$ \rule[-5pt]{0pt}{0pt} \\ \hline
\phi_6^2 & $\{\frac16,\frac16,\frac56,\frac56\}$ &
$\{\frac16,\frac56,\frac76,\frac{11}{6}\}$,
$\{\frac16,\frac16,\frac{11}{6},\frac{11}{6}\}$,
$\{\frac56,\frac56,\frac76,\frac76\}$\rule[-5pt]{0pt}{0pt} \\ \hline
\phi_8 & $\{\frac18,\frac38,\frac58,\frac78\}$ &
\rule[-5pt]{0pt}{0pt}$\{\frac18,\frac38,\frac{13}{8},\frac{15}{8}\}$,
$\{\frac18,\frac58,\frac{11}{8},\frac{15}{8}\}$,
$\{\frac38,\frac78,\frac98,\frac{13}{8}\}$,
$\{\frac58,\frac78,\frac98,\frac{11}{8}\}$ \rule[-5pt]{0pt}{0pt} \\ \hline
\phi_{10} & $\{\frac{1}{10},\frac{3}{10},\frac{7}{10},\frac{9}{10}\}$
& \rule[-5pt]{0pt}{0pt}$\{\frac{7}{10},\frac{9}{10},\frac{11}{10},\frac{13}{10}\}$,
$\{\frac{1}{10},\frac{3}{10},\frac{17}{10},\frac{19}{10}\}$,
$\{\frac{1}{10},\frac{7}{10},\frac{13}{10},\frac{19}{10}\}$,
$\{\frac{3}{10},\frac{9}{10},\frac{11}{10},\frac{17}{10}\}$
\rule[-5pt]{0pt}{0pt} \\ \hline
\phi_{12} & $\{\frac{1}{12},\frac{5}{12},\frac{7}{12},\frac{11}{12}\}$
& \rule[-5pt]{0pt}{0pt}$\{\frac{1}{12},\frac{5}{12},\frac{19}{12},\frac{23}{12}\}$,
$\{\frac{1}{12},\frac{7}{12},\frac{17}{12},\frac{23}{12}\}$,
$\{\frac{5}{12},\frac{11}{12},\frac{13}{12},\frac{19}{12}\}$,
$\{\frac{7}{12},\frac{11}{12},\frac{13}{12},\frac{17}{12}\}$
\rule[-5pt]{0pt}{0pt} \\ \hline
\end{tabular}
\end{center}
\caption{Possible characteristic polynomials and spectra}
\label{tab:spectra}
\end{table}

One step up from the hypergeometric equations are the Calabi--Yau
equations with $k=2$. In this case the most general form is
\begin{equation}
\label{eq:four}
\begin{split}
  D_{4,2}^{\text{gen}} &= \theta^4 -c z (A \theta^4 + 2A\theta^3 +
  (B+A) \theta^2 + B\theta + C) \\
  &\qquad - d z^2 (v \theta + u) (x \theta + w) (x \theta + 2x-w)
  (v \theta + 2v-u).
\end{split}
\end{equation}
The first and the last term are again restricted because of the
required monodromies (Conditions~\ref{en:condMUM}
and~\ref{en:condspecinf}), where again we suppose that the coefficient
of $\theta^4$ is 1. The only difference is that for the
spectrum at $z=\infty$ we now use $s=2$ instead of $s=1$. This is just
based on the experimental fact that for $k=2$ most Calabi--Yau
equations turn out to have $s=2$. In fact our database contains one
exception ($s=\frac{3}{2}$). For larger $k$ we observe that several
values of $s$ occur (based on the equations in our database).

Condition~\ref{en:condAZ} forces the middle term (with $z$) to be of
the given form, as can be checked by computing the coefficients
$a_i(z)$ for a general polynomial $P_1(\theta)$.  In practice in most
cases the middle term factors and the operator is of the form
\begin{equation}
\label{eq:fourf}
\begin{split}
  D_{4,2}^{\text{fact}} &= \theta^4 -c z (v \theta + u)
    (v \theta +v-u) (A \theta^2 + A \theta + B) \\
  &\qquad - d z^2 (v \theta + u) (x \theta + w) (x \theta + 2x-w)
  (v \theta + 2v-u).
\end{split}
\end{equation}
Using this form reduces the number of parameters by one which makes
the computer search a lot faster.

As mentioned above, we also search for equations of lower order ($n=2$
or $n=3$) because they can be used as factors in the Hadamard product
to find interesting fourth order equations. We have searched for
operators of the following forms (see \cite{AZ})
\begin{align*}
  D_{2,2}^{\text{had}} &= \theta^2 - c(A\theta^2+A\theta+B) -
  d(\theta+1)^2 \\
\intertext{and}
  D_{3,2}^{\text{had}} &= \theta^3 -c(2\theta+1)(A\theta^2+A\theta+B)-
  d(\theta+1)^3.
\end{align*}

The special forms of the operator discussed above allow us to do a
computer search for Calabi--Yau equations. For fourth order equations
we use the following two step process.

\begin{step}
 A program written in C searches a certain range of the
  parameters in the expressions for the operators discussed above
  (e.g., $A$, $B$, $C$, $c$, $d$, $\frac{u}{v}$, and $\frac{w}{x}$ for
  $D_{4,2}^{\text{gen}}$). The C program only checks for an integral
  solution to the recursion relation. If an integral power series
  solution exists, the corresponding values of the parameters are
  written to a Maple array.
\end{step}

\begin{step}
 A Maple program reads the candidate parameter combinations and
  computes the instanton numbers. These need not be integral, because
  we can still multiply with the degree. Therefore, we compute the
  lowest common multiple of denominators of the first 15 and the first
  20 instanton numbers. If they are equal, we consider this to be
  acceptable. In practice we observe that in these cases the lowest
  common multiple is also very small (below $100$). We also check if
  the differential operator factors as a differential operator
  (using the \texttt{DEtools} package). If does, we remove the
  operator from our list.
\end{step}

For second and third order equations we cannot directly define the
instanton numbers. Instead, we check if using them as factors in a
Hadamard product yields a fourth order equation of Calabi--Yau type.

The first step is the most time consuming one. So we use a few tricks
to save some time. First note that replacing $z$ by $\lambda z$
multiplies the coefficients $A_n$ of the power series solution by
$\lambda^n$. In terms of the parameters this amounts to multiplying
$c$ and $d$ by $\lambda$ and $\lambda^2$ respectively. Therefore we
can put $c=1$ and compute the first $N$ coefficients of the power
series. In general the coefficients will be rational. We will require
the denominators to be a product of powers of a limited number of
small primes. We can then compute the scaling factor to make the
first $N$ coefficients integral. Using the rescaled parameters $c$ and
$d$ we can compute the first $M$ ($M>N$) coefficients and check that
they are integral. We always used $N=50$ and $M=1000$, but it never
happened that with the rescaled parameters we found non integral
coefficients. Note that we do not \emph{prove} that all the
coefficients of the power series solution are integral, but as we
check the first 1000 coefficients, we are pretty confident that they
are.

In this way we did extensive searches using up to a dozen computers.
The results are summarized in \pref{tab:compsearches}. The last column
contains the number of solutions found by the C program. For the
fourth order equations we also indicate between parentheses the number
of solutions that survive the second step of the searching process and
are contained in our database. For the spectra searched we refer to
\pref{tab:spectra}. The star indicates that we restrict to the subset
of the spectra from \sref{tab:spectra} labelled by a star.

The difference between the programs \texttt{search4} and
\texttt{search4q} is that the former is completely general, whereas
the latter assumes that apart from $z=0$ and $z=\infty$ there is just
one rational critical point. This yields big restrictions on the
parameters and allows us to compute $d$. Therefore this program is
much faster and we were able to search a much larger parameter
range. Of course it also means that we cannot find any equations which
do not satisfy this extra condition.

We do not claim that our search has been exhaustive. First of all the
parameter ranges are rather arbitrary. It is to be expected that
enlarging these ranges we will find more equations. Furthermore, the
expressions \pref{eq:four} and \pref{eq:fourf} that we use are not
completely general. Both the assumption that the coefficient
of $\theta^4$ is 1 and the restriction to $s=2$ are arbitrary. As
mentioned above, we already have one example with
$s=\frac{3}{2}$. There are also 32 equations where the coefficient of
$\theta^4$ is bigger than 1. However, in all these cases $k$ is at
least 4. All in all, this is still very much work in progress. We are
still adding equations to our list every now and then and we keep
enhancing the information on the equations in our database.

\begin{table}
\begin{center}
\setlength{\extrarowheight}{2pt}
\begin{tabular}{|>{\ttfamily}l|l|p{5cm}|l|p{8mm}|}
\hline
\multicolumn{1}{|l|}{Program} & Operator             & Parameter range
& Spectra & \# \\ \hline
search2  & $D_{2,2}^{\text{had}}$    & $1\le A,B \le 500$,
$d=2^p3^q5^r7^s$ with $0\le p \le 20$, $0\le q \le 6$, $0\le r \le 5$,
$0 \le s \le 2$ & --- & 126 \\ \hline
search3  & $D_{3,2}^{\text{had}}$    & $1\le A,B \le 500$, $d=2^p3^q5^r7^s$ with $0\le p \le 20$, $0\le q \le 6$, $0\le r \le 5$,
$0 \le s \le 2$ & --- & 139 \\ \hline
search4  & $D_{4,2}^{\text{gen}}$  & $1 \le A,B,C \le 200$, $d=2^p3^q5^r7^s$ with $0\le p \le 15$, $0\le q \le 5$, $0\le r \le 4$,
$0 \le s \le 2$ & $s=2^*$ & 66 (15) \\ \hline
search4f & $D_{4,2}^{\text{fact}}$ & $1 \le A,B \le 100$, $d=2^p3^q5^r7^s$ with $0\le p \le 15$, $0\le q \le 6$, $0\le r \le 5$,
$0 \le s \le 1$ & $s=2$ & 364 (58) \\ \hline
search4q & $D_{4,2}^{\text{gen}}$  & $1\le A,B,C \le 1000$, $d$ is
computed & $s=2$ & 1321 (30) \\ \hline
\end{tabular}
\end{center}
\caption{Computer searches}
\label{tab:compsearches}
\end{table}

\subsection{How to sum for $k=-n$ to $k=-1$}
\label{subsec:howtosum}
As already mentioned, the work \cite{AZ} explains
manipulating binomial expressions to deduce new examples
of Calabi--Yau equations. This is certainly a ``hypergeometric''
machinery, but, surprisingly, it works for certain unusual
series $y_0=\sum_{n=0}^\infty A_nz^n$ as well, namely,
when coefficients $A_n$ involve not only binomials but
also harmonic sums and consist of several summations.

In some formulas for $A_n$ (particular cases are \#211 and \#264)
it is necessary to sum also for $k=-n$ to $k=-1$
to get the correct coefficients of the solution $y_0$
to the differential equation. This requires expressing
binomial coefficients near $k=-\bar k$,
where $\bar k$~ is a positive integer:
\begin{align*}
\binom nk
&\approx\binom n{-\bar k-\varepsilon}
=\frac{n!}{\Gamma(1-\bar k-\varepsilon)\Gamma(1+n+\bar k+\varepsilon)}
\\
&\;\phantom{\approx\binom n{-\bar k-\varepsilon}}
=(-1)^{\bar k}\bar k^{-1}\binom{n+\bar k}n^{-1}\varepsilon+O(\varepsilon^2),
\displaybreak[2]\\
\binom{2k}k
&\approx\binom{-2\bar k-2\varepsilon}{-\bar k-\varepsilon}
=\frac{\Gamma(1-2\bar k-2\varepsilon)}{\Gamma(1-\bar k-\varepsilon)^2}
=\bar k^{-1}\binom{2\bar k}{\bar k}^{-1}\varepsilon+O(\varepsilon^2),
\displaybreak[2]\\
\binom{n+k}n
&\approx\binom{n-\bar k-\varepsilon}n
=\frac{\Gamma(1+n-\bar k-\varepsilon)}{n!\Gamma(1-\bar k-\varepsilon)}
=(-1)^{\bar k}{\bar k}^{-1}\binom n{\bar k}^{-1}\varepsilon
+O(\varepsilon^2),
\end{align*}
and so on.

The differential equation for case \#211 was found
in two ways: first, as the reflection of case \#210
at infinity and, secondly, by taking Maple's
\texttt{Zeilberger} on the following binomial sum:
\[
\text{``}A_n\text{''}
=\binom{2n}n^4\sum_k(-1)^{n+k}\binom nk^2\binom{2k}k
\binom{4n-2k}{2n-k}\binom{n+k}n^{-2}\binom{2n}k^{-1}.
\]
Usual summing for $k=0$ to $k=n$ does not give
the correct coefficient of $y_0$. To remove the
latter quotation one should sum for $k=-n$ to~$n$.
So we have to consider also negative $k.$
The 'difficult' (i.e., negative) part of the
summand is
\[
\binom nk^2\binom{2k}k\binom{n+k}n^{-2}\binom{2n}k^{-1},
\]
which near $k=-\bar k$ becomes (after simple reduction)
\[
(-1)^{\bar k}\binom{n+\bar k}n^{-2}\binom{2\bar k}{\bar k}^{-1}
\binom n{\bar k}^2\binom{2n+\bar k}{2n}
+O(\varepsilon),
\]
thus giving us the formula indicated below in Table~A.

In case \#264 we use \texttt{Zeilberger} on
\begin{align*}
\text{``}A_n\text{''}
&=16^{-n}\binom{2n}n^2\sum_k(n-2k)\binom nk\binom{2k}k
\binom{2n-2k}{n-k}\binom{2n+2k}{n+k}^2
\\ &\qquad\times
\binom{4n-2k}{2n-k}^2\binom{2n}k^{-1}\binom{2n}{n-k}^{-1}
\end{align*}
and proceed, as before, for negative~$k$ to obtain
for the `difficult' part
\[
\binom nk\binom{2k}k\binom{2n}k^{-1}
\]
near $k=-\bar k$ the following result:
\[
\bar k^{-1}\binom{n+\bar k}n^{-1}\binom{2\bar k}{\bar k}^{-1}
\binom{2n+\bar k}{2n}\varepsilon
+O(\varepsilon^2).
\]
Therefore, for negative $k$, we get the sum
(this formula was found by C.~Krattenthaler~\cite{Kr})
\begin{align*}
\frac{\d\text{``}A_n\text{''}}{\d\varepsilon}
&=16^{-n}\binom{2n}n^2\sum_{k=1}^n\frac{n+2k}k%
\binom{2n+k}{2n}\binom{2n+2k}{n+k}
\\ &\qquad\times
\binom{2n-2k}{n-k}^2\binom{4n+2k}{2n+k}^2
\binom{n+k}n^{-1}\binom{2k}k^{-1},
\end{align*}
where we replace $\bar k=-k$ by~$k$.
Finally, adding the derivative of ``$A_n$'' for positive $k$, we deduce
our formula in Table~A for case \#264.

\section{Electronic database of Calabi--Yau equations}
All equations we found are also available in electronic form. The web
interface to our database is at:
\[
  \text{\hlink{\dburl}}
\]
The database includes some extra information that is not contained in
the tables in this paper. Furthermore we hope to be able to keep it up
to date so that it can serve as a list of all known Calabi--Yau
equations. Therefore we encourage anybody who finds new equations of
Calabi--Yau type to send them to us so that we can include them into
our database. Comments and suggestions are also welcome. The database
and the web interface are still a work in progress.

Information about the use of the database can be found at the web
address mentioned above. A small warning for people comparing the
table in this article and the database: for equations past \#180 the
numbering of the equations in the database is different from the one
used in this article. However, by searching for \texttt{almkvist[n]}
in the source field, where $n$ is the number used in this article one
can easily find the equation in the database.

Part of the extra information that we do not provide in this article,
but is contained in the database, is related to numerical computations
of the monodromies around the singular points of the differential
equations. In many cases we succeeded in guessing the exact
monodromies using the approximate monodromies obtained numerically.
In turn the exact monodromies allowed us to compute some geometric
invariants of potential underlying Calabi--Yau manifolds. These, in
combination with further information obtained from the monodromies,
enabled us to compute elliptic invariants in many cases. This work is
discussed in detail in~\cite{ES}.

\begin{remark}
We feel that our table in \pref{app:CYlist} is not perfectly
organized, and we plan to systemize it in the nearest future.
On the other hand, the papers \cite{AZ}, \cite{ES},
as well as the electronic database have crucial references
to the absolute numbering in \pref{app:CYlist}.
Therefore, we have decided to keep the present
numeration at least in the first arXiv-version
of this work.
\end{remark}

\subsection*{Acknowledgments}
We would like to thank frankly C.~Krattenthaler who has
provided~\cite{Kr} us with several explicit formulas
for $A_n$ appearing now in Table~A. The work of the
fourth-named author was partially supported by grant
no.~03-01-00359 of the Russian Foundation for Basic Research,
but this was many years ago.

Since the appearance of the first version of the database in 2005,
we benefit a lot from valuable comments and numerous corrections.
We thank V.~Batyrev, M.~Bogner, D.~Broadhurst, T.~Coates, T.~Guttmann, M.~Kreuzer,
S.~Reiter, H.~Verrill, D.~Zagier and many others for their feedback and discussions.

\newpage

\appendix
\cleardoublepage
\section{Table of Calabi--Yau equations}
\label{app:CYlist}

The first column is reserved for numeration of cases.
The second column of the following table contains a 4th-order linear
differential operator~$D$ (by means of $\theta=z\frac{\d}{\d z}$).
The third column indicates coefficients
$\{A_n\}_{n=0,1,2,\dots}$ of the analytic solution
$y_0(z)=\sum_{n=0}^\infty A_nz^n$ to the MUM differential equation
$Dy=0$ normalized by the condition $y_0(0)=A_0=1$.

See \cite{Mo}, \cite{KT1} for cases \#1, \#2;
\cite{LT} for cases \#3--6;
\cite{KT1} cases \#7, \#8;
\cite{KT2} for cases \#10--14;
\cite{BS1} for cases \#15--23, and
\cite{BS2} for cases \#24--28.

\medskip
\hbox to\hsize{\hss\vbox{\offinterlineskip
\halign to120mm{\strut\tabskip=100pt minus 100pt
\strut\vrule\vphantom{\vrule height11.2pt}#&\hbox to5mm{\hfil$#$\,}&%
\vrule#&\hbox to75mm{\hfil$#$\hfil}&%
\vrule#&\hbox to39mm{\hfil$\dsize#$\hfil}&%
\vrule#\tabskip=0pt\cr\noalign{\hrule}
& \# &%
&\multispan3\hss differential operator $D$ and
coefficients $A_n$, $n=0,1,2,\dots$\hss &\cr
\noalign{\hrule\vskip1pt\hrule}
& \eqnno{1} &%
& D=\theta^4-5^5z\bigl(\theta+\frac15\bigr)\bigl(\theta+\frac25\bigr)
\bigl(\theta+\frac35\bigr)\bigl(\theta+\frac45\bigr) &%
& A_n=\frac{(5n)!}{n!^5} \ph &\cr
\noalign{\hrule}
& \eqnno{2} &%
& D=\theta^4-8\cdot10^5z
\bigl(\theta+\frac1{10}\bigr)\bigl(\theta+\frac3{10}\bigr)
\bigl(\theta+\frac7{10}\bigr)\bigl(\theta+\frac9{10}\bigr) &%
& A_n=\frac{(10n)!}{n!^3(2n)!(5n)!} \ph &\cr
\noalign{\hrule}
& \eqnno{3} &%
& D=\theta^4-256z\bigl(\theta+\frac12\bigr)^4 &%
& A_n=\binom{2n}n^4 \ph &\cr
\noalign{\hrule}
& \eqnno{4} &%
& D=\theta^4-3^6z\bigl(\theta+\frac13\bigr)^2\bigl(\theta+\frac23\bigr)^2 &%
& A_n=\biggl(\frac{(3n)!}{n!^3}\biggr)^2 &\cr
\noalign{\hrule}
& \eqnno{5} &%
& D=\theta^4-432z\bigl(\theta+\frac12\bigr)^2
\bigl(\theta+\frac13\bigr)\bigl(\theta+\frac23\bigr) &%
& A_n=\binom{2n}n^2\frac{(3n)!}{n!^3} \ph &\cr
\noalign{\hrule}
& \eqnno{6} &%
& D=\theta^4-2^{10}z\bigl(\theta+\frac12\bigr)^2
\bigl(\theta+\frac14\bigr)\bigl(\theta+\frac34\bigr) &%
& A_n=\binom{2n}n\frac{(4n)!}{n!^4} \ph &\cr
\noalign{\hrule}
& \eqnno{7} &%
& D=\theta^4-2^{16}z\bigl(\theta+\frac18\bigr)\bigl(\theta+\frac38\bigr)
\bigl(\theta+\frac58\bigr)\bigl(\theta+\frac78\bigr) &%
& A_n=\frac{(8n)!}{n!^4(4n)!} \ph &\cr
\noalign{\hrule}
& \eqnno{8} &%
& D=\theta^4-11664z\bigl(\theta+\frac16\bigr)\bigl(\theta+\frac13\bigr)
\bigl(\theta+\frac23\bigr)\bigl(\theta+\frac56\bigr) &%
& A_n=\frac{(6n)!}{n!^4(2n)!} \ph &\cr
\noalign{\hrule}
& \eqnno{9} &%
& D=\theta^4-12^6z\bigl(\theta+\frac1{12}\bigr)\bigl(\theta+\frac5{12}\bigr)
\bigl(\theta+\frac7{12}\bigr)\bigl(\theta+\frac{11}{12}\bigr) &%
& A_n=\binom{2n}n\frac{(12n)!}{n!^2(4n)!(6n)!} \ph &\cr
\noalign{\hrule}
& \eqnno{10} &%
& D=\theta^4-2^{12}z\bigl(\theta+\frac14\bigr)^2\bigl(\theta+\frac34\bigr)^2 &%
& A_n=\biggl(\frac{(4n)!}{n!^2(2n)!}\biggr)^2 &\cr
\noalign{\hrule}
& \eqnno{11} &%
& D=\theta^4-12^3z\bigl(\theta+\frac14\bigr)\bigl(\theta+\frac34\bigr)
\bigl(\theta+\frac13\bigr)\bigl(\theta+\frac23\bigr) &%
& A_n=\binom{3n}n\frac{(4n)!}{n!^4} \ph &\cr
\noalign{\hrule}
& \eqnno{12} &%
& D=\theta^4-2^{10}\cdot3^3z
\bigl(\theta+\frac14\bigr)\bigl(\theta+\frac34\bigr)
\bigl(\theta+\frac16\bigr)\bigl(\theta+\frac56\bigr) &%
& A_n=\binom{4n}n\frac{(6n)!}{n!^2(2n)!^2} \ph &\cr
\noalign{\hrule}
& \eqnno{13} &%
& D=\theta^4-2^8\cdot3^6z
\bigl(\theta+\frac16\bigr)^2\bigl(\theta+\frac56\bigr)^2 &%
& A_n=\biggl(\frac{(6n)!}{n!(2n)!(3n)!}\biggr)^2 &\cr
\noalign{\hrule}
& \eqnno{14} &%
& D=\theta^4-2^8\cdot3^3z\bigl(\theta+\frac12\bigr)^2
\bigl(\theta+\frac16\bigr)\bigl(\theta+\frac56\bigr) &%
& A_n=\binom{2n}n\frac{(6n)!}{n!^3(3n)!} \ph &\cr
\noalign{\hrule}
}}\hss}

\newpage

\hbox to\hsize{\hss\vbox{\offinterlineskip
\halign to120mm{\strut\tabskip=100pt minus 100pt
\strut\vrule#&\hbox to5.5mm{\hss$#$\hss}&%
\vrule#&\hbox to114mm{\hfil$\dsize#$\hfil}&%
\vrule#\tabskip=0pt\cr\noalign{\hrule}
& \# &%
& \text{differential operator $D$ and coefficients $A_n$, $n=0,1,2,\dots$} &\cr
\noalign{\hrule\vskip1pt\hrule}
& \eqnno{15} &&
\aligned
\\[-12pt]
&
D=\theta^4-3z(3\theta+1)(3\theta+2)(7\theta^2+7\theta+2)
\\[-2pt] &\;
-72z^2(3\theta+1)(3\theta+2)(3\theta+4)(3\theta+5)
\\[-8pt]
\endaligned &\cr
& &\span\hrulefill&\cr
& &&
A_n=\frac{(3n)!}{n!^3}\sum_{k=0}^n\dbinom nk^3
&\cr
\noalign{\hrule}
& \eqnno{16} &&
\aligned
\\[-12pt]
&
D=\theta^4-4z(2\theta+1)^2(5\theta^2+5\theta+2)
+2^8z^2(\theta+1)^2(2\theta+1)(2\theta+3)
\\[-8pt]
\endaligned &\cr
& &\span\hrulefill&\cr
& &&
A_n
=\binom{2n}n\sum_{j+k+l+m=n}
\biggl(\frac{n!}{j!\,k!\,l!\,m!}\biggr)^2
=\binom{2n}n\sum_{k=0}^n
\binom nk^2\binom{2k}k\binom{2n-2k}{n-k}
&\cr
\noalign{\hrule}
& \eqnno{17} &&
\aligned
\\[-12pt]
&
D=25\theta^4-15z(51\theta^4+84\theta^3+72\theta^2+30\theta+5)
\\[-2pt] &\;
+6z^2(531\theta^4+828\theta^3+541\theta^2+155\theta+15)
\\[-2pt] &\;
-54z^3(423\theta^4+2160\theta^3+4399\theta^2+3795\theta+1170)
\\[-2pt] &\;
+3^5z^4(279\theta^4+1368\theta^3+2270\theta^2+1586\theta+402)
-3^{10}z^5(\theta+1)^4
\\[-8pt]
\endaligned &\cr
& &\span\hrulefill&\cr
& &&
A_n=\sum_{j+k+l=n}
\biggl(\frac{n!}{j!\,k!\,l!}\biggr)^3
&\cr
\noalign{\hrule}
& \eqnno{18} &&
\aligned
\\[-12pt]
&
D=\theta^4-4z(2\theta+1)^2(3\theta^2+3\theta+1)
-4z^2(4\theta+2)(4\theta+3)(4\theta+5)(4\theta+6)
\\[-8pt]
\endaligned &\cr
& &\span\hrulefill&\cr
& &&
A_n=\binom{2n}n\sum_{k=0}^n\binom nk^4
&\cr
\noalign{\hrule}
& \eqnno{19} &&
\aligned
\\[-12pt]
&
D=23^2\theta^4
-23z(921\theta^4+2046\theta^3+1644\theta^2+621\theta+92)
\\[-2pt] &\;
-z^2(380851\theta^4+1328584\theta^3+1772673\theta^2+1033528\theta+221168)
\\[-2pt] &\;
-2z^3(475861\theta^4+1310172\theta^3+1028791\theta^2+208932\theta-27232)
\\[-2pt] &\;
-2^2\cdot17z^4(8873\theta^4+14020\theta^3+5139\theta^2-1664\theta-976)
\\[-2pt] &\;
+2^3\cdot3\cdot17^2z^5(\theta+1)^2(3\theta+2)(3\theta+4)
\\[-8pt]
\endaligned &\cr
& &\span\hrulefill&\cr
& &&
A_n=\sum_{k=0}^n\binom nk^3\binom{n+k}n\binom{2n-k}n
&\cr
\noalign{\hrule}
& \eqnno{20} &&
\aligned
\\[-12pt]
&
D=\theta^4
-3z(48\theta^4+60\theta^3+53\theta^2+23\theta+4)
\\[-2pt] &\;
+9z^2(873\theta^4+1980\theta^3+2319\theta^2+1344\theta+304)
\\[-2pt] &\;
-2\cdot3^4z^3(1269\theta^4+3888\theta^3+5259\theta^2+3348\theta+800)
\\[-2pt] &\;
+2^2\cdot3^6z^4(891\theta^4+3240\theta^3+4653\theta^2+2952\theta+688)
\\[-2pt] &\;
-2^3\cdot3^{11}z^5(\theta+1)^2(3\theta+2)(3\theta+4)
\\[-8pt]
\endaligned &\cr
& &\span\hrulefill&\cr
& &&
A_n=\sum_{k=0}^n\binom nk\frac{(3k)!}{k!^3}
\,\frac{(3n-3k)!}{(n-k)!^3}
&\cr
\noalign{\hrule}
}}\hss}

\newpage

\hbox to\hsize{\hss\vbox{\offinterlineskip
\halign to120mm{\strut\tabskip=100pt minus 100pt
\strut\vrule#&\hbox to5.5mm{\hss$#$\hss}&%
\vrule#&\hbox to114mm{\hfil$\dsize#$\hfil}&%
\vrule#\tabskip=0pt\cr\noalign{\hrule}
& \# &%
& \text{differential operator $D$ and coefficients $A_n$, $n=0,1,2,\dots$} &\cr
\noalign{\hrule\vskip1pt\hrule}
& \eqnno{21} &&
\aligned
\\[-12pt]
&
D=25\theta^4
-20z(36\theta^4+84\theta^3+72\theta^2+30\theta+5)
\\[-2pt] &\;
-2^4z^2(181\theta^4+268\theta^3+71\theta^2-70\theta-35)
\\[-2pt] &\;
+2^8z^3(\theta+1)(37\theta^3+248\theta^2+375\theta+165)
\\[-2pt] &\;
+2^{10}z^4(39\theta^4+198\theta^3+331\theta^2+232\theta+59)
+2^{15}z^5(\theta+1)^4
\\[-8pt]
\endaligned &\cr
& &\span\hrulefill&\cr
& &&
A_n=\sum_{k=0}^n\binom nk^3\binom{2k}k\binom{2n-2k}{n-k}
&\cr
\noalign{\hrule}
& \eqnno{22} &&
\aligned
\\[-12pt]
&
D=49\theta^4
-7z(155\theta^4+286\theta^3+234\theta^2+91\theta+14)
\\[-2pt] &\;
-z^2(16105\theta^4+68044\theta^3+102261\theta^2+66094\theta+15736)
\\[-2pt] &\;
+2^3z^3(2625\theta^4+8589\theta^3+9071\theta^2+3759\theta+476)
\\[-2pt] &\;
-2^4z^4(465\theta^4+1266\theta^3+1439\theta^2+806\theta+184)
+2^9z^5(\theta+1)^4
\\[-8pt]
\endaligned &\cr
& &\span\hrulefill&\cr
& &&
A_n=\sum_{k=0}^n\binom nk^5
&\cr
\noalign{\hrule}
& \eqnno{23} &&
\aligned
\\[-12pt]
&
D=9\theta^4
-12z(64\theta^4+80\theta^3+73\theta^2+33\theta+6)
\\[-2pt] &\;
+2^7z^2(194\theta^4+440\theta^3+527\theta^2+315\theta+75)
\\[-2pt] &\;
-2^{12}z^3(94\theta^4+288\theta^3+397\theta^2+261\theta+66)
\\[-2pt] &\;
+2^{17}z^4(22\theta^4+80\theta^3+117\theta^2+77\theta+19)
-2^{23}z^5(\theta+1)^4
\\[-8pt]
\endaligned &\cr
& &\span\hrulefill&\cr
& &&
A_n=\sum_{k=0}^n\binom nk\binom{2k}k^2\binom{2n-2k}{n-k}^2
&\cr
\noalign{\hrule}
& \eqnno{24} &&
\aligned
\\[-12pt]
&
D=\theta^4-3z(3\theta+2)(3\theta+1)(11\theta^2+11\theta+3)
\\[-2pt] &\;
-9z^2(3\theta+5)(3\theta+2)(3\theta+4)(3\theta+1)
\\[-8pt]
\endaligned &\cr
& &\span\hrulefill&\cr
& &&
A_n=\frac{(3n)!}{n!^3}\sum_{k=0}^n\binom nk^2\binom{n+k}k
&\cr
\noalign{\hrule}
& \eqnno{25} &&
\aligned
\\[-12pt]
&
D=\theta^4-4z(2\theta+1)^2(11\theta^2+11\theta+3)
-16z^2(2\theta+3)^2(1+2\theta)^2
\\[-8pt]
\endaligned &\cr
& &\span\hrulefill&\cr
& &&
A_n=\binom{2n}n^2\sum_{k=0}^n\binom nk^2\binom{n+k}k
&\cr
\noalign{\hrule}
& \eqnno{26} &&
\aligned
\\[-12pt]
&
D=\theta^4-2z(2\theta+1)^2(13\theta^2+13\theta+4)
\\[-2pt] &\;
-12z^2(2\theta+3)(2\theta+1)(3\theta+2)(3\theta+4)
\\[-8pt]
\endaligned &\cr
& &\span\hrulefill&\cr
& &&
\aligned
A_n
&=\binom{2n}n\sum_k
\binom nk^2\binom{n+k}k\binom{2k}n
\\[-2pt]
&=\binom{2n}n\sum_{k,l}\binom nk^2\binom nl^2
\binom{k+l}n
\\[2pt]
\endaligned &\cr
\noalign{\hrule}
}}\hss}

\newpage

\hbox to\hsize{\hss\vbox{\offinterlineskip
\halign to120mm{\strut\tabskip=100pt minus 100pt
\strut\vrule#&\hbox to5.5mm{\hss$#$\hss}&%
\vrule#&\hbox to114mm{\hfil$\dsize#$\hfil}&%
\vrule#\tabskip=0pt\cr\noalign{\hrule}
& \# &%
& \text{differential operator $D$ and coefficients $A_n$, $n=0,1,2,\dots$} &\cr
\noalign{\hrule\vskip1pt\hrule}
& \eqnno{27} &&
\aligned
\\[-12pt]
&
D=9\theta^4
-3z(173\theta^4+340\theta^3+272\theta^2+102\theta+15)
\\[-2pt] &\;
-2z^2(1129\theta^4+5032\theta^3+7597\theta^2+4773\theta+1083)
\\[-2pt] &\;
+2z^3(843\theta^4+2628\theta^3+2353\theta^2+675\theta+6)
\\[-2pt] &\;
-z^4(295\theta^4+608\theta^3+478\theta^2+174\theta+26)
+z^5(\theta+1)^4
\\[-8pt]
\endaligned &\cr
& &\span\hrulefill&\cr
& &&
A_n=\sum_{k,l}\binom nk^2\binom nl^2\binom{k+l}n\binom{2n-k}n
&\cr
\noalign{\hrule}
& \eqnno{28} &&
\aligned
\\[-12pt]
&
D=\theta^4-z(65\theta^4+130\theta^3+105\theta^2+40\theta+6)
\\[-2pt] &\;
+4z^2(4\theta+3)(4\theta+5)(\theta+1)^2
\\[-8pt]
\endaligned &\cr
& &\span\hrulefill&\cr
& &&
A_n=\sum_{k,l,m}\binom nk\binom nl\binom mk\binom ml
\binom{k+l}k\binom nm^2
=\sum_{k,l}\binom nk^2\binom nl^2\binom{k+l}n^2
&\cr
\noalign{\hrule}
& \eqnno{29} &&
\aligned
\\[-12pt]
&
D=\theta^4-2z(2\theta+1)^2(17\theta^2+17\theta+5)
\\[-2pt] &\;
+4z^2(2\theta+3)(2\theta+1)(\theta+1)^2
\\[-8pt]
\endaligned &\cr
& &\span\hrulefill&\cr
& &&
A_n=\binom{2n}n\sum_{k=0}^n\binom nk^2\binom{n+k}k^2
&\cr
\noalign{\hrule}
}}\hss}

\begin{remark}
The original formulas for $A_n$ in cases \#24, 25, 27--29
were simplified with the help of the binomial identity
$$
\sum_j\binom nj\binom kj=\binom{n+k}k.
$$
\end{remark}

\newpage

The origin of the following new cases is explained
in the introductory part and \cite{AZ,ES}.
In particular, cases \#\eqnno{30} and \#\eqnno{31} are obtained
in \cite{AZ}, Section~6; case \#\eqnno{32} (the pullback of the differential
equation related to~$\zeta(4)$) is the subject of \cite{AZ}, Section~4;
case \#\eqnno{130} was communicated to us by H.~Verrill~\cite{Ve},
case \#\eqnnol{186}{247} (except for an explicit formula
for~$A_n$'s) by Tj{\o}tta~\cite{Tj},
and case \#\eqnnoltmp{366}{426} by T.~Guttmann~\cite{Gu}.
Entries \#\eqnnoltmp{369}{429}--\#\eqnnoltmp{372}{432}
correspond to the so-called \emph{Hurwitz product},
$(A_n)\circ(B_n)=\sum_k\binom nkA_kB_{n-k}$
(actually, the Hurwitz square) and are due to M.~Bogner.
Equations \#\eqnnoltmp{374}{434}, \#\eqnnoltmp{375}{435}, and \#\eqnnoltmp{376}{436}
were found by P.~Metelitsyn.
Equations \#\eqnnoltmp{386}{446}--\eqnnoltmp{393}{453} were communicated
to us by M.~Bogner and S.~Reiter~\cite{BR},
while equations \#\eqnnoltmp{402}{462}--\eqnnoltmp{404}{464} are due to T.~Coates~\cite{Co}.

\hbox to\hsize{\hss\vbox{\offinterlineskip
\halign to120mm{\strut\tabskip=100pt minus 100pt
\strut\vrule#&\hbox to5.5mm{\hss$#$\hss}&%
\vrule#&\hbox to114mm{\hfil$\dsize#$\hfil}&%
\vrule#\tabskip=0pt\cr\noalign{\hrule}
& \# &%
& \text{differential operator $D$ and coefficients $A_n$, $n=0,1,2,\dots$} &\cr
\noalign{\hrule\vskip1pt\hrule}
& \eqnnol{\hbox{$3^*$}}{35} &&
\aligned
\\[-12pt]
&
D=\theta^4
-2^4z(2\theta+1)^2(2\theta^2+2\theta+1)
+2^{10}z^2(\theta+1)^2(2\theta+1)(2\theta+3)
\\[-8pt]
\endaligned &\cr
& &\span\hrulefill&\cr
& &&
A_n
=16^n\binom{2n}n\sum_k\binom{-1/2}k^2\binom{-1/2}{n-k}^2
=\binom{2n}n\sum_{k=0}^n\binom{2k}k^2\binom{2n-2k}{n-k}^2
&\cr
\noalign{\hrule}
& \eqnnol{\hbox{$4^*$}}{186} &&
\aligned
\\[-12pt]
&
D=\theta^4-6z(2\theta+1)^2(9\theta^2+9\theta+5)
+2^2\cdot3^6z^2(\theta+1)^2(2\theta+1)(2\theta+3)
\\[-8pt]
\endaligned &\cr
& &\span\hrulefill&\cr
& &&
A_n=27^n\binom{2n}n\sum_k\binom{-1/3}k^2\binom{-2/3}{n-k}^2
&\cr
\noalign{\hrule}
& \eqnnol{\hbox{$4^{**}$}}{46} &&
\aligned
\\[-12pt]
&
D=\theta^4-6z(2\theta+1)^2(9\theta^2+9\theta+4)
\\[-2pt] &\;
+2^2\cdot3^4z^2(2\theta+1)(2\theta+3)(3\theta+2)(3\theta+4)
\\[-8pt]
\endaligned &\cr
& &\span\hrulefill&\cr
& &&
\aligned
\\[-12pt]
A_n
&=27^n\binom{2n}n\sum_k\binom{-1/3}k\binom{-2/3}k
\binom{-1/3}{n-k}\binom{-2/3}{n-k}
\\[-2pt]
&=\binom{2n}n\sum_{k=0}^n\frac{(3k)!}{k!^3}\,\frac{(3n-3k)!}{(n-k)!^3}
\vspace{1.5pt}
\endaligned &\cr
\noalign{\hrule}
& \eqnnol{\hbox{$6^*$}}{220} &&
\aligned
\\[-12pt]
&
D=\theta^4-2^5z(2\theta+1)^2(2\theta^2+2\theta+1)
\\[-2pt] &\;
+2^8z^2(2\theta+1)(2\theta+3)(4\theta+3)(4\theta+5)
\\[-8pt]
\endaligned &\cr
& &\span\hrulefill&\cr
& &&
A_n=32^n\binom{2n}n\sum_k\binom{-1/2}k\binom{-1/4}k
\binom{-1/2}{n-k}\binom{-3/4}{n-k}
&\cr
\noalign{\hrule}
& \eqnnol{\hbox{$7^*$}}{192} &&
\aligned
\\[-12pt]
&
D=\theta^4-2^4z(2\theta+1)^2(32\theta^2+32\theta+19)
\\[-2pt] &\;
+2^{14}z^2(2\theta+1)(2\theta+3)(4\theta+3)(4\theta+5)
\\[-8pt]
\endaligned &\cr
& &\span\hrulefill&\cr
& &&
A_n=256^n\binom{2n}n\sum_k\binom{-1/8}k\binom{-3/8}k
\binom{-5/8}{n-k}\binom{-7/8}{n-k}
&\cr
\noalign{\hrule}
& \eqnnol{\hbox{$7^{**}$}}{189} &&
\aligned
\\[-12pt]
&
D=\theta^4-2^4z(2\theta+1)^2(32\theta^2+32\theta+13)
+2^{16}z^2(2\theta+1)^2(2\theta+3)^2
\\[-8pt]
\endaligned &\cr
& &\span\hrulefill&\cr
& &&
A_n=256^n\binom{2n}n\sum_k\binom{-1/8}k\binom{-5/8}k
\binom{-3/8}{n-k}\binom{-7/8}{n-k}
&\cr
\noalign{\hrule}
}}\hss}

\newpage

\hbox to\hsize{\hss\vbox{\offinterlineskip
\halign to120mm{\strut\tabskip=100pt minus 100pt
\strut\vrule#&\hbox to6.5mm{\hss$#$\hss}&%
\vrule#&\hbox to113mm{\hfil$\dsize#$\hfil}&%
\vrule#\tabskip=0pt\cr\noalign{\hrule}
& \# &%
& \text{differential operator $D$ and coefficients $A_n$, $n=0,1,2,\dots$} &\cr
\noalign{\hrule\vskip1pt\hrule}
& \eqnnol{\hbox{$8^*$}}{193} &&
\aligned
\\[-12pt]
&
D=\theta^4-12z(2\theta+1)^2(18\theta^2+18\theta+11)
\\[-2pt] &\;
+2^4\cdot3^4z^2(2\theta+1)(2\theta+3)(6\theta+5)(6\theta+7)
\\[-8pt]
\endaligned &\cr
& &\span\hrulefill&\cr
& &&
A_n=108^n\binom{2n}n\sum_k\binom{-1/6}k\binom{-1/3}k
\binom{-5/6}{n-k}\binom{-2/3}{n-k}
&\cr
\noalign{\hrule}
& \eqnnol{\hbox{$8^{**}$}}{141} &&
\aligned
\\[-12pt]
&
D=\theta^4-12z(2\theta+1)^2(18\theta^2+18\theta+7)
+2^4\cdot3^6z^2(2\theta+1)^2(2\theta+3)^2
\\[-8pt]
\endaligned &\cr
& &\span\hrulefill&\cr
& &&
A_n=108^n\binom{2n}n\sum_k\binom{-1/6}k\binom{-2/3}k
\binom{-5/6}{n-k}\binom{-1/3}{n-k}
&\cr
\noalign{\hrule}
& \eqnnol{\hbox{$9^*$}}{191} &&
\aligned
\\[-12pt]
&
D=\theta^4-48z(2\theta+1)^2(72\theta^2+72\theta+41)
\\[-2pt] &\;
+2^{14}\cdot3^4z^2(2\theta+1)(2\theta+3)(3\theta+2)(3\theta+4)
\\[-8pt]
\endaligned &\cr
& &\span\hrulefill&\cr
& &&
A_n=1728^n\binom{2n}n\sum_k\binom{-1/12}k\binom{-5/12}k
\binom{-7/12}{n-k}\binom{-11/12}{n-k}
&\cr
\noalign{\hrule}
& \eqnnol{\hbox{$9^{**}$}}{190} &&
\aligned
\\[-12pt]
&
D=\theta^4-48z(2\theta+1)^2(72\theta^2+72\theta+31)
+2^{12}\cdot3^6z^2(2\theta+1)^2(2\theta+3)^2
\\[-8pt]
\endaligned &\cr
& &\span\hrulefill&\cr
& &&
A_n=432^n\binom{2n}n^2\sum_k(-1)^k\binom{-5/6}k\binom{-1/6}{n-k}^2
&\cr
\noalign{\hrule}
& \eqnnol{\hbox{$10^*$}}{185} &&
\aligned
\\[-12pt]
&
D=\theta^4-2^4z(2\theta+1)^2(8\theta^2+8\theta+5)
+2^{14}z^2(\theta+1)^2(2\theta+1)(2\theta+3)
\\[-8pt]
\endaligned &\cr
& &\span\hrulefill&\cr
& &&
A_n=64^n\binom{2n}n\sum_k\binom{-1/4}k^2\binom{-3/4}{n-k}^2
&\cr
\noalign{\hrule}
& \eqnnol{\hbox{$10^{**}$}}{111} &&
\aligned
\\[-12pt]
&
D=\theta^4-2^4z(2\theta+1)^2(8\theta^2+8\theta+3)
+2^{12}z^2(2\theta+1)^2(2\theta+3)^2
\\[-8pt]
\endaligned &\cr
& &\span\hrulefill&\cr
& &&
\aligned
\\[-12pt]
A_n
&=64^n\binom{2n}n\sum_k\binom{-1/4}k\binom{-3/4}k
\binom{-3/4}{n-k}\binom{-1/4}{n-k}
\\[-2pt]
&=\binom{2n}n\sum_{k=0}^n\frac{(4k)!}{k!^2(2k)!}
\,\frac{(4n-4k)!}{(n-k)!^2(2n-2k)!}
\\[2pt]
\endaligned &\cr
\noalign{\hrule}
}}\hss}

\newpage

\hbox to\hsize{\hss\vbox{\offinterlineskip
\halign to120mm{\strut\tabskip=100pt minus 100pt
\strut\vrule#&\hbox to6.5mm{\hss$#$\hss}&%
\vrule#&\hbox to113mm{\hfil$\dsize#$\hfil}&%
\vrule#\tabskip=0pt\cr\noalign{\hrule}
& \# &%
& \text{differential operator $D$ and coefficients $A_n$, $n=0,1,2,\dots$} &\cr
\noalign{\hrule\vskip1pt\hrule}
& \eqnnol{\hbox{$13^*$}}{188} &&
\aligned
\\[-12pt]
&
D=\theta^4-48z(2\theta+1)^2(18\theta^2+18\theta+13)
\\[-2pt] &\;
+2^{10}\cdot3^6z^2(\theta+1)^2(2\theta+1)(2\theta+3)
\\[-8pt]
\endaligned &\cr
& &\span\hrulefill&\cr
& &&
A_n=432^n\binom{2n}n\sum_k\binom{-1/6}k^2\binom{-5/6}{n-k}^2
&\cr
\noalign{\hrule}
& \eqnnol{\hbox{$13^{**}$}}{47} &&
\aligned
\\[-12pt]
&
D=\theta^4-48z(2\theta+1)^2(18\theta^2+18\theta+5)
\\[-2pt] &\;
+2^{10}\cdot3^4z^2(2\theta+1)(2\theta+3)(3\theta+1)(3\theta+5)
\\[-8pt]
\endaligned &\cr
& &\span\hrulefill&\cr
& &&
\aligned
\\[-12pt]
A_n
&=432^n\binom{2n}n\sum_k\binom{-1/6}k\binom{-5/6}k
\binom{-5/6}{n-k}\binom{-1/6}{n-k}
\\[-2pt]
&=\binom{2n}n\sum_{k=0}^n\frac{(6k)!}{k!(2k)!(3k)!}
\,\frac{(6n-6k)!}{(n-k)!(2n-k)!(3n-3k)!}
\vspace{1.5pt}
\endaligned &\cr
\noalign{\hrule}
& \eqnnol{\hbox{$\wh1$}}{206} &&
\aligned
\\[-12pt]
&
D=\theta^4-5z(10000\theta^4+12500\theta^3+9500\theta^2+3250\theta+399)
\\[-2pt] &\;
+5^8z^2(2400\theta^4+6000\theta^3+6290\theta^2+2800\theta+399)
\\[-2pt] &\;
-4\cdot 5^{14}z^3(4\theta+3)(80\theta^3+240\theta^2+221\theta+42)
\\[-2pt] &\;
+5^{20}z^4(4\theta+1)(4\theta+3)(4\theta+7)(4\theta+9)
\\[-8pt]
\endaligned &\cr
& &\span\hrulefill&\cr
& &&
\gathered
\text{the pullback of the 5th-order differential equation $D'y=0$, where}
\\[-8pt]
\endgathered &\cr
& &\span\hrulefill&\cr
& &&
\aligned
\\[-12pt]
&
D'=\theta^5
-10z(2\theta+1)(3125\theta^4+6250\theta^3+7875\theta^2+4750\theta+1226)
\\[-2pt] &\;
+2^45^6z^2(\theta+1)(6250\theta^4+25000\theta^3+4550\theta^2+41000\theta+14851)
\\[-2pt] &\;
-2^65^{14}z^3(\theta+1)(\theta+2)(2\theta+3)(25\theta^2+75\theta+82)
\\[-2pt] &\;
+2^85^{19}z^4(\theta+1)(\theta+2)(\theta+3)(25\theta^2+100\theta+113)
\\[-2pt] &\;
-2^95^{25}z^5(\theta+1)(\theta+2)(\theta+3)(\theta+4)(2\theta+5)
\\[-8pt]
\endaligned &\cr
& &\span\hrulefill&\cr
& &&
A_n'=(4\cdot 5^5)^n\sum_k(-1)^k\binom nk
\frac{(1/2)_k(1/5)_k(2/5)_k(3/5)_k(4/5)_k}{k!^5}
&\cr
\noalign{\hrule}
}}\hss}

\newpage

\hbox to\hsize{\hss\vbox{\offinterlineskip
\halign to120mm{\strut\tabskip=100pt minus 100pt
\strut\vrule#&\hbox to6.5mm{\hss$#$\hss}&%
\vrule#&\hbox to113mm{\hfil$\dsize#$\hfil}&%
\vrule#\tabskip=0pt\cr\noalign{\hrule}
& \# &%
& \text{differential operator $D$ and coefficients $A_n$, $n=0,1,2,\dots$} &\cr
\noalign{\hrule\vskip1pt\hrule}
& \eqnnol{\hbox{$\wh2$}}{207} &&
\aligned
\\[-12pt]
&
D=\theta^4
-2^45z(160000\theta^4+200000\theta^3+154000\theta^2+54000\theta+7189)
\\[-2pt] &\;
+2^{14}5^8z^2(9600\theta^4+24000\theta^3+25320\theta^2+1400\theta+1669)
\\[-2pt] &\;
-2^{24}5^{14}z^3(4\theta+3)(320\theta^3+960\theta^2+888\theta+171)
\\[-2pt] &\;
+2^{32}5^{20}z^4(4\theta+1)(4\theta+3)(4\theta+7)(4\theta+9)
\\[-8pt]
\endaligned &\cr
& &\span\hrulefill&\cr
& &&
\gathered
\text{the pullback of the 5th-order differential equation $D'y=0$, where}
\\[-8pt]
\endgathered &\cr
& &\span\hrulefill&\cr
& &&
\aligned
\\[-12pt]
&
D'=\theta^5
\\[-2pt] &\;
-2^55z(2\theta+1)(50000\theta^4+100000\theta^3+127000\theta^2+77000\theta+19811)
\\[-2pt] &\;
+2^{16}5^6z^2(\theta+1)(100000\theta^4+400000\theta^3+731000\theta^2+662000\theta+240811)
\\[-2pt] &\;
-2^{28}5^{14}z^3(\theta+1)(\theta+2)(2\theta+3)(100\theta^2+300\theta+331)
\\[-2pt] &\;
+2^{39}5^{19}z^4(\theta+1)(\theta+2)(\theta+3)(50\theta^2+200\theta+227)
\\[-2pt] &\;
-2^{49}5^{25}z^5(\theta+1)(\theta+2)(\theta+3)(\theta+4)(2\theta+5)
\\[-8pt]
\endaligned &\cr
& &\span\hrulefill&\cr
& &&
A_n'=3200000^n\sum_k(-1)^k\binom nk
\frac{(1/2)_k(1/10)_k(3/10)_k(7/10)_k(9/10)_k}{k!^5}
&\cr
\noalign{\hrule}
& \eqnnol{\hbox{$\wh3$}}{208} &&
\aligned
\\[-12pt]
&
D=\theta^4-2^4z(4\theta+1)(64\theta^3+64\theta^2+44\theta+9)
\\[-2pt] &\;
+2^{14}z^2(4\theta+1)(96\theta^3+216\theta^2+196\theta+61)
\\[-2pt] &\;
-2^{24}z^3(4\theta+1)(4\theta+3)(16\theta^2+44\theta+33)
\\[-2pt] &\;
+2^{32}z^4(4\theta+1)(4\theta+3)(4\theta+7)(4\theta+9)
\\[-8pt]
\endaligned &\cr
& &\span\hrulefill&\cr
& &&
\gathered
\text{the pullback of the 5th-order differential equation $D'y=0$, where}
\\[-8pt]
\endgathered &\cr
& &\span\hrulefill&\cr
& &&
\aligned
\\[-12pt]
&
D'=\theta^5
-2^5z(2\theta+1)(80\theta^4+160\theta^3+200\theta^2+120\theta+31)
\\[-2pt] &\;
+5\cdot 2^{16}z^2(\theta+1)(32\theta^4+128\theta^3+232\theta^2+208\theta+75)
\\[-2pt] &\;
-5\cdot 2^{28}z^3(\theta+1)(\theta+2)(2\theta+3)(4\theta^2+12\theta+13)
\\[-2pt] &\;
+5\cdot 2^{37}z^4(\theta+1)(\theta+2)(\theta+3)(2\theta^2+8\theta+9)
\\[-2pt] &\;
-2^{49}z^5(\theta+1)(\theta+2)(\theta+3)(\theta+4)(2\theta+5)
\\[-8pt]
\endaligned &\cr
& &\span\hrulefill&\cr
& &&
A_n'=1024^n\sum_k(-1)^k\binom nk\frac{(1/2)_k^5}{k!^5}
&\cr
\noalign{\hrule}
& \eqnnol{\hbox{$\wh4$}}{209} &&
\aligned
\\[-12pt]
&
D=\theta^4-3^2z(1296\theta^4+1620\theta^3+1224\theta^2+414\theta+49)
\\[-2pt] &\;
+3^8z^2(7776\theta^4+19440\theta^3+20322\theta^2+9000\theta+1267)
\\[-2pt] &\;
-4\cdot 3^{16}z^3(4\theta+3)(144\theta^3+432\theta^2+397\theta+75)
\\[-2pt] &\;
+3^{24}z^4(4\theta+1)(4\theta+3)(4\theta+7)(4\theta+9)
\\[-8pt]
\endaligned &\cr
& &\span\hrulefill&\cr
& &&
\gathered
\text{the pullback of the 5th-order differential equation $D'y=0$, where}
\\[-8pt]
\endgathered &\cr
& &\span\hrulefill&\cr
& &&
\aligned
\\[-12pt]
&
D'=\theta^5
-2\cdot 3^2z(2\theta+1)(405\theta^4+810\theta^3+1017\theta^2+612\theta+158)
\\[-2pt] &\;
+2^43^8z^2(\theta+1)(810\theta^4+3240\theta^3+5886\theta^2+5292\theta+1913)
\\[-2pt] &\;
-2^63^{17}z^3(\theta+1)(\theta+2)(2\theta+3)(15\theta^2+45\theta+49)
\\[-2pt] &\;
+2^83^{22}z^4(\theta+1)(\theta+2)(\theta+3)(45\theta^2+180\theta+203)
\\[-2pt] &\;
-2^93^{30}z^5(\theta+1)(\theta+2)(\theta+3)(\theta+4)(2\theta+5)
\\[-8pt]
\endaligned &\cr
& &\span\hrulefill&\cr
& &&
A_n'=(4\cdot 3^6)^n\sum_k(-1)^k\binom nk
\frac{(1/2)_k(1/3)_k^2(2/3)_k^2}{k!^5}
&\cr
\noalign{\hrule}
}}\hss}

\newpage

\hbox to\hsize{\hss\vbox{\offinterlineskip
\halign to120mm{\strut\tabskip=100pt minus 100pt
\strut\vrule#&\hbox to6.5mm{\hss$#$\hss}&%
\vrule#&\hbox to113mm{\hfil$\dsize#$\hfil}&%
\vrule#\tabskip=0pt\cr\noalign{\hrule}
& \# &%
& \text{differential operator $D$ and coefficients $A_n$, $n=0,1,2,\dots$} &\cr
\noalign{\hrule\vskip1pt\hrule}
& \eqnnol{\hbox{$\wh5$}}{210} &&
\aligned
\\[-12pt]
&
D=\theta^4-2^23z(576\theta^4+720\theta^3+542\theta^2+182\theta+21)
\\[-2pt] &\;
+2^43^2z^2(124416\theta^4+311040\theta^3+324576\theta^2+143280\theta+20017)
\\[-2pt] &\;
-2^{11}3^7z^3(4\theta+3)(1152\theta^3+3456\theta^2+3172\theta+597)
\\[-2pt] &\;
+2^{16}3^{12}z^4(4\theta+1)(4\theta+3)(4\theta+7)(4\theta+9)
\\[-8pt]
\endaligned &\cr
& &\span\hrulefill&\cr
& &&
\gathered
\text{the pullback of the 5th-order differential equation $D'y=0$, where}
\\[-8pt]
\endgathered &\cr
& &\span\hrulefill&\cr
& &&
\aligned
\\[-12pt]
&
D'=\theta^5
-2^33z(2\theta+1)(180\theta^4+360\theta^3+451\theta^2+271\theta+70)
\\[-2pt] &\;
+2^{10}3^4z^2(\theta+1)(360\theta^4+1440\theta^3+2613\theta^2+2346\theta+847)
\\[-2pt] &\;
-2^{15}3^8z^3(\theta+1)(\theta+2)(2\theta+3)(120\theta^2+360\theta+391)
\\[-2pt] &\;
+2^{22}3^{10}z^4(\theta+1)(\theta+2)(\theta+3)(180\theta^2+720\theta
+811)
\\[-2pt] &\;
-2^{29}3^{15}z^5(\theta+1)(\theta+2)(\theta+3)(\theta+4)(2\theta+5)
\\[-8pt]
\endaligned &\cr
& &\span\hrulefill&\cr
& &&
A_n'=1728^n\sum_k(-1)^k\binom nk
\frac{(1/2)_k^3(1/3)_k(2/3)_k}{k!^5}
&\cr
\noalign{\hrule}
& \eqnnol{\hbox{$\wh6$}}{211} &&
\aligned
\\[-12pt]
&
D=\theta^4-2^4z(1024\theta^4+1280\theta^3+968\theta^2+328\theta+39)
\\[-2pt] &\;
+2^{12}z^2(24576\theta^4+61440\theta^3+64256\theta^2+28480\theta+4017)
\\[-2pt] &\;
-2^{27}z^3(4\theta+3)(512\theta^3+1536\theta^2+1412\theta+267)
\\[-2pt] &\;
+2^{40}z^4(4\theta+1)(4\theta+3)(4\theta+7)(4\theta+9)
\\[-8pt]
\endaligned &\cr
& &\span\hrulefill&\cr
& &&
\gathered
\text{the pullback of the 5th-order differential equation $D'y=0$, where}
\\[-8pt]
\endgathered &\cr
& &\span\hrulefill&\cr
& &&
\aligned
\\[-12pt]
&
D'=\theta^5
-2^5z(2\theta+1)(320\theta^4+640\theta^3+804\theta^2+484\theta+125)
\\[-2pt] &\;
+2^{18}z^2(\theta+1)(640\theta^4+2560\theta^3+4652\theta^2+4184\theta+1513)
\\[-2pt] &\;
-2^{31}z^3(\theta+1)(\theta+2)(2\theta+3)(160\theta^2+480\theta+523)
\\[-2pt] &\;
+2^{44}z^4(\theta+1)(\theta+2)(\theta+3)(80\theta^2+320\theta+361)
\\[-2pt] &\;
-2^{54}z^5(\theta+1)(\theta+2)(\theta+3)(\theta+4)(2\theta+5)
\\[-8pt]
\endaligned &\cr
& &\span\hrulefill&\cr
& &&
A_n'=2^{12n}\sum_k(-1)^k\binom nk
\frac{(1/2)_k^3(1/4)_k(3/4)_k}{k!^5}
&\cr
\noalign{\hrule}
& \eqnnol{\hbox{$\wh7$}}{212} &&
\aligned
\\[-12pt]
&
D=\theta^4
-2^4z(49152\theta^4+81920\theta^3+62720\theta^2+21760\theta+2793)
\\[-2pt] &\;
+2^{26}z^2(4\theta+3)(768\theta^3+1984\theta^2+1572\theta+291)
\\[-2pt] &\;
-2^{46}z^3(4\theta+1)(4\theta+3)(4\theta+7)(4\theta+9)
\\[-8pt]
\endaligned &\cr
& &\span\hrulefill&\cr
& &&
\gathered
\text{the pullback of the 5th-order differential equation $D'y=0$, where}
\\[-8pt]
\endgathered &\cr
& &\span\hrulefill&\cr
& &&
\aligned
\\[-12pt]
&
D'=\theta^5
-2^5z(2\theta+1)(20480\theta^4+40960\theta^3+51840\theta^2+31360\theta+8087)
\\[-2pt] &\;
+2^{24}z^2(\theta+1)(40960\theta^4+163840\theta^3+298880\theta^2+270080\theta+98071)
\\[-2pt] &\;
-2^{48}z^3(\theta+1)(\theta+2)(2\theta+3)(64\theta^2+192\theta+211)
\\[-2pt] &\;
+2^{67}z^4(\theta+1)(\theta+2)(\theta+3)(32\theta^2+128\theta+145)
\\[-2pt] &\;
-2^{89}z^5(\theta+1)(\theta+2)(\theta+3)(\theta+4)(2\theta+5)
\\[-8pt]
\endaligned &\cr
& &\span\hrulefill&\cr
& &&
A_n'=2^{18n}\sum_k(-1)^k\binom nk
\frac{(1/2)_k(1/8)_k(3/8)_k(5/8)_k(7/8)_k}{k!^5}
&\cr
\noalign{\hrule}
}}\hss}

\newpage

\hbox to\hsize{\hss\vbox{\offinterlineskip
\halign to120mm{\strut\tabskip=100pt minus 100pt
\strut\vrule#&\hbox to6.5mm{\hss$#$\hss}&%
\vrule#&\hbox to113mm{\hfil$\dsize#$\hfil}&%
\vrule#\tabskip=0pt\cr\noalign{\hrule}
& \# &%
& \text{differential operator $D$ and coefficients $A_n$, $n=0,1,2,\dots$} &\cr
\noalign{\hrule\vskip1pt\hrule}
& \eqnnol{\hbox{$\wh8$}}{213} &&
\aligned
\\[-12pt]
&
D=\theta^4-2^23^2z(3888\theta^4+6480\theta^3+4950\theta^2+1710\theta+217)
\\[-2pt] &\;
+2^43^{10}z^2(4\theta+3)(1728\theta^3+4464\theta^2+3532\theta+651)
\\[-2pt] &\;
-2^{10}3^{18}z^3(4\theta+1)(4\theta+3)(4\theta+7)(4\theta+9)
\\[-8pt]
\endaligned &\cr
& &\span\hrulefill&\cr
& &&
\gathered
\text{the pullback of the 5th-order differential equation $D'y=0$, where}
\\[-8pt]
\endgathered &\cr
& &\span\hrulefill&\cr
& &&
\aligned
\\[-12pt]
&
D'=\theta^5-2^33^2z(1620\theta^4+3240\theta^3+4095\theta^2+2475\theta+638)
\\[-2pt] &\;
+2^{10}3^8z^2(\theta+1)(3240\theta^4+12960\theta^3+23625\theta^2+21350\theta+7739)
\\[-2pt] &\;
-5\cdot 2^{15}3^{17}z^3(\theta+1)(\theta+2)(2\theta+3)(24\theta^2+72\theta+79)
\\[-2pt] &\;
+5\cdot 2^{22}3^{22}z^4(\theta+1)(\theta+2)(\theta+3)(36\theta^2+144\theta+163)
\\[-2pt] &\;
-2^{29}3^{30}z^5(\theta+1)(\theta+2)(\theta+3)(\theta+4)(2\theta+5)
\\[-8pt]
\endaligned &\cr
& &\span\hrulefill&\cr
& &&
A_n'=(4\cdot 11664)^n\sum_k(-1)^k\binom nk
\frac{(1/2)_k(1/3)_k(2/3)_k(1/6)_k(5/6)_k}{k!^5}
&\cr
\noalign{\hrule}
& \eqnnol{\hbox{$\wh9$}}{214} &&
\aligned
\\[-12pt]
&
D=\theta^4
-2^43^2z(248832\theta^4+41472\theta^3+318528\theta^3+11116\theta+14497)
\\[-2pt] &\;
+2^{22}3^{10}z^2(4\theta+3)(432\theta^3+1116\theta^2+886\theta+165)
\\[-2pt] &\;
-2^{34}3^{18}z^3(4\theta+1)(4\theta+3)(4\theta+7)(4\theta+9)
\\[-8pt]
\endaligned &\cr
& &\span\hrulefill&\cr
& &&
\gathered
\text{the pullback of the 5th-order differential equation $D'y=0$, where}
\\[-8pt]
\endgathered &\cr
& &\span\hrulefill&\cr
& &&
\aligned
\\[-12pt]
&
D'=\theta^5
-2^53^2z(2\theta+1)(103680\theta^4+207360\theta^3+262944\theta^2
\\[-2pt] &\;\quad
+159264\theta+41087)
\\[-2pt] &\;
+2^{20}3^{17}z^2(\theta+1)(207360\theta^4+829440\theta^3+1514592\theta^2
\\[-2pt] &\;\quad
+1370304\theta+498143)
\\[-2pt] &\;
-2^{38}3^{17}z^3(\theta+1)(\theta+2)(2\theta+3)(240\theta^2+720\theta+793)
\\[-2pt] &\;
+2^{53}3^{22}z^4(\theta+1)(\theta+2)(\theta+3)(360\theta^2+1440\theta+1633)
\\[-2pt] &\;
-2^{69}3^{30}z^5(\theta+1)(\theta+2)(\theta+3)(\theta+4)(2\theta+5)
\\[-8pt]
\endaligned &\cr
& &\span\hrulefill&\cr
& &&
A_n'=3456^{2n}\sum_k(-1)^k\binom nk
\frac{(1/2)_k(1/12)_k(5/12)_k(7/12)_k(11/12)k}{k!^5}
&\cr
\noalign{\hrule}
& \eqnnol{\hbox{$\wh{10}$}}{215} &&
\aligned
\\[-12pt]
&
D=\theta^4-2^4z(3072\theta^4+5120\theta^3+3904\theta^2+1344\theta+169)
\\[-2pt] &\;
+2^{23}z^2(4\theta+3)(24\theta^3+62\theta^2+49\theta+9)
\\[-2pt] &\;
-2^{34}z^3(4\theta+1)(4\theta+3)(4\theta+7)(4\theta+9)
\\[-8pt]
\endaligned &\cr
& &\span\hrulefill&\cr
& &&
\gathered
\text{the pullback of the 5th-order differential equation $D'y=0$, where}
\\[-8pt]
\endgathered &\cr
& &\span\hrulefill&\cr
& &&
\aligned
\\[-12pt]
&
D'=\theta^5
-2^5z(2\theta+1)(1280\theta^4+2560\theta^3+3232\theta^2+1952\theta+503)
\\[-2pt] &\;
+2^{20}z^2(\theta+1)(2560\theta^4+10240\theta^3+18656\theta^2+16832\theta+6103)
\\[-2pt] &\;
-2^{38}z^3(\theta+1)(\theta+2)(2\theta+3)(80\theta^2+240\theta+263)
\\[-2pt] &\;
+2^{53}z^4(\theta+1)(\theta+2)(\theta+3)(40\theta^2+160\theta+181)
\\[-2pt] &\;
-2^{69}z^5(\theta+1)(\theta+2)(\theta+3)(\theta+4)(2\theta+5)
\\[-8pt]
\endaligned &\cr
& &\span\hrulefill&\cr
& &&
A_n'=2^{14n}\sum_k(-1)^k\binom nk
\frac{(1/2)_k(1/4)_k^2(3/4)_k^2}{k!^5}
&\cr
\noalign{\hrule}
}}\hss}

\newpage

\hbox to\hsize{\hss\vbox{\offinterlineskip
\halign to120mm{\strut\tabskip=100pt minus 100pt
\strut\vrule#&\hbox to6.5mm{\hss$#$\hss}&%
\vrule#&\hbox to113mm{\hfil$\dsize#$\hfil}&%
\vrule#\tabskip=0pt\cr\noalign{\hrule}
& \# &%
& \text{differential operator $D$ and coefficients $A_n$, $n=0,1,2,\dots$} &\cr
\noalign{\hrule\vskip1pt\hrule}
& \eqnnol{\hbox{$\wh{11}$}}{216} &&
\aligned
\\[-12pt]
&
D=\theta^4-2^23z(2304\theta^4+2880\theta^3+2186\theta^2+746\theta+91)
\\[-2pt] &\;
+2^43^3z^2(1990656\theta^4+4976640\theta^3+5213952\theta^2+2318400\theta+329497)
\\[-2pt] &\;
-2^{15}3^7z^3(4\theta+3)(4608\theta^3+13824\theta^2+12724\theta+2415)
\\[-2pt] &\;
+2^{24}3^{12}z^4(4\theta+1)(4\theta+3)(4\theta+7)(4\theta+9)
\\[-8pt]
\endaligned &\cr
& &\span\hrulefill&\cr
& &&
\gathered
\text{the pullback of the 5th-order differential equation $D'y=0$, where}
\\[-8pt]
\endgathered &\cr
& &\span\hrulefill&\cr
& &&
\aligned
\\[-12pt]
&
D'=\theta^5
-2^33z(2\theta+1)(720\theta^4+1440\theta^3+1813\theta^2+1093\theta+282)
\\[-2pt] &\;
+2^{12}3^5z^2(\theta+1)(480\theta^4+1920\theta^3+3493\theta^2+3146\theta+1139)
\\[-2pt] &\;
-2^{19}3^8z^3(\theta+1)(\theta+2)(2\theta+3)(480\theta^2+1440\theta+1573)
\\[-2pt] &\;
+2^{28}3^{10}z^4(\theta+1)(\theta+2)(\theta+3)(720\theta^2+2880\theta+3253)
\\[-2pt] &\;
-2^{39}3^{15}z^5(\theta+1)(\theta+2)(\theta+3)(\theta+4(2\theta+5)
\\[-8pt]
\endaligned &\cr
& &\span\hrulefill&\cr
& &&
A_n'=(4\cdot 12^3)^n\sum_k(-1)^k\binom nk
\frac{(1/2)_k(1/3)_k(2/3)_k(1/4)_k(3/4)_k}{k!^5}
&\cr
\noalign{\hrule}
& \eqnnol{\hbox{$\wh{12}$}}{217} &&
\aligned
\\[-12pt]
&
D=\theta^4-2^43z(9216\theta^4+11520\theta^3+8840\theta^23080\theta+403)
\\[-2pt] &\;
+2^{12}3^2z^2(1990656\theta^4+4976640\theta^3+5241600\theta^2+2352960\theta+3420049)
\\[-2pt] &\;
-2^{27}3^7z^3(4\theta+3)(4608\theta^3+13824\theta^2+12772\theta+2451)
\\[-2pt] &\;
+2^{40}3^{12}z^4(4\theta+1)(4\theta+3)(4\theta+7)(4\theta+9)
\\[-8pt]
\endaligned &\cr
& &\span\hrulefill&\cr
& &&
\gathered
\text{the pullback of the 5th-order differential equation $D'y=0$, where}
\\[-8pt]
\endgathered &\cr
& &\span\hrulefill&\cr
& &&
\aligned
\\[-12pt]
&
D'=\theta^5
-2^5z(2\theta+1)(2880\theta^4+5760\theta^3+7300\theta^2+4420\theta+1137)
\\[-2pt] &\;
+2^{18}3^5z^2(\theta+1)(1920\theta^4+7680\theta^3+14020\theta^2+12680\theta+4607)
\\[-2pt] &\;
-5\cdot 2^{31}3^8z^3(\theta+1)(\theta+2)(2\theta+3)(96\theta^2+288\theta+317)
\\[-2pt] &\;
+5\cdot 2^{44}3^{10}z^4(\theta+1)(\theta+2)(\theta+3)(144\theta^2+576\theta+653)
\\[-2pt] &\;
-2^{59}3^{15}z^5(\theta+1)(\theta+2)(\theta+3)(\theta+4)(2\theta+5)
\\[-8pt]
\endaligned &\cr
& &\span\hrulefill&\cr
& &&
A_n'=(2^{12}3^3)^n\sum_k(-1)^k\binom nk
\frac{(1/2)_k(1/4)_k(3/4)_k(1/6)_k(5/6)_k}{k!^5}
&\cr
\noalign{\hrule}
& \eqnnol{\hbox{$\wh{13}$}}{218} &&
\aligned
\\[-12pt]
&
D=\theta^4
-2^43^2z(20736\theta^4+25920\theta^3+20016\theta^2+7056\theta+961)
\\[-2pt] &\;
+2^{14}3^8z^2(31104\theta^4+77760\theta^3+82152\theta^2+37080\theta+5461)
\\[-2pt] &\;
-2^{24}3^{16}z^3(4\theta+3)(576\theta^3+1728\theta^2+1600\theta+309)
\\[-2pt] &\;
+2^{32}3^{24}z^4(4\theta+1)(4\theta+3)(4\theta+7)(4\theta+9)
\\[-8pt]
\endaligned &\cr
& &\span\hrulefill&\cr
& &&
\gathered
\text{the pullback of the 5th-order differential equation $D'y=0$, where}
\\[-8pt]
\endgathered &\cr
& &\span\hrulefill&\cr
& &&
\aligned
\\[-12pt]
&
D'=\theta^5
-2^53^2z(2\theta+1)(6480\theta^4+12960\theta^3+16488\theta^2+10008\theta+2567)
\\[-2pt] &\;
+2^{16}3^8z^2(\theta+1)(12960\theta^4+51840\theta^3+94824\theta^2+85968\theta+31295)
\\[-2pt] &\;
-2^{28}3^{17}z^3(\theta+1)(\theta+2)(2\theta+3)(60\theta^2+180\theta+199)
\\[-2pt] &\;
+2^{39}3^{22}z^4(\theta+1)(\theta+2)(\theta+3)(90\theta^2+360\theta+409)
\\[-2pt] &\;
-2^{49}3^{30}z^5(\theta+1)(\theta+2)(\theta+3)(\theta+4)(2\theta+5)
\\[-8pt]
\endaligned &\cr
& &\span\hrulefill&\cr
& &&
A_n'=(2^{10}3^6)^n\sum_k(-1)^k\binom nk
\frac{(1/2)_k(1/6)_k^2(5/6)_k^2}{k!^5}
&\cr
\noalign{\hrule}
}}\hss}

\newpage

\hbox to\hsize{\hss\vbox{\offinterlineskip
\halign to120mm{\strut\tabskip=100pt minus 100pt
\strut\vrule#&\hbox to6.5mm{\hss$#$\hss}&%
\vrule#&\hbox to113mm{\hfil$\dsize#$\hfil}&%
\vrule#\tabskip=0pt\cr\noalign{\hrule}
& \# &%
& \text{differential operator $D$ and coefficients $A_n$, $n=0,1,2,\dots$} &\cr
\noalign{\hrule\vskip1pt\hrule}
& \eqnnol{\hbox{$\wh{14}$}}{219} &&
\aligned
\\[-12pt]
&
D=\theta^4-2^43z(2304\theta^4+2880\theta^3+2192\theta^2+752\theta+93)
\\[-2pt] &\;
+2^{14}3^4z^2(31104\theta^4+77760\theta^3+81576\theta^2+36360\theta+5197)
\\[-2pt] &\;
-2^{24}3^7z^3(4\theta+3)(576\theta^3+1728\theta^2+1592\theta+303)
\\[-2pt] &\;
+2^{32}3^{12}z^4(4\theta+1)(4\theta+3)(4\theta+7)(4\theta+9)
\\[-8pt]
\endaligned &\cr
& &\span\hrulefill&\cr
& &&
\gathered
\text{the pullback of the 5th-order differential equation $D'y=0$, where}
\\[-8pt]
\endgathered &\cr
& &\span\hrulefill&\cr
& &&
\aligned
\\[-12pt]
&
D'=\theta^5
-2^53z(2\theta+1)(720\theta^4+1440\theta^3+1816\theta^2+1096\theta+283)
\\[-2pt] &\;
+2^{16}3^4z^2(\theta+1)(1440\theta^4+5760\theta^3+10488\theta^2+9456\theta+3427)
\\[-2pt] &\;
-2^{28}3^8z^3(\theta+1)(\theta+2)(2\theta+3)(60\theta^2+180\theta+197)
\\[-2pt] &\;
+2^{39}3^{10}z^4(\theta+1)(\theta+2)(\theta+3)(90\theta^2+360\theta+407)
\\[-2pt] &\;
-2^{49}3^{15}z^5(\theta+1)(\theta+2)(\theta+3)(\theta+4)(2\theta+5)
\\[-8pt]
\endaligned &\cr
& &\span\hrulefill&\cr
& &&
A_n'=(2^{10}3^3)^n\sum_k(-1)^k\binom nk
\frac{(1/2)_k^3(1/6)_k(5/6)_k}{k!^5}
&\cr
\noalign{\hrule}
& \eqnno{30} &&
\aligned
\\[-12pt]
&
D=\theta^4
-2^4z(4\theta+1)(4\theta+3)(8\theta^2+8\theta+3)
\\[-2pt] &\;
+2^{12}z^2(4\theta+1)(4\theta+3)(4\theta+5)(4\theta+7)
\\[-8pt]
\endaligned &\cr
& &\span\hrulefill&\cr
& &&
A_n
=\frac{(4n)!}{n!^2(2n)!}
\sum_{k=0}^n2^{2(n-k)}\binom{2k}k^2\binom{2n-2k}{n-k}
\; \text{(see \cite{AZ}, \thetag{6.7})}
&\cr
\noalign{\hrule}
& \eqnno{31} &&
\aligned
\\[-12pt]
&
D=\theta^4
-2^4z(4\theta+1)(32\theta^3+40\theta^2+28\theta+7)
\\[-2pt] &\;
+2^{12}z^2(4\theta+1)(4\theta+3)^2(4\theta+5)
\\[-8pt]
\endaligned &\cr
& &\span\hrulefill&\cr
& &&
A_n=A_n^{(-)} \; \text{(see \cite{AZ}, \thetag{6.8})}
&\cr
\noalign{\hrule}
}}\hss}

\newpage

\hbox to\hsize{\hss\vbox{\offinterlineskip
\halign to120mm{\strut\tabskip=100pt minus 100pt
\strut\vrule#&\hbox to5.5mm{\hss$#$\hss}&%
\vrule#&\hbox to114mm{\hfil$\dsize#$\hfil}&%
\vrule#\tabskip=0pt\cr\noalign{\hrule}
& \# &%
& \text{differential operator $D$ and coefficients $A_n$, $n=0,1,2,\dots$} &\cr
\noalign{\hrule\vskip1pt\hrule}
& \eqnno{32} &&
\aligned
\\[-12pt]
&
D=\theta^4
-3z(360\theta^4+90\theta^3+27\theta^2-18\theta-11)
\\[-2pt] &\;
+3z^2(145764\theta^4+72882\theta^3+24899\theta^2+6094\theta+13224)
\\[-2pt] &\;
-3^4z^3(970920\theta^4+728190\theta^3+279069\theta^2+130766\theta-16005)
\\[-2pt] &\;
+3^2z^4(587866086\theta^4+587866086\theta^3+249671835\theta^2
\\[-2pt] &\;\quad
+51395166\theta+4547776)
\\[-2pt] &\;
+3^5z^5(8738280\theta^4+10922850\theta^3+5667111\theta^2+606852\theta-464210)
\\[-2pt] &\;
+3^5z^6(1311876\theta^4+1967814\theta^3+1171557\theta^2+284652\theta+134263)
\\[-2pt] &\;
+3^8z^7(3240\theta^4+5670\theta^3+3753\theta^2+1260\theta+116)
\\[-2pt] &\;
+3^8z^8(3\theta+1)^2(3\theta+2)^2
\\[-8pt]
\endaligned &\cr
& &\span\hrulefill&\cr
& &&
\gathered
\text{the pullback of the 5th-order differential equation $D'y=0$, where}
\\[-8pt]
\endgathered &\cr
& &\span\hrulefill&\cr
& &&
\aligned
\\[-12pt]
&
D'=\theta^5
-3z(2\theta+1)(3\theta^2+3\theta+1)(15\theta^2+15\theta+4)
\\[-2pt] &\;
-3z^2(\theta+1)^3(3\theta+2)(3\theta+4)
\\[-8pt]
\endaligned &\cr
& &\span\hrulefill&\cr
& &&
A_n'=\sum_{k,l}\binom nk^2\binom nl^2\binom{n+k}n\binom{n+l}n\binom{k+l}n
&\cr
\noalign{\hrule}
& \eqnno{33} &&
\aligned
\\[-12pt]
&
D=\theta^4
-2^2z(324\theta^4+456\theta^3+321\theta^2+93\theta+10)
\\[-2pt] &\;
+2^9z^2(584\theta^4+584\theta^3+4\theta^2-71\theta-13)
\\[-2pt] &\;
-2^{16}z^3(324\theta^4+192\theta^3+123\theta^2+48\theta+7)
+2^{24}z^4(2\theta+1)^4
\\[-8pt]
\endaligned &\cr
& &\span\hrulefill&\cr
& &&
A_n=\binom{2n}n^2\sum_{k=0}^{2n}\binom{2n}k^3
&\cr
\noalign{\hrule}
& \eqnno{34} &&
\aligned
\\[-12pt]
&
D=\theta^4
-z(35\theta^4+70\theta^3+63\theta^2+28\theta+5)
\\[-2pt] &\;
+z^2(\theta+1)^2(259\theta^2+518\theta+285)
-225z^3(\theta+1)^2(\theta+2)^2
\\[-8pt]
\endaligned &\cr
& &\span\hrulefill&\cr
& &&
A_n=\sum_{i+j+k+l+m=n}\biggl(\frac{n!}{i!j!k!l!m!}\biggr)^2
&\cr
\noalign{\hrule}
& \eqnno{35} &&
\aligned
\\[-12pt]
&
D=\theta^4
-2^2\cdot3z(192\theta^4+240\theta^3+191\theta^2+71\theta+10)
\\[-2pt] &\;
+2^7\cdot3^2z^2(1746\theta^4+3960\theta^3+4323\theta^2+2247\theta+395)
\\[-2pt] &\;
-2^{12}\cdot3^4z^3(2538\theta^4+7776\theta^3+9915\theta^2+5643\theta+1030)
\\[-2pt] &\;
+2^{17}\cdot3^6z^4(1782\theta^4+6480\theta^3+8793\theta^2+4905\theta+875)
\\[-2pt] &\;
-2^{23}\cdot3^{11}z^5(\theta+1)^2(3\theta+1)(3\theta+5)
\\[-8pt]
\endaligned &\cr
& &\span\hrulefill&\cr
& &&
A_n=\sum_{k=0}^n\binom nk\frac{(6k)!}{k!(2k)!(3k)!}
\,\frac{(6n-6k)!}{(n-k)!(2n-2k)!(3n-3k)!}
&\cr
\noalign{\hrule}
& \eqnno{36} &&
\aligned
\\[-12pt]
&
D=\theta^4
-2^4z(2\theta+1)^2(3\theta^2+3\theta+1)
+2^9z^2(2\theta+1)^2(2\theta+3)^2
\\[-8pt]
\endaligned &\cr
& &\span\hrulefill&\cr
& &&
A_n=\binom{2n}n^2\sum_{k=0}^n
\binom nk\binom{2k}k\binom{2n-2k}{n-k}
&\cr
\noalign{\hrule}
}}\hss}

\newpage

\hbox to\hsize{\hss\vbox{\offinterlineskip
\halign to120mm{\strut\tabskip=100pt minus 100pt
\strut\vrule#&\hbox to5.5mm{\hss$#$\hss}&%
\vrule#&\hbox to114mm{\hfil$\dsize#$\hfil}&%
\vrule#\tabskip=0pt\cr\noalign{\hrule}
& \# &%
& \text{differential operator $D$ and coefficients $A_n$, $n=0,1,2,\dots$} &\cr
\noalign{\hrule\vskip1pt\hrule}
& \eqnno{37} &&
\aligned
\\[-12pt]
&
D=\theta^4
-4z(1280\theta^4+320\theta^3+94\theta^2-66\theta-39)
\\[-2pt] &\;
+2^4z^2(679936\theta^4+339968\theta^3+114304\theta^2-19008\theta+11601)
\\[-2pt] &\;
-2^{13}z^3(1515520\theta^4+1136640\theta^3+430656\theta^2+6144\theta+1329)
\\[-2pt] &\;
+2^{20}z^4(7868416\theta^4+7868416\theta^3+3318656\theta^2+156096\theta-118695)
\\[-2pt] &\;
-2^{31}z^5(1515520\theta^4+1894400\theta^3+880576\theta^2+62592\theta-4599)
\\[-2pt] &\;
+2^{40}z^6(679936\theta^4+1019904\theta^3+518016\theta^2+64704\theta+19713)
\\[-2pt] &\;
-2^{57}z^7\theta(640\theta^3+1120\theta^2+617\theta+105)
\\[-2pt] &\;
+2^{68}z^8\theta(\theta+1)(4\theta+1)(4\theta+3)
\\[-8pt]
\endaligned &\cr
& &\span\hrulefill&\cr
& &&
\gathered
\text{the pullback of the 5th-order differential equation $D'y=0$, where}
\\[-8pt]
\endgathered &\cr
& &\span\hrulefill&\cr
& &&
\aligned
\\[-12pt]
&
D'=\theta^5
-8z(2\theta+1)(4\theta+1)(4\theta+3)(5\theta^2+5\theta+2)
\\ \vspace{-1pt} &\;
+2^{10}z^2(\theta+1)(4\theta+1)(4\theta+3)(4\theta+5)(4\theta+7)
\\[-8pt]
\endaligned &\cr
& &\span\hrulefill&\cr
& &&
A_n'=\frac{(4n)!}{n!^2(2n)!}\sum_{k=0}^n
\binom nk^2\binom{2k}k\binom{2n-2k}{n-k}
&\cr
\noalign{\hrule}
& \eqnno{38} &&
\aligned
\\[-12pt]
&
D=\theta^4
-2^4z(4\theta+1)(4\theta+3)(3\theta^2+3\theta+1)
\\[-2pt] &\;
+2^9z^2(4\theta+1)(4\theta+3)(4\theta+5)(4\theta+7)
\\[-8pt]
\endaligned &\cr
& &\span\hrulefill&\cr
& &&
A_n=\frac{(4n)!}{n!^2(2n)!}\sum_{k=0}^n
\binom nk\binom{2k}k\binom{2n-2k}{n-k}
&\cr
\noalign{\hrule}
& \eqnno{39} &&
\aligned
\\[-12pt]
&
D=\theta^4
-4z(320\theta^4+80\theta^3+26\theta^2-14\theta-9)
\\[-2pt] &\;
+2^4z^2(42496\theta^4+21248\theta^3+7808\theta^2-560\theta+825)
\\[-2pt] &\;
-2^{12}z^3(47360\theta^4+35520\theta^3+14568\theta^2+1212\theta+141)
\\[-2pt] &\;
+2^{19}z^4(61472\theta^4+61472\theta^3+27848\theta^2+2992\theta-783)
\\[-2pt] &\;
-2^{27}z^5(23680\theta^4+29600\theta^3+14684\theta^2+1858\theta-15)
\\[-2pt] &\;
+2^{32}z^6(42496\theta^4+63744\theta^3+34368\theta^2+6000\theta+1293)
\\[-2pt] &\;
-2^{45}z^7\theta(160\theta^3+280\theta^2+163\theta+35)
+2^{54}z^8\theta(\theta+1)(2\theta+1)^2
\\[-8pt]
\endaligned &\cr
& &\span\hrulefill&\cr
& &&
\gathered
\text{the pullback of the 5th-order differential equation $D'y=0$, where}
\\[-8pt]
\endgathered &\cr
& &\span\hrulefill&\cr
& &&
\aligned
\\[-12pt]
&
D'=\theta^5
-8z(2\theta+1)^3(5\theta^2+5\theta+2)
+2^{10}z^2(\theta+1)(2\theta+1)^2(2\theta+3)^2
\\[-8pt]
\endaligned &\cr
& &\span\hrulefill&\cr
& &&
A_n'=\binom{2n}n^2\sum_{k=0}^n\binom nk^2\binom{2k}k\binom{2n-2k}{n-k}
&\cr
\noalign{\hrule}
& \eqnno{40} &&
\aligned
\\[-12pt]
&
D=\theta^4
-2^4z(64\theta^4+32\theta^3+12\theta^2-4\theta-3)
\\[-2pt] &\;
+2^{13}z^2(48\theta^4+48\theta^3+22\theta^2+2\theta+3)
\\[-2pt] &\;
-2^{22}z^3\theta(4\theta+3)(4\theta^2+3\theta+1)
+2^{30}z^4\theta(\theta+1)(2\theta+1)^2
\\[-8pt]
\endaligned &\cr
& &\span\hrulefill&\cr
& &&
\gathered
\text{the pullback of the 5th-order differential equation $D'y=0$, where}
\\[-8pt]
\endgathered &\cr
& &\span\hrulefill&\cr
& &&
\aligned
\\[-12pt]
&
D'=\theta^5
-2^5z(2\theta+1)^3(2\theta^2+2\theta+1)
+2^{12}z^2(\theta+1)(2\theta+1)^2(2\theta+3)^2
\\[-8pt]
\endaligned &\cr
& &\span\hrulefill&\cr
& &&
A_n'=\binom{2n}n^2\sum_{k=0}^n\binom{2k}k^2\binom{2n-2k}{n-k}^2
&\cr
\noalign{\hrule}
}}\hss}

\newpage

\hbox to\hsize{\hss\vbox{\offinterlineskip
\halign to120mm{\strut\tabskip=100pt minus 100pt
\strut\vrule#&\hbox to5.5mm{\hss$#$\hss}&%
\vrule#&\hbox to114mm{\hfil$\dsize#$\hfil}&%
\vrule#\tabskip=0pt\cr\noalign{\hrule}
& \# &%
& \text{differential operator $D$ and coefficients $A_n$, $n=0,1,2,\dots$} &\cr
\noalign{\hrule\vskip1pt\hrule}
& \eqnno{41} &&
\aligned
\\[-12pt]
&
D=\theta^4-2z(2\theta+1)^2(7\theta^2+7\theta+3)
+324z^2(\theta+1)^2(2\theta+1)(2\theta+3)
\\[-8pt]
\endaligned &\cr
& &\span\hrulefill&\cr
& &&
A_n=\binom{2n}n\sum_{k=0}^n(-1)^k3^{n-3k}\binom n{3k}
\binom{n+k}k\frac{(3k)!}{k!^3}
&\cr
\noalign{\hrule}
& \eqnno{42} &&
\aligned
\\[-12pt]
&
D=\theta^4-8z(2\theta+1)^2(3\theta^2+3\theta+1)
+64z^2(\theta+1)^2(2\theta+1)(2\theta+3)
\\[-8pt]
\endaligned &\cr
& &\span\hrulefill&\cr
& &&
A_n=\binom{2n}n\sum_k\binom nk^2\binom{2k}n^2
&\cr
\noalign{\hrule}
& \eqnno{43} &&
\aligned
\\[-12pt]
&
D=\theta^4
-2^4z(256\theta^4+128\theta^3+44\theta^2-20\theta-13)
\\[-2pt] &\;
+2^{14}z^2(384\theta^4+384\theta^3+164\theta^2+4\theta+23)
\\[-2pt] &\;
-2^{26}z^3\theta(64\theta^3+96\theta^2+49\theta+9)
+2^{36}z^4\theta(\theta+1)(4\theta+1)(4\theta+3)
\\[-8pt]
\endaligned &\cr
& &\span\hrulefill&\cr
& &&
\gathered
\text{the pullback of the 5th-order differential equation $D'y=0$, where}
\\[-8pt]
\endgathered &\cr
& &\span\hrulefill&\cr
& &&
\aligned
\\[-12pt]
&
D'=\theta^5
-2^5z(2\theta+1)(4\theta+1)(4\theta+3)(2\theta^2+2\theta+1)
\\[-2pt] &\;
+2^{12}z^2(\theta+1)(4\theta+1)(4\theta+3)(4\theta+5)(4\theta+7)
\\[-8pt]
\endaligned &\cr
& &\span\hrulefill&\cr
& &&
A_n'=\frac{(4n)!}{n!^2(2n)!}\sum_{k=0}^n
\binom{2k}k^2\binom{2n-2k}{n-k}^2
&\cr
\noalign{\hrule}
& \eqnno{44} &&
\aligned
\\[-12pt]
&
D=\theta^4
-2^2z(544\theta^4+136\theta^3+37\theta^2-31\theta-18)
\\[-2pt] &\;
+2^6z^2(27760\theta^4+13880\theta^3+4355\theta^2+823\theta+2472)
\\[-2pt] &\;
-2^{10}z^3(630496\theta^4+472872\theta^3+168105\theta^2+75333\theta-9924)
\\[-2pt] &\;
+2^{14}z^4(5400856\theta^4+5400856\theta^3+2146159\theta^2+393923\theta+22032)
\\[-2pt] &\;
-2^{18}z^5(630496\theta^4+788120\theta^3+365135\theta^2-5981\theta-55194)
\\[-2pt] &\;
+2^{22}z^6(27760\theta^4+41640\theta^3+21705\theta^2+2373\theta+1752)
\\[-2pt] &\;
-2^{26}z^7\theta(544\theta^3+952\theta^2+547\theta+119)
+2^{30}z^8\theta(\theta+1)(2\theta+1)^2
\\[-8pt]
\endaligned &\cr
& &\span\hrulefill&\cr
& &&
\gathered
\text{the pullback of the 5th-order differential equation $D'y=0$, where}
\\[-8pt]
\endgathered &\cr
& &\span\hrulefill&\cr
& &&
\aligned
\\[-12pt]
&
D'=\theta^5
-4z(2\theta+1)^3(17\theta^2+17\theta+5)
+16z^2(\theta+1)(2\theta+1)^2(2\theta+3)^2
\\[-8pt]
\endaligned &\cr
& &\span\hrulefill&\cr
& &&
A_n'=\binom{2n}n^2\sum_{k=0}^n\binom nk^2\binom{n+k}k^2
&\cr
\noalign{\hrule}
& \eqnno{45} &&
\aligned
\\[-12pt]
&
D=\theta^4
-4z(2\theta+1)^2(7\theta^2+7\theta+2)
-2^7z^2(2\theta+1)^2(2\theta+3)^2
\\[-8pt]
\endaligned &\cr
& &\span\hrulefill&\cr
& &&
A_n=\binom{2n}n^2\sum_{k=0}^n\binom nk^3
&\cr
\noalign{\hrule}
& \eqnno{46} &&
\aligned
\\[-12pt]
&
D=\theta^4
-6z(2\theta+1)^2(9\theta^2+9\theta+4)
\\[-2pt] &\;
+324z^2(2\theta+1)(2\theta+3)(3\theta+2)(3\theta+4)
\\[-8pt]
\endaligned &\cr
& &\span\hrulefill&\cr
& &&
A_n=\binom{2n}n\sum_{k=0}^n
\frac{(3k)!}{k!^3}\,\frac{(3n-3k)!}{(n-k)!^3}
&\cr
\noalign{\hrule}
}}\hss}

\newpage

\hbox to\hsize{\hss\vbox{\offinterlineskip
\halign to120mm{\strut\tabskip=100pt minus 100pt
\strut\vrule#&\hbox to5.5mm{\hss$#$\hss}&%
\vrule#&\hbox to114mm{\hfil$\dsize#$\hfil}&%
\vrule#\tabskip=0pt\cr\noalign{\hrule}
& \# &%
& \text{differential operator $D$ and coefficients $A_n$, $n=0,1,2,\dots$} &\cr
\noalign{\hrule\vskip1pt\hrule}
& \eqnno{47} &&
\aligned
\\[-12pt]
&
D=\theta^4
-2^43z(2\theta+1)^2(18\theta^2+18\theta+5)
\\[-2pt] &\;
+2^{10}3^4z^2(2\theta+1)(2\theta+3)(3\theta+1)(3\theta+5)
\\[-8pt]
\endaligned &\cr
& &\span\hrulefill&\cr
& &&
A_n=\binom{2n}n\sum_{k=0}^n\frac{(6k)!}{k!(2k)!(3k)!}
\,\frac{(6n-6k)!}{(n-k)!(2n-2k)!(3n-3k)!}
&\cr
\noalign{\hrule}
& \eqnno{48} &&
\aligned
\\[-12pt]
&
D=\theta^4
-12z(3\theta+1)(3\theta+2)(3\theta^2+3\theta+1)
\\[-2pt] &\;
+288z^2(3\theta+1)(3\theta+2)(3\theta+4)(3\theta+5)
\\[-8pt]
\endaligned &\cr
& &\span\hrulefill&\cr
& &&
A_n=\frac{(3n)!}{n!^3}\sum_{k=0}^n
\binom nk\binom{2k}k\binom{2n-2k}{n-k}
&\cr
\noalign{\hrule}
& \eqnno{49} &&
\aligned
\\[-12pt]
&
D=\theta^4
-2^2\cdot3z(144\theta^4+72\theta^3+26\theta^2-10\theta-7)
\\[-2pt] &\;
+2^4\cdot3^4z^2(864\theta^4+864\theta^3+384\theta^2+24\theta+53)
\\[-2pt] &\;
-2^{11}\cdot3^8z^3\theta(24\theta^3+36\theta^2+19\theta+4)
\\[-2pt] &\;
+2^{16}\cdot3^{10}z^4\theta(\theta+1)(3\theta+1)(3\theta+2)
\\[-8pt]
\endaligned &\cr
& &\span\hrulefill&\cr
& &&
\gathered
\text{the pullback of the 5th-order differential equation $D'y=0$, where}
\\[-8pt]
\endgathered &\cr
& &\span\hrulefill&\cr
& &&
\aligned
\\[-12pt]
&
D'=\theta^5
-24z(2\theta+1)(3\theta+1)(3\theta+2)(2\theta^2+2\theta+1)
\\[-2pt] &\;
+2304z^2(\theta+1)(3\theta+1)(3\theta+2)(3\theta+4)(3\theta+5)
\\[-8pt]
\endaligned &\cr
& &\span\hrulefill&\cr
& &&
A_n'=\frac{(3n)!}{n!^3}\sum_{k=0}^n
\binom{2k}k^2\binom{2n-2k}{n-k}^2
&\cr
\noalign{\hrule}
& \eqnno{50} &&
\aligned
\\[-12pt]
&
D=\theta^4
-3z(720\theta^4+180\theta^3+56\theta^2-34\theta-21)
\\[-2pt] &\;
+3^4z^2(23904\theta^4+11952\theta^3+4226\theta^2-472\theta+439)
\\[-2pt] &\;
-3^7z^3(426240\theta^4+319680\theta^3+126672\theta^2+6828\theta+868)
\\[-2pt] &\;
+3^{10}z^4(4425984\theta^4+4425984\theta^3+1943584\theta^2+158704\theta-61035)
\\[-2pt] &\;
-2^7\cdot3^{13}z^5(213120\theta^4+266400\theta^3+128456\theta^2+13202\theta-364)
\\[-2pt] &\;
+2^{13}\cdot3^{16}z^6(11952\theta^4+17928\theta^3+9417\theta^2+1443\theta+356)
\\[-2pt] &\;
-2^{20}\cdot3^{19}z^7\theta(180\theta^3+315\theta^2+179\theta+35)
\\[-2pt] &\;
+2^{24}\cdot3^{22}z^8\theta(\theta+1)(3\theta+1)(3\theta+2)
\\[-8pt]
\endaligned &\cr
& &\span\hrulefill&\cr
& &&
\gathered
\text{the pullback of the 5th-order differential equation $D'y=0$, where}
\\[-8pt]
\endgathered &\cr
& &\span\hrulefill&\cr
& &&
\aligned
\\[-12pt]
&
D'=\theta^5
-6z(2\theta+1)(3\theta+1)(3\theta+2)(5\theta^2+5\theta+2)
\\[-2pt] &\;
+576z^2(\theta+1)(3\theta+1)(3\theta+2)(3\theta+4)(3\theta+5)
\\[-8pt]
\endaligned &\cr
& &\span\hrulefill&\cr
& &&
\aligned
\\[-12pt]
A_n'
&=\frac{(3n)!}{n!^3}\sum_{j+k+l+m=n}\biggl(\frac{n!}{j!\,k!\,l!\,m!}\biggr)^2
=\frac{(3n)!}{n!^3}\sum_{k=0}^n\binom nk^2
\binom{2k}k\binom{2n-2k}{n-k}
\\[-2pt]
&=4^{-n}\frac{(3n)!}{n!^3}\binom{2n}n^3
\sum_{k=0}^n\binom nk^4\binom{2n}{2k}^{-3}
\vspace{1.5pt}
\endaligned &\cr
\noalign{\hrule}
& \eqnno{51} &&
\aligned
\\[-12pt]
&
D=\theta^4
-4z(4\theta+1)(4\theta+3)(11\theta^2+11\theta+3)
\\[-2pt] &\;
-16z^2(4\theta+1)(4\theta+3)(4\theta+5)(4\theta+7)
\\[-8pt]
\endaligned &\cr
& &\span\hrulefill&\cr
& &&
A_n=\frac{(4n)!}{n!^2(2n)!}\sum_{k=0}^n\binom nk^2\binom{n+k}k
&\cr
\noalign{\hrule}
}}\hss}

\newpage

\hbox to\hsize{\hss\vbox{\offinterlineskip
\halign to120mm{\strut\tabskip=100pt minus 100pt
\strut\vrule#&\hbox to5.5mm{\hss$#$\hss}&%
\vrule#&\hbox to114mm{\hfil$\dsize#$\hfil}&%
\vrule#\tabskip=0pt\cr\noalign{\hrule}
& \# &%
& \text{differential operator $D$ and coefficients $A_n$, $n=0,1,2,\dots$} &\cr
\noalign{\hrule\vskip1pt\hrule}
& \eqnno{52} &&
\aligned
\\[-12pt]
&
D=\theta^4
-2^2z(2176\theta^4+544\theta^3+131\theta^2-141\theta-78)
\\[-2pt] &\;
+2^7z^2(222080\theta^4+111040\theta^3+31370\theta^2+3690\theta+19557)
\\[-2pt] &\;
-2^{14}z^3(2521984\theta^4+1891488\theta^3+613311\theta^2+261807\theta-37362)
\\[-2pt] &\;
+2^{20}z^4(21603424\theta^4+21603424\theta^3+7909529\theta^2+1234089\theta+58191)
\\[-2pt] &\;
-2^{26}z^5(2521984\theta^4+3152480\theta^3+1362025\theta^2-102855\theta-233532)
\\[-2pt] &\;
+2^{31}z^6(222080\theta^4+333120\theta^3+163230\theta^2+9150\theta+13317)
\\[-2pt] &\;
-2^{38}z^7\theta(2176\theta^3+3808\theta^2+2069\theta+357)
\\[-2pt] &\;
+2^{44}z^8\theta(\theta+1)(4\theta+1)(4\theta+3)
\\[-8pt]
\endaligned &\cr
& &\span\hrulefill&\cr
& &&
\gathered
\text{the pullback of the 5th-order differential equation $D'y=0$, where}
\\[-8pt]
\endgathered &\cr
& &\span\hrulefill&\cr
& &&
\aligned
\\[-12pt]
&
D'=\theta^5
-4z(2\theta+1)(4\theta+1)(4\theta+3)(17\theta^2+17\theta+5)
\\[-2pt] &\;
+16z^2(\theta+1)(4\theta+1)(4\theta+3)(4\theta+5)(4\theta+7)
\\[-8pt]
\endaligned &\cr
& &\span\hrulefill&\cr
& &&
A_n'=\frac{(4n)!}{n!^2(2n)!}\sum_{k=0}^n\binom nk^2\binom{n+k}k^2
&\cr
\noalign{\hrule}
& \eqnno{53} &&
\aligned
\\[-12pt]
&
D=\theta^4
-3z(1224\theta^4+306\theta^3+79\theta^2-74\theta-42)
\\[-2pt] &\;
+3^4z^2(62460\theta^4+31230\theta^3+9365\theta^2+1490\theta+5534)
\\[-2pt] &\;
-3^7z^3(1418616\theta^4+1063962\theta^3+363459\theta^2+159618\theta-21788)
\\[-2pt] &\;
+3^{10}z^4(12151926\theta^4+12151926\theta^3+4660081\theta^2+800926\theta+41364)
\\[-2pt] &\;
-3^{13}z^5(1418616\theta^4+1773270\theta^3+796925\theta^2-33190\theta-127418)
\\[-2pt] &\;
+3^{16}z^6(62460\theta^4+93690\theta^3+47535\theta^2+4110\theta+3854)
\\[-2pt] &\;
-3^{19}z^7\theta(1224\theta^3+2142\theta^2+1201\theta+238)
\\[-2pt] &\;
+3^{22}z^8\theta(\theta+1)(3\theta+1)(3\theta+2)
\\[-8pt]
\endaligned &\cr
& &\span\hrulefill&\cr
& &&
\gathered
\text{the pullback of the 5th-order differential equation $D'y=0$, where}
\\[-8pt]
\endgathered &\cr
& &\span\hrulefill&\cr
& &&
\aligned
\\[-12pt]
&
D'=\theta^5
-3z(2\theta+1)(3\theta+1)(3\theta+2)(17\theta^2+17\theta+5)
\\[-2pt] &\;
+9z^2(\theta+1)(3\theta+1)(3\theta+2)(3\theta+4)(3\theta+5)
\\[-8pt]
\endaligned &\cr
& &\span\hrulefill&\cr
& &&
A_n'=\frac{(3n)!}{n!^3}\sum_{k=0}^n\binom nk^2\binom{n+k}k^2
&\cr
\noalign{\hrule}
& \eqnno{54} &&
\aligned
\\[-12pt]
&
D=\theta^4
-2z(6\theta^4+10\theta^3+10\theta^2+5\theta+1)
\\[-2pt] &\;
+4z^2(1039\theta^4+4146\theta^3+5707\theta^2+3127\theta+592)
\\[-2pt] &\;
-8z^3(\theta+1)(4116\theta^3+18512\theta^2+27768\theta+13889)
\\[-2pt] &\;
+16z^4(\theta+1)(\theta+2)(6159\theta^2+24631\theta+25143)
\\[-2pt] &\;
-64z^5(\theta+1)(\theta+2)(\theta+3)(2051\theta+5127)
\\[-2pt] &\;
+65600z^6(\theta+1)(\theta+2)(\theta+3)(\theta+4)
\\[-8pt]
\endaligned &\cr
& &\span\hrulefill&\cr
& &&
A_n=\sum_k(-1)^k2^{n-2k}\binom n{2k}\biggl(\frac{(4k)!}{k!^2(2k)!}\biggr)^2
&\cr
\noalign{\hrule}
}}\hss}

\newpage

\hbox to\hsize{\hss\vbox{\offinterlineskip
\halign to120mm{\strut\tabskip=100pt minus 100pt
\strut\vrule#&\hbox to5.5mm{\hss$#$\hss}&%
\vrule#&\hbox to114mm{\hfil$\dsize#$\hfil}&%
\vrule#\tabskip=0pt\cr\noalign{\hrule}
& \# &%
& \text{differential operator $D$ and coefficients $A_n$, $n=0,1,2,\dots$} &\cr
\noalign{\hrule\vskip1pt\hrule}
& \eqnno{55} &&
\aligned
\\[-12pt]
&
D=9\theta^4
-12z(208\theta^4+224\theta^3+163\theta^2+51\theta+6)
\\[-2pt] &\;
+2^9z^2(32\theta^4-928\theta^3-1606\theta^2-837\theta-141)
\\[-2pt] &\;
+2^{16}z^3(144\theta^4+576\theta^3+467\theta^2+144\theta+15)
-2^{24}z^4(2\theta+1)^4
\\[-8pt]
\endaligned &\cr
& &\span\hrulefill&\cr
& &&
A_n=\binom{2n}n^2\sum_{k=0}^n\binom nk^2\binom{2n}{2k}
&\cr
\noalign{\hrule}
& \eqnno{56} &&
\aligned
\\[-12pt]
&
D=\theta^4
-2^4z(22\theta^4+8\theta^3+9\theta^2+5\theta+1)
\\[-2pt] &\;
+2^9z^2(94\theta^4+88\theta^3+97\theta^2+45\theta+8)
\\[-2pt] &\;
-2^{14}z^3(194\theta^4+336\theta^3+371\theta^2+195\theta+41)
\\[-2pt] &\;
+2^{19}\cdot3z^4(64\theta^4+176\theta^3+217\theta^2+129\theta+30)
-2^{27}\cdot3^2z^5(\theta+1)^4
\\[-8pt]
\endaligned &\cr
& &\span\hrulefill&\cr
& &&
A_n=\sum_{k=0}^n\binom nk^{-1}\binom{2k}k^3\binom{2n-2k}{n-k}^3
&\cr
\noalign{\hrule}
& \eqnno{57} &&
\aligned
\\[-12pt]
&
D=\theta^4
-2^2z(8\theta^4+10\theta^3+10\theta^2+5\theta+1)
\\[-2pt] &\;
+2^4z^2(28\theta^4+70\theta^3+95\theta^2+65\theta+18)
\\[-2pt] &\;
-2^7z^3(28\theta^4+105\theta^3+180\theta^2+150\theta+49)
\\[-2pt] &\;
+2^8z^4(799\theta^4+6182\theta^3+15629\theta^2+13660\theta+3856)
\\[-2pt] &\;
-2^{11}z^5(\theta+1)(1486\theta^3+11082\theta^2+25388\theta+17397)
\\[-2pt] &\;
+2^{13}z^6(\theta+1)(\theta+2)(2201\theta^2+13185\theta+18827)
\\[-2pt] &\;
-2^{16}z^7(\theta+1)(\theta+2)(\theta+3)(731\theta+2557)
\\[-2pt] &\;
+2^{17}\cdot5\cdot73z^8(\theta+1)(\theta+2)(\theta+3)(\theta+4)
\\[-8pt]
\endaligned &\cr
& &\span\hrulefill&\cr
& &&
A_n=\sum_k(-1)^k4^{n-4k}\binom n{4k}
\biggl(\frac{(6k)!}{k!(2k)!(3k)!}\biggr)^2
&\cr
\noalign{\hrule}
& \eqnno{58} &&
\aligned
\\[-12pt]
&
D=\theta^4
-4z(2\theta+1)^2(10\theta^2+10\theta+3)
+144z^2(2\theta+1)^2(2\theta+3)^2
\\[-8pt]
\endaligned &\cr
& &\span\hrulefill&\cr
& &&
A_n=\binom{2n}n^2\sum_{k=0}^n\binom nk^2\binom{2k}k
&\cr
\noalign{\hrule}
& \eqnno{59} &&
\aligned
\\[-12pt]
&
D=49\theta^4
-2z(1799\theta^4+3640\theta^3+3045\theta^2+1225\theta+196)
\\[-2pt] &\;
+2^2z^2(13497\theta^4+55536\theta^3+81222\theta^2+50337\theta+11396)
\\[-2pt] &\;
-2^3z^3(17201\theta^4+114996\theta^3+248466\theta^2+202629\theta+55412)
\\[-2pt] &\;
-2^4z^4(5762\theta^4+29668\theta^3+48150\theta^2+31741\theta+7412)
\\[-2pt] &\;
-2^5\cdot3z^5(4\theta+5)(3\theta+2)(3\theta+4)(4\theta+3)
\\[-8pt]
\endaligned &\cr
& &\span\hrulefill&\cr
& &&
A_n=\sum_{k=0}^n\binom nk\binom{2n-k}n
\binom{n+k}k\binom{2k}k\binom{2n-2k}{n-k}
&\cr
\noalign{\hrule}
}}\hss}

\newpage

\hbox to\hsize{\hss\vbox{\offinterlineskip
\halign to120mm{\strut\tabskip=100pt minus 100pt
\strut\vrule#&\hbox to5.5mm{\hss$#$\hss}&%
\vrule#&\hbox to114mm{\hfil$\dsize#$\hfil}&%
\vrule#\tabskip=0pt\cr\noalign{\hrule}
& \# &%
& \text{differential operator $D$ and coefficients $A_n$, $n=0,1,2,\dots$} &\cr
\noalign{\hrule\vskip1pt\hrule}
& \eqnno{60} &&
\aligned
\\[-12pt]
&
D=\theta^4
-2z(248\theta^4+62\theta^3+23\theta^2-8\theta-6)
\\[-2pt] &\;
+2^2z^2(24792\theta^4+12396\theta^3+5153\theta^2+704\theta+1126)
\\[-2pt] &\;
-2^2z^3(2549440\theta^4+1912080\theta^3+882790\theta^2+230890\theta+21527)
\\[-2pt] &\;
+2^3z^4(71646752\theta^4+71646752\theta^3+36508992\theta^2+8176202\theta-785881)
\\[-2pt] &\;
-2^4\cdot3z^5(367119360\theta^4+458899200\theta^3+256824520\theta^2
\\[-2pt] &\;\quad
+54290940\theta+1905253)
\\[-2pt] &\;
+2^4\cdot3^2z^6(2056347648\theta^4+3084521472\theta^3+1880419968\theta^2
\\[-2pt] &\;\quad
+499802112\theta+140357665)
\\[-2pt] &\;
-2^9\cdot3^5z^7(20570112\theta^4+35997696\theta^3+23709888\theta^2
\\[-2pt] &\;\quad
+7437024\theta+640447)
\\[-2pt] &\;
+2^{12}\cdot3^8z^8(24\theta+5)(24\theta+11)(24\theta+13)(24\theta+19)
\\[-8pt]
\endaligned &\cr
& &\span\hrulefill&\cr
& &&
\gathered
\text{the pullback of the 5th-order differential equation $D'y=0$, where}
\\[-8pt]
\endgathered &\cr
& &\span\hrulefill&\cr
& &&
\aligned
\\[-12pt]
&
D'=\theta^5
-2z(2\theta+1)(31\theta^4+62\theta^3+54\theta^2+23\theta+4)
\\[-2pt] &\;
+12z^2(\theta+1)(3\theta+2)(3\theta+4)(4\theta+3)(4\theta+5)
\\[-8pt]
\endaligned &\cr
& &\span\hrulefill&\cr
& &&
A_n'=\sum_{k=0}^n\binom nk^2\binom{2n-k}n
\binom{n+k}k\binom{2k}k\binom{2n-2k}{n-k}
&\cr
\noalign{\hrule}
& \eqnno{61} &&
\aligned
\\[-12pt]
&
D=\theta^4
-2^43^2z(6\theta+1)(6\theta+5)(72\theta^2+72\theta+31)
\\[-2pt] &\;
+2^{12}3^8z^2(6\theta+1)(6\theta+5)(6\theta+7)(6\theta+11)
\\[-8pt]
\endaligned &\cr
& &\span\hrulefill&\cr
& &&
A_n=\frac{(6n)!}{n!\,(2n)!\,(3n)!}
\cdot432^n\sum_k(-1)^k\binom{-5/6}k\binom{-1/6}{n-k}^2
&\cr
\noalign{\hrule}
& \eqnno{62} &&
\aligned
\\[-12pt]
&
D=\theta^4
-12z(6\theta+1)(6\theta+5)(7\theta^2+7\theta+2)
\\[-2pt] &\;
+1152z^2(6\theta+1)(6\theta+5)(6\theta+7)(6\theta+11)
\\[-8pt]
\endaligned &\cr
& &\span\hrulefill&\cr
& &&
A_n=\frac{(6n)!}{n!(2n)!(3n)!}\sum_{k=0}^n\binom nk^3
&\cr
\noalign{\hrule}
& \eqnno{63} &&
\aligned
\\[-12pt]
&
D=\theta^4
-12z(6\theta+1)(6\theta+5)(11\theta^2+11\theta+3)
\\[-2pt] &\;
+144z^2(6\theta+1)(6\theta+5)(6\theta+7)(6\theta+11)
\\[-8pt]
\endaligned &\cr
& &\span\hrulefill&\cr
& &&
A_n=\frac{(6n)!}{n!(2n)!(3n)!}\sum_{k=0}^n\binom nk^2\binom{n+k}k
&\cr
\noalign{\hrule}
}}\hss}

\newpage

\hbox to\hsize{\hss\vbox{\offinterlineskip
\halign to120mm{\strut\tabskip=100pt minus 100pt
\strut\vrule#&\hbox to5.5mm{\hss$#$\hss}&%
\vrule#&\hbox to114mm{\hfil$\dsize#$\hfil}&%
\vrule#\tabskip=0pt\cr\noalign{\hrule}
& \# &%
& \text{differential operator $D$ and coefficients $A_n$, $n=0,1,2,\dots$} &\cr
\noalign{\hrule\vskip1pt\hrule}
& \eqnno{64} &&
\aligned
\\[-12pt]
&
D=\theta^4
-12z(6\theta+1)(6\theta+5)(10\theta^2+10\theta+3)
\\[-2pt] &\;
+1296z^2(6\theta+1)(6\theta+5)(6\theta+7)(6\theta+11)
\\[-8pt]
\endaligned &\cr
& &\span\hrulefill&\cr
& &&
A_n=\frac{(6n)!}{n!(2n)!(3n)!}\sum_{k=0}^n\binom nk^2\binom{2k}k
&\cr
\noalign{\hrule}
& \eqnno{65} &&
\aligned
\\[-12pt]
&
D=\theta^4
-48z(6\theta+1)(6\theta+5)(3\theta^2+3\theta+1)
\\[-2pt] &\;
+4608z^2(6\theta+1)(6\theta+5)(6\theta+7)(6\theta+11)
\\[-8pt]
\endaligned &\cr
& &\span\hrulefill&\cr
& &&
A_n=\frac{(6n)!}{n!(2n)!(3n)!}\sum_{k=0}^n\binom nk
\binom{2k}k\binom{2n-2k}{n-k}
&\cr
\noalign{\hrule}
& \eqnno{66} &&
\aligned
\\[-12pt]
&
D=\theta^4
-2^2\cdot3z(2880\theta^4+720\theta^3+194\theta^2-166\theta-93)
\\[-2pt] &\;
+2^4\cdot3^4z^2(382464\theta^4+191232\theta^3+59648\theta^2-15088\theta+5833)
\\[-2pt] &\;
-2^{12}\cdot3^7z^3(426240\theta^4+319680\theta^3+113352\theta^2-5412\theta-305)
\\[-2pt] &\;
+2^{19}\cdot3^{10}z^4(553248\theta^4+553248\theta^3+219896\theta^2-1432\theta-9327)
\\[-2pt] &\;
-2^{28}\cdot3^{13}z^5(106560\theta^4+133200\theta^3+58678\theta^2+1321\theta-518)
\\[-2pt] &\;
+2^{32}\cdot3^{16}z^6(382464\theta^4+573696\theta^3+277440\theta^2+22704\theta+10669)
\\[-2pt] &\;
-2^{45}\cdot3^{19}z^7\theta(1440\theta^3+2520\theta^2+1327\theta+175)
\\[-2pt] &\;
+2^{54}\cdot3^{22}z^8\theta(\theta+1)(6\theta+1)(6\theta+5)
\\[-8pt]
\endaligned &\cr
& &\span\hrulefill&\cr
& &&
\gathered
\text{the pullback of the 5th-order differential equation $D'y=0$, where}
\\[-8pt]
\endgathered &\cr
& &\span\hrulefill&\cr
& &&
\aligned
\\[-12pt]
&
D'=\theta^5
-24z(2\theta+1)(6\theta+1)(6\theta+5)(5\theta^2+5\theta+2)
\\[-2pt] &\;
+9216z^2(\theta+1)(6\theta+1)(6\theta+5)(6\theta+7)(6\theta+11)
\\[-8pt]
\endaligned &\cr
& &\span\hrulefill&\cr
& &&
A_n'=\frac{(6n)!}{n!(2n)!(3n)!}\sum_{k=0}^n\binom nk^2
\binom{2k}k\binom{2n-2k}{n-k}
&\cr
\noalign{\hrule}
& \eqnno{67} &&
\aligned
\\[-12pt]
&
D=\theta^4
-2^4\cdot3z(576\theta^4+288\theta^3+92\theta^2-52\theta-31)
\\[-2pt] &\;
+2^{13}\cdot3^4z^2(432\theta^4+432\theta^3+174\theta^2-6\theta+25)
\\[-2pt] &\;
-2^{22}\cdot3^8z^3\theta(4\theta+1)(12\theta^2+15\theta+5)
\\[-2pt] &\;
+2^{30}\cdot3^{10}z^4\theta(\theta+1)(6\theta+1)(6\theta+5)
\\[-8pt]
\endaligned &\cr
& &\span\hrulefill&\cr
& &&
\gathered
\text{the pullback of the 5th-order differential equation $D'y=0$, where}
\\[-8pt]
\endgathered &\cr
& &\span\hrulefill&\cr
& &&
\aligned
\\[-12pt]
&
D'=\theta^5
-96z(2\theta+1)(6\theta+1)(6\theta+5)(3\theta^2+3\theta+1)
\\[-2pt] &\;
+36864z^2(\theta+1)(6\theta+1)(6\theta+5)(6\theta+7)(6\theta+11)
\\[-8pt]
\endaligned &\cr
& &\span\hrulefill&\cr
& &&
A_n'=\frac{(6n)!}{n!(2n)!(3n)!}\sum_{k=0}^n
\binom{2k}k^2\binom{2n-2k}{n-k}^2
&\cr
\noalign{\hrule}
& \eqnno{68} &&
\aligned
\\[-12pt]
&
D=\theta^4
-4z(4\theta+1)(4\theta+3)(7\theta^2+7\theta+2)
\\[-2pt] &\;
-128z^2(4\theta+1)(4\theta+3)(4\theta+5)(4\theta+7)
\\[-8pt]
\endaligned &\cr
& &\span\hrulefill&\cr
& &&
A_n=\frac{(4n)!}{n!^2(2n)!}\sum_{k=0}^n\binom nk^3
&\cr
\noalign{\hrule}
}}\hss}

\newpage

\hbox to\hsize{\hss\vbox{\offinterlineskip
\halign to120mm{\strut\tabskip=100pt minus 100pt
\strut\vrule#&\hbox to5.5mm{\hss$#$\hss}&%
\vrule#&\hbox to114mm{\hfil$\dsize#$\hfil}&%
\vrule#\tabskip=0pt\cr\noalign{\hrule}
& \# &%
& \text{differential operator $D$ and coefficients $A_n$, $n=0,1,2,\dots$} &\cr
\noalign{\hrule\vskip1pt\hrule}
& \eqnno{69} &&
\aligned
\\[-12pt]
&
D=\theta^4
-4z(4\theta+1)(4\theta+3)(10\theta^2+10\theta+3)
\\[-2pt] &\;
+144z^2(4\theta+1)(4\theta+3)(4\theta+5)(4\theta+7)
\\[-8pt]
\endaligned &\cr
& &\span\hrulefill&\cr
& &&
A_n=\frac{(4n)!}{n!^2(2n)!}\sum_{k=0}^n\binom nk^2\binom{2k}k
&\cr
\noalign{\hrule}
& \eqnno{70} &&
\aligned
\\[-12pt]
&
D=\theta^4
-3z(3\theta+1)(3\theta+2)(10\theta^2+10\theta+3)
\\[-2pt] &\;
+81z^2(3\theta+1)(3\theta+2)(3\theta+4)(3\theta+5)
\\[-8pt]
\endaligned &\cr
& &\span\hrulefill&\cr
& &&
A_n=\frac{(3n)!}{n!^3}\sum_{k=0}^n\binom nk^2\binom{2k}k
&\cr
\noalign{\hrule}
& \eqnno{71} &&
\aligned
\\[-12pt]
&
D=\theta^4
+2^4z(39\theta^4-42\theta^3-29\theta^2-8\theta-1)
\\[-2pt] &\;
+2^{11}z^2\theta(37\theta^3-137\theta^2-10\theta-1)
\\[-2pt] &\;
-2^{16}z^3(181\theta^4+456\theta^3+353\theta^2+132\theta+19)
\\[-2pt] &\;
-2^{23}5z^4(36\theta^4+60\theta^3+36\theta^2+6\theta-1)
+2^{30}5^2z^5(\theta+1)^4
\\[-8pt]
\endaligned &\cr
& &\span\hrulefill&\cr
& &&
A_n=2^{-n}\sum_{k=0}^n
\binom{2k}k^4\binom{2n-2k}{n-k}^4\binom nk^{-3}
&\cr
\noalign{\hrule}
& \eqnno{72} &&
\aligned
\\[-12pt]
&
D=\theta^4
-2z(136\theta^4+260\theta^3+198\theta^2+68\theta+9)
\\[-2pt] &\;
+2^2z^2(1048\theta^4+2584\theta^3+3102\theta^2+2188\theta+703)
\\[-2pt] &\;
-2^4z^3(1552\theta^4+4632\theta^3+6552\theta^2+4614\theta+1089)
\\[-2pt] &\;
+2^4z^4(4112\theta^4+14368\theta^3+21528\theta^2+14344\theta+3585)
-2^{16}z^5(\theta+1)^4
\\[-8pt]
\endaligned &\cr
& &\span\hrulefill&\cr
& &&
A_n=\binom{2n}n\sum_{k=0}^n\binom nk^2\binom{2k}k^3
\binom{2n}{2k}^{-1}
&\cr
\noalign{\hrule}
& \eqnno{73} &&
\aligned
\\[-12pt]
&
D=\theta^4
-18z(42\theta^4+60\theta^3+45\theta^2+15\theta+2)
\\[-2pt] &\;
+3\cdot 18^2z^2(180\theta^4+432\theta^3+453\theta^2+222\theta+40)
\\[-2pt] &\;
-3\cdot 18^4z^3(2\theta+1)^2(13\theta^2+29\theta+20)
+18^6z^4(2\theta+1)^2(2\theta+3)^2
\\[-8pt]
\endaligned &\cr
& &\span\hrulefill&\cr
& &&
A_n=\binom{2n}n\sum_{k=0}^n\frac{(3k)!}{k!^3}
\frac{(3n-3k)!}{(n-k)!^3}\binom{2n}k\binom nk^{-1}
&\cr
\noalign{\hrule}
& \eqnno{74} &&
\aligned
\\[-12pt]
&
D=\theta^4
-2\cdot3z(99\theta^4+36\theta^3+39\theta^2+21\theta+4)
\\[-2pt] &\;
+2^2\cdot3^2z^2(3807\theta^4+3564\theta^3+3798\theta^2+1683\theta+284)
\\[-2pt] &\;
-2^3\cdot3^5z^3(7857\theta^4+13608\theta^3+14562\theta^2+7317\theta+1444)
\\[-2pt] &\;
+2^4\cdot3^9z^4(2592\theta^4+7128\theta^3+8550\theta^2+4851\theta+1052)
\\[-2pt] &\;
-2^5\cdot3^{13}z^5(3\theta+2)(3\theta+4)(6\theta+5)(6\theta+7)
\\[-8pt]
\endaligned &\cr
& &\span\hrulefill&\cr
& &&
A_n=\sum_{k=0}^n\frac{(3k)!}{k!^3}\frac{(3n-3k)!}{(n-k)!^3}
\binom{2k}k\binom{2n-2k}{n-k}\binom nk^{-1}
&\cr
\noalign{\hrule}
& \eqnno{75} &&
\aligned
\\[-12pt]
&
D=\theta^4
-10z(2\theta+1)^2(26\theta^2+26\theta+5)
\\[-2pt] &\;
+4z^2(2\theta+1)(2\theta+3)(774\theta^2+1548\theta+823)
\\[-2pt] &\;
-3088z^3(2\theta+1)(2\theta+3)^2(2\theta+5)
\\[-2pt] &\;
+4112z^4(2\theta+1)(2\theta+3)(2\theta+5)(2\theta+7)
\\[-8pt]
\endaligned &\cr
& &\span\hrulefill&\cr
& &&
A_n=\binom{2n}n\sum_{k=0}^n\binom nk\frac{(4k)!}{k!^4}
&\cr
\noalign{\hrule}
}}\hss}

\newpage

\hbox to\hsize{\hss\vbox{\offinterlineskip
\halign to120mm{\strut\tabskip=100pt minus 100pt
\strut\vrule#&\hbox to5.5mm{\hss$#$\hss}&%
\vrule#&\hbox to114mm{\hfil$\dsize#$\hfil}&%
\vrule#\tabskip=0pt\cr\noalign{\hrule}
& \# &%
& \text{differential operator $D$ and coefficients $A_n$, $n=0,1,2,\dots$} &\cr
\noalign{\hrule\vskip1pt\hrule}
& \eqnno{76} &&
\aligned
\\[-12pt]
&
D=\theta^4
-z(1029\theta^4+2058\theta^3+1482\theta^2+453\theta+49)
\\[-2pt] &\;
+z^2(\theta+1)^2(4106\theta^2+8212\theta+5007)
\\[-2pt] &\;
-z^3(\theta+1)(\theta+2)(6154\theta^2+18462\theta+14809)
\\[-2pt] &\;
+4101z^4(\theta+1)(\theta+2)^2(\theta+3)
\\[-2pt] &\;
-1025z^5(\theta+1)(\theta+2)(\theta+3)(\theta+4)
\\[-8pt]
\endaligned &\cr
& &\span\hrulefill&\cr
& &&
A_n=\sum_{k=0}^n\binom nk\binom{2k}k\frac{(4k)!}{k!^4}
&\cr
\noalign{\hrule}
& \eqnno{77} &&
\aligned
\\[-12pt]
&
D=\theta^4
-z(16416\theta^4+4116\theta^3+906\theta^2-1152\theta-625)
\\[-2pt] &\;
+z^2(101122496\theta^4+50758192\theta^3+13326174\theta^2+773908\theta+8758171)
\\[-2pt] &\;
-z^3(277299334656\theta^4+209187183168\theta^3+63952323504\theta^2
\\[-2pt] &\;\quad
+25856143860\theta-4222453166)
\\[-2pt] &\;
+z^4(286996730758656\theta^4+290319472045824\theta^3+102144306403616\theta^2
\\[-2pt] &\;\quad
+13734541538384\theta+29090804825)
\\[-2pt] &\;
-2^4\cdot3\cdot5^2\cdot41z^5(4\theta-1)(23108277888\theta^3+52094682096\theta^2
\\[-2pt] &\;\quad
+44583010636\theta+13882711785)
\\[-2pt] &\;
+2^5\cdot5^4\cdot41^2z^6(4\theta-1)(4\theta+3)(50561248\theta^2+126501592\theta+87298415)
\\[-2pt] &\;
-2^8\cdot3\cdot5^6\cdot41^3z^7(4\theta-1)(4\theta+3)(4\theta+7)(1368\theta+2395)
\\[-2pt] &\;
+2^8\cdot5^8\cdot41^4z^8(4\theta-1)(4\theta+3)(4\theta+7)(4\theta+11)
\\[-8pt]
\endaligned &\cr
& &\span\hrulefill&\cr
& &&
\gathered
\text{the pullback of the 5th-order differential equation $D'y=0$, where}
\\[-8pt]
\endgathered &\cr
& &\span\hrulefill&\cr
& &&
\aligned
\\[-12pt]
&
D'=\theta^5
-2(2\theta+1)(1029\theta^4+2058\theta^3+1482\theta^2+453\theta+49)
\\[-2pt] &\;
+4(\theta+1)(2\theta+1)(2\theta+3)(4106\theta^2+8212\theta+5007)
\\[-2pt] &\;
-8(2\theta+1)(2\theta+3)(2\theta+5)(6154\theta^2+18462\theta+14809)
\\[-2pt] &\;
+48\cdot1367(\theta+2)(2\theta+1)(2\theta+3)(2\theta+5)(2\theta+7)
\\[-2pt] &\;
-2^5\cdot5^2\cdot41(2\theta+1)(2\theta+3)(2\theta+5)(2\theta+7)(2\theta+9)
\\[-8pt]
\endaligned &\cr
& &\span\hrulefill&\cr
& &&
A_n'=\binom{2n}n\sum_{k=0}^n\binom nk
\binom{2k}k\frac{(4k)!}{k!^4}
&\cr
\noalign{\hrule}
}}\hss}

\newpage

\hbox to\hsize{\hss\vbox{\offinterlineskip
\halign to120mm{\strut\tabskip=100pt minus 100pt
\strut\vrule#&\hbox to5.5mm{\hss$#$\hss}&%
\vrule#&\hbox to114mm{\hfil$\dsize#$\hfil}&%
\vrule#\tabskip=0pt\cr\noalign{\hrule}
& \# &%
& \text{differential operator $D$ and coefficients $A_n$, $n=0,1,2,\dots$} &\cr
\noalign{\hrule\vskip1pt\hrule}
& \eqnno{78} &&
\aligned
\\[-12pt]
&
D=\theta^4
-z(16392\theta^4+4102\theta^3+900\theta^2-1151\theta-624)
\\[-2pt] &\;
+z^2(100778012\theta^4+50454570\theta^3+13189539\theta^2+803211\theta+8784048)
\\[-2pt] &\;
-z^3(275482230840\theta^4+207014719614\theta^3+62883681150\theta^2
\\[-2pt] &\;\quad
+25936450350\theta-4109625168)
\\[-2pt] &\;
+z^4(282850876768326\theta^4+283952402645202\theta^3+98665243460597\theta^2
\\[-2pt] &\;\quad
+13790192460644\theta+614356880784)
\\[-2pt] &\;
-17\cdot241z^5\theta(275482230840\theta^3+620238065874\theta^2
\\[-2pt] &\;\quad
+530195542680\theta+165017893373)
\\[-2pt] &\;
+17^2\cdot241^2z^6\theta(\theta+1)(100778012\theta^2+252010594\theta+173831131)
\\[-2pt] &\;
-2\cdot17^3\cdot241^3z^7\theta(\theta+1)(\theta+2)(8196\theta+14345)
\\[-2pt] &\;
+17^4\cdot241^4z^8\theta(\theta+1)(\theta+2)(\theta+3)
\\[-8pt]
\endaligned &\cr
& &\span\hrulefill&\cr
& &&
\gathered
\text{the pullback of the 5th-order differential equation $D'y=0$, where}
\\[-8pt]
\endgathered &\cr
& &\span\hrulefill&\cr
& &&
\aligned
\\[-12pt]
&
D'=\theta^5
-z(2\theta+1)(2051\theta^4+4102\theta^3+2951\theta^2+900\theta+97)
\\[-2pt] &\;
+z^2(\theta+1)(20495\theta^4+81980\theta^3+132457\theta^2+100954\theta+30175)
\\[-2pt] &\;
-2z^3(\theta+1)(\theta+2)(2\theta+3)(10245\theta^2+30735\theta+27983)
\\[-2pt] &\;
+z^4(\theta+1)(\theta+2)(\theta+3)(40975\theta^2+163900\theta+173889)
\\[-2pt] &\;
-10243z^5(\theta+1)(\theta+2)(\theta+3)(\theta+4)(2\theta+5)
\\[-2pt] &\;
+4097z^6(\theta+1)(\theta+2)(\theta+3)(\theta+4)(\theta+5)
\\[-8pt]
\endaligned &\cr
& &\span\hrulefill&\cr
& &&
A_n'=\sum_{k=0}^n\binom nk\binom{2k}k^2\frac{(4k)!}{k!^4}
&\cr
\noalign{\hrule}
& \eqnno{79} &&
\aligned
\\[-12pt]
&
D=\theta^4
-z(3130\theta^4+6260\theta^3+4385\theta^2+1255\theta+121)
\\[-2pt] &\;
+15z^2(\theta+1)^2(834\theta^2+1668\theta+1001)
\\[-2pt] &\;
-5z^3(\theta+1)(\theta+2)(3752\theta^2+11256\theta+9005)
\\[-2pt] &\;
+12505z^4(\theta+1)(\theta+2)^2(\theta+3)
\\[-2pt] &\;
-3126z^5(\theta+1)(\theta+2)(\theta+3)(\theta+4)
\\[-8pt]
\endaligned &\cr
& &\span\hrulefill&\cr
& &&
\gathered
\text{same $K(q)$ as in \#1}
\\[-8pt]
\endgathered &\cr
& &\span\hrulefill&\cr
& &&
A_n=\sum_{k=0}^n\binom nk\frac{(5k)!}{k!^5}
&\cr
\noalign{\hrule}
}}\hss}

\newpage

\hbox to\hsize{\hss\vbox{\offinterlineskip
\halign to120mm{\strut\tabskip=100pt minus 100pt
\strut\vrule#&\hbox to5.5mm{\hss$#$\hss}&%
\vrule#&\hbox to114mm{\hfil$\dsize#$\hfil}&%
\vrule#\tabskip=0pt\cr\noalign{\hrule}
& \# &%
& \text{differential operator $D$ and coefficients $A_n$, $n=0,1,2,\dots$} &\cr
\noalign{\hrule\vskip1pt\hrule}
& \eqnno{80} &&
\aligned
\\[-12pt]
&
D=\theta^4
-2z(25016\theta^4+6260\theta^3+1255\theta^2-1875\theta-998)
\\[-2pt] &\;
+2^2z^2(234725112\theta^4+117512640\theta^3+28486025\theta^2-10620\theta+20240412)
\\[-2pt] &\;
-2^4z^3(489688550224\theta^4+367970437920\theta^3+104723790940\theta^2
\\[-2pt] &\;\quad
+41242253625\theta-7232306723)
\\[-2pt] &\;
+2^6z^4(383914650187780\theta^4+385383011813200\theta^3+126338801127175\theta^2
\\[-2pt] &\;\quad
+14531449104300\theta+193021814409)
\\[-2pt] &\;
-16\cdot24\cdot521z^5(4\theta-1)(489688550224\theta^3+1102503263256\theta^2
\\[-2pt] &\;\quad
+937734079014\theta+289728842991)
\\[-2pt] &\;
+4\cdot24^2\cdot521^2z^6(4\theta-1)(4\theta+3)(234725112\theta^2
\\[-2pt] &\;\quad
+586962864\theta+404126249)
\\[-2pt] &\;
-16\cdot24^3\cdot521^3z^7(4\theta-1)(4\theta+3)(4\theta+7)(3127\theta+5473)
\\[-2pt] &\;
+24^4\cdot521^4z^8(4\theta-1)(4\theta+3)(4\theta+7)(4\theta+11)
\\[-8pt]
\endaligned &\cr
& &\span\hrulefill&\cr
& &&
\gathered
\text{the pullback of the 5th-order differential equation $D'y=0$, where}
\\[-8pt]
\endgathered &\cr
& &\span\hrulefill&\cr
& &&
\aligned
\\[-12pt]
&
D'=\theta^5
-2z(2\theta+1)(3130\theta^4+6260\theta^3+4385\theta^2+1255\theta+121)
\\[-2pt] &\;
+60z^2(\theta+1)(2\theta+1)(2\theta+3)(834\theta^2+1668\theta+1001)
\\[-2pt] &\;
-40z^3(2\theta+1)(2\theta+3)(2\theta+5)(3752\theta^2+11256\theta+9005)
\\[-2pt] &\;
+2^4\cdot5\cdot41\cdot61z^4(\theta+2)(2\theta+1)(2\theta+3)(2\theta+5)(2\theta+7)
\\[-2pt] &\;
-2^6\cdot3\cdot521z^5(2\theta+1)(2\theta+3)(2\theta+5)(2\theta+7)(2\theta+9)
\\[-8pt]
\endaligned &\cr
& &\span\hrulefill&\cr
& &&
A_n'=\binom{2n}n\sum_{k=0}^n\binom nk\frac{(5k)!}{k!^5}
&\cr
\noalign{\hrule}
}}\hss}

\newpage

\hbox to\hsize{\hss\vbox{\offinterlineskip
\halign to120mm{\strut\tabskip=100pt minus 100pt
\strut\vrule#&\hbox to5.5mm{\hss$#$\hss}&%
\vrule#&\hbox to114mm{\hfil$\dsize#$\hfil}&%
\vrule#\tabskip=0pt\cr\noalign{\hrule}
& \# &%
& \text{differential operator $D$ and coefficients $A_n$, $n=0,1,2,\dots$} &\cr
\noalign{\hrule\vskip1pt\hrule}
& \eqnno{81} &&
\aligned
\\[-12pt]
&
D=\theta^4
-z(50008\theta^4+12506\theta^3+2504\theta^2-3749\theta-1995)
\\[-2pt] &\;
+z^2(937850028\theta^4+469125042\theta^3+113531285\theta^2+50011\theta+81040890)
\\[-2pt] &\;
-z^3(7818126050056\theta^4+5867345737626\theta^3+1665767387626\theta^2
\\[-2pt] &\;\quad
+660699781296\theta-114681662235)
\\[-2pt] &\;
+z^4(24453139064250070\theta^4+24484407817250210\theta^3
\\[-2pt] &\;\quad
+7988410528587745\theta^2+932188857443850\theta+28724348623005)
\\[-2pt] &\;
-3^3\cdot463z^5\theta(7818126050056\theta^3+17594534812710\theta^2
\\[-2pt] &\;\quad
+14959207512780\theta+4621083875125)
\\[-2pt] &\;
+3^6\cdot463^2z^6\theta(\theta+1)(937850028\theta^2+2344825098\theta+1614171341)
\\[-2pt] &\;
-2\cdot3^9\cdot463^3z^7\theta(\theta+1)(\theta+2)(25004\theta+43759)
\\[-2pt] &\;
+3^{12}\cdot463^4z^8\theta(\theta+1)(\theta+2)(\theta+3)
\\[-8pt]
\endaligned &\cr
& &\span\hrulefill&\cr
& &&
\gathered
\text{same $K(q)$ as in \#80}
\\
\text{the pullback of the 5th-order differential equation $D'y=0$, where}
\\[-8pt]
\endgathered &\cr
& &\span\hrulefill&\cr
& &&
\aligned
\\[-12pt]
&
D'=\theta^5
-z(2\theta+1)(6253\theta^4+12506\theta^3+8757\theta^2+2504\theta+241)
\\[-2pt] &\;
+z^2(\theta+1)(62515\theta^4+250060\theta^3+402605\theta^2+305090\theta+90511)
\\[-2pt] &\;
-70z^3(\theta+1)(\theta+2)(2\theta+3)(893\theta^2+2679\theta+2429)
\\[-2pt] &\;
+5z^4(\theta+1)(\theta+2)(\theta+3)(25003\theta^2+100012\theta+106013)
\\[-2pt] &\;
-31253z^5(\theta+1)(\theta+2)(\theta+3)(\theta+4)(2\theta+5)
\\[-2pt] &\;
+3^3\cdot463z^6(\theta+1)(\theta+2)(\theta+3)(\theta+4)(\theta+5)
\\[-8pt]
\endaligned &\cr
& &\span\hrulefill&\cr
& &&
A_n'=\sum_{k=0}^n\binom nk\binom{2k}k\frac{(5k)!}{k!^5}
&\cr
\noalign{\hrule}
& \eqnno{82} &&
\aligned
\\[-12pt]
&
D=\theta^4
-z(186631\theta^4+46662\theta^3+8428\theta^2-14903\theta-7812)
\\[-2pt] &\;
+z^2(13061813781\theta^4+6531700068\theta^3+1452108415\theta^2
\\[-2pt] &\;\quad
-105618806\theta+1116194184)
\\[-2pt] &\;
-z^3(406305132943139\theta^4+304784361622746\theta^3
\\[-2pt] &\;\quad
+80477825393615\theta^2+30347763996276\theta-5895787146228)
\\[-2pt] &\;
+13\cdot37\cdot97z^4\theta(101592610503011\theta^3+101629619768568\theta^2
\\[-2pt] &\;\quad
+31089416550198\theta+2838662698240)
\\[-2pt] &\;
-13^2\cdot37^2\cdot97^2z^5\theta(\theta+1)(6531186837\theta^2+10886104581\theta+5020240741)
\\[-2pt] &\;
+13^3\cdot37^3\cdot97^3z^6\theta(\theta+1)(\theta+2)(139975\theta+186639)
\\[-2pt] &\;
-13^4\cdot37^4\cdot97^4z^7\theta(\theta+1)(\theta+2)(\theta+3)
\\[-8pt]
\endaligned &\cr
& &\span\hrulefill&\cr
& &&
\gathered
\text{the pullback of the 5th-order differential equation $D'y=0$, where}
\\[-8pt]
\endgathered &\cr
& &\span\hrulefill&\cr
& &&
\aligned
\\[-12pt]
&
D'=\theta^5
-7z(2\theta+1)(3333\theta^4+6666\theta^3+4537\theta^2+1204\theta+103)
\\[-2pt] &\;
+z^2(\theta+1)(233295\theta^4+933180\theta^3+1496985\theta^2+1127610\theta+331951)
\\[-2pt] &\;
-10z^3(\theta+1)(\theta+2)(\theta+3)(23329\theta^2+69987\theta+63183)
\\[-2pt] &\;
+5z^4(\theta+1)(\theta+2)(\theta+3)(93315\theta^2+373260\theta+395293)
\\[-2pt] &\;
-3\cdot59\cdot659z^5(\theta+1)(\theta+2)(\theta+3)(\theta+4)(2\theta+5)
\\[-2pt] &\;
+13\cdot37\cdot97z^6(\theta+1)(\theta+2)(\theta+3)(\theta+4)(\theta+5)
\\[-8pt]
\endaligned &\cr
& &\span\hrulefill&\cr
& &&
A_n'=\sum_{k=0}^n\binom nk\frac{(6k)!}{k!^6}
&\cr
\noalign{\hrule}
}}\hss}

\newpage

\hbox to\hsize{\hss\vbox{\offinterlineskip
\halign to120mm{\strut\tabskip=100pt minus 100pt
\strut\vrule#&\hbox to5.5mm{\hss$#$\hss}&%
\vrule#&\hbox to114mm{\hfil$\dsize#$\hfil}&%
\vrule#\tabskip=0pt\cr\noalign{\hrule}
& \# &%
& \text{differential operator $D$ and coefficients $A_n$, $n=0,1,2,\dots$} &\cr
\noalign{\hrule\vskip1pt\hrule}
& \eqnno{83} &&
\aligned
\\[-12pt]
&
D=\theta^4
-2^4z(88\theta^4+32\theta^3+33\theta^2+17\theta+3)
\\[-2pt] &\;
+2^9z^2(1504\theta^4+1408\theta^3+1436\theta^2+596\theta+93)
\\[-2pt] &\;
-2^{18}z^3(776\theta^4+1344\theta^3+1381\theta^2+651\theta+117)
\\[-2pt] &\;
+3\cdot2^{23}z^4(2\theta+1)(512\theta^3+1152\theta^2+1054\theta+339)
\\[-2pt] &\;
-9\cdot2^{31}z^5(2\theta+1)(2\theta+3)(4\theta+3)(4\theta+5)
\\[-8pt]
\endaligned &\cr
& &\span\hrulefill&\cr
& &&
A_n=\binom{2n}n\sum_{k=0}^n\binom nk\binom{2n}{2k}^{-1}
\frac{(4k)!}{k!^2(2k)!}\frac{(4n-4k)!}{(n-k)!^2(2n-2k)!}
&\cr
\noalign{\hrule}
& \eqnno{84} &&
\aligned
\\[-12pt]
&
D=\theta^4
-4z(32\theta^4+64\theta^3+63\theta^2+31\theta+6)
\\[-2pt] &\;
+256z^2(\theta+1)^2(4\theta+3)(4\theta+5)
\\[-8pt]
\endaligned &\cr
& &\span\hrulefill&\cr
& &&
A_n=\sum_{k=0}^n\binom nk^2\binom{2n}{2k}^{-1}
\frac{(4k)!}{k!^2(2k)!}\frac{(4n-4k)!}{(n-k)!^2(2n-2k)!}
&\cr
\noalign{\hrule}
& \eqnno{85} &&
\aligned
\\[-12pt]
&
D=\theta^4
-2z(2\theta+1)^2(866\theta^2+866\theta+121)
\\[-2pt] &\;
+4z^2(2\theta+1)(2\theta+3)(5190\theta^2+10380\theta+5431)
\\[-2pt] &\;
-2^4\cdot1297z^3(2\theta+1)(2\theta+3)^2(2\theta+5)
\\[-2pt] &\;
+2^4\cdot7\cdot13\cdot19z^4(2\theta+1)(2\theta+3)(2\theta+5)(2\theta+7)
\\[-8pt]
\endaligned &\cr
& &\span\hrulefill&\cr
& &&
A_n=\binom{2n}n\sum_{k=0}^n\binom nk\binom{2k}k
\frac{(6k)!}{k!(2k)!(3k)!}
&\cr
\noalign{\hrule}
& \eqnno{86} &&
\aligned
\\[-12pt]
&
D=\theta^4
-z(6917\theta^4+13834\theta^3+9610\theta^2+2693\theta+241)
\\[-2pt] &\;
+z^2(\theta+1)^2(27658\theta^2+55316\theta+33039)
\\[-2pt] &\;
-z^3(\theta+1)(\theta+2)(41482\theta^2+124446\theta+99481)
\\[-2pt] &\;
+27653z^4(\theta+1)(\theta+2)^2(\theta+3)
\\[-2pt] &\;
-31\cdot223z^5(\theta+1)(\theta+2)(\theta+3)(\theta+4)
\\[-8pt]
\endaligned &\cr
& &\span\hrulefill&\cr
& &&
\gathered
\text{same $K(q)$ as in \#85}
\\[-8pt]
\endgathered &\cr
& &\span\hrulefill&\cr
& &&
A_n=\sum_{k=0}^n\binom nk\binom{2k}k\frac{(6k)!}{k!^3(3k)!}
&\cr
\noalign{\hrule}
& \eqnno{87} &&
\aligned
\\[-12pt]
&
D=\theta^4
-z(3145\theta^4+6282\theta^3+4403\theta^2+1262\theta+122)
\\[-2pt] &\;
+4z^2(2\theta+1)(6270\theta^3+18804\theta^2+20058\theta+7523)
\\[-2pt] &\;
-8z^3(2\theta+1)(2\theta+3)(9395\theta^2+28181\theta+22544)
\\[-2pt] &\;
+16z^4(2\theta+1)(2\theta+3)(2\theta+5)(6260\theta+12519)
\\[-2pt] &\;
-2^4\cdot3\cdot7\cdot149z^5(2\theta+1)(2\theta+3)(2\theta+5)(2\theta+7)
\\[-8pt]
\endaligned &\cr
& &\span\hrulefill&\cr
& &&
A_n=\binom{2n}n\sum_{k=0}^n\binom nk\frac{(5k)!}{k!^3(2k)!}
&\cr
\noalign{\hrule}
}}\hss}

\newpage

\hbox to\hsize{\hss\vbox{\offinterlineskip
\halign to120mm{\strut\tabskip=100pt minus 100pt
\strut\vrule#&\hbox to5.5mm{\hss$#$\hss}&%
\vrule#&\hbox to114mm{\hfil$\dsize#$\hfil}&%
\vrule#\tabskip=0pt\cr\noalign{\hrule}
& \# &%
& \text{differential operator $D$ and coefficients $A_n$, $n=0,1,2,\dots$} &\cr
\noalign{\hrule\vskip1pt\hrule}
& \eqnno{88} &&
\aligned
\\[-12pt]
&
D=\theta^4
-z(110624\theta^4+27668\theta^3+5386\theta^2-8448\theta-4465)
\\[-2pt] &\;
+z^2(4589568448\theta^4+2296111664\theta^3+542921054\theta^2
\\[-2pt] &\;\quad
-10713580\theta+395878043)
\\[-2pt] &\;
-z^3(84647953763840\theta^4+63541034834496\theta^3+17687257564080\theta^2
\\[-2pt] &\;\quad
+6901304092148\theta-1219151757870)
\\[-2pt] &\;
+z^4(586017416803534336\theta^4+587032971970653952\theta^3
\\[-2pt] &\;\quad
+188288528023623968\theta^2+20443050826683472\theta+526999484846681)
\\[-2pt] &\;
-2^4\cdot31\cdot223z^5(4\theta-1)(21161988440960\theta^3+47628241370064\theta^2
\\[-2pt] &\;\quad
+40437124384228\theta+12464650057531)
\\[-2pt] &\;
+2^5\cdot31^2\cdot223^2z^6(4\theta-1)(4\theta+3)(2294784224\theta^2
\\[-2pt] &\;\quad
+5737624280\theta+3947717231)
\\[-2pt] &\;
-2^8\cdot31^3\cdot223^3z^7(4\theta-1)(4\theta+3)(4\theta+7)(27656\theta+48401)
\\[-2pt] &\;
+2^8\cdot31^4\cdot223^4z^8(4\theta-1)(4\theta+3)(4\theta+7)(4\theta+11)
\\[-8pt]
\endaligned &\cr
& &\span\hrulefill&\cr
& &&
\gathered
\text{the pullback of the 5th-order differential equation $D'y=0$, where}
\\[-8pt]
\endgathered &\cr
& &\span\hrulefill&\cr
& &&
\aligned
\\[-12pt]
&
D'=\theta^5
-2z(2\theta+1)(6917\theta^4+13834\theta^3+9610\theta^2+2693\theta+241)
\\[-2pt] &\;
+4z^2(\theta+1)(2\theta+1)(2\theta+3)(27658\theta^2+55316\theta+33039)
\\[-2pt] &\;
-8z^3(2\theta+1)(2\theta+3)(2\theta+5)(41482\theta^2+124446\theta+99481)
\\[-2pt] &\;
+442448z^4(\theta+2)(2\theta+1)(2\theta+3)(2\theta+5)(2\theta+7)
\\[-2pt] &\;
-221216z^5(2\theta+1)(2\theta+3)(2\theta+5)(2\theta+7)(2\theta+9)
\\[-8pt]
\endaligned &\cr
& &\span\hrulefill&\cr
& &&
A_n'=\binom{2n}n\sum_{k=0}^n\binom nk\binom{2k}k\frac{(6k)!}{k!^3(3k)!}
&\cr
\noalign{\hrule}
}}\hss}

\newpage

\hbox to\hsize{\hss\vbox{\offinterlineskip
\halign to120mm{\strut\tabskip=100pt minus 100pt
\strut\vrule#&\hbox to5.5mm{\hss$#$\hss}&%
\vrule#&\hbox to114mm{\hfil$\dsize#$\hfil}&%
\vrule#\tabskip=0pt\cr\noalign{\hrule}
& \# &%
& \text{differential operator $D$ and coefficients $A_n$, $n=0,1,2,\dots$} &\cr
\noalign{\hrule\vskip1pt\hrule}
& \eqnno{89} &&
\aligned
\\[-12pt]
&
D=\theta^4
-z(110600\theta^4+27654\theta^3+5380\theta^2-8447\theta-4464)
\\[-2pt] &\;
+z^2(4587245596\theta^4+2294065194\theta^3+542010659\theta^2
\\[-2pt] &\;\quad
-10506997\theta+396053424)
\\[-2pt] &\;
-z^3(84565362438200\theta^4+63442370368638\theta^3+17639445591678\theta^2
\\[-2pt] &\;\quad
+6905443579438\theta-1214113199184)
\\[-2pt] &\;
+z^4(584748316984283206\theta^4+585086560085496018\theta^3
\\[-2pt] &\;\quad
+187240905834613493\theta^2+20468485628797540\theta+703267501469328)
\\[-2pt] &\;
-43\cdot643z^5\theta(84565362438200\theta^3+190290414025938\theta^2
\\[-2pt] &\;\quad
+161530840617496\theta+49787404171133)
\\[-2pt] &\;
+43^2\cdot643^2z^6\theta(\theta+1)(4587245596\theta^2
+11468556386\theta+7890274651)
\\[-2pt] &\;
-2\cdot43^3\cdot643^3z^7\theta(\theta+1)(\theta+2)(55300\theta+96777)
\\[-2pt] &\;
+43^4\cdot643^4z^8\theta(\theta+1)(\theta+2)(\theta+3)
\\[-8pt]
\endaligned &\cr
& &\span\hrulefill&\cr
& &&
\gathered
\text{same $K(q)$ as in \#88}
\\
\text{the pullback of the 5th-order differential equation $D'y=0$, where}
\\[-8pt]
\endgathered &\cr
& &\span\hrulefill&\cr
& &&
\aligned
\\[-12pt]
&
D'=\theta^5
-z(2\theta+1)(13827\theta^4+27654\theta^3+19207\theta^2+5380\theta+481)
\\[-2pt] &\;
+z^2(\theta+1)(138255\theta^4+553020\theta^3+889449\theta^2+672858\theta+199135)
\\[-2pt] &\;
-2z^3(\theta+1)(\theta+2)(2\theta+3)(69125\theta^2+207375\theta+187791)
\\[-2pt] &\;
+z^4(\theta+1)(\theta+2)(\theta+3)(276495\theta^2+1105980\theta+1172033)
\\[-2pt] &\;
-69123z^5(\theta+1)(\theta+2)(\theta+3)(\theta+4)(2\theta+5)
\\[-2pt] &\;
+27649z^6(\theta+1)(\theta+2)(\theta+3)(\theta+4)(\theta+5)
\\[-8pt]
\endaligned &\cr
& &\span\hrulefill&\cr
& &&
A_n'=\sum_{k=0}^n\binom nk
\binom{2k}k^2\frac{(6k)!}{k!^3(3k)!}
&\cr
\noalign{\hrule}
& \eqnno{90} &&
\aligned
\\[-12pt]
&
D=\theta^4
-2z(2\theta+1)^2(56\theta^2+56\theta+13)
\\[-2pt] &\;
+20z^2(2\theta+1)(2\theta+3)(66\theta^2+132\theta+71)
\\[-2pt] &\;
-1312z^3(2\theta+1)(2\theta+3)^2(2\theta+5)
\\[-2pt] &\;
+1744z^4(2\theta+1)(2\theta+3)(2\theta+5)(2\theta+7)
\\[-8pt]
\endaligned &\cr
& &\span\hrulefill&\cr
& &&
A_n=\binom{2n}n\sum_{k=0}^n\binom nk\binom{2k}k\frac{(3k)!}{k!^3}
&\cr
\noalign{\hrule}
& \eqnno{91} &&
\aligned
\\[-12pt]
&
D=\theta^4
-z(437\theta^4+874\theta^3+646\theta^2209\theta+25)
\\[-2pt] &\;
+z^2(\theta+1)^2(1738\theta^2+3476\theta+2151)
\\[-2pt] &\;
-z^3(\theta+1)(\theta+2)(2602\theta^2+7806\theta+6277)
\\[-2pt] &\;
+1733z^4(\theta+1)(\theta+2)^2(\theta+3)
\\[-2pt] &\;
-433z^5(\theta+1)(\theta+2)(\theta+3)(\theta+4)
\\[-8pt]
\endaligned &\cr
& &\span\hrulefill&\cr
& &&
A_n=\sum_{k=0}^n\binom nk\binom{2k}k^2\frac{(3k)!}{k!^3}
&\cr
\noalign{\hrule}
}}\hss}

\newpage

\hbox to\hsize{\hss\vbox{\offinterlineskip
\halign to120mm{\strut\tabskip=100pt minus 100pt
\strut\vrule#&\hbox to5.5mm{\hss$#$\hss}&%
\vrule#&\hbox to114mm{\hfil$\dsize#$\hfil}&%
\vrule#\tabskip=0pt\cr\noalign{\hrule}
& \# &%
& \text{differential operator $D$ and coefficients $A_n$, $n=0,1,2,\dots$} &\cr
\noalign{\hrule\vskip1pt\hrule}
& \eqnno{92} &&
\aligned
\\[-12pt]
&
D=\theta^4
-2^2z(12\theta^4+16\theta^3+14\theta^2+6\theta+1)
\\[-2pt] &\;
+2^4z^2(789\theta^4+3076\theta^3+4167\theta^2+2222\theta+385)
\\[-2pt] &\;
-2^7z^3(2\theta+1)(1498\theta^3+6701\theta^2+9987\theta+4949)
\\[-2pt] &\;
+2^8z^4(2\theta+1)(2\theta+3)(4434\theta^2+1769\theta+18003)
\\[-2pt] &\;
-2^{11}z^5(2\theta+1)(2\theta+3)(2\theta+5)(1470\theta+3671)
\\[-2pt] &\;
+2^{12}\cdot733z^6(2\theta+1)(2\theta+3)(2\theta+5)(2\theta+7)
\\[-8pt]
\endaligned &\cr
& &\span\hrulefill&\cr
& &&
A_n=\binom{2n}n\sum_k(-1)^k2^{n-2k}\binom n{2k}\binom{2k}k^2\binom{6k}{2k}
&\cr
\noalign{\hrule}
& \eqnno{93} &&
\aligned
\\[-12pt]
&
D=\theta^4
-z(448\theta^4+872\theta^3+648\theta^2+212\theta+26)
\\[-2pt] &\;
+z^2(7008\theta^4+17376\theta^3+20664\theta^2+14496\theta+4652)
\\[-2pt] &\;
-z^3(41728\theta^4+124800\theta^3+175488\theta^2+122784\theta+28576)
\\[-2pt] &\;
+z^4(110848\theta^4+387584\theta^3+577920\theta^2+381824\theta+94096)
\\[-2pt] &\;
-3072z^5(\theta+1)^2(6\theta+5)(6\theta+7)
\\[-8pt]
\endaligned &\cr
& &\span\hrulefill&\cr
& &&
A_n=\binom{2n}n\sum_{k=0}^n\binom nk^2
\binom{2k}k\binom{2n}{2k}^{-1}\frac{(3k)!}{k!^3}
&\cr
\noalign{\hrule}
& \eqnno{94} &&
\aligned
\\[-12pt]
&
D=\theta^4
-z(27656\theta^4+6918\theta^3+1420\theta^2-2039\theta-1092)
\\[-2pt] &\;
+z^2(286848028\theta^4+143534634\theta^3+35472587\theta^2+608435\theta+24800580)
\\[-2pt] &\;
-z^3(1322624277560\theta^4+993115489662\theta^3+287226757326\theta^2
\\[-2pt] &\;\quad
+115129154686\theta-19751450796)
\\[-2pt] &\;
+z^4(2289130534388806\theta^4+2294419884217554\theta^3+762036995515973\theta^2
\\[-2pt] &\;\quad
+93972178515892\theta+2628910994892)
\\[-2pt] &\;
-31\cdot223z^5\theta(1322624277560\theta^3+2977051906002\theta^2
\\[-2pt] &\;\quad
+2534884699840\theta+784655821733)
\\[-2pt] &\;
+31^2\cdot223^2z^6\theta(\theta+1)(286848028\theta^2+717230690\theta+494009059)
\\[-2pt] &\;
-2\cdot31^3\cdot223^3z^7\theta(\theta+1)(\theta+2)(13828\theta+24201)
\\[-2pt] &\;
+31^4\cdot223^4z^8\theta(\theta+1)(\theta+2)(\theta+3)
\\[-8pt]
\endaligned &\cr
& &\span\hrulefill&\cr
& &&
\gathered
\text{the pullback of the 5th-order differential equation $D'y=0$, where}
\\[-8pt]
\endgathered &\cr
& &\span\hrulefill&\cr
& &&
\aligned
\\[-12pt]
&
D'=\theta^5
-z(2\theta+1)(3459\theta^4+6918\theta^3+4879\theta^2+1420\theta+145)
\\[-2pt] &\;
+z^2(\theta+1)(34575\theta^4+138300\theta^3+222873\theta^2+169146\theta+50287)
\\[-2pt] &\;
-2z^3(\theta+1)(\theta+2)(2\theta+3)(17285\theta^2+51855\theta+47067)
\\[-2pt] &\;
+z^4(\theta+1)(\theta+2)(\theta+3)(69135\theta^2+276540\theta+293201)
\\[-2pt] &\;
-17283z^5(\theta+1)(\theta+2)(\theta+3)(\theta+4)(2\theta+5)
\\[-2pt] &\;
+6913z^6(\theta+1)(\theta+2)(\theta+3)(\theta+4)(\theta+5)
\\[-8pt]
\endaligned &\cr
& &\span\hrulefill&\cr
& &&
A_n'=\sum_{k=0}^n\binom nk\binom{2k}k
\frac{(3k)!}{k!^3}\frac{(4k)!}{k!^2(2k)!}
&\cr
\noalign{\hrule}
}}\hss}

\newpage

\hbox to\hsize{\hss\vbox{\offinterlineskip
\halign to120mm{\strut\tabskip=100pt minus 100pt
\strut\vrule#&\hbox to5.5mm{\hss$#$\hss}&%
\vrule#&\hbox to114mm{\hfil$\dsize#$\hfil}&%
\vrule#\tabskip=0pt\cr\noalign{\hrule}
& \# &%
& \text{differential operator $D$ and coefficients $A_n$, $n=0,1,2,\dots$} &\cr
\noalign{\hrule\vskip1pt\hrule}
& \eqnno{95} &&
\aligned
\\[-12pt]
&
D=\theta^4
-z(1733\theta^4+3466\theta^3+2446\theta^2+713\theta+73)
\\[-2pt] &\;
+z^2(\theta+1)^2(6922\theta^2+13844\theta+8343)
\\[-2pt] &\;
-z^3(\theta+1)(\theta+2)(10378\theta^2+31134\theta+24925)
\\[-2pt] &\;
+6917z^4(\theta+1)(\theta+2)^2(\theta+3)
\\[-2pt] &\;
-1729z^5(\theta+1)(\theta+2)(\theta+3)(\theta+4)
\\[-8pt]
\endaligned &\cr
& &\span\hrulefill&\cr
& &&
A_n=\sum_{k=0}^n\binom nk\frac{(3k)!}{k!^3}\frac{(4k)!}{k!^2(2k)!}
&\cr
\noalign{\hrule}
& \eqnno{96} &&
\aligned
\\[-12pt]
&
D=\theta^4
-z(1040\theta^4+2056\theta^3+1484\theta^2+456\theta+50)
\\[-2pt] &\;
+4z^2(4120\theta^4+10264\theta^3+12062\theta^2+8396\theta+2687)
\\[-2pt] &\;
-16z^3(6160\theta^4+18456\theta^3+25752\theta^2+17862\theta+4081)
\\[-2pt] &\;
+16z^4(16400\theta^4+57376\theta^3+85016\theta^2+55560\theta+13441)
\\[-2pt] &\;
-2^{14}z^5(\theta+1)^2(4\theta+3)(4\theta+5)
\\[-8pt]
\endaligned &\cr
& &\span\hrulefill&\cr
& &&
A_n=\binom{2n}n\sum_{k=0}^n\binom nk^2\binom{2n}{2k}^{-1}\frac{(4k)!}{k!^4}
&\cr
\noalign{\hrule}
& \eqnno{97} &&
\aligned
\\[-12pt]
&
D=\theta^4
-z(16400\theta^4+4100\theta^3+902\theta^2-1148\theta-623)
\\[-2pt] &\;
+z^2(100925536\theta^4+50462768\theta^3+13232182\theta^2+847628\theta+8796151)
\\[-2pt] &\;
-z^3(276490092800\theta^4+207367569600\theta^3+63242834064\theta^2
\\[-2pt] &\;\quad
+26362467316\theta-3821932666)
\\[-2pt] &\;
+z^4(285882691092736\theta^4+285882691092736\theta^3+100029707616352\theta^2
\\[-2pt] &\;\quad
+15451954237456\theta+1121739103233)
\\[-2pt] &\;
-2^{11}z^5(2211920742400\theta^4+2764900928000\theta^3+2006667953392\theta^2
\\[-2pt] &\;\quad
+1165909352584\theta+332723569061)
\\[-2pt] &\;
+2^{20}z^6(25836937216\theta^4+38755405824\theta^3+33098847360\theta^2
\\[-2pt] &\;\quad
+13914543616\theta-569939335)
\\[-2pt] &\;
-2^{33}z^7(8396800\theta^4+14694400\theta^3+12330304\theta^2+3444000\theta+461113)
\\[-2pt] &\;
+2^{44}z^8(8\theta+3)^2(8\theta+5)^2
\\[-8pt]
\endaligned &\cr
& &\span\hrulefill&\cr
& &&
\gathered
\text{the pullback of the 5th-order differential equation $D'y=0$, where}
\\[-8pt]
\endgathered &\cr
& &\span\hrulefill&\cr
& &&
\aligned
\\[-12pt]
&
D'=\theta^5
-2z(2058\theta^5+5125\theta^4+5002\theta^3+2378\theta^2+549\theta+49)
\\[-2pt] &\;
+2^2z^2(20520\theta^5+61480\theta^4+90950\theta^3+89906\theta^2+55877\theta+15545)
\\[-2pt] &\;
-2^3z^3(82000\theta^5+286840\theta^4+507060\theta^3+532058\theta^2+277120\theta+38601)
\\[-2pt] &\;
+2^4z^4(163920\theta^5+655520\theta^4+1259720\theta^3+1286224\theta^2
\\[-2pt] &\;\quad
+611097\theta+122882)
\\[-2pt] &\;
-2^5z^5(163872\theta^5+737360\theta^4+1464400\theta^3+1446952\theta^2
\\[-2pt] &\;\quad
+709642\theta+138241)
\\[-2pt] &\;
+2^{18}z^6(\theta+1)^3(4\theta+3)(4\theta+5)
\\[-8pt]
\endaligned &\cr
& &\span\hrulefill&\cr
& &&
A_n'=\binom{2n}n\sum_{k=0}^n\binom nk^2\binom{2k}k^2
\binom{2n}{2k}^{-1}\frac{(4k)!}{k!^2(2k)!}
&\cr
\noalign{\hrule}
}}\hss}

\newpage

\hbox to\hsize{\hss\vbox{\offinterlineskip
\halign to120mm{\strut\tabskip=100pt minus 100pt
\strut\vrule#&\hbox to5.5mm{\hss$#$\hss}&%
\vrule#&\hbox to114mm{\hfil$\dsize#$\hfil}&%
\vrule#\tabskip=0pt\cr\noalign{\hrule}
& \# &%
& \text{differential operator $D$ and coefficients $A_n$, $n=0,1,2,\dots$} &\cr
\noalign{\hrule\vskip1pt\hrule}
& \eqnno{98} &&
\aligned
\\[-12pt]
&
D=\theta^4
-z(6928\theta^4+1732\theta^3+414\theta^2-452\theta-251)
\\[-2pt] &\;
+z^2(18026592\theta^4+9013296\theta^3+2541558\theta^2+295932\theta+1580631)
\\[-2pt] &\;
-z^3(20926439680\theta^4+15694829760\theta^3+5119446672\theta^2
\\[-2pt] &\;\quad
+2221999476\theta-262743618)
\\[-2pt] &\;
+z^4(9248048087296\theta^4+9248048087296\theta^3+3470622495840\theta^2
\\[-2pt] &\;\quad
+668590646032\theta+65658526849)
\\[-2pt] &\;
-2^4\cdot24z^5(376675914240\theta^4+470844892800\theta^3+348726894896\theta^2
\\[-2pt] &\;\quad
+202514342736\theta+57000462817)
\\[-2pt] &\;
+2^6\cdot24^2z^6(23362463232\theta^4+35043694848\theta^3+30328003584\theta^2
\\[-2pt] &\;\quad
+12869932032\theta-426340673)
\\[-2pt] &\;
-2^4\cdot24^5z^7(17957376\theta^4+31425408\theta^3+26673840\theta^2
\\[-2pt] &\;\quad
+7638120\theta+1060093)
\\[-2pt] &\;
+24^8z^8(12\theta+5)^2(12\theta+7)^2
\\[-8pt]
\endaligned
&\cr
& &\span\hrulefill&\cr
& &&
\gathered
\text{the pullback of the 5th-order differential equation $D'y=0$, where}
\\[-8pt]
\endgathered &\cr
& &\span\hrulefill&\cr
& &&
\aligned
\\[-12pt]
&
D'=\theta^5
-2z(874\theta^5+2165\theta^4+2146\theta^3+1054\theta^2+257\theta+25)
\\[-2pt] &\;
+2^2z^2(8680\theta^5+25960\theta^4+38710\theta^3+38498\theta^2+23981\theta+6689)
\\[-2pt] &\;
-2^3z^3(34640\theta^5+121080\theta^4+215220\theta^3+227034\theta^2+119160\theta+17085)
\\[-2pt] &\;
+2^4z^4(69200\theta^5+276640\theta^4+533960\theta^3+548432\theta^2+263401\theta+53762)
\\[-2pt] &\;
-2^5z^5(69152\theta^5+311120\theta^4+620240\theta^3+616936\theta^2+305866\theta+60481)
\\[-2pt] &\;
+2^{14}\cdot3z^6(\theta+1)^3(6\theta+5)(6\theta+7)
\\[-8pt]
\endaligned &\cr
& &\span\hrulefill&\cr
& &&
A_n'=\binom{2n}n\sum_{k=0}^n\binom nk^2
\binom{2k}k^2\binom{2n}{2k}^{-1}\frac{(3k)!}{k!^3}
&\cr
\noalign{\hrule}
& \eqnno{99} &&
\aligned
\\[-12pt]
&
D=13^2\theta^4
-z(59397\theta^4+117546\theta^3+86827\theta^2+28054\theta+3380)
\\[-2pt] &\;
+2^4z^2(6386\theta^4-1774\theta^3-17898\theta^2-11596\theta-2119)
\\[-2pt] &\;
+2^8z^3(67\theta^4+1248\theta^3+1091\theta^2+312\theta+26)
-2^{12}z^4(2\theta+1)^4
\\[-8pt]
\endaligned &\cr
& &\span\hrulefill&\cr
& &&
\aligned
A_n
&=\binom{2n}n^2\sum_{k=0}^n\binom nk^2\binom{2n+k}n
\\[-2pt]
&=\sum_{k,l}\binom nk\binom nl\binom{n+k}n\binom{n+l}n
\binom{2n+l}n\binom n{l-k}
\endaligned &\cr
\noalign{\hrule}
& \eqnno{100} &&
\aligned
\\[-12pt]
&
D=\theta^4
-z(73\theta^4+98\theta^3+77\theta^2+28\theta+4)
\\[-2pt] &\;
+z^2(520\theta^4-1040\theta^3-2904\theta^2-2048\theta-480)
\\[-2pt] &\;
+2^6z^3(65\theta^4+390\theta^3+417\theta^2+180\theta+28)
\\[-2pt] &\;
-2^9z^4(73\theta^4+194\theta^3+221\theta^2+124\theta+28)
+2^{15}z^5(\theta+1)^4
\\[-8pt]
\endaligned &\cr
& &\span\hrulefill&\cr
& &&
A_n=\Biggl\{\sum_{k=0}^n\binom nk^3\Biggr\}^2
&\cr
\noalign{\hrule}
}}\hss}

\newpage

\hbox to\hsize{\hss\vbox{\offinterlineskip
\halign to120mm{\strut\tabskip=100pt minus 100pt
\strut\vrule#&\hbox to5.5mm{\hss$#$\hss}&%
\vrule#&\hbox to114mm{\hfil$\dsize#$\hfil}&%
\vrule#\tabskip=0pt\cr\noalign{\hrule}
& \# &%
& \text{differential operator $D$ and coefficients $A_n$, $n=0,1,2,\dots$} &\cr
\noalign{\hrule\vskip1pt\hrule}
& \eqnno{101} &&
\aligned
\\[-12pt]
&
D=\theta^4
-z(124\theta^4+242\theta^3+187\theta^2+66\theta+9)
\\[-2pt] &\;
+z^2(123\theta^4-246\theta^3-787\theta^2-554\theta-124)
\\[-2pt] &\;
+z^3(123\theta^4+738\theta^3+689\theta^2+210\theta+12)
\\[-2pt] &\;
-z^4(124\theta^4+254\theta^3+205\theta^2+78\theta+12)
+z^5(\theta+1)^4
\\[-8pt]
\endaligned &\cr
& &\span\hrulefill&\cr
& &&
\aligned
A_n
&=\Biggl\{\sum_{k=0}^n\binom nk^2\binom{n+k}k\Biggr\}^2
\\[-2pt]
&=\sum_{k,l}\binom nk\binom nl\binom{n+k}n\binom{n+l}n
\binom{n+l-k}n\binom n{l-k}
\\[2pt]
\endaligned &\cr
\noalign{\hrule}
& \eqnno{102} &&
\aligned
\\[-12pt]
&
D=\theta^4
-z(7\theta^2+7\theta+2)(11\theta^2+11\theta+3)
\\[-2pt] &\;
-z^2(1049\theta^4+4100\theta^3+5689\theta^2+3178\theta+640)
\\[-2pt] &\;
+2^3z^3(77\theta^4-462\theta^3-1420\theta^2-1053\theta-252)
\\[-2pt] &\;
+2^4z^4(1041\theta^4+2082\theta^3-1406\theta^2-2447\theta-746)
\\[-2pt] &\;
+2^6z^5(77\theta^4+770\theta^3+428\theta^2-93\theta-80)
\\[-2pt] &\;
-2^6z^6(1049\theta^4+96\theta^3-317\theta^2+96\theta+100)
\\[-2pt] &\;
-2^9z^7(7\theta^2+7\theta+2)(11\theta^2+11\theta+3)
+2^{12}z^8(\theta+1)^4
\\[-8pt]
\endaligned &\cr
& &\span\hrulefill&\cr
& &&
A_n=\sum_{k=0}^n\binom nk^3\sum_{k=0}^n\binom nk^2\binom{n+k}k
&\cr
\noalign{\hrule}
& \eqnno{103} &&
\aligned
\\[-12pt]
&
D=\theta^4
-z(73\theta^4+200\theta^3+160\theta^2+60\theta+9)
\\[-2pt] &\;
-z^2(738\theta^4-1476\theta^3-5274\theta^2-3816\theta-918)
\\[-2pt] &\;
+z^3(6642\theta^4+39852\theta^3+32238\theta^2+5832\theta-1458)
\\[-2pt] &\;
+3^6z^4(73\theta^4+92\theta^3-2\theta^2-48\theta-18)
-3^{10}z^5(\theta+1)^4
\\[-8pt]
\endaligned &\cr
& &\span\hrulefill&\cr
& &&
A_n=\Biggl\{\sum_{k=0}^n\binom nk^2\binom{2k}k\Biggr\}^2
&\cr
\noalign{\hrule}
& \eqnno{104} &&
\aligned
\\[-12pt]
&
D=\theta^4
-z(10\theta^2+10\theta+3)(7\theta^2+7\theta+2)
\\[-2pt] &\;
-z^2(71\theta^4+1148\theta^3+1591\theta^2+886\theta+192)
\\[-2pt] &\;
-2^3\cdot3^2z^3(70\theta^4-420\theta^3-1289\theta^2-963\theta-240)
\\[-2pt] &\;
-2^4\cdot3^2z^4(143\theta^4+286\theta^3-1138\theta^2-1281\theta-414)
\\[-2pt] &\;
+2^6\cdot3^4z^5(70\theta^4+700\theta^3+391\theta^2-75\theta-76)
\\[-2pt] &\;
+2^6\cdot3^4z^6(-71\theta^4+864\theta^3+1427\theta^2+864\theta+180)
\\[-2pt] &\;
+2^9\cdot3^6z^7(10\theta^2+10\theta+3)(7\theta^2+7\theta+2)
+2^{12}\cdot3^8z^8(\theta+1)^4
\\[-8pt]
\endaligned &\cr
& &\span\hrulefill&\cr
& &&
A_n=\sum_{k=0}^n\binom nk^3\cdot\sum_{k=0}^n\binom nk^2\binom{2k}k
&\cr
\noalign{\hrule}
& \eqnno{105} &&
\aligned
\\[-12pt]
&
D=\theta^4
-4z(3\theta^2+3\theta+1)(7\theta^2+7\theta+1)
\\[-2pt] &\;
+2^5z^2(45\theta^4+84\theta^3+117\theta^2+66\theta+12)
\\[-2pt] &\;
-2^{10}z^3(21\theta^4-126\theta^3-386\theta^2-291\theta-76)
\\[-2pt] &\;
+2^{14}z^4(37\theta^4+74\theta^3+50\theta^2+13\theta+6)
\\[-2pt] &\;
+2^{18}z^5(21\theta^4+210\theta^3+118\theta^2-19\theta-24)
\\[-2pt] &\;
+3\cdot 2^{21}z^6(15\theta^4+32\theta^3+45\theta^2+32\theta+8)
\\[-2pt] &\;
+2^{26}z^7(3\theta^2+3\theta+1)(7\theta^{}+7\theta+2)
+2^{32}z^8(\theta+1)^4
\\[-8pt]
\endaligned &\cr
& &\span\hrulefill&\cr
& &&
A_n=\sum_{k=0}^n\binom nk^3\cdot\sum_{k=0}^n\binom nk
\binom{2k}k\binom{2n-2k}{n-k}
&\cr
\noalign{\hrule}
}}\hss}

\newpage

\hbox to\hsize{\hss\vbox{\offinterlineskip
\halign to120mm{\strut\tabskip=100pt minus 100pt
\strut\vrule#&\hbox to5.5mm{\hss$#$\hss}&%
\vrule#&\hbox to114mm{\hfil$\dsize#$\hfil}&%
\vrule#\tabskip=0pt\cr\noalign{\hrule}
& \# &%
& \text{differential operator $D$ and coefficients $A_n$, $n=0,1,2,\dots$} &\cr
\noalign{\hrule\vskip1pt\hrule}
& \eqnno{106} &&
\aligned
\\[-12pt]
&
D=\theta^4
-4z(3\theta^2+3\theta+1)(11\theta^2+11\theta+3)
\\[-2pt] &\;
+2^4z^2(241\theta^4+940\theta^3+1303\theta^2+726\theta+145)
\\[-2pt] &\;
-2^7z^3(33\theta^4-198\theta^3-607\theta^2-456\theta-117)
\\[-2pt] &\;
+2^{10}z^4(239\theta^4+478\theta^3-322\theta^2-561\theta-169)
\\[-2pt] &\;
+2^{12}z^5(33\theta^4+330\theta^3+185\theta^2-32\theta-37)
\\[-2pt] &\;
+2^{14}z^6(241\theta^4+24\theta^3-71\theta^2+24\theta+23)
\\[-2pt] &\;
+2^{17}z^7(3\theta^2+3\theta+1)(11\theta^2+11\theta+3)
+2^{20}z^8(\theta+1)^4
\\[-8pt]
\endaligned &\cr
& &\span\hrulefill&\cr
& &&
A_n=\sum_{k=0}^n\binom nk^2\binom{n+k}k
\cdot\sum_{k=0}^n\binom nk\binom{2k}k\binom{2n-2k}{n-k}
&\cr
\noalign{\hrule}
& \eqnno{107} &&
\aligned
\\[-12pt]
&
D=\theta^4
-2^4z(3\theta^4+18\theta^3+15\theta^2+6\theta+1)
\\[-2pt] &\;
-2^9z^2(5\theta^4-10\theta^3-45\theta^2-34\theta-9)
\\[-2pt] &\;
+2^{14}z^3(5\theta^4+30\theta^3+15\theta^2-6\theta-5)
\\[-2pt] &\;
+2^{19}z^4(3\theta^4-6\theta^3-21\theta^2-18\theta-5)
-2^{25}z^5(\theta+1)^4
\\[-8pt]
\endaligned &\cr
& &\span\hrulefill&\cr
& &&
A_n=\Biggl\{\sum_{k=0}^n\binom nk\binom{2k}k\binom{2n-2k}{n-k}\Biggr\}^2
&\cr
\noalign{\hrule}
& \eqnno{108} &&
\aligned
\\[-12pt]
&
D=\theta^4
-2z(6\theta^4+10\theta^3+10\theta^2+5\theta+1)
\\[-2pt] &\;
+2^2z^2(46671\theta^4+186674\theta^3+238539\theta^2+103735\theta+14416)
\\[-2pt] &\;
-2^3z^3(\theta+1)(186644\theta^3+839888\theta^2+1223544\theta+575489)
\\[-2pt] &\;
+2^4z^4(\theta+1)(\theta+2)(279951\theta^2+1119788\theta+1124983)
\\[-2pt] &\;
-2^6z^5(\theta+1)(\theta+2)(\theta+3)(93315\theta+233287)
\\[-2pt] &\;
+2^6\cdot13\cdot37\cdot97z^6(\theta+1)(\theta+2)(\theta+3)(\theta+4)
\\[-8pt]
\endaligned &\cr
& &\span\hrulefill&\cr
& &&
A_n=\sum_k(-1)^k2^{n-2k}\binom n{2k}
\biggl(\frac{(6k)!}{k!(2k)!(3k)!}\biggr)^2
&\cr
\noalign{\hrule}
& \eqnno{109} &&
\aligned
\\[-12pt]
&
D=49\theta^4
-2\cdot3z(8904\theta^4+17556\theta^3+12453\theta^2+3675\theta+392)
\\[-2pt] &\;
+2^2\cdot3z^2(43704\theta^4+38088\theta^3-25757\theta^2-20608\theta-3360)
\\[-2pt] &\;
-2^4\cdot3^3z^3(2736\theta^4-1512\theta^3-1672\theta^2-357\theta-14)
\\[-2pt] &\;
-2^6\cdot3^5z^4(2\theta+1)^2(3\theta+1)(3\theta+2)
\\[-8pt]
\endaligned &\cr
& &\span\hrulefill&\cr
& &&
A_n=\binom{2n}n^2\sum_{k=0}^{2n}\binom{n+k}k\binom{2n}k^2
&\cr
\noalign{\hrule}
& \eqnno{110} &&
\aligned
\\[-12pt]
&
D=\theta^4
-12z(3\theta+1)(3\theta+2)(8\theta^2+8\theta+3)
\\[-2pt] &\;
+2304z^2(3\theta+1)(3\theta+2)(3\theta+4)(3\theta+5)
\\[-8pt]
\endaligned &\cr
& &\span\hrulefill&\cr
& &&
A_n=\frac{(3n)!}{n!^3}\sum_{k=0}^n4^{n-k}\binom{2k}k^2\binom{2n-2k}{n-k}
&\cr
\noalign{\hrule}
& \eqnno{111} &&
\aligned
\\[-12pt]
&
D=\theta^4
-16z(2\theta+1)^2(8\theta^2+8\theta+3)
\\[-2pt] &\;
+4096z^2(2\theta+1)^2(2\theta+3)^2
\\[-8pt]
\endaligned &\cr
& &\span\hrulefill&\cr
& &&
A_n=\binom{2n}n^2\sum_{k=0}^n4^{n-k}\binom{2k}k^2\binom{2n-2k}{n-k}
&\cr
\noalign{\hrule}
}}\hss}

\newpage

\hbox to\hsize{\hss\vbox{\offinterlineskip
\halign to120mm{\strut\tabskip=100pt minus 100pt
\strut\vrule#&\hbox to5.5mm{\hss$#$\hss}&%
\vrule#&\hbox to114mm{\hfil$\dsize#$\hfil}&%
\vrule#\tabskip=0pt\cr\noalign{\hrule}
& \# &%
& \text{differential operator $D$ and coefficients $A_n$, $n=0,1,2,\dots$} &\cr
\noalign{\hrule\vskip1pt\hrule}
& \eqnno{112} &&
\aligned
\\[-12pt]
&
D=\theta^4
-48z(6\theta+1)(6\theta+5)(8\theta^2+8\theta+3)
\\[-2pt] &\;
+36864z^2(6\theta+1)(6\theta+5)(6\theta+7)(6\theta+11)
\\[-8pt]
\endaligned &\cr
& &\span\hrulefill&\cr
& &&
A_n=\frac{(6n)!}{n!(2n)!(3n)!}
\sum_{k=0}^n4^{n-k}\binom{2k}k^2\binom{2n-2k}{n-k}
&\cr
\noalign{\hrule}
& \eqnno{113} &&
\aligned
\\[-12pt]
&
D=\theta^4
-z(10\theta^2+10\theta+3)(11\theta^2+11\theta+3)
\\[-2pt] &\;
+z^2(1025\theta^4+3992\theta^3+5533\theta^2+3082\theta+615)
\\[-2pt] &\;
+3^2z^3(-110\theta^4+660\theta^3+2027\theta^2+1509\theta+369)
\\[-2pt] &\;
+3^2z^4(2032\theta^4+4064\theta^3-2726\theta^2-4758\theta-1431)
\\[-2pt] &\;
+3^4z^5(110\theta^4+1100\theta^3+613\theta^2-125\theta-117)
\\[-2pt] &\;
+3^4z^6(1025\theta^4+108\theta^3-293\theta^2+108\theta+99)
\\[-2pt] &\;
+3^6z^7(10\theta^2+10\theta+3)(11\theta^2+11\theta+3)
+3^8z^8(\theta+1)^4
\\[-8pt]
\endaligned &\cr
& &\span\hrulefill&\cr
& &&
A_n=\sum_{k=0}^n\binom nk^2\binom{n+k}k
\cdot\sum_{k=0}^n\binom nk^2\binom{2k}k
&\cr
\noalign{\hrule}
& \eqnno{114} &&
\aligned
\\[-12pt]
&
D=\theta^4
-4z(7\theta^2+7\theta+2)(8\theta^2+8\theta+3)
\\[-2pt] &\;
+2^7z^2(98\theta^4+200\theta^3+274\theta^2+148\theta+23)
\\[-2pt] &\;
+2^{13}z^3(-56\theta^4+336\theta^3+1027\theta^2+783\theta+216)
\\[-2pt] &\;
+2^{19}z^4(82\theta^4+164\theta^3+84\theta^2+2\theta+11)
\\[-2pt] &\;
+2^{24}z^5(56\theta^4+560\theta^3+317\theta^2-39\theta-68)
\\[-2pt] &\;
+2^{29}z^6(98\theta^4+192\theta^3+262\theta^2+192\theta+47)
\\[-2pt] &\;
+2^{35}z^7(7\theta^2+7\theta+2)(8\theta^2+8\theta+3)
+2^{44}z^8(\theta+1)^4
\\[-8pt]
\endaligned &\cr
& &\span\hrulefill&\cr
& &&
A_n=\sum_{k=0}^n\binom nk^3
\cdot\sum_{k=0}^n4^{n-k}\binom{2k}k^2\binom{2n-2k}{n-k}
&\cr
\noalign{\hrule}
& \eqnno{115} &&
\aligned
\\[-12pt]
&
D=\theta^4
-2^4z(16\theta^4+128\theta^3+112\theta^2+48\theta+9)
\\[-2pt] &\;
+2^{12}z^2(-32\theta^4+64\theta^3+304\theta^2+240\theta+71)
\\[-2pt] &\;
+2^{20}z^3(32\theta^4+192\theta^3+80\theta^2-48\theta-39)
\\[-2pt] &\;
+2^{28}z^4(16\theta^4-64\theta^3-176\theta^2-144\theta-39)
-2^{40}z^5(\theta+1)^4
\\[-8pt]
\endaligned &\cr
& &\span\hrulefill&\cr
& &&
A_n=\biggl\{\sum_{k=0}^n4^{n-k}\binom{2k}k^2\binom{2n-2k}{n-k}\biggr\}^2
&\cr
\noalign{\hrule}
& \eqnno{116} &&
\aligned
\\[-12pt]
&
D=\theta^4
-2^5z(10\theta^4+26\theta^3+20\theta^2+7\theta+1)
\\[-2pt] &\;
+2^8z^2(52\theta^4+472\theta^3+832\theta^2+492\theta+103)
\\[-2pt] &\;
+2^{16}z^3(14\theta^4+12\theta^3-96\theta^2-105\theta-29)
\\[-2pt] &\;
-2^{18}z^4(2\theta+1)(56\theta^3+468\theta^2+646\theta+249)
\\[-2pt] &\;
-2^{24}z^5(2\theta+1)(2\theta+3)(4\theta+3)(4\theta+5)
\\[-8pt]
\endaligned &\cr
& &\span\hrulefill&\cr
& &&
A_n=\binom{2n}n^2\sum_{k=0}^n4^{n-k}\binom nk^2
\binom{n+k}k\binom{2k}k\binom{2n}{2k}^{-1}
&\cr
\noalign{\hrule}
}}\hss}

\newpage

\hbox to\hsize{\hss\vbox{\offinterlineskip
\halign to120mm{\strut\tabskip=100pt minus 100pt
\strut\vrule#&\hbox to5.5mm{\hss$#$\hss}&%
\vrule#&\hbox to114mm{\hfil$\dsize#$\hfil}&%
\vrule#\tabskip=0pt\cr\noalign{\hrule}
& \# &%
& \text{differential operator $D$ and coefficients $A_n$, $n=0,1,2,\dots$} &\cr
\noalign{\hrule\vskip1pt\hrule}
& \eqnno{117} &&
\aligned
\\[-12pt]
&
D=3^2\theta^4
+12z(256\theta^4+176\theta^3+133\theta^2+45\theta+6)
\\[-2pt] &\;
+2^7z^2(2588\theta^4+1952\theta^3+584\theta^2+15\theta-15)
\\[-2pt] &\;
+2^{12}z^3(3183\theta^4+2466\theta^3+1801\theta^2+711\theta+111)
\\[-2pt] &\;
+2^{17}\cdot7z^4(134\theta^4+250\theta^3+180\theta^2+55\theta+5)
\\[-2pt] &\;
-2^{22}\cdot7^2z^5(\theta+1)^4
\\[-8pt]
\endaligned &\cr
& &\span\hrulefill&\cr
& &&
\text{a formula for $A_n$ is not known}
&\cr
\noalign{\hrule}
& \eqnno{118} &&
\aligned
\\[-12pt]
&
D=\theta^4
-z(465\theta^4+594\theta^3+431\theta^2+134\theta+16)
\\[-2pt] &\;
+2^4z^2(2625\theta^4+1911\theta^3-946\theta^2-884\theta-176)
\\[-2pt] &\;
+2^6z^3(-16105\theta^4+3624\theta^3+5241\theta^2+1284\theta+36)
\\[-2pt] &\;
-2^{11}\cdot7z^4(155\theta^4+334\theta^3+306\theta^2+139\theta+26)
+2^{16}\cdot7^2z^5(\theta+1)^4
\\[-8pt]
\endaligned &\cr
& &\span\hrulefill&\cr
& &&
\text{a formula for $A_n$ is not known}
&\cr
\noalign{\hrule}
& \eqnno{119} &&
\aligned
\\[-12pt]
&
D=9\theta^4
-12z(256\theta^4+320\theta^3+271\theta^2+111\theta+18)
\\[-2pt] &\;
+2^7z^2(3104\theta^4+7040\theta^3+8012\theta^2+4452\theta+927)
\\[-2pt] &\;
-2^{15}z^3(752\theta^4+2304\theta^3+3042\theta^2+1854\theta+405)
\\[-2pt] &\;
+2^{21}z^4(2\theta+1)(176\theta^3+552\theta^2+622\theta+231)
\\[-2pt] &\;
-2^{31}z^5(\theta+1)^2(2\theta+1)(2\theta+3)
\\[-8pt]
\endaligned &\cr
& &\span\hrulefill&\cr
& &&
A_n=\sum_{k=0}^n\binom nk\frac{(4k)!}{k!^2(2k)!}
\,\frac{(4n-4k)!}{(n-k)!^2(2n-2k)!^2}
&\cr
\noalign{\hrule}
& \eqnno{120} &&
\aligned
\\[-12pt]
&
D=\theta^4
-4z(8\theta^2+8\theta+3)(10\theta^2+10\theta+3)
\\[-2pt] &\;
+2^4z^2(1600\theta^4+8128\theta^3+11408\theta^2+6560\theta+1473)
\\[-2pt] &\;
+2^{10}\cdot3^2z^3(80\theta^4-480\theta^3-1466\theta^2-1122\theta-315)
\\[-2pt] &\;
-2^{13}\cdot3^2z^4(1744\theta^4+3488\theta^3-4256\theta^2-6000\theta-2079)
\\[-2pt] &\;
+2^{18}\cdot3^4z^5(80\theta^4+800\theta^3+454\theta^2-50\theta-99)
\\[-2pt] &\;
+2^{20}\cdot3^4z^6(1600\theta^4-1728\theta^3-3376\theta^2-1728\theta-207)
\\[-2pt] &\;
-2^{26}\cdot3^6z^7(8\theta^2+8\theta+3)(10\theta^2+10\theta+3)
+2^{32}\cdot3^8z^8(\theta+1)^4
\\[-8pt]
\endaligned &\cr
& &\span\hrulefill&\cr
& &&
A_n=\sum_{k=0}^n\binom nk^2\binom{2k}k
\cdot\sum_{k=0}^n4^{n-k}\binom{2k}k^2\binom{2n-2k}{n-k}
&\cr
\noalign{\hrule}
& \eqnno{121} &&
\aligned
\\[-12pt]
&
D=\theta^4
-4z(8\theta^2+8\theta+3)(11\theta^2+11\theta+3)
\\[-2pt] &\;
+2^4z^2(1936\theta^4+7552\theta^3+10464\theta^2+5824\theta+1159)
\\[-2pt] &\;
+2^{10}z^3(-88\theta^4+528\theta^3+1615\theta^2+1227\theta+333)
\\[-2pt] &\;
+2^{13}z^4(1920\theta^4+3840\theta^3-2592\theta^2-4512\theta-1353)
\\[-2pt] &\;
+2^{18}z^5(88\theta^4+880\theta^3+497\theta^2-67\theta-105)
\\[-2pt] &\;
+2^{20}z^6(1936\theta^4+192\theta^3-576\theta^2+192\theta+183)
\\[-2pt] &\;
+2^{26}z^7(8\theta^2+8\theta+3)(11\theta^2+11\theta+3)
+2^{30}z^8(\theta+1)^4
\\[-8pt]
\endaligned &\cr
& &\span\hrulefill&\cr
& &&
A_n=\sum_{k=0}^n\binom nk^2\binom{n+k}k
\cdot\sum_{k=0}^n4^{n-k}\binom{2k}k^2\binom{2n-2k}{n-k}
&\cr
\noalign{\hrule}
}}\hss}

\newpage

\hbox to\hsize{\hss\vbox{\offinterlineskip
\halign to120mm{\strut\tabskip=100pt minus 100pt
\strut\vrule#&\hbox to5.5mm{\hss$#$\hss}&%
\vrule#&\hbox to114mm{\hfil$\dsize#$\hfil}&%
\vrule#\tabskip=0pt\cr\noalign{\hrule}
& \# &%
& \text{differential operator $D$ and coefficients $A_n$, $n=0,1,2,\dots$} &\cr
\noalign{\hrule\vskip1pt\hrule}
& \eqnno{122} &&
\aligned
\\[-12pt]
&
D=\theta^4
-2^4z(3\theta^2+3\theta+1)(8\theta^2+8\theta+3)
\\[-2pt] &\;
+2^9z^2(72\theta^4+480\theta^3+680\theta^2+400\theta+97)
\\[-2pt] &\;
+2^{17}z^3(24\theta^4-144\theta^3-439\theta^2-339\theta-99)
\\[-2pt] &\;
+2^{23}z^4(-88\theta^4-176\theta^3+320\theta^2+408\theta+151)
\\[-2pt] &\;
+2^{30}z^5(24\theta^4+240\theta^3+137\theta^2-11\theta-31)
\\[-2pt] &\;
+2^{35}z^6(72\theta^4-192\theta^3-328\theta^2-192\theta-31)
\\[-2pt] &\;
-2^{43}z^7(3\theta^2+3\theta+1)(8\theta^2+8\theta+3)
+2^{52}z^8(\theta+1)^4
\\[-8pt]
\endaligned &\cr
& &\span\hrulefill&\cr
& &&
A_n=\sum_{k=0}^n\binom nk\binom{2k}k\binom{2n-2k}{n-k}
\cdot\sum_{k=0}^n4^{n-k}\binom{2k}k^2\binom{2n-2k}{n-k}
&\cr
\noalign{\hrule}
& \eqnno{123} &&
\aligned
\\[-12pt]
&
D=\theta^4
-4z(3\theta^2+3\theta+1)(10\theta^2+10\theta+3)
\\[-2pt] &\;
+2^4z^2(209\theta^4+1052\theta^3+1471\theta^2+838\theta+183)
\\[-2pt] &\;
+2^7\cdot3^2z^3(30\theta^4-180\theta^3-551\theta^2-417\theta-111)
\\[-2pt] &\;
-2^{10}\cdot3^2z^4(227\theta^4+454\theta^3-550\theta^2-777\theta-261)
\\[-2pt] &\;
+2^{12}\cdot3^4z^5(30\theta^4+300\theta^3+169\theta^2-25\theta-35)
\\[-2pt] &\;
+2^{14}\cdot3^4z^6(209\theta^4-216\theta^3-431\theta^2-216\theta-27)
\\[-2pt] &\;
-2^{17}\cdot3^6z^7(3\theta^2+3\theta+1)(10\theta^2+10\theta+3)
+2^{20}\cdot3^8z^8(\theta+1)^4
\\[-8pt]
\endaligned &\cr
& &\span\hrulefill&\cr
& &&
A_n=\sum_{k=0}^n\binom{n+k}k^2\binom{2k}k
\cdot\sum_{k=0}^n4^{n-k}\binom nk^2\binom{2k}k
&\cr
\noalign{\hrule}
& \eqnno{124} &&
\aligned
\\[-12pt]
&
D=61^2\theta^4
-61z(3029\theta^4+5572\theta^3+4677\theta^2+1891\theta+305)
\\[-2pt] &\;
+z^2(1215215\theta^4+3428132\theta^3+4267228\theta^2+2572675\theta+611586)
\\[-2pt] &\;
-3^4z^3(39370\theta^4+140178\theta^3+206807\theta^2+142191\theta+37332)
\\[-2pt] &\;
+3^8z^4(566\theta^4+2230\theta^3+3356\theta^2+2241\theta+558)
-3^{13}z^5(\theta+1)^4
\\[-8pt]
\endaligned &\cr
& &\span\hrulefill&\cr
& &&
A_n=\sum_{k,l}\binom nk^2\binom nl\binom kl\binom{k+l}k\binom{2n-k-l}{n-k}
&\cr
\noalign{\hrule}
& \eqnno{125} &&
\aligned
\\[-12pt]
&
D=\theta^4
-z(11669\theta^4+23338\theta^3+15886\theta^2+4217\theta+361)
\\[-2pt] &\;
+z^2(\theta+1)^2(4666\theta^2+93332\theta+55095)
\\[-2pt] &\;
-z^3(\theta+1)(\theta+2)(69994\theta^2+209982\theta+167533)
\\[-2pt] &\;
+46661z^4(\theta+1)(\theta+2)^2(\theta+3)
\\[-2pt] &\;
-11665z^5(\theta+1)(\theta+2)(\theta+3)(\theta+4)
\\[-8pt]
\endaligned &\cr
& &\span\hrulefill&\cr
& &&
A_n=\sum_{k=0}^n\binom nk\frac{(6k)!}{k!^4(2k)!}
&\cr
\noalign{\hrule}
}}\hss}

\newpage

\hbox to\hsize{\hss\vbox{\offinterlineskip
\halign to120mm{\strut\tabskip=100pt minus 100pt
\strut\vrule#&\hbox to5.5mm{\hss$#$\hss}&%
\vrule#&\hbox to114mm{\hfil$\dsize#$\hfil}&%
\vrule#\tabskip=0pt\cr\noalign{\hrule}
& \# &%
& \text{differential operator $D$ and coefficients $A_n$, $n=0,1,2,\dots$} &\cr
\noalign{\hrule\vskip1pt\hrule}
& \eqnno{126} &&
\aligned
\\[-12pt]
&
D=\theta^4
-z(186656\theta^4+46676\theta^3+8434\theta^2-14904\theta-7813)
\\[-2pt] &\;
+z^2(13065919936\theta^4+6535199792\theta^3+1453643198\theta^2
\\[-2pt] &\;\quad
-105990748\theta+1115889227)
\\[-2pt] &\;
-z^3(406553346039296\theta^4+305071791609408\theta^3+80614299041712\theta^2
\\[-2pt] &\;\quad
+30335098270388\theta-5908917544614)
\\[-2pt] &\;
+z^4(4746509269839709696\theta^4+4751387282863861504\theta^3
\\[-2pt] &\;\quad
+1455612685803281696\theta^2+132330243679138384\theta
\\[-2pt] &\;\quad
-845722995408775)
\\[-2pt] &\;
-2^4\cdot5\cdot2333z^5(4\theta-1)(101638336509824\theta^3
+228725452667088\theta^2
\\[-2pt] &\;\quad
+193462702609348\theta+59317557294355)
\\[-2pt] &\;
+2^5\cdot5^2\cdot2333^2z^6(4\theta-1)(4\theta+3)(6532959968\theta^2
\\[-2pt] &\;\quad
+16333519832\theta+11222488415)
\\[-2pt] &\;
-2^8\cdot5^3\cdot2333^3z^7(4\theta-1)(4\theta+3)(4\theta+7)(46664\theta+81665)
\\[-2pt] &\;
+2^8\cdot5^4\cdot2333^4z^8(4\theta-1)(4\theta+3)(4\theta+7)(4\theta+11)
\\[-8pt]
\endaligned &\cr
& &\span\hrulefill&\cr
& &&
\gathered
\text{same $K(q)$ as in \#82}
\\
\text{the pullback of the 5th-order differential equation $D'y=0$, where}
\\[-8pt]
\endgathered &\cr
& &\span\hrulefill&\cr
& &&
\aligned
\\[-12pt]
&
D'=\theta^5
-2z(2\theta+1)(11669\theta^4+23338\theta^3+15886\theta^2+4217\theta+361)
\\[-2pt] &\;
+4z^2(\theta+1)(2\theta+1)(2\theta+3)(4666\theta^2+93332\theta+55095)
\\[-2pt] &\;
-8z^3(2\theta+1)(2\theta+3)(2\theta+5)(69994\theta^2+209982\theta+167533)
\\[-2pt] &\;
+746576z^4(\theta+2)(2\theta+1)(2\theta+3)(2\theta+5)(2\theta+7)
\\[-2pt] &\;
-373280z^5(2\theta+1)(2\theta+3)(2\theta+5)(2\theta+7)(2\theta+9)
\endaligned &\cr
& &\span\hrulefill&\cr
& &&
A_n'=\binom{2n}n\sum_{k=0}^n\binom nk\frac{(6k)!}{k!^4(2k)!}
&\cr
\noalign{\hrule}
}}\hss}

\newpage

\hbox to\hsize{\hss\vbox{\offinterlineskip
\halign to120mm{\strut\tabskip=100pt minus 100pt
\strut\vrule#&\hbox to5.5mm{\hss$#$\hss}&%
\vrule#&\hbox to114mm{\hfil$\dsize#$\hfil}&%
\vrule#\tabskip=0pt\cr\noalign{\hrule}
& \# &%
& \text{differential operator $D$ and coefficients $A_n$, $n=0,1,2,\dots$} &\cr
\noalign{\hrule\vskip1pt\hrule}
& \eqnno{127} &&
\aligned
\\[-12pt]
&
D=\theta^4
-z(186640\theta^4+46660\theta^3+8430\theta^2-14900\theta-7811)
\\[-2pt] &\;
+z^2(13063680096\theta^4+6531840048\theta^3+1452582774\theta^2
\\[-2pt] &\;\quad
-105150948\theta+1116311607)
\\[-2pt] &\;
-z^3(406448815694080\theta^4+304836611770560\theta^3
+80522819145360\theta^2
\\[-2pt] &\;\quad
+30399976916340\theta-5858109970482)
\\[-2pt] &\;
+z^4(4744882429422797056\theta^4+4744882429422797056\theta^3
\\[-2pt] &\;\quad
+1452477726672583776\theta^2+134807094643795216\theta
\\[-2pt] &\;\quad
+727318797282433)
\\[-2pt] &\;
-2^7\cdot3^2z^5(65844708142440960\theta^4+82305885178051200\theta^3
\\[-2pt] &\;\quad
+57162742271218800\theta^2+33149511110202480\theta+9595488553244821)
\\[-2pt] &\;
+2^{12}\cdot3^8z^6(16930529404416\theta^4+25395794106624\theta^3
+21045353440896\theta^2
\\[-2pt] &\;\quad
+8640860120064\theta-514138329761)
\\[-2pt] &\;
-2^{19}\cdot3^{14}z^7(483770880\theta^4+846599040\theta^3
+692048880\theta^2
\\[-2pt] &\;\quad
+182253960\theta+22046293)
\\[-2pt] &\;
+2^{24}\cdot3^{22}z^8(4\theta+1)(4\theta+3)(12\theta+5)(12\theta+7)
\\[-8pt]
\endaligned &\cr
& &\span\hrulefill&\cr
& &&
\gathered
\text{same $K(q)$ as in \#82}
\\
\text{the pullback of the 5th-order differential equation $D'y=0$, where}
\\[-8pt]
\endgathered &\cr
& &\span\hrulefill&\cr
& &&
\aligned
\\[-12pt]
&
D'=\theta^5
+z(46676\theta^5+116650\theta^4+110180\theta^3+48620\theta^2+9874\theta+722)
\\[-2pt] &\;
-z^2(93328\theta^5+2799520\theta^4+4069720\theta^3+3965960\theta^2
\\[-2pt] &\;\quad
+2452340\theta+678212)
\\[-2pt] &\;
+z^3(7465600\theta^5+26128320\theta^4+45620640\theta^3+47278800\theta^2
\\[-2pt] &\;\quad
+24206400\theta+3155112)
\\[-2pt] &\;
-2^4z^4(1866320\theta^5+7465120\theta^4+14204360\theta^3+14307920\theta^2
\\[-2pt] &\;\quad
+6634825\theta+1290242)
\\[-2pt] &\;
+2^5z^5(1866272\theta^5+8398160\theta^4+16537040\theta^3+16096360\theta^2
\\[-2pt] &\;\quad
+7704010\theta+1451521)
\\[-2pt] &\;
-4608z^6(\theta+1)(3\theta+2)(3\theta+4)(6\theta+5)(6\theta+7)
\endaligned &\cr
& &\span\hrulefill&\cr
& &&
A_n'=\sum_{k=0}^n\binom{2n-2k}{n-k}\frac{(6k)!}{k!^6}
&\cr
\noalign{\hrule}
& \eqnno{128} &&
\aligned
\\[-12pt]
&
D=\theta^4
-z(3145\theta^4+6282\theta^3+4403\theta^2+1262\theta+122)
\\[-2pt] &\;
+4z^2(2\theta+1)(6270\theta^3+18804\theta^2+20058\theta+7523)
\\[-2pt] &\;
-8z^3(2\theta+1)(2\theta+3)(9395\theta^2+28181\theta+22544)
\\[-2pt] &\;
+16z^4(2\theta+1)(2\theta+3)(2\theta+5)(6260\theta+2519)
\\[-2pt] &\;
-50064z^5(2\theta+1)(2\theta+3)(2\theta+5)(2\theta+7)
\\[-8pt]
\endaligned &\cr
& &\span\hrulefill&\cr
& &&
\gathered
\text{same $K(q)$ as in \#1}
\\[-8pt]
\endgathered &\cr
& &\span\hrulefill&\cr
& &&
A_n=\binom{2n}n\sum_{k=0}^n\binom nk\frac{(5k)!}{k!^3(2k)!}
&\cr
\noalign{\hrule}
}}\hss}

\newpage

\hbox to\hsize{\hss\vbox{\offinterlineskip
\halign to120mm{\strut\tabskip=100pt minus 100pt
\strut\vrule#&\hbox to5.5mm{\hss$#$\hss}&%
\vrule#&\hbox to114mm{\hfil$\dsize#$\hfil}&%
\vrule#\tabskip=0pt\cr\noalign{\hrule}
& \# &%
& \text{differential operator $D$ and coefficients $A_n$, $n=0,1,2,\dots$} &\cr
\noalign{\hrule\vskip1pt\hrule}
& \eqnno{129} &&
\aligned
\\[-12pt]
&
D=\theta^4
-2\cdot3z\theta(4\theta+1)(6\theta^2+1)
+2^2\cdot3(3z)^2\theta(92\theta^3+46\theta^2+33\theta+9)
\\[-2pt] &\;
-2^4(3z)^3\theta(1012\theta^3+759\theta^2+583\theta+171)
\\[-2pt] &\;
+2^4\cdot7(3z)^4\theta(1518\theta^3+1518\theta^2+1243\theta+388)
\\[-2pt] &\;
-2^5\cdot3\cdot7(3z)^5\theta(2024\theta^3+2530\theta^2+2200\theta+725)
\\[-2pt] &\;
+2^5(3z)^6(269200\theta^4+403800\theta^3+371615\theta^2+128340\theta-434)
\\[-2pt] &\;
-2^6(3z)^7(692352\theta^4+1211616\theta^3+1176528\theta^2+422499\theta-6076)
\\[-2pt] &\;
+2^7(3z)^8(1472166\theta^4+2944332\theta^3+3008574\theta^2+1113993\theta-39494)
\\[-2pt] &\;
-2^8(3z)^9(2621536\theta^4+5898456\theta^3+6326976\theta^2+2395311\theta-157976)
\\[-2pt] &\;
+2^9\cdot7(3z)^{10}(563856\theta^4+1409640\theta^3+1583800\theta^2+608575\theta-62062)
\\[-2pt] &\;
-2^{10}\cdot7(3z)^{11}(722976\theta^4+1988184\theta^3+2335176\theta^2+906129\theta-124124)
\\[-2pt] &\;
+2^8(3z)^{12}(44454688\theta^4+133364064\theta^3+163448032\theta^2
\\[-2pt] &\;\quad
+63949536\theta-10415783)
\\[-2pt] &\;
-2^{12}(3z)^{13}(5247024\theta^4+17052828\theta^3+21770956\theta^2
\\[-2pt] &\;\quad
+8611822\theta-1478855)
\\[-2pt] &\;
+2^{12}(3z)^{14}(8546736\theta^4+29913576\theta^3+39718716\theta^2
\\[-2pt] &\;\quad
+15991626\theta-2532173)
\\[-2pt] &\;
-2^{14}(3z)^{15}(3006608\theta^4+11274780\theta^3+15546100\theta^2
\\[-2pt] &\;\quad
+6429885\theta-794437)
\\[-2pt] &\;
+2^{13}(3z)^{16}(7307784\theta^4+29231136\theta^3+41794464\theta^2
\\[-2pt] &\;\quad
+17919876\theta-1365581)
\\[-2pt] &\;
-2^{16}(3z)^{17}(956304\theta^4+4064292\theta^3+6017616\theta^2+2691333\theta-83545)
\\[-2pt] &\;
+2^{17}(3z)^{18}(428800\theta^4+1929600\theta^3+2954715\theta^2+1380690\theta-2429)
\\[-2pt] &\;
-2^{18}(3z)^{19}(163104\theta^4+774744\theta^3+1225444\theta^2+596629\theta+4256)
\\[-2pt] &\;
+2^{16}(3z)^{20}(414336\theta^4+2071680\theta^3+3381120\theta^2+1704960\theta+5957)
\\[-2pt] &\;
-2^{21}(3z)^{21}\theta(6688\theta^3+35112\theta^2+59068\theta+30597)
\\[-2pt] &\;
+2^{23}(3z)^{22}\theta(\theta+1)(672\theta^2+3024\theta+3379)
\\[-2pt] &\;
-2^{27}\cdot3(3z)^{23}\theta(\theta+1)(\theta+2)(4\theta+11)
+2^{28}(3z)^{24}\theta(\theta+1)(\theta+2)(\theta+3)
\\[-8pt]
\endaligned &\cr
& &\span\hrulefill&\cr
& &&
\gathered
\text{the pullback of the 5th-order differential equation $D'y=0$, where}
\\[-8pt]
\endgathered &\cr
& &\span\hrulefill&\cr
& &&
\aligned
\\[-12pt]
&
D'=\theta^5
-6z\theta(\theta^2+\theta+1)(3\theta^2+3\theta+1)
\\[-2pt] &\;
+(6z)^2(\theta+1)(15\theta^4+60\theta^3+105\theta^2+90\theta+31)
\\[-2pt] &\;
-10(6z)^3(\theta+1)(\theta+2)(2\theta+3)(\theta^2+3\theta+3)
\\[-2pt] &\;
+5(6z)^4(\theta+1)(\theta+2)(\theta+3)(3\theta^2+12\theta+13)
\\[-2pt] &\;
-3(6z)^5(\theta+1)(\theta+2)(\theta+3)(\theta+4)(2\theta+5)
\\[-2pt] &\;
+2(6z)^6(\theta+1)(\theta+2)(\theta+3)(\theta+4)(\theta+5)
\endaligned &\cr
& &\span\hrulefill&\cr
& &&
A_n'=\sum_{k=0}^n(-1)^k6^{n-6k}\binom n{6k}\frac{(6k)!}{k!^6}
&\cr
\noalign{\hrule}
}}\hss}

\newpage

\hbox to\hsize{\hss\vbox{\offinterlineskip
\halign to120mm{\strut\tabskip=100pt minus 100pt
\strut\vrule#&\hbox to5.5mm{\hss$#$\hss}&%
\vrule#&\hbox to114mm{\hfil$\dsize#$\hfil}&%
\vrule#\tabskip=0pt\cr\noalign{\hrule}
& \# &%
& \text{differential operator $D$ and coefficients $A_n$, $n=0,1,2,\dots$} &\cr
\noalign{\hrule\vskip1pt\hrule}
& \eqnno{130} &&
\aligned
\\[-12pt]
&
D=\theta^4
-z(224\theta^4+56\theta^3+28\theta^2-3)
\\[-2pt] &\;
+z^2(21952\theta^4+10976\theta^3+5998\theta^2+1020\theta+261)
\\[-2pt] &\;
-z^3(1238528\theta^4+928896\theta^3+551280\theta^2+139800\theta+29316)
\\[-2pt] &\;
+z^4(44574208\theta^4+44574208\theta^3+28575008\theta^2+8415456\theta+1270689)
\\[-2pt] &\;
-2^4z^5(67010048\theta^4+83762560\theta^3+57732992\theta^2+18330336\theta+2276415)
\\[-2pt] &\;
+2^5z^6(549092864\theta^4+823639296\theta^3+607886736\theta^2
\\[-2pt] &\;\quad
+202347456\theta+21851403)
\\[-2pt] &\;
-2^8z^7(767502848\theta^4+1343129984\theta^3+1057767088\theta^2
\\[-2pt] &\;\quad
+363949656\theta+34797501)
\\[-2pt] &\;
+2^8z^8(5778481408\theta^4+11556962816\theta^3+9682036960\theta^2
\\[-2pt] &\;\quad
+3432553440\theta+350349201)
\\[-2pt] &\;
-2^{13}\cdot3^2z^9(98832896\theta^4+222374016\theta^3+197614416\theta^2
\\[-2pt] &\;\quad
+73072472\theta+9376113)
\\[-2pt] &\;
+2^{16}\cdot3^4z^{10}(4204032\theta^4+10510080\theta^3+9877952\theta^2+3874976\theta+592815)
\\[-2pt] &\;
-2^{21}\cdot3^6z^{11}(4\theta+1)(6272\theta^3+15680\theta^2+13172\theta+3783)
\\[-2pt] &\;
+2^{24}\cdot3^8z^{12}(4\theta+1)(4\theta+3)^2(4\theta+5)
\\[-8pt]
\endaligned &\cr
& &\span\hrulefill&\cr
& &&
\gathered
\text{the pullback of the 5th-order differential equation $D'y=0$, where}
\\[-8pt]
\endgathered &\cr
& &\span\hrulefill&\cr
& &&
\aligned
\\[-12pt]
&
D'=\theta^5
-2z(2\theta+1)(14\theta^4+28\theta^3+28\theta^2+14\theta+3)
\\[-2pt] &\;
+4z^2(\theta+1)^3(196\theta^2+392\theta+255)
-1152z^3(\theta+1)^2(\theta+2)^2(2\theta+3)
\endaligned &\cr
& &\span\hrulefill&\cr
& &&
A_n'=\sum_{i+j+k+l+m+s=n}\biggl(\frac{n!}{i!j!k!l!m!s!}\biggr)^2
&\cr
\noalign{\hrule}
}}\hss}

\newpage

\hbox to\hsize{\hss\vbox{\offinterlineskip
\halign to120mm{\strut\tabskip=100pt minus 100pt
\strut\vrule#&\hbox to5.5mm{\hss$#$\hss}&%
\vrule#&\hbox to114mm{\hfil$\dsize#$\hfil}&%
\vrule#\tabskip=0pt\cr\noalign{\hrule}
& \# &%
& \text{differential operator $D$ and coefficients $A_n$, $n=0,1,2,\dots$} &\cr
\noalign{\hrule\vskip1pt\hrule}
& \eqnno{131} &&
\aligned
\\[-12pt]
&
D=\theta^4
-5z\theta(24\theta^3+6\theta^2+4\theta+1)+\dotsb
\\[-2pt] &\;
+3^4\cdot5^{24}z^{24}\theta(\theta+1)(\theta+2)(\theta+4)
\\[-8pt]
\endaligned &\cr
& &\span\hrulefill&\cr
& &&
\gathered
\text{the pullback of the 5th-order differential equation $D'y=0$, where}
\\[-8pt]
\endgathered &\cr
& &\span\hrulefill&\cr
& &&
\aligned
\\[-12pt]
&
D'=\theta^5
-5z(2\theta+1)(\theta^2+\theta+1)(3\theta^2+3\theta+1)
\\[-2pt] &\;
+5^2z^2(\theta+1)(15\theta^4+60\theta^3+105\theta^2+90\theta+31)
\\[-2pt] &\;
-2\cdot5^4z^3(\theta+1)(\theta+2)(2\theta+3)(\theta^2+3\theta+3)
\\[-2pt] &\;
+5^5z^4(\theta+1)(\theta+2)(\theta+3)(3\theta^2+12\theta+12)
\\[-2pt] &\;
-5^5z^5(\theta+1)(\theta+2)(\theta+3)(\theta+4)(2\theta+5)
\\[-2pt] &\;
-3\cdot5^6z^6(\theta+1)(\theta+2)(\theta+3)(\theta+4)(\theta+5)
\endaligned &\cr
& &\span\hrulefill&\cr
& &&
A_n'=\sum_{k=0}^n(-1)^k5^{n-5k}\binom{2k}k\binom n{5k}\frac{(5k)!}{k!^5}
&\cr
\noalign{\hrule}
& \eqnno{132} &&
\aligned
\\[-12pt]
&
D=\theta^4
-5z(80\theta^4+20\theta^3+10\theta^2-1)+\dotsb
\\[-2pt] &\;
+2^{30}\cdot5^{20}z^{16}\theta(\theta+2)(2\theta+1)(2\theta+3)
\\[-8pt]
\endaligned &\cr
& &\span\hrulefill&\cr
& &&
\gathered
\text{the pullback of the 5th-order differential equation $D'y=0$, where}
\\[-8pt]
\endgathered &\cr
& &\span\hrulefill&\cr
& &&
\aligned
\\[-12pt]
&
D'=\theta^5
-10z(2\theta+1)(5\theta^4+10\theta^3+10\theta^2+5\theta+1)
\\[-2pt] &\;
+500z^2(\theta+1)(2\theta+1)(2\theta+3)(2\theta^2+4\theta+3)
\\[-2pt] &\;
-5000z^3(2\theta+1)(2\theta+3)(2\theta+5)(2\theta^2+6\theta+5)
\\[-2pt] &\;
+50000z^4(\theta+2)(2\theta+1)(2\theta+3)(2\theta+5)(2\theta+7)
\endaligned &\cr
& &\span\hrulefill&\cr
& &&
A_n'=\binom{2n}n\sum_{k=0}^n(-1)^k5^{n-5k}\binom n{5k}\frac{(5k)!}{k!^5}
&\cr
\noalign{\hrule}
&\multispan3\hss\text{cases \#\eqnno{133}--\eqnno{143} are described
in \cite{AZ}, Section 7}\hss&\cr
\noalign{\hrule}
}}\hss}

\newpage

\hbox to\hsize{\hss\vbox{\offinterlineskip
\halign to120mm{\strut\tabskip=100pt minus 100pt
\strut\vrule#&\hbox to5.5mm{\hss$#$\hss}&%
\vrule#&\hbox to114mm{\hfil$\dsize#$\hfil}&%
\vrule#\tabskip=0pt\cr\noalign{\hrule}
& \# &%
& \text{differential operator $D$ and coefficients $A_n$, $n=0,1,2,\dots$} &\cr
\noalign{\hrule\vskip1pt\hrule}
& \eqnno{144} &&
\aligned
\\[-12pt]
&
D=\theta^4
-z(73\theta^4+578\theta^3+493\theta^2+204\theta+36)
\\[-2pt] &\;
+2^3\cdot3^2z^2(-145\theta^4+290\theta^3+1387\theta^2+1072\theta+300)
\\[-2pt] &\;
+2^6\cdot3^4z^3(145\theta^4+870\theta^3+353\theta^2-252\theta-180)
\\[-2pt] &\;
-2^9\cdot3^6z^4(-73\theta^4+286\theta^3+803\theta^2+660\theta+180)
-2^{15}\cdot3^{10}z^5(\theta+1)^4
\\[-8pt]
\endaligned &\cr
& &\span\hrulefill&\cr
& &&
\sum_{n=0}^\infty A_nz^n=\text{(g)}*\text{(g)}
\qquad\text{(see \cite{AZ}, Section 7)}
&\cr
\noalign{\hrule}
& \eqnno{145} &&
\aligned
\\[-12pt]
&
D=\theta^4
-3^2z(81\theta^4+648\theta^3+576\theta^2+252\theta+49)
\\[-2pt] &\;
+2\cdot3^8z^2(-81\theta^4+162\theta^3+765\theta^2+612\theta+187)
\\[-2pt] &\;
+2\cdot3^{14}z^3(81\theta^4+486\theta^3+207\theta^2-108\theta-97)
\\[-2pt] &\;
+3^{20}z^4(81\theta^4-324\theta^3-882\theta^2-720\theta-194)
-3^{30}z^5(\theta+1)^4
\\[-8pt]
\endaligned &\cr
& &\span\hrulefill&\cr
& &&
\sum_{n=0}^\infty A_nz^n=\text{(h)}*\text{(h)}
\qquad\text{(see \cite{AZ}, Section 7)}
&\cr
\noalign{\hrule}
& \eqnno{146} &&
\aligned
\\[-12pt]
&
D=\theta^4
-4z(5\theta^4+10\theta^3+10\theta^2+5\theta+1)
\\[-2pt] &\;
+80z^2(\theta+1)^2(2\theta^2+4\theta+3)
-320z^3(\theta+1)(\theta+2)(2\theta^2+6\theta+5)
\\[-2pt] &\;
+2304z^4(\theta+1)(\theta+2)^2(\theta+3)
-5120z^5(\theta+1)(\theta+2)(\theta+3)(\theta+4)
\\[-8pt]
\endaligned &\cr
& &\span\hrulefill&\cr
& &&
A_n=\sum_{k=0}^n(-1)^k4^{n-4k}\binom{2k}k\binom n{4k}\frac{(4k)!}{k!^4}
&\cr
\noalign{\hrule}
& \eqnno{147} &&
\aligned
\\[-12pt]
&
D=\theta^4
-8z(2\theta+1)^2(2\theta^2+2\theta+1)
\\[-2pt] &\;
+64z^2(2\theta+1)(2\theta+3)(6\theta^2+12\theta+7)
\\[-2pt] &\;
-1024z^3(2\theta+1)(2\theta+3)^2(2\theta+5)
\\[-2pt] &\;
+8192z^4(2\theta+1)(2\theta+3)(2\theta+5)(2\theta+7)
\\[-8pt]
\endaligned &\cr
& &\span\hrulefill&\cr
& &&
A_n=\binom{2n}n\sum_{k=0}^n(-1)^k4^{n-4k}\binom n{4k}\frac{(4k)!}{k!^4}
&\cr
\noalign{\hrule}
& \eqnno{148} &&
\aligned
\\[-12pt]
&
D=\theta^4
-5z(5\theta^4+10\theta^3+10\theta^2+5\theta+1)
+125z^2(\theta+1)^2(2\theta^2+4\theta+3)
\\[-2pt] &\;
-625z^3(\theta+1)(\theta+2)(2\theta^2+6\theta+5)
-3125z^4(\theta+1)(\theta+2)^2(\theta+3)
\\[-8pt]
\endaligned &\cr
& &\span\hrulefill&\cr
& &&
A_n=\sum_{k=0}^n(-1)^k5^{n-5k}\binom n{5k}\frac{(5k)!}{k!^5}
&\cr
\noalign{\hrule}
}}\hss}

\newpage

\hbox to\hsize{\hss\vbox{\offinterlineskip
\halign to120mm{\strut\tabskip=100pt minus 100pt
\strut\vrule#&\hbox to5.5mm{\hss$#$\hss}&%
\vrule#&\hbox to114mm{\hfil$\dsize#$\hfil}&%
\vrule#\tabskip=0pt\cr\noalign{\hrule}
& \# &%
& \text{differential operator $D$ and coefficients $A_n$, $n=0,1,2,\dots$} &\cr
\noalign{\hrule\vskip1pt\hrule}
& \eqnno{149} &&
\aligned
\\[-12pt]
&
D=\theta^4
-2^23z(4896\theta^4+1224\theta^3+265\theta^2-347\theta-186)
\\[-2pt] &\;
+2^63^4z^2(249840\theta^4+124920\theta^3+32255\theta^2+1619\theta+21824)
\\[-2pt] &\;
-2^{10}3^7z^3(5674464\theta^4+4255848\theta^3
+1276509\theta^2+519897\theta-79028)
\\[-2pt] &\;
+2^{14}3^{10}z^4(48607704\theta^4+486077040\theta^3+16615003\theta^2
\\[-2pt] &\;\quad
+2178895\theta+94752)
\\[-2pt] &\;
-2^{18}3^{13}z^5(5674464\theta^4+7093080\theta^3+2892155\theta^2
-369553\theta-546818)
\\[-2pt] &\;
+2^{22}3^{16}z^6(249840\theta^4+374760\theta^3+174525\theta^2
+1689\theta+14384)
\\[-2pt] &\;
-2^{26}3^{19}z^7\theta(4896\theta^3+8568\theta^3+4447\theta+595)
\\[-2pt] &\;
+2^{30}3^{22}z^8\theta(\theta+1)(6\theta+1)(6\theta+5)
\\[-8pt]
\endaligned &\cr
& &\span\hrulefill&\cr
& &&
\gathered
\text{the pullback of the 5th-order differential equation $D'y=0$, where}
\\[-8pt]
\endgathered &\cr
& &\span\hrulefill&\cr
& &&
\aligned
\\[-12pt]
&
D'=\theta^5
-12z(2\theta+1)(6\theta+1)(6\theta+5)(17\theta^2+17\theta+5)
\\[-2pt] &\;
+144z^2(\theta+1)(6\theta+1)(6\theta+5)(6\theta+7)(6\theta+11)
\endaligned &\cr
& &\span\hrulefill&\cr
& &&
A_n'=\frac{(6n)!}{n!(2n)!(3n)!}\sum_{k=0}^n\binom nk^2\binom{n+k}k^2
&\cr
\noalign{\hrule}
& \eqnno{150} &&
\aligned
\\[-12pt]
&
D=\theta^4
-2^2z(224\theta^4+56\theta^3+19\theta^2-9\theta-6)
\\[-2pt] &\;
+2^6z^2(6000\theta^4+3000\theta^3+1139\theta^2-537\theta-420)
\\[-2pt] &\;
-2^{10}z^3(98336\theta^4+73752\theta^3+31007\theta^2-5469\theta+1524)
\\[-2pt] &\;
+2^{14}z^4(1073176\theta^4+1073176\theta^3+495055\theta^2+34851\theta+112248)
\\[-2pt] &\;
-2^{18}3^4z^5(98336\theta^4+122920\theta^3+61737\theta^2+13173\theta+6210)
\\[-2pt] &\;
+2^{18}3^8z^6(96000\theta^4+144000\theta^3+78224\theta^2+19792\theta-576)
\\[-2pt] &\;
-2^{26}3^{12}z^7\theta(224\theta^3+392\theta^2+229\theta+49)
+2^{30}3^{16}z^8\theta(\theta+1)(2\theta+1)^2
\\[-8pt]
\endaligned &\cr
& &\span\hrulefill&\cr
& &&
\gathered
\text{the pullback of the 5th-order differential equation $D'y=0$, where}
\\[-8pt]
\endgathered &\cr
& &\span\hrulefill&\cr
& &&
\aligned
\\[-12pt]
&
D'=\theta^5
-4z(2\theta+1)^3(7\theta^2+7\theta+3)
+1296z^2(\theta+1)(2\theta+1)^2(2\theta+3)^2
\endaligned &\cr
& &\span\hrulefill&\cr
& &&
A_n'=\binom{2n}n^2\sum_{k=0}^n(-1)^k3^{n-3k}\binom n{3k}
\binom{n+k}k\frac{(3k)!}{k!^3}
&\cr
\noalign{\hrule}
}}\hss}

\newpage

\hbox to\hsize{\hss\vbox{\offinterlineskip
\halign to120mm{\strut\tabskip=100pt minus 100pt
\strut\vrule#&\hbox to5.5mm{\hss$#$\hss}&%
\vrule#&\hbox to114mm{\hfil$\dsize#$\hfil}&%
\vrule#\tabskip=0pt\cr\noalign{\hrule}
& \# &%
& \text{differential operator $D$ and coefficients $A_n$, $n=0,1,2,\dots$} &\cr
\noalign{\hrule\vskip1pt\hrule}
& \eqnno{151} &&
\aligned
\\[-12pt]
&
D=\theta^4
-3z(504\theta^4+126\theta^3+41\theta^2-22\theta-14)
\\[-2pt] &\;
+3^2z^2(121500\theta^4+60750\theta^3+22221\theta^2-11790\theta-8762)
\\[-2pt] &\;
-3^5z^3(1991304\theta^4+1493478\theta^3+607149\theta^2-133506\theta+26716)
\\[-2pt] &\;
+3^8z^4(21731814\theta^4+21731814\theta^3+9723033\theta^2+378126\theta+2244916)
\\[-2pt] &\;
-3^{15}z^5(1991304\theta^4+2489130\theta^3+1215603\theta^2+230166\theta+125122)
\\[-2pt] &\;
+3^{22}z^6(121500\theta^4+182250\theta^3+96471\theta^2+22446\theta-698)
\\[-2pt] &\;
-3^{31}z^7\theta(504\theta^3+882\theta^2+503\theta+98)
+3^{38}z^8\theta(\theta+1)(3\theta+1)(3\theta+2)
\\[-8pt]
\endaligned &\cr
& &\span\hrulefill&\cr
& &&
\gathered
\text{the pullback of the 5th-order differential equation $D'y=0$, where}
\\[-8pt]
\endgathered &\cr
& &\span\hrulefill&\cr
& &&
\aligned
\\[-12pt]
&
D'=\theta^5
-3z(2\theta+1)(3\theta+1)(3\theta+2)(7\theta^2+7\theta+3)
\\[-2pt] &\;
+729z^2(\theta+1)(3\theta+1)(3\theta+2)(3\theta+4)(3\theta+5)
\endaligned &\cr
& &\span\hrulefill&\cr
& &&
A_n'=\frac{(3n)!}{n!^3}\sum_{k=0}^n(-1)^k3^{n-3k}\binom n{3k}
\binom{n+k}k\frac{(3k)!}{k!^3}
&\cr
\noalign{\hrule}
& \eqnno{152} &&
\aligned
\\[-12pt]
&
D=\theta^4
-2^2z(896\theta^4+224\theta^3+69\theta^2-43\theta-26)
\\[-2pt] &\;
+2^7z^2(48000\theta^4+24000\theta^3+8362\theta^2-5110\theta-3589)
\\[-2pt] &\;
-2^{14}z^3(393344\theta^4+295008\theta^3+114809\theta^2-31991\theta+4246)
\\[-2pt] &\;
+2^{20}z^4(4292704\theta^4+4292704\theta^3+1846073\theta^2-6199\theta+436401)
\\[-2pt] &\;
-2^{26}3^4z^5(393344\theta^4+491680\theta^3+231583\theta^2+36431\theta+24552)
\\[-2pt] &\;
+2^{31}3^8z^6(48000\theta^4+72000\theta^3+36862\theta^2+7582\theta-261)
\\[-2pt] &\;
-2^{38}3^{12}z^7\theta(896\theta^3+1568\theta^2+867\theta+147)
\\[-2pt] &\;
+2^{44}3^{16}z^8\theta(\theta+1)(4\theta+1)(4\theta+3)
\\[-8pt]
\endaligned &\cr
& &\span\hrulefill&\cr
& &&
\gathered
\text{the pullback of the 5th-order differential equation $D'y=0$, where}
\\[-8pt]
\endgathered &\cr
& &\span\hrulefill&\cr
& &&
\aligned
\\[-12pt]
&
D'=\theta^5
-4z(2\theta+1)(4\theta+1)(4\theta+3)(7\theta^2+7\theta+3)
\\[-2pt] &\;
+1296z^2(\theta+1)(4\theta+1)(4\theta+3)(4\theta+5)(4\theta+7)
\endaligned &\cr
& &\span\hrulefill&\cr
& &&
A_n'=\frac{(4n)!}{n!^2(2n)!}\sum_{k=0}^n(-1)^k3^{n-3k}\binom n{3k}
\binom{n+k}k\frac{(3k)!}{k!^3}
&\cr
\noalign{\hrule}
}}\hss}

\newpage

\hbox to\hsize{\hss\vbox{\offinterlineskip
\halign to120mm{\strut\tabskip=100pt minus 100pt
\strut\vrule#&\hbox to5.5mm{\hss$#$\hss}&%
\vrule#&\hbox to114mm{\hfil$\dsize#$\hfil}&%
\vrule#\tabskip=0pt\cr\noalign{\hrule}
& \# &%
& \text{differential operator $D$ and coefficients $A_n$, $n=0,1,2,\dots$} &\cr
\noalign{\hrule\vskip1pt\hrule}
& \eqnno{153} &&
\aligned
\\[-12pt]
&
D=\theta^4
-2^2\cdot3z(2016\theta^4+504\theta^3+143\theta^2-109\theta-62)
\\[-2pt] &\;
+2^63^2z^2(486000\theta^4+243000\theta^3+78759\theta^2-58149\theta-38156)
\\[-2pt] &\;
-2^{10}3^5z^3(7965216\theta^4+5973912\theta^3+2179683\theta^2
-807129\theta+56452)
\\[-2pt] &\;
+2^{14}3^8z^4(86927256\theta^4+86927256\theta^3+35270163\theta^2
\\[-2pt] &\;\quad
-2418777\theta+8633800)
\\[-2pt] &\;
-2^{18}3^{15}z^5(7965216\theta^4+9956520\theta^3+4447557\theta^2
+481617\theta+492250)
\\[-2pt] &\;
+2^{22}3^{21}z^6(1458000\theta^4+2187000\theta^3+1066527\theta^2
+175635\theta-7332)
\\[-2pt] &\;
-2^{26}3^{31}z^7\theta(2016\theta^3+3528\theta^2+1865\theta+245)
\\[-2pt] &\;
+2^{30}3^{38}z^8\theta(\theta+1)(6\theta+1)(6\theta+5)
\\[-8pt]
\endaligned &\cr
& &\span\hrulefill&\cr
& &&
\gathered
\text{the pullback of the 5th-order differential equation $D'y=0$, where}
\\[-8pt]
\endgathered &\cr
& &\span\hrulefill&\cr
& &&
\aligned
\\[-12pt]
&
D'=\theta^5
-12z(2\theta+1)(6\theta+1)(6\theta+5)(7\theta^2+7\theta+3)
\\[-2pt] &\;
+11664z^2(\theta+1)(6\theta+1)(6\theta+5)(6\theta+7)(6\theta+11)
\endaligned &\cr
& &\span\hrulefill&\cr
& &&
A_n'=\frac{(6n)!}{n!(2n)!(3n)!}\sum_{k=0}^n(-1)^k3^{n-3k}
\binom n{3k}\binom{n+k}k\frac{(3k)!}{k!^3}
&\cr
\noalign{\hrule}
& \eqnno{154} &&
\aligned
\\[-12pt]
&
D=\theta^4
+2^23z(36\theta^4-144\theta^3-99\theta^2-27\theta-3)
\\[-2pt] &\;
+2^43^5z^2(312\theta^2+156\theta+25)
\\[-2pt] &\;
+2^{10}3^8z^3(12\theta^4+72\theta^3+57\theta^2+18\theta+2)
+2^{12}3^{12}z^4(2\theta+1)^4
\\[-8pt]
\endaligned &\cr
& &\span\hrulefill&\cr
& &&
A_n=\binom{2n}n^2\sum_{k=0}^n(-1)^k3^{2n-3k}\binom{2n}{3k}\frac{(3k)!}{k!^3}
&\cr
\noalign{\hrule}
& \eqnno{155} &&
\aligned
\\[-12pt]
&
D=\theta^4
-2^4z(256\theta^4+2048\theta^3+1856\theta^2+832\theta+169)
\\[-2pt] &\;
-2^{16}z^2(512\theta^4-1024\theta^3-4800\theta^2-3904\theta-1239)
\\[-2pt] &\;
+2^{28}z^3(512\theta^4+3072\theta^3+1344\theta^2-576\theta-599)
\\[-2pt] &\;
+2^{40}z^4(256\theta^4-1024\theta^3-2752\theta^2-2240\theta-599)
-2^{60}z^5(\theta+1)^4
\\[-8pt]
\endaligned &\cr
& &\span\hrulefill&\cr
& &&
A_n=\biggl\{64^n\sum_k(-1)^k\binom{-3/4}k\binom{-1/4}{n-k}^2\biggr\}^2
&\cr
\noalign{\hrule}
& \eqnno{156} &&
\aligned
\\[-12pt]
&
D=\theta^4
-6z(2\theta+1)^2(2\theta^2+2\theta+1)
\\[-2pt] &\;
+36z^2(2\theta+1)(2\theta+3)(6\theta^2+12\theta+7)
\\[-2pt] &\;
-3888z^4(2\theta+1)(2\theta+3)(2\theta+5)(2\theta+7)
\\[-8pt]
\endaligned &\cr
& &\span\hrulefill&\cr
& &&
A_n=\binom{2n}n\sum_{k=0}^n(-1)^k3^{n-3k}\binom n{3k}
\binom{2k}k\frac{(3k)!}{k!^3}
&\cr
\noalign{\hrule}
& \eqnno{157} &&
\aligned
\\[-12pt]
&
D=\theta^4
-3z(5\theta^4+10\theta^3+10\theta^2+3\theta+1)
+45z^2(\theta+1)^2(2\theta^2+4\theta+3)
\\[-2pt] &\;
+27z^3(\theta+1)(\theta+2)(6\theta^2+18\theta+11)
-2187z^4(\theta+1)(\theta+2)^2(\theta+3)
\\[-2pt] &\;
+3645z^5(\theta+1)(\theta+2)(\theta+3)(\theta+4)
\\[-8pt]
\endaligned &\cr
& &\span\hrulefill&\cr
& &&
A_n=\sum_{k=0}^n(-1)^k3^{n-3k}\binom n{3k}\binom{2k}k^2\frac{(3k)!}{k!^3}
&\cr
\noalign{\hrule}
}}\hss}

\newpage

\hbox to\hsize{\hss\vbox{\offinterlineskip
\halign to120mm{\strut\tabskip=100pt minus 100pt
\strut\vrule#&\hbox to5.5mm{\hss$#$\hss}&%
\vrule#&\hbox to114mm{\hfil$\dsize#$\hfil}&%
\vrule#\tabskip=0pt\cr\noalign{\hrule}
& \# &%
& \text{differential operator $D$ and coefficients $A_n$, $n=0,1,2,\dots$} &\cr
\noalign{\hrule\vskip1pt\hrule}
& \eqnno{158} &&
\aligned
\\[-12pt]
&
D=\theta^4
-2\cdot 3z(12\theta^4+16\theta^3+14\theta^2+6\theta+1)
\\[-2pt] &\;
+2^23^2z^2(60\theta^4+160\theta^3+198\theta^2+116\theta+25)
\\[-2pt] &\;
-2^23^3z^3(2\theta+1)(106\theta^3+263\theta^2+231\theta+86)
\\[-2pt] &\;
-2^43^4z^4(2\theta+1)(2\theta+3)(21\theta^2+178\theta+249)
\\[-2pt] &\;
+2^43^5z^5(2\theta+1)(2\theta+3)(2\theta+5)(114\theta+355)
\\[-2pt] &\;
-2^63^623z^6(2\theta+1)(2\theta+3)(2\theta+5)(2\theta+7)
\\[-8pt]
\endaligned &\cr
& &\span\hrulefill&\cr
& &&
A_n=\binom{2n}n\sum_{k=0}^n(-1)^k3^{n-3k}\binom n{3k}
\binom{3k}k\frac{(3k)!}{k!^3}
&\cr
\noalign{\hrule}
& \eqnno{159} &&
\aligned
\\[-12pt]
&
D=\theta^4
-2\cdot 5z(14\theta^4+16\theta^3+14\theta^2+6\theta+1)
\\[-2pt] &\;
+2^25^2z^2(84\theta^4+192\theta^3+226\theta^2+128\theta+27)
\\[-2pt] &\;
-2^45^3z^3(140\theta^4+480\theta^3+710\theta^2+480\theta+114)
\\[-2pt] &\;
+2^45^5z^4(2\theta+1)(56\theta^3+228\theta^2+342\theta+181)
\\[-2pt] &\;
-2^45^5z^5(2\theta+1)(2\theta+3)(332\theta^2+1216\theta+1189)
\\[-2pt] &\;
+2^75^6z^6(2\theta+1)(2\theta+3)(2\theta+5)(26\theta+57)
\\[-2pt] &\;
-2^83\cdot 5^7z^7(2\theta+1)(2\theta+3)(2\theta+5)(2\theta+7)
\\[-8pt]
\endaligned &\cr
& &\span\hrulefill&\cr
& &&
A_n=\binom{2n}n\sum_{k=0}^n(-1)^k5^{n-5k}\binom n{5k}
\binom{2k}k^{-1}\frac{(5k)!}{k!^5}
&\cr
\noalign{\hrule}
& \eqnno{160} &&
\aligned
\\[-12pt]
&
D=\theta^4
-3z(3\theta^2+3\theta+1)(7\theta^2+7\theta+3)
\\[-2pt] &\;
+3^2z^2(171\theta^4+396\theta^3+555\theta^2+318\theta+64)
\\[-2pt] &\;
+2^33^4z^3(-21\theta^4+126\theta^3+386\theta^2+291\theta+76)
\\[-2pt] &\;
+2^43^5z^4(147\theta^4+294\theta^3+102\theta^2-45\theta-14)
\\[-2pt] &\;
+2^63^7z^5(21\theta^4+210\theta^3+118\theta^2-19\theta-24)
\\[-2pt] &\;
+2^63^8z^6(171\theta^4+288\theta^3+393\theta^2+288\theta+76)
\\[-2pt] &\;
+2^93^{10}z^7(3\theta^2+3\theta+1)(7\theta^2+7\theta+3)
+2^{12}3^{12}z^8(\theta+1)^4
\\[-8pt]
\endaligned &\cr
& &\span\hrulefill&\cr
& &&
A_n=\sum_{k=0}^n\binom nk^3\sum_{k=0}^n(-1)^k3^{n-3k}
\binom n{3k}\frac{(3k)!}{k!^3}
&\cr
\noalign{\hrule}
& \eqnno{161} &&
\aligned
\\[-12pt]
&
D=\theta^4
-3z(3\theta^2+3\theta+1)(11\theta^2+11\theta+3)
\\[-2pt] &\;
+3^2z^2(366\theta^4+1428\theta^3+1980\theta^2+1104\theta+221)
\\[-2pt] &\;
+3^4z^3(-33\theta^4+198\theta^3+607\theta^2+456\theta+117)
\\[-2pt] &\;
+3^5z^4(726\theta^4+1452\theta^3-978\theta^2-1704\theta-515)
\\[-2pt] &\;
+3^7z^5(33\theta^4+330\theta^3+185\theta^2-32\theta-37)
\\[-2pt] &\;
+3^8z^6(366\theta^4+36\theta^3-108\theta^2+36\theta+35)
\\[-2pt] &\;
+3^{10}z^7(3\theta^2+3\theta+1)(11\theta^2+11\theta+3)
+3^{12}z^8(\theta+1)^4
\\[-8pt]
\endaligned &\cr
& &\span\hrulefill&\cr
& &&
A_n=\sum_{k=0}^n\binom nk^2\binom{n+k}k\sum_{k=0}^n(-1)^k3^{n-3k}
\binom n{3k}\frac{(3k)!}{k!^3}
&\cr
\noalign{\hrule}
}}\hss}

\newpage

\hbox to\hsize{\hss\vbox{\offinterlineskip
\halign to120mm{\strut\tabskip=100pt minus 100pt
\strut\vrule#&\hbox to5.5mm{\hss$#$\hss}&%
\vrule#&\hbox to114mm{\hfil$\dsize#$\hfil}&%
\vrule#\tabskip=0pt\cr\noalign{\hrule}
& \# &%
& \text{differential operator $D$ and coefficients $A_n$, $n=0,1,2,\dots$} &\cr
\noalign{\hrule\vskip1pt\hrule}
& \eqnno{162} &&
\aligned
\\[-12pt]
&
D=\theta^4
-3z(3\theta^2+3\theta+1)(10\theta^2+10\theta+3)
\\[-2pt] &\;
+3^3z^2(91\theta^4+472\theta^3+659\theta^2+374\theta+81)
\\[-2pt] &\;
+3^6z^3(30\theta^4-180\theta^3-551\theta^2-417\theta-111)
\\[-2pt] &\;
+3^8z^4(-200\theta^4-400\theta^3+514\theta^2+714\theta+237)
\\[-2pt] &\;
+3^{11}z^5(30\theta^4+300\theta^3+169\theta^2-25\theta-35)
\\[-2pt] &\;
+3^{13}z^6(91\theta^4-108\theta^3-211\theta^2-108\theta-15)
\\[-2pt] &\;
-3^{16}z^7(3\theta^2+3\theta+1)(10\theta^2+10\theta+3)
+3^{20}z^8(\theta+1)^4
\\[-8pt]
\endaligned &\cr
& &\span\hrulefill&\cr
& &&
A_n=\sum_{k=0}^n\binom nk^2\binom{2k}k\sum_{k=0}^n(-1)^k3^{n-3k}
\binom n{3k}\frac{(3k)!}{k!^3}
&\cr
\noalign{\hrule}
& \eqnno{163} &&
\aligned
\\[-12pt]
&
D=\theta^4
-2^23z(3\theta^2+3\theta+1)^2
\\[-2pt] &\;
+2^43^2z^2(21\theta^4+156\theta^3+219\theta^2+126\theta+29)
\\[-2pt] &\;
+2^73^4z^3(3\theta^2+3\theta+1)(3\theta^2-21\theta-35)
\\[-2pt] &\;
+2^{10}3^5z^4(-27\theta^4-54\theta^3+114\theta^2+141\theta+49)
\\[-2pt] &\;
+2^{12}3^7z^5(3\theta^2+3\theta+1)(3\theta^2+27\theta-11)
\\[-2pt] &\;
+2^{14}3^8z^6(21\theta^4-72\theta^3-123\theta^2-72\theta-13)
\\[-2pt] &\;
-2^{17}3^{10}z^7(3\theta^2+3\theta+1)^2
+2^{20}3^{12}z^8(\theta+1)^4
\\[-8pt]
\endaligned &\cr
& &\span\hrulefill&\cr
& &&
A_n=\sum_{k=0}^n\binom nk\binom{2k}k\binom{2n-2k}{n-k}
\cdot\sum_{k=0}^n(-1)^k3^{n-3k}\binom n{3k}\frac{(3k)!}{k!^3}
&\cr
\noalign{\hrule}
& \eqnno{164} &&
\aligned
\\[-12pt]
&
D=\theta^4
-2^23z(3\theta^2+3\theta+1)(8\theta^2+8\theta+3)
\\[-2pt] &\;
+2^43^2z^2(144\theta^4+1152\theta^3+1632\theta^2+960\theta+235)
\\[-2pt] &\;
+2^{10}3^4z^3(24\theta^4-144\theta^3-439\theta^2-339\theta-99)
\\[-2pt] &\;
+2^{13}3^5z^4(-192\theta^4-384\theta^3+864\theta^2+1056\theta+389)
\\[-2pt] &\;
+2^{18}3^7z^5(24\theta^4+240\theta^3+137\theta^2-11\theta-31)
\\[-2pt] &\;
+2^{20}3^8z^6(144\theta^4-576\theta^3-960\theta^2-576\theta-101)
\\[-2pt] &\;
-2^{26}3^{10}z^7(3\theta^2+3\theta+1)(8\theta^2+8\theta+3)
+2^{32}3^{12}z^8(\theta+1)^4
\\[-8pt]
\endaligned &\cr
& &\span\hrulefill&\cr
& &&
A_n=\sum_{k=0}^n4^{n-k}\binom{2k}k^2\binom{2n-2k}{n-k}
\cdot\sum_{k=0}^n(-1)^k3^{n-3k}\binom n{3k}\frac{(3k)!}{k!^3}
&\cr
\noalign{\hrule}
& \eqnno{165} &&
\aligned
\\[-12pt]
&
D=\theta^4
-3^2z(3\theta+1)(6\theta^2+3\theta+1)
\\[-2pt] &\;
+3^5z^2(-3\theta^4+6\theta^3+39\theta^2+30\theta+8)
\\[-2pt] &\;
+3^8z^3(3\theta^4+18\theta^3-3\theta^2-18\theta-8)
\\[-2pt] &\;
-3^{11}z^4(3\theta+2)(6\theta^2+9\theta+4)
-3^{15}z^5(\theta+1)^4
\\[-8pt]
\endaligned &\cr
& &\span\hrulefill&\cr
& &&
A_n=\Biggl\{\sum_{k=0}^n(-1)^k3^{n-3k}\binom n{3k}
\frac{(3k)!}{k!^3}\Biggr\}^2
&\cr
\noalign{\hrule}
}}\hss}

\newpage

\hbox to\hsize{\hss\vbox{\offinterlineskip
\halign to120mm{\strut\tabskip=100pt minus 100pt
\strut\vrule#&\hbox to5.5mm{\hss$#$\hss}&%
\vrule#&\hbox to114mm{\hfil$\dsize#$\hfil}&%
\vrule#\tabskip=0pt\cr\noalign{\hrule}
& \# &%
& \text{differential operator $D$ and coefficients $A_n$, $n=0,1,2,\dots$} &\cr
\noalign{\hrule\vskip1pt\hrule}
& \eqnno{166} &&
\aligned
\\[-12pt]
&
D=\theta^4
-2^43^2z(1296\theta^4+10368\theta^3+9648\theta^2+4464\theta+961)
\\[-2pt] &\;
-2^{12}3^8z^2(2592\theta^4-5184\theta^3-24048\theta^2-20016\theta-6671)
\\[-2pt] &\;
+2^{20}3^{14}z^3(2592\theta^4+15552\theta^3+7056\theta^2-2160\theta-2927)
\\[-2pt] &\;
+2^{28}3^{20}z^4(1296\theta^4-5184\theta^3-13680\theta^2-11088\theta-2927)
\\[-2pt] &\;
-2^{40}3^{30}z^5(\theta+1)^4
\\[-8pt]
\endaligned &\cr
& &\span\hrulefill&\cr
& &&
A_n=\biggl\{432^n\sum_k(-1)^k\binom{-5/6}k\binom{-1/6}{n-k}^2\biggr\}^2
&\cr
\noalign{\hrule}
& \eqnno{167} &&
\aligned
\\[-12pt]
&
D=\theta^4
-3z(7\theta^2+7\theta+2)(18\theta^2+18\theta+7)
\\[-2pt] &\;
+3^2z^2(3969\theta^4+8100\theta^3+11025\theta^2+5850\theta+832)
\\[-2pt] &\;
+2^33^7z^3(-126\theta^4+756\theta^3+2309\theta^2+1767\theta+496)
\\[-2pt] &\;
+2^43^8z^4(3321\theta^4+6642\theta^3+3330\theta^2+9\theta+526)
\\[-2pt] &\;
+2^63^{13}z^5(126\theta^4+1260\theta^3+715\theta^2-79\theta-156)
\\[-2pt] &\;
+2^63^{14}z^6(3969\theta^4+7776\theta^3+10593\theta^2+7776\theta+1876)
\\[-2pt] &\;
+2^93^{19}z^7(7\theta^2+7\theta+2)(18\theta^2+18\theta+7)
+2^{12}3^{24}z^8(\theta+1)^4
\\[-8pt]
\endaligned &\cr
& &\span\hrulefill&\cr
& &&
\sum_{n=0}^\infty A_nz^n=\text{(a)}*\text{(h)}
\qquad\text{(see \cite{AZ}, Section 7)}
&\cr
\noalign{\hrule}
& \eqnno{168} &&
\aligned
\\[-12pt]
&
D=\theta^4
-3z(11\theta^2+11\theta+3)(18\theta^2+18\theta+7)
\\[-2pt] &\;
+3^2z^2(9801\theta^4+38232\theta^3+52965\theta^2+29466\theta+5855)
\\[-2pt] &\;
+3^7z^3(-198\theta^4+1188\theta^3+3631\theta^2+2769\theta+765)
\\[-2pt] &\;
+3^8z^4(19440\theta^4+38880\theta^3-26262\theta^2-45702\theta-13679)
\\[-2pt] &\;
+3^{13}z^5(198\theta^4+1980\theta^3+1121\theta^2-137\theta-241)
\\[-2pt] &\;
+3^{14}z^6(9801\theta^4+972\theta^3-2925\theta^2+972\theta+923)
\\[-2pt] &\;
+3^{19}z^7(11\theta^2+11\theta+3)(18\theta^2+18\theta+7)
+3^{24}z^8(\theta+1)^4
\\[-8pt]
\endaligned &\cr
& &\span\hrulefill&\cr
& &&
\sum_{n=0}^\infty A_nz^n=\text{(b)}*\text{(h)}
\qquad\text{(see \cite{AZ}, Section 7)}
&\cr
\noalign{\hrule}
}}\hss}

\newpage

\hbox to\hsize{\hss\vbox{\offinterlineskip
\halign to120mm{\strut\tabskip=100pt minus 100pt
\strut\vrule#&\hbox to5.5mm{\hss$#$\hss}&%
\vrule#&\hbox to114mm{\hfil$\dsize#$\hfil}&%
\vrule#\tabskip=0pt\cr\noalign{\hrule}
& \# &%
& \text{differential operator $D$ and coefficients $A_n$, $n=0,1,2,\dots$} &\cr
\noalign{\hrule\vskip1pt\hrule}
& \eqnno{169} &&
\aligned
\\[-12pt]
&
D=\theta^4
-3z(10\theta^2+10\theta+3)(18\theta^2+18\theta+7)
\\[-2pt] &\;
+3^4z^2(900\theta^2+4592\theta^3+6426\theta^2+3708\theta+841)
\\[-2pt] &\;
+3^9z^3(180\theta^4-1080\theta^3-3296\theta^2-2532\theta-723)
\\[-2pt] &\;
+3^{12}z^4(-1962\theta^4-3924\theta^3+4770\theta^2+6732\theta+2359)
\\[-2pt] &\;
+3^{17}z^5(180\theta^4+1800\theta^3+1024\theta^2-100\theta-227)
\\[-2pt] &\;
+3^{20}z^6(900\theta^4-972\theta^3-1890\theta^2-972\theta-113)
\\[-2pt] &\;
-3^{25}z^7(10\theta^2+10\theta+3)(18\theta^2+18\theta+7)
+3^{32}z^8(\theta+1)^4
\\[-8pt]
\endaligned &\cr
& &\span\hrulefill&\cr
& &&
\sum_{n=0}^\infty A_nz^n=\text{(c)}*\text{(h)}
\qquad\text{(see \cite{AZ}, Section 7)}
&\cr
\noalign{\hrule}
& \eqnno{170} &&
\aligned
\\[-12pt]
&
D=\theta^4
-2^23z(3\theta^2+3\theta+1)(18\theta^2+18\theta+7)
\\[-2pt] &\;
+2^43^2z^2(729\theta^4+4860\theta^3+6903\theta^2+4086\theta+1007)
\\[-2pt] &\;
+2^73^7z^3(54\theta^4-324\theta^3-987\theta^2-765\theta-227)
\\[-2pt] &\;
+2^{10}3^8z^4(-891\theta^4-1782\theta^3+3222\theta^2+4113\theta+1549)
\\[-2pt] &\;
+2^{12}3^{13}z^5(54\theta^4+540\theta^3+309\theta^2-21\theta-71)
\\[-2pt] &\;
+2^{14}3^{14}z^6(729\theta^4-1944\theta^3-3303\theta^2-1944\theta-307)
\\[-2pt] &\;
-2^{17}3^{19}z^7(3\theta^2+3\theta+1)(18\theta^2+18\theta+7)
+2^{20}3^{24}z^8(\theta+1)^4
\\[-8pt]
\endaligned &\cr
& &\span\hrulefill&\cr
& &&
\sum_{n=0}^\infty A_nz^n=\text{(d)}*\text{(h)}
\qquad\text{(see \cite{AZ}, Section 7)}
&\cr
\noalign{\hrule}
& \eqnno{171} &&
\aligned
\\[-12pt]
&
D=\theta^4
-2^23z(72\theta^4+288\theta^3+254\theta^2+110\theta+21)
\\[-2pt] &\;
+2^43^2z^2(-1296\theta^4+15552\theta^3+33696\theta^2+23472\theta+6625)
\\[-2pt] &\;
+2^93^5z^3(2592\theta^4+5184\theta^3-9000\theta^2-11592\theta-4499)
\\[-2pt] &\;
+2^{12}3^8z^4(-1296\theta^4-20736\theta^3-20736\theta^2-7920\theta+1)
\\[-2pt] &\;
+2^{18}3^{13}z^5(-72\theta^4+178\theta^2+178\theta+51)
+2^{24}3^{18}z^6(\theta+1)^4
\\[-8pt]
\endaligned &\cr
& &\span\hrulefill&\cr
& &&
\sum_{n=0}^\infty A_nz^n=\text{(e)}*\text{(h)}
\qquad\text{(see \cite{AZ}, Section 7)}
&\cr
\noalign{\hrule}
}}\hss}

\newpage

\hbox to\hsize{\hss\vbox{\offinterlineskip
\halign to120mm{\strut\tabskip=100pt minus 100pt
\strut\vrule#&\hbox to5.5mm{\hss$#$\hss}&%
\vrule#&\hbox to114mm{\hfil$\dsize#$\hfil}&%
\vrule#\tabskip=0pt\cr\noalign{\hrule}
& \# &%
& \text{differential operator $D$ and coefficients $A_n$, $n=0,1,2,\dots$} &\cr
\noalign{\hrule\vskip1pt\hrule}
& \eqnno{172} &&
\aligned
\\[-12pt]
&
D=\theta^4
-3^2z(3\theta^2+3\theta+1)(18\theta^2+18\theta+7)
\\[-2pt] &\;
+3^5z^2(243\theta^4+1944\theta^3+2763\theta^2+1638\theta+409)
\\[-2pt] &\;
+3^{11}z^3(54\theta^4-324\theta^3-987\theta^2-765\theta-227)
\\[-2pt] &\;
+3^{14}z^4(-648\theta^4-1296\theta^3+2898\theta^2+3456\theta+1333)
\\[-2pt] &\;
+3^{20}z^5(54\theta^4+540\theta^3+309\theta^2-21\theta-7)
\\[-2pt] &\;
+3^{23}z^6(243\theta^4-972\theta^3-1611\theta^2-972\theta-167)
\\[-2pt] &\;
-3^{29}z^7(3\theta^2+3\theta+1)(18\theta^2+18\theta+7)
+3^{36}z^8(\theta+1)^4
\\[-8pt]
\endaligned &\cr
& &\span\hrulefill&\cr
& &&
\sum_{n=0}^\infty A_nz^n=\text{(f)}*\text{(h)}
\qquad\text{(see \cite{AZ}, Section 7)}
&\cr
\noalign{\hrule}
& \eqnno{173} &&
\aligned
\\[-12pt]
&
D=\theta^4
-z(7\theta^2+7\theta+2)(17\theta^2+17\theta+6)
\\[-2pt] &\;
+2^6z^2(55\theta^4+112\theta^3+155\theta^2+86\theta+15)
\\[-2pt] &\;
+2^63^2z^3(-119\theta^4+714\theta^3+2185\theta^2+1656\theta+444)
\\[-2pt] &\;
+2^{12}3^2z^4(92\theta^4+184\theta^3+98\theta^2+6\theta+9)
\\[-2pt] &\;
+2^{12}3^4z^5(119\theta^4+1190\theta^3+671\theta^2-96\theta-140)
\\[-2pt] &\;
+2^{18}3^4z^6(55\theta^4+108\theta^3+149\theta^2+108\theta+27)
\\[-2pt] &\;
+2^{18}3^6z^7(7\theta^2+7\theta+2)(17\theta^2+17\theta+6)
+2^{24}3^8z^8(\theta+1)^4
\\[-8pt]
\endaligned &\cr
& &\span\hrulefill&\cr
& &&
\sum_{n=0}^\infty A_nz^n=\text{(a)}*\text{(g)}
\qquad\text{(see \cite{AZ}, Section 7)}
&\cr
\noalign{\hrule}
& \eqnno{174} &&
\aligned
\\[-12pt]
&
D=\theta^4
-z(11\theta^2+11\theta+3)(17\theta^2+17\theta+6)
\\[-2pt] &\;
+z^2(8711\theta^4+33980\theta^3+47095\theta^2+26230\theta+5232)
\\[-2pt] &\;
+2^33^2z^3(-187\theta^4+1122\theta^3+3436\theta^2+2595\theta+684)
\\[-2pt] &\;
+2^43^2z^4(8639\theta^4+17278\theta^3-11650\theta^2-20289\theta-6102)
\\[-2pt] &\;
+2^63^4z^5(187\theta^4+1870\theta^3+1052\theta^2-163\theta-216)
\\[-2pt] &\;
+2^63^4z^6(8711\theta^4+864\theta^3-2579\theta^2+864\theta+828)
\\[-2pt] &\;
+2^93^6z^7(11\theta^2+11\theta+3)(17\theta^2+17\theta+6)
+2^{12}3^8z^8(\theta+1)^4
\\[-8pt]
\endaligned &\cr
& &\span\hrulefill&\cr
& &&
\sum_{n=0}^\infty A_nz^n=\text{(b)}*\text{(g)}
\qquad\text{(see \cite{AZ}, Section 7)}
&\cr
\noalign{\hrule}
& \eqnno{175} &&
\aligned
\\[-12pt]
&
D=\theta^4
-z(10\theta^2+10\theta+3)(17\theta^2+17\theta+6)
\\[-2pt] &\;
+3^4z^2(89\theta^4+452\theta^3+633\theta^2+362\theta+80)
\\[-2pt] &\;
+2^33^4z^3(170\theta^4-1020\theta^3-3119\theta^2-2373\theta-648)
\\[-2pt] &\;
+2^43^8z^4(-97\theta^4-194\theta^3+238\theta^2+335\theta+114)
\\[-2pt] &\;
+2^63^8z^5(170\theta^4+1700\theta^3+961\theta^2-125\theta-204)
\\[-2pt] &\;
+2^63^{12}z^6(89\theta^4-96\theta^3-189\theta^2-96\theta-12)
\\[-2pt] &\;
-2^93^{12}z^7(10\theta^2+10\theta+3)(17\theta^2+17\theta+6)
+2^{12}3^{16}z^8(\theta+1)^4
\\[-8pt]
\endaligned &\cr
& &\span\hrulefill&\cr
& &&
\sum_{n=0}^\infty A_nz^n=\text{(c)}*\text{(g)}
\qquad\text{(see \cite{AZ}, Section 7)}
&\cr
\noalign{\hrule}
}}\hss}

\newpage

\hbox to\hsize{\hss\vbox{\offinterlineskip
\halign to120mm{\strut\tabskip=100pt minus 100pt
\strut\vrule#&\hbox to5.5mm{\hss$#$\hss}&%
\vrule#&\hbox to114mm{\hfil$\dsize#$\hfil}&%
\vrule#\tabskip=0pt\cr\noalign{\hrule}
& \# &%
& \text{differential operator $D$ and coefficients $A_n$, $n=0,1,2,\dots$} &\cr
\noalign{\hrule\vskip1pt\hrule}
& \eqnno{176} &&
\aligned
\\[-12pt]
&
D=\theta^4
-2^2z(3\theta^2+3\theta+1)(17\theta^2+17\theta+6)
\\[-2pt] &\;
+2^5z^2(325\theta^4+2164\theta^3+3053\theta^2+1778\theta+420)
\\[-2pt] &\;
+2^{10}3^2z^3(51\theta^4-306\theta^3-934\theta^2-717\theta-204)
\\[-2pt] &\;
+2^{14}3^2z^4(-397\theta^4-794\theta^3+1454\theta^2+1851\theta+666)
\\[-2pt] &\;
+2^{18}3^4z^5(51\theta^4+510\theta^3+290\theta^2-29\theta-64)
\\[-2pt] &\;
+2^{21}3^4z^6(325\theta^4-864\theta^3-1489\theta^2-864\theta-144)
\\[-2pt] &\;
-2^{26}3^6z^7(3\theta^2+3\theta+1)(17\theta^2+17\theta+6)
+2^{32}3^8z^8(\theta+1)^4
\\[-8pt]
\endaligned &\cr
& &\span\hrulefill&\cr
& &&
\sum_{n=0}^\infty A_nz^n=\text{(d)}*\text{(g)}
\qquad\text{(see \cite{AZ}, Section 7)}
&\cr
\noalign{\hrule}
& \eqnno{177} &&
\aligned
\\[-12pt]
&
D=\theta^4
-2^2z(8\theta^2+8\theta+3)(17\theta^2+17\theta+6)
\\[-2pt] &\;
+2^7z^2(578\theta^4+4040\theta^3+5746\theta^2+3412\theta+849)
\\[-2pt] &\;
+2^{13}3^2z^3(136\theta^4-816\theta^3-2485\theta^2-1929\theta-576)
\\[-2pt] &\;
+2^{19}3^2z^4(-722\theta^4-1444\theta^3+2764\theta^2+3486\theta+1323)
\\[-2pt] &\;
+2^{24}3^4z^5(136\theta^4+1360\theta^3+779\theta^2-49\theta-180)
\\[-2pt] &\;
+2^{29}3^4z^6(578\theta^4-1728\theta^3-2906\theta^2-1728\theta-279)
\\[-2pt] &\;
-2^{35}3^6z^7(8\theta^2+8\theta+3)(17\theta^2+17\theta+6)
+2^{44}3^8z^8(\theta+1)^4
\\[-8pt]
\endaligned &\cr
& &\span\hrulefill&\cr
& &&
\sum_{n=0}^\infty A_nz^n=\text{(e)}*\text{(g)}
\qquad\text{(see \cite{AZ}, Section 7)}
&\cr
\noalign{\hrule}
& \eqnno{178} &&
\aligned
\\[-12pt]
&
D=\theta^4
-3z(3\theta^2+3\theta+1)(17\theta^2+17\theta+6)
\\[-2pt] &\;
+3^3z^2(217\theta^4+1732\theta^3+2441\theta^2+1418\theta+336)
\\[-2pt] &\;
+2^33^6z^3(51\theta^4-306\theta^3-934\theta^2-717\theta-204)
\\[-2pt] &\;
+2^43^8z^4(-289\theta^4-578\theta^3+1310\theta^2+1599\theta+570)
\\[-2pt] &\;
+2^63^{11}z^5(51\theta^4+510\theta^3+290\theta^2-29\theta-64)
\\[-2pt] &\;
+2^63^{13}z^6(217\theta^4-864\theta^3-1453\theta^2-864\theta-156)
\\[-2pt] &\;
-2^93^{16}z^7(3\theta^2+3\theta+1)(17\theta^2+17\theta+6)
+2^{12}3^{20}z^8(\theta+1)^4
\\[-8pt]
\endaligned &\cr
& &\span\hrulefill&\cr
& &&
\sum_{n=0}^\infty A_nz^n=\text{(f)}*\text{(g)}
\qquad\text{(see \cite{AZ}, Section 7)}
&\cr
\noalign{\hrule}
& \eqnno{179} &&
\aligned
\\[-12pt]
&
D=\theta^4
-3z(17\theta^2+17\theta+6)(18\theta^2+18\theta+7)
\\[-2pt] &\;
+3^4z^2(2601\theta^4+18180\theta^3+25929\theta^2+15498\theta+3920)
\\[-2pt] &\;
+2^33^9z^3(306\theta^4-1836\theta^3-5587\theta^2-4353\theta-1320)
\\[-2pt] &\;
+2^43^{12}z^4(-3249\theta^4-6498\theta^3+12366\theta^2+15615\theta+6034)
\\[-2pt] &\;
+2^63^{17}z^5(306\theta^4+3060\theta^3+1757\theta^2-89\theta-412)
\\[-2pt] &\;
+2^63^{20}z^6(2601\theta^4-7776\theta^3-13005\theta^2-7776\theta-1228)
\\[-2pt] &\;
-2^93^{25}z^7(17\theta^2+17\theta+6)(18\theta^2+18\theta+7)
+2^{12}3^{32}z^8(\theta+1)^4
\\[-8pt]
\endaligned &\cr
& &\span\hrulefill&\cr
& &&
\sum_{n=0}^\infty A_nz^n=\text{(g)}*\text{(h)}
\qquad\text{(see \cite{AZ}, Section 7)}
&\cr
\noalign{\hrule}
}}\hss}

\newpage

\hbox to\hsize{\hss\vbox{\offinterlineskip
\halign to120mm{\strut\tabskip=100pt minus 100pt
\strut\vrule#&\hbox to6.5mm{\hss$#$\hss}&%
\vrule#&\hbox to113mm{\hfil$\dsize#$\hfil}&%
\vrule#\tabskip=0pt\cr\noalign{\hrule}
& \# &%
& \text{differential operator $D$ and coefficients $A_n$, $n=0,1,2,\dots$} &\cr
\noalign{\hrule\vskip1pt\hrule}
& \eqnno{180} &&
\aligned
\\[-12pt]
&
D=\theta^4-2^43z(198\theta^4+72\theta^3+69\theta^2+33\theta+5)
\\[-2pt] &\;
+2^93^2z^2(7614\theta^4+7128\theta^3+6813\theta^2+2529\theta+340)
\\[-2pt] &\;
-2^{14}3^5z^3(15714\theta^4+27216\theta^3+26343\theta^2+11151\theta+1685)
\\[-2pt] &\;
+2^{19}3^9z^4(3\theta+1)(3\theta+2)(576\theta^2+1008\theta+605)
\\[-2pt] &\;
-2^{27}3^{13}z^5(3\theta+1)(3\theta+2)(3\theta+4)(3\theta+5)
\\[-8pt]
\endaligned &\cr
& &\span\hrulefill&\cr
& &&
A_n=\binom{2n}n\sum_k\binom nk\binom{2n}{2k}^{-1}
\frac{(6k)!}{k!(2k)!(3k)!}\,
\frac{(6n-6k)!}{(n-k)!(2n-2k)!(3n-3k)!}
&\cr
\noalign{\hrule}
& \eqnno{181} &&
\aligned
\\[-12pt]
&
D=\theta^4
-18z(324\theta^4+648\theta^3+765\theta^2+441\theta+97)
\\[-2pt] &\;
+236196z^2(\theta+1)^2(6\theta+5)(6\theta+7)
\\[-8pt]
\endaligned &\cr
& &\span\hrulefill&\cr
& &&
\text{a formula for $A_n$ is not known}
&\cr
\noalign{\hrule}
& \eqnno{182} &&
\aligned
\\[-12pt]
&
D=\theta^4
-z(43\theta^4+86\theta^3+77\theta^2+34\theta+6)
\\[-2pt] &\;
+12z^2(\theta+1)^2(6\theta+5)(6\theta+7)
\\[-8pt]
\endaligned &\cr
& &\span\hrulefill&\cr
& &&
\text{a formula for $A_n$ is not known}
&\cr
\noalign{\hrule}
& \eqnno{183} &&
\aligned
\\[-12pt]
&
D=\theta^4
-4z(2\theta+1)^2(7\theta^2+7\theta+3)
\\[-2pt] &\;
+48z^2(2\theta+1)(2\theta+3)(4\theta+3)(4\theta+5)
\\[-8pt]
\endaligned &\cr
& &\span\hrulefill&\cr
& &&
A_0=1, \;
A_n=3\binom{2n}n\sum_{k=0}^{[n/3]}(-1)^k\frac{n-2k}{2n-3k}\binom nk
\binom{2k}k\binom{2n-2k}{n+k}\binom{2n-3k}n
&\cr
\noalign{\hrule}
& \eqnnol{184}{187} &&
\aligned
\\[-12pt]
&
D=\theta^4
-2z(2\theta+1)^2(11\theta^2+11\theta+5)
\\[-2pt] &\;
+500z^2(\theta+1)^2(2\theta+1)(2\theta+3)
\\[-8pt]
\endaligned &\cr
& &\span\hrulefill&\cr
& &&
A_0=1, \quad
A_n=5\binom{2n}n\sum_{k=0}^{[n/5]}(-1)^k\frac{n-2k}{4n-5k}\binom nk^3\binom{4n-5k}{3n}
&\cr
\noalign{\hrule}
& \eqnnol{185}{184} &&
\aligned
\\[-12pt]
&
D=\theta^4
-6z(2\theta+1)^2(3\theta^2+3\theta+1)
\\[-2pt] &\;
-108z^2(\theta+1)^2(2\theta+1)(2\theta+3)
\\[-8pt]
\endaligned &\cr
& &\span\hrulefill&\cr
& &&
A_n=\binom{2n}n\sum_{k,l}\binom nk\binom nl
\binom{k+l}k^2\binom n{k+l}
&\cr
\noalign{\hrule}
& \eqnnol{186}{247} &&
\aligned
\\[-12pt]
&
D=19^2\theta^4-19z(700\theta^4+1238\theta^3+999\theta^2+380\theta+57)
\\[-2pt] &\;
-z^2(64745\theta^4+368006\theta^3+609133\theta^2+412756\theta+102258)
\\[-2pt] &\;
+3^3z^3(6397\theta^4+12198\theta^3-11923\theta^2-27360\theta-11286)
\\[-2pt] &\;
+3^6z^4(64\theta^4+1154\theta^3+2425\theta^2+1848\theta+486)
\\[-2pt] &\;
-3^{11}z^5(\theta+1)^4
\\[-8pt]
\endaligned &\cr
& &\span\hrulefill&\cr
& &&
A_n=\sum_{k,l}\binom nk^2\binom nl\binom{k+l}k
\binom{2n-k-l}n\binom{2k}{n-l}
&\cr
\noalign{\hrule}
& \eqnnol{187}{248} &&
\aligned
\\[-12pt]
&
D=\theta^4+z(-64\theta^4+898\theta^3+653\theta^2+204\theta+27)
\\[-2pt] &\;
+3^2z^2(-6397\theta^4-13390\theta^3+10135\theta^2+7492\theta+1850)
\\[-2pt] &\;
+3^4z^3(64745\theta^4-109026\theta^3-106415\theta^2-39528\theta-4626)
\\[-2pt] &\;
+3^919z^4(700\theta^4+1562\theta^3+1485\theta^2+704\theta+138)
\\[-2pt] &\;
-3^{14}19^2z^5(\theta+1)^4
\\[-8pt]
\endaligned &\cr
& &\span\hrulefill&\cr
& &&
\gathered
\text{the reflection of \#186 at infinity}
\\[-8pt]
\endgathered &\cr
& &\span\hrulefill&\cr
& &&
\text{a formula for $A_n$ is not known}
&\cr
\noalign{\hrule}
}}\hss}

\newpage

\hbox to\hsize{\hss\vbox{\offinterlineskip
\halign to120mm{\strut\tabskip=100pt minus 100pt
\strut\vrule#&\hbox to6.5mm{\hss$#$\hss}&%
\vrule#&\hbox to113mm{\hfil$\dsize#$\hfil}&%
\vrule#\tabskip=0pt\cr\noalign{\hrule}
& \# &%
& \text{differential operator $D$ and coefficients $A_n$, $n=0,1,2,\dots$} &\cr
\noalign{\hrule\vskip1pt\hrule}
& \eqnnol{188}{245} &&
\aligned
\\[-12pt]
&
D=\theta^4-2z(280\theta^4+70\theta^3+28\theta^2-7\theta-6)
\\[-2pt] &\;
+2^2z^2(33544\theta^4+16772\theta^3+7451\theta^2+223\theta+354)
\\[-2pt] &\;
-2^5z^3(562360\theta^4+421770\theta^3+206211\theta^2+30144\theta+6888)
\\[-2pt] &\;
+2^6z^4(23219644\theta^4+23219644\theta^3+12399919\theta^2+2343719\theta+158340)
\\[-2pt] &\;
-2^{10}z^5(76599320\theta^4+95749150\theta^3+55494982\theta^2+11674319\theta+584280)
\\[-2pt] &\;
+2^{11}z^6(1300767032\theta^4+1951150548\theta^3+1220689587\theta^2
\\[-2pt] &\;\quad
+274877463\theta+6995856)
\\[-2pt] &\;
-2^{14}z^7(3491380760\theta^4+6109916330\theta^3+4107826961\theta^2
\\[-2pt] &\;\quad
+935352670\theta-69427776)
\\[-2pt] &\;
+2^{15}z^8(22736902622\theta^4+45473805244\theta^3+32745416803\theta^2
\\[-2pt] &\;\quad
+7035711131\theta-1066887690)
\\[-2pt] &\;
-2^{17}15^2z^9(188347832\theta^4+423782622\theta^3+326067072\theta^2
\\[-2pt] &\;\quad
+67245937\theta-10221330)
\\[-2pt] &\;
+2^{18}15^4z^{10}(1735944\theta^4+4339860\theta^3+3554495\theta^2+820435\theta+24846)
\\[-2pt] &\;
-2^{21}15^6z^{11}\theta(2\theta+1)(1036\theta^2+2331\theta+1308)
\\[-2pt] &\;
+2^{22}15^8z^{12}\theta(\theta+1)(2\theta+1)(2\theta+3)
\\[-8pt]
\endaligned &\cr
& &\span\hrulefill&\cr
& &&
\gathered
\text{the pullback of the 5th-order differential equation $D'y=0$, where}
\\[-8pt]
\endgathered &\cr
& &\span\hrulefill&\cr
& &&
\aligned
\\[-12pt]
&
D'=\theta^5
-2z(2\theta+1)(35\theta^4+70\theta^3+63\theta^2+28\theta+5)
\\[-2pt] &\;
+4z^2(\theta+1)(2\theta+1)(2\theta+3)(259\theta^2+518\theta+285)
\\[-2pt] &\;
-1800z^3(\theta+1)(\theta+2)(2\theta+1)(2\theta+3)(2\theta+5)
\\[-8pt]
\endaligned &\cr
& &\span\hrulefill&\cr
& &&
A_n'=\binom{2n}n\sum_{i+j+k+l+m=n}
\biggl(\frac{n!}{i!j!k!l!m!}\biggr)^2
&\cr
\noalign{\hrule}
& \eqnnol{189}{246} &&
\aligned
\\[-12pt]
&
D=\theta^4-z(1040\theta^4+260\theta^3+80\theta^2-50\theta-31)
\\[-2pt] &\;
+z^2(409696\theta^4+204848\theta^3+71714\theta^2+15832\theta+34347)
\\[-2pt] &\;
-2^2z^3(18374720\theta^4+13781040\theta^3+5417700\theta^2+2280015\theta-228513)
\\[-2pt] &\;
+z^4(5406720256\theta^4+5406720256\theta^3+2370020896\theta^2
\\[-2pt] &\;\quad
+470740976\theta-45885347)
\\[-2pt] &\;
-2^7z^5(587991040\theta^4+734988800\theta^3+370894640\theta^2
\\[-2pt] &\;\quad
+26874280\theta-35246511)
\\[-2pt] &\;
+2^{12}z^6(104882176\theta^4+157323264\theta^3+88826496\theta^2
\\[-2pt] &\;\quad
+16395264\theta+8271969)
\\[-2pt] &\;
-2^{21}z^7(532480\theta^4+931840\theta^3+577600\theta^2+160160\theta+10273)
\\[-2pt] &\;
+2^{28}z^8(8\theta+1)(8\theta+3)(8\theta+5)(8\theta+7)
\\[-8pt]
\endaligned &\cr
& &\span\hrulefill&\cr
& &&
\gathered
\text{the pullback of the 5th-order differential equation $D'y=0$, where}
\\[-8pt]
\endgathered &\cr
& &\span\hrulefill&\cr
& &&
\aligned
\\[-12pt]
&
D'=\theta^5-2z(2\theta+1)(65\theta^4+130\theta^3+105\theta^2+40\theta+6)
\\[-2pt] &\;
+16z^2(\theta+1)(2\theta+1)(2\theta+3)(4\theta+3)(4\theta+5)
\\[-8pt]
\endaligned &\cr
& &\span\hrulefill&\cr
& &&
A_n'=\binom{2n}n\sum_{k,l}\binom nk^2\binom nl^2\binom{k+l}n^2
&\cr
\noalign{\hrule}
}}\hss}

\newpage

\hbox to\hsize{\hss\vbox{\offinterlineskip
\halign to120mm{\strut\tabskip=100pt minus 100pt
\strut\vrule#&\hbox to6.5mm{\hss$#$\hss}&%
\vrule#&\hbox to113mm{\hfil$\dsize#$\hfil}&%
\vrule#\tabskip=0pt\cr\noalign{\hrule}
& \# &%
& \text{differential operator $D$ and coefficients $A_n$, $n=0,1,2,\dots$} &\cr
\noalign{\hrule\vskip1pt\hrule}
& \eqnnol{190}{286} &&
\aligned
\\[-12pt]
&
D=\theta^4-z(4112\theta^4+1028\theta^3+262\theta^2-252\theta-143)
\\[-2pt] &\;
+z^2(6357088\theta^4+3178544\theta^3+949302\theta^2+146956\theta+559863)
\\[-2pt] &\;
-2z^3(2198012032\theta^4+1648509024\theta^3+567945288\theta^2
\\[-2pt] &\;\quad
+254177786\theta-23879293)
\\[-2pt] &\;
+z^4(1168836133120\theta^4+1168836133120\theta^3+465472077920\theta^2
\\[-2pt] &\;\quad
+105421717520\theta+12297022465)
\\[-2pt] &\;
-2^{13}z^5(2198012032\theta^4+2747515040\theta^3+2068510800\theta^2
\\[-2pt] &\;\quad
+1199997450\theta+333111035)
\\[-2pt] &\;
+2^{24}z^6(6357088\theta^4+9535632\theta^3+8338518\theta^2+3563032\theta-95477)
\\[-2pt] &\;
-2^{32}z^7(65792\theta^4+115136\theta^3+98608\theta^2+28784\theta+4107)
\\[-2pt] &\;
+2^{44}z^8(2\theta+1)^4
\\[-8pt]
\endaligned &\cr
& &\span\hrulefill&\cr
& &&
\gathered
\text{the pullback of the 5th-order differential equation $D'y=0$, where}
\\[-8pt]
\endgathered &\cr
& &\span\hrulefill&\cr
& &&
\aligned
\\[-12pt]
&
D'=\theta^5-2z(522\theta^5+1285\theta^4+1290\theta^3+650\theta^2+165\theta+17)
\\[-2pt] &\;
+5\cdot 2^2z^2(1032\theta^5+3080\theta^4+4622\theta^3+4618\theta^2+2881\theta+805)
\\[-2pt] &\;
-5\cdot 2^3z^3(4112\theta^5+14360\theta^4+25636\theta^3+27154\theta^2+14336\theta+2101)
\\[-2pt] &\;
+5\cdot 2^4z^4(8208\theta^5+32800\theta^4+63528\theta^3+65552\theta^2+31749\theta+6554)
\\[-2pt] &\;
-2^5z^5(40992\theta^5+184400\theta^4+368720\theta^3+368680\theta^2+184330\theta+36865)
\\[-2pt] &\;
-2^{20}z^6(\theta+1)^5
\\[-8pt]
\endaligned &\cr
& &\span\hrulefill&\cr
& &&
A_n'=\sum_k\binom{2n-2k}{n-k}\binom{2k}k^5
&\cr
\noalign{\hrule}
}}\hss}

\newpage

\hbox to\hsize{\hss\vbox{\offinterlineskip
\halign to120mm{\strut\tabskip=100pt minus 100pt
\strut\vrule#&\hbox to6.5mm{\hss$#$\hss}&%
\vrule#&\hbox to113mm{\hfil$\dsize#$\hfil}&%
\vrule#\tabskip=0pt\cr\noalign{\hrule}
& \# &%
& \text{differential operator $D$ and coefficients $A_n$, $n=0,1,2,\dots$} &\cr
\noalign{\hrule\vskip1pt\hrule}
& \eqnnol{191}{287} &&
\aligned
\\[-12pt]
&
D=\theta^4-z(27664\theta^4+6916\theta^3+1422\theta^2-2036\theta-1091)
\\[-2pt] &\;
+z^2(287096928\theta^4+143548464\theta^3+35543478\theta^2+682140\theta+24820311)
\\[-2pt] &\;
-2z^3(662746448000\theta^4+497059836000\theta^3+144112489416\theta^2
\\[-2pt] &\;\quad
+58159349850\theta-9469320873)
\\[-2pt] &\;
+z^4(2303683700982016\theta^4+2303683700982016\theta^3+768404767144032\theta^2
\\[-2pt] &\;\quad
+101811826189072\theta+5044177631041)
\\[-2pt] &\;
-2^73z^5(95435488512000\theta^4+119294360640000\theta^3+85218257140784\theta^2
\\[-2pt] &\;\quad
+49484972267112\theta+14178690000383)
\\[-2pt] &\;
+2^{12}3^2z^6(5953241899008\theta^4+8929862848512\theta^3+7544293361280\theta^2
\\[-2pt] &\;\quad
+3145351380480\theta-150788399663)
\\[-2pt] &\;
-2^{21}3^5z^7(1147281408\theta^4+2007742464\theta^3+1668935232\theta^2
\\[-2pt] &\;\quad
+456621984\theta+58834537)
\\[-2pt] &\;
+2^{28}3^8z^8(24\theta+7)(24\theta+11)(24\theta+13)(24\theta+17)
\\[-8pt]
\endaligned &\cr
& &\span\hrulefill&\cr
& &&
\gathered
\text{the pullback of the 5th-order differential equation $D'y=0$, where}
\\[-8pt]
\endgathered &\cr
& &\span\hrulefill&\cr
& &&
\aligned
\\[-12pt]
&
D'=\theta^5-2z(3466\theta^5+8645\theta^4+8338\theta^3+3862\theta^2+857\theta+73)
\\[-2pt] &\;
+2^2z^2(34600\theta^5+103720\theta^4+152470\theta^3+149954\theta^2+93053\theta+25841)
\\[-2pt] &\;
-2^3z^3(138320\theta^5+483960\theta^4+851700\theta^3+889722\theta^2+460680\theta+62757)
\\[-2pt] &\;
+2^4z^4(276560\theta^5+1106080\theta^4+2117960\theta^3+2152016\theta^2+1013833)
\\[-2pt] &\;
-2^5z^5(276512\theta^5+1244240\theta^4+2463440\theta^3+2420968\theta^2
\\[-2pt] &\;\quad
+1177258\theta+226801)
\\[-2pt] &\;
+3\cdot 2^{12}z^6(\theta+1)(4\theta+3)(4\theta+5)(6\theta+5)(6\theta+7)
\\[-8pt]
\endaligned &\cr
& &\span\hrulefill&\cr
& &&
A_n'=\sum_k\binom{2n-2k}{n-k}\frac{(3k)!\,(4k)!}{k!^7}
&\cr
\noalign{\hrule}
}}\hss}

\newpage

\hbox to\hsize{\hss\vbox{\offinterlineskip
\halign to120mm{\strut\tabskip=100pt minus 100pt
\strut\vrule#&\hbox to6.5mm{\hss$#$\hss}&%
\vrule#&\hbox to113mm{\hfil$\dsize#$\hfil}&%
\vrule#\tabskip=0pt\cr\noalign{\hrule}
& \# &%
& \text{differential operator $D$ and coefficients $A_n$, $n=0,1,2,\dots$} &\cr
\noalign{\hrule\vskip1pt\hrule}
& \eqnnol{192}{288} &&
\aligned
\\[-12pt]
&
D=\theta^4-z(442384\theta^4+110596\theta^3+18054\theta^2-37244\theta-19343)
\\[-2pt] &\;
+3z^2(24463540256\theta^4+12231770128\theta^3+2506663954\theta^2
\\[-2pt] &\;\quad
-374614268\theta+2063881213)
\\[-2pt] &\;
-2z^3(2705798010765440\theta^4+2029348508074080\theta^3+500517564438600\theta^2
\\[-2pt] &\;\quad
+178669097012730\theta-39466352094909)
\\[-2pt] &\;
+z^4(149673916892499149056\theta^4+149673916892499149056\theta^3
\\[-2pt] &\;\quad
+43170889056861388896\theta^2+2898231542674452496\theta
\\[-2pt] &\;\quad
-195585894794578943)
\\[-2pt] &\;
-2^{11}3z^5(389634913550223360\theta^4+487043641937779200\theta^3
\\[-2pt] &\;\quad
+331475300236995440\theta^2+191938784333046600\theta
\\[-2pt] &\;\quad
+55636145439400469)
\\[-2pt] &\;
+2^{20}3^2z^6(1521827912245248\theta^4+2282741868367872\theta^3
\\[-2pt] &\;\quad
+1865287735567488\theta^2+756921798415872\theta-52441938125015)
\\[-2pt] &\;
-2^{33}3^5z^7(18346549248\theta^4+32106461184\theta^3+25927048512\theta^2
\\[-2pt] &\;\quad
+6633105696\theta+751215433)
\\[-2pt] &\;
+2^{44}3^8z^8(24\theta+5)(24\theta+11)(24\theta+13)(24\theta+19)
\\[-8pt]
\endaligned &\cr
& &\span\hrulefill&\cr
& &&
\gathered
\text{the pullback of the 5th-order differential equation $D'y=0$, where}
\\[-8pt]
\endgathered &\cr
& &\span\hrulefill&\cr
& &&
\aligned
\\[-12pt]
&
D'=\theta^5-2z(55306\theta^5+138245\theta^4+128650\theta^3+54730\theta^2+10469\theta+721)
\\[-2pt] &\;
+2^2z^2(553000\theta^5+1658920\theta^4+2392390\theta^3+2315570\theta^2
\\[-2pt] &\;\quad
+1429445\theta+394553)
\\[-2pt] &\;
-2^3z^3(2211920\theta^5+7741560\theta^4+13440180\theta^3+13847130\theta^2
\\[-2pt] &\;\quad
+7036800\theta+889929)
\\[-2pt] &\;
+2^4z^4(4423760\theta^5+17694880\theta^4+33515720\theta^3+33546320\theta^2
\\[-2pt] &\;\quad
+15386905\theta+2949122)
\\[-2pt] &\;
-2^5z^5(4423712\theta^5+19906640\theta^4+39045200\theta^3+37739560\theta^2
\\[-2pt] &\;\quad
+17863690\theta+3317761)
\\[-2pt] &\;
+2^{18}3z^6(\theta+1)(3\theta+2)(3\theta+4)(4\theta+3)(4\theta+5)
\\[-8pt]
\endaligned &\cr
& &\span\hrulefill&\cr
& &&
A_n'=\sum_k\binom{2k}k\binom{2n-2k}{n-k}\frac{(4k)!}{k!^2(2k)!}
\,\frac{(6k)!}{k!(2k)!(3k)!}
&\cr
\noalign{\hrule}
& \eqnnol{193}{254} &&
\aligned
\\[-12pt]
&
D=7^2\theta^4-7z(1135\theta^4+2204\theta^3+1683\theta^2+581\theta+77)
\\[-2pt] &\;
+z^2(28723\theta^4+40708\theta^3+13260\theta^2-1337\theta-896)
\\[-2pt] &\;
-z^3(32126\theta^4+38514\theta^3+26511\theta^2+10731\theta+1806)
\\[-2pt] &\;
+7\cdot 11z^4(130\theta^4+254\theta^3+192\theta^2+65\theta+8)
+11^2z^5(\theta+1)^4
\\[-8pt]
\endaligned &\cr
& &\span\hrulefill&\cr
& &&
A_n=\sum_{i,j}\binom ni^2\binom nj^2\binom{i+j}j\binom{n+i+j}n
&\cr
\noalign{\hrule}
}}\hss}

\newpage

\hbox to\hsize{\hss\vbox{\offinterlineskip
\halign to120mm{\strut\tabskip=100pt minus 100pt
\strut\vrule#&\hbox to6.5mm{\hss$#$\hss}&%
\vrule#&\hbox to113mm{\hfil$\dsize#$\hfil}&%
\vrule#\tabskip=0pt\cr\noalign{\hrule}
& \# &%
& \text{differential operator $D$ and coefficients $A_n$, $n=0,1,2,\dots$} &\cr
\noalign{\hrule\vskip1pt\hrule}
& \eqnnol{194}{255} &&
\aligned
\\[-12pt]
&
D=17^2\theta^4-17z(1465\theta^4+2768\theta^3+2200\theta^2+816\theta+119)
\\[-2pt] &\;
+2z^2(62015\theta^4+131582\theta^3+125017\theta^2+65926\theta+15300)
\\[-2pt] &\;
-2\cdot 3^3z^3(4325\theta^4+10914\theta^3+12803\theta^2+7446\theta+1700)
\\[-2pt] &\;
+3^6z^4(265\theta^4+836\theta^3+1118\theta^2+700\theta+168)
-3^{10}z^5(\theta+1)^4
\\[-8pt]
\endaligned &\cr
& &\span\hrulefill&\cr
& &&
A_n=\sum_{i,j}\binom ni^2\binom nj^2\binom{i+j}j^2
&\cr
\noalign{\hrule}
& \eqnnol{195}{256} &&
\aligned
\\[-12pt]
&
D=29^2\theta^4-29z(3026\theta^4+5848\theta^3+4577\theta^2+1653\theta+232)
\\[-2pt] &\;
+z^2(258647\theta^4+424220\theta^3+239159\theta^2+57768\theta+5568)
\\[-2pt] &\;
-z^3(272743\theta^4+532614\theta^3+581647\theta^2+336864\theta+76560)
\\[-2pt] &\;
+2^217z^4(1922\theta^4+6193\theta^3+8121\theta^2+4894\theta+1112)
\\[-2pt] &\;
-2^23\cdot 17^2z^5(\theta+1)^2(3\theta+2)(3\theta+4)
\\[-8pt]
\endaligned &\cr
& &\span\hrulefill&\cr
& &&
A_n=\sum_{i,j}\binom ni^2\binom nj^2\binom{i+j}j\binom{n+i}n
&\cr
\noalign{\hrule}
& \eqnnol{196}{257} &&
\aligned
\\[-12pt]
&
D=47^2\theta^4-47z(2489\theta^4+4984\theta^3+4043\theta^2+1551\theta+235)
\\[-2pt] &\;
-z^2(208867\theta^4+790072\theta^3+1135848\theta^2+701851\theta+161022)
\\[-2pt] &\;
+z^3(37085\theta^4+637644\theta^3+383912\theta^2+149319\theta+38352)
\\[-2pt] &\;
+z^4(291161\theta^4-511820\theta^3-4424049\theta^2-5161283\theta-1770676)
\\[-2pt] &\;
-z^5(406192\theta^4+749482\theta^3+750755\theta^2+260936\theta-2151)
\\[-2pt] &\;
+3^3z^6(5305\theta^4+90750\theta^3+152551\theta^2+91194\theta+17914)
\\[-2pt] &\;
+2\cdot 3^6z^7(106\theta^4+230\theta^3+197\theta^2+82\theta+15)
-2^23^{10}z^8(\theta+1)^4
\\[-8pt]
\endaligned &\cr
& &\span\hrulefill&\cr
& &&
A_n=\sum_{i,j}\binom ni^2\binom nj^2\binom{i+j}j\binom{n+i-j}n
&\cr
\noalign{\hrule}
& \eqnnol{197}{258} &&
\aligned
\\[-12pt]
&
D=13^2\theta^4-13^2z(41\theta^4+82\theta^3+67\theta^2+26\theta+4)
\\[-2pt] &\;
-13\cdot 2^3z^2(471\theta^4+1788\theta^3+2555\theta^2+1534\theta+338)
\\[-2pt] &\;
+13\cdot 2^6z^3(251\theta^4+1014\theta^3+1798\theta^2+1413\theta+405)
\\[-2pt] &\;
+2^9z^4(749\theta^4+436\theta^3-4908\theta^2-6266\theta-2145)
\\[-2pt] &\;
-2^{12}z^5(379\theta^4+1270\theta^3+967\theta^2-42\theta-178)
\\[-2pt] &\;
+2^{15}z^6(-9\theta^4+156\theta^3+273\theta^2-156\theta-28)
\\[-2pt] &\;
+2^{18}z^7(13\theta^4+26\theta^3+20\theta^2+7\theta+1)
-2^{21}z^8(\theta+1)^4
\\[-8pt]
\endaligned &\cr
& &\span\hrulefill&\cr
& &&
A_n=\sum_{i,j}\binom ni^2\binom nj^2\binom ji\binom{i+j}j
&\cr
\noalign{\hrule}
}}\hss}

\newpage

\hbox to\hsize{\hss\vbox{\offinterlineskip
\halign to120mm{\strut\tabskip=100pt minus 100pt
\strut\vrule#&\hbox to6.5mm{\hss$#$\hss}&%
\vrule#&\hbox to113mm{\hfil$\dsize#$\hfil}&%
\vrule#\tabskip=0pt\cr\noalign{\hrule}
& \# &%
& \text{differential operator $D$ and coefficients $A_n$, $n=0,1,2,\dots$} &\cr
\noalign{\hrule\vskip1pt\hrule}
& \eqnnol{198}{259} &&
\aligned
\\[-12pt]
&
D=11^2\theta^4-7\cdot 11z(130\theta^4+266\theta^3+210\theta^2+77\theta+11)
\\[-2pt] &\;
-z^2(32126\theta^4+89990\theta^3+103725\theta^2+55253\theta+11198)
\\[-2pt] &\;
+z^3(28723\theta^4+74184\theta^3+63474\theta^2+20625\theta+1716)
\\[-2pt] &\;
-7z^4(1135\theta^4+2336\theta^3+1881\theta^2+713\theta+110)
+7^2z^5(\theta+1)^4
\\[-8pt]
\endaligned &\cr
& &\span\hrulefill&\cr
& &&
\gathered
\text{the reflection of \#193 at infinity}
\\[-8pt]
\endgathered &\cr
& &\span\hrulefill&\cr
& &&
A_n=(-1)^n\sum_{i,j}\binom ni^2\binom nj^2\binom{i+j}j\binom{2n-i}n
&\cr
\noalign{\hrule}
& \eqnnol{199}{260} &&
\aligned
\\[-12pt]
&
D=\theta^4-z(265\theta^4+224\theta^3+200\theta^2+88\theta+15)
\\[-2pt] &\;
+2\cdot 3z^2(4325\theta^4+6386\theta^3+6011\theta^2+2718\theta+468)
\\[-2pt] &\;
-2\cdot 3^2z^3(62015\theta^4+116478\theta^3+102361\theta^2+37422\theta+4824)
\\[-2pt] &\;
+17\cdot 3^6z^4(1465\theta^4+3092\theta^3+2686\theta^2+1140\theta+200)
\\[-2pt] &\;
-17^23^{10}z^5(\theta+1)^4
\\[-8pt]
\endaligned &\cr
& &\span\hrulefill&\cr
& &&
\gathered
\text{the reflection of \#194 at infinity}
\\[-8pt]
\endgathered &\cr
& &\span\hrulefill&\cr
& &&
\text{a formula for $A_n$ is not known}
&\cr
\noalign{\hrule}
& \eqnnol{200}{261} &&
\aligned
\\[-12pt]
&
D=2^2\theta^4-2z(106\theta^4+194\theta^3+143\theta^2+46\theta+6)
\\[-2pt] &\;
+3z^2(-5305\theta^4+6953\theta^3+87869\theta^2+37122\theta+6174)
\\[-2pt] &\;
+3^2z^3(406192\theta^4+875286\theta^3+939461\theta^2+616896\theta+144378)
\\[-2pt] &\;
+3^6z^4(-291161\theta^4-1676464\theta^3+1141623\theta^2+986711\theta+230461)
\\[-2pt] &\;
-3^{10}z^5(370857\theta^4+845784\theta^3+696122\theta^2+189001\theta+6158)
\\[-2pt] &\;
+3^{14}z^6(208867\theta^4+45396\theta^3+18834\theta^2+35097\theta+13814)
\\[-2pt] &\;
+47\cdot 3^{18}z^7(2489\theta^4+4972\theta^3+4025\theta^2+1539\theta+232)
\\[-2pt] &\;
-47^23^{22}z^8(\theta+1)^4
\\[-8pt]
\endaligned &\cr
& &\span\hrulefill&\cr
& &&
\gathered
\text{the reflection of \#196 at infinity}
\\[-8pt]
\endgathered &\cr
& &\span\hrulefill&\cr
& &&
\text{a formula for $A_n$ is not known}
&\cr
\noalign{\hrule}
& \eqnnol{201}{262} &&
\aligned
\\[-12pt]
&
D=\theta^4-2^4z(13\theta^4+26\theta^3+20\theta^2+7\theta+1)
\\[-2pt] &\;
+2^8z^2(9\theta^4+192\theta^3+249\theta^2+114\theta+20)
\\[-2pt] &\;
+2^{12}z^3(379\theta^4+246\theta^3-569\theta^2-318\theta-60)
\\[-2pt] &\;
+2^{16}z^4(-749\theta^4-2560\theta^3+1722\theta^2+1862\theta+474)
\\[-2pt] &\;
-13\cdot 2^{20}z^5(251\theta^4-10\theta^3+262\theta^2+145\theta+27)
\\[-2pt] &\;
+13\cdot 2^{24}z^6(471\theta^4+96\theta^3+17\theta^2+96\theta+42)
\\[-2pt] &\;
+13^22^{28}z^7(41\theta^4+82\theta^3+67\theta^2+26\theta+4)
-13^22^{35}z^8(\theta+1)^4
\\[-8pt]
\endaligned &\cr
& &\span\hrulefill&\cr
& &&
\gathered
\text{the reflection of \#197 at infinity}
\\[-8pt]
\endgathered &\cr
& &\span\hrulefill&\cr
& &&
\text{a formula for $A_n$ is not known}
&\cr
\noalign{\hrule}
}}\hss}

\newpage

\hbox to\hsize{\hss\vbox{\offinterlineskip
\halign to120mm{\strut\tabskip=100pt minus 100pt
\strut\vrule#&\hbox to6.5mm{\hss$#$\hss}&%
\vrule#&\hbox to113mm{\hfil$\dsize#$\hfil}&%
\vrule#\tabskip=0pt\cr\noalign{\hrule}
& \# &%
& \text{differential operator $D$ and coefficients $A_n$, $n=0,1,2,\dots$} &\cr
\noalign{\hrule\vskip1pt\hrule}
& \eqnnol{202}{272} &&
\aligned
\\[-12pt]
&
D=19^2\theta^4-19z(1370\theta^4+2620\theta^3+2089\theta^2+778\theta+114)
\\[-2pt] &\;
+z^2(39521\theta^4-3916\theta^3-106779\theta^2-95266\theta-25384)
\\[-2pt] &\;
+2^3z^3(1649\theta^4+19779\theta^3+29667\theta^2+17613\theta+3876)
\\[-2pt] &\;
-5\cdot 2^4z^4(\theta+1)(499\theta^3+1411\theta^2+1378\theta+456)
+5^22^9z^5(\theta+1)^4
\\[-8pt]
\endaligned &\cr
& &\span\hrulefill&\cr
& &&
A_n=\sum_{i,j}\binom ni^2\binom nj^2\binom{i+j}j\binom{n+i-j}{n-j}
&\cr
\noalign{\hrule}
& \eqnnol{203}{273} &&
\aligned
\\[-12pt]
&
D=5^2\theta^4-5z\theta(499\theta^3+86\theta^2+53\theta+10)
\\[-2pt] &\;
+2^4z^2(1649\theta^4-13183\theta^3-19776\theta^2-11020\theta-2200)
\\[-2pt] &\;
+2^6z^3(39521\theta^4+162000\theta^3+142095\theta^2+51540\theta+6540)
\\[-2pt] &\;
-19\cdot 2^{11}z^4(1370\theta^4+2860\theta^3+2449\theta^2+1019\theta+174)
\\[-2pt] &\;
+19^22^{16}z^5(\theta+1)^4
\\[-8pt]
\endaligned &\cr
& &\span\hrulefill&\cr
& &&
\gathered
\text{the reflection of \#202 at infinity}
\\[-8pt]
\endgathered &\cr
& &\span\hrulefill&\cr
& &&
\text{a formula for $A_n$ is not known}
&\cr
\noalign{\hrule}
& \eqnnol{204}{200} &&
\aligned
\\[-12pt]
&
D=\theta^4-16z(128\theta^4+256\theta^3+304\theta^2+176\theta+39)
+2^{20}z^2(\theta+1)^4
\\[-8pt]
\endaligned &\cr
& &\span\hrulefill&\cr
& &&
\gathered
\text{self-dual at infinity}
\\[-8pt]
\endgathered &\cr
& &\span\hrulefill&\cr
& &&
\text{a formula for $A_n$ is not known}
&\cr
\noalign{\hrule}
& \eqnno{205} &&
\aligned
\\[-12pt]
&
D=\theta^4-z(59\theta^4+118\theta^3+105\theta^2+46\theta+8)
\\[-2pt] &\;
+96z^2(\theta+1)^2(3\theta+2)(3\theta+4)
\\[-8pt]
\endaligned &\cr
& &\span\hrulefill&\cr
& &&
A_0=1, \quad
A_n=4\sum_{k=0}^{[n/4]}\frac{n-2k}{3n-4k}\binom nk^2\binom{2k}k\binom{2n-2k}{n-k}\binom{3n-4k}{2n}
&\cr
\noalign{\hrule}
& \eqnnoltmp{206}{244} &&
\aligned
\\[-12pt]
&
D=\theta^4-2^2z\theta(\theta+1)(2\theta+1)^2
\\[-2pt] &\;
-2^5z^2(2\theta+1)(2\theta+3)(11\theta^2+22\theta+12)
\\[-2pt] &\;
-2^4\cdot3\cdot5^2z^3(2\theta+1)(2\theta+3)^2(2\theta+5)
\\[-2pt] &\;
-2^8\cdot19z^4(2\theta+1)(2\theta+3)(2\theta+5)(2\theta+7)
\\[-8pt]
\endaligned &\cr
& &\span\hrulefill&\cr
& &&
\gathered
\text{$A_n$ is the constant term of $S^{2n}$
(Batyrev \#14.16326; Kreuzer $X_{44,92}^{65}$),}
\\[-2pt]
\aligned
\\[-13pt]
\text{where}\;
S&=x+y+x+t+\frac z{xy}+\frac{zt}y
+\frac y{xz}+\frac t{xy}+\frac y{xt}
\\[-2pt] &\quad
+\frac1{x^2z}+\frac1{x^2y}+\frac1{x^2t}
+\frac1{x^3zt}+\frac y{x^2zt}
\endaligned
\endgathered &\cr
\noalign{\hrule}
}}\hss}

\newpage

\hbox to\hsize{\hss\vbox{\offinterlineskip
\halign to120mm{\strut\tabskip=100pt minus 100pt
\strut\vrule#&\hbox to6.5mm{\hss$#$\hss}&%
\vrule#&\hbox to113mm{\hfil$\dsize#$\hfil}&%
\vrule#\tabskip=0pt\cr\noalign{\hrule}
& \# &%
& \text{differential operator $D$ and coefficients $A_n$, $n=0,1,2,\dots$} &\cr
\noalign{\hrule\vskip1pt\hrule}
& \eqnnol{207}{225} &&
\aligned
\\[-12pt]
&
D=\theta^4+2^4z(-1072\theta^4+17824\theta^3+10888\theta^2+1976\theta+145)
\\[-2pt] &\;
+2^{17}z^2(-51088\theta^4-116368\theta^3+45264\theta^2+14228\theta+1397)
\\[-2pt] &\;
+13\cdot 2^{28}z^3(73104\theta^4+1536\theta^3-488\theta^2+384\theta+97)
\\[-2pt] &\;
-13^22^{44}z^4(2\theta+1)^4
\\[-8pt]
\endaligned &\cr
& &\span\hrulefill&\cr
& &&
\gathered
\text{the reflection of \#99 at infinity}
\\[-8pt]
\endgathered &\cr
& &\span\hrulefill&\cr
& &&
\text{a formula for $A_n$ is not known}
&\cr
\noalign{\hrule}
& \eqnnol{208}{289} &&
\aligned
\\[-12pt]
&
D=7^2\theta^4-14z(1056\theta^4+1884\theta^3+1397\theta^2+455\theta+56)
\\[-2pt] &\;
+2^2z^2(68280\theta^4+41016\theta^3-67611\theta^2-54348\theta-10752)
\\[-2pt] &\;
+2^4z^3(-53312\theta^4+162120\theta^3+195172\theta^2+78561\theta+11130)
\\[-2pt] &\;
-19\cdot 2^6z^4(2\theta+1)^2(1189\theta^2+2533\theta+1646)
\\[-2pt] &\;
+19^22^{11}z^5(2\theta+1)^2(2\theta+3)^2
\\[-8pt]
\endaligned &\cr
& &\span\hrulefill&\cr
& &&
A_n=\binom{2n}n^2\sum_k\binom nk^2\binom{n+2k}n
&\cr
\noalign{\hrule}
& \eqnnol{209}{290} &&
\aligned
\\[-12pt]
&
D=17^2\theta^4-34z(1902\theta^4+3708\theta^3+2789\theta^2+935\theta+119)
\\[-2pt] &\;
+2^2z^2(62408\theta^4+68576\theta^3-10029\theta^2-24106\theta-5661)
\\[-2pt] &\;
-2^2z^3(66180\theta^4+33048\theta^3+20785\theta^2+17799\theta+4794)
\\[-2pt] &\;
+2^7z^4(2\theta+1)(196\theta^3+498\theta^2+487\theta+169)
\\[-2pt] &\;
-2^{12}z^5(\theta+1)^2(2\theta+1)(2\theta+3)
\\[-8pt]
\endaligned &\cr
& &\span\hrulefill&\cr
& &&
A_n=\binom{2n}n\sum_k\binom nk^2\binom{n+k}n\binom{n+2k}n
&\cr
\noalign{\hrule}
& \eqnnol{210}{263} &&
\aligned
\\[-12pt]
&
D=5^2\theta^4
-20z(688\theta^4+1352\theta^3+981\theta^2+305\theta+35)
\\[-2pt] &\;
+2^4z^2(5856\theta^4+7008\theta^3+96\theta^2-1260\theta-265)
\\[-2pt] &\;
+2^{10}z^3(176\theta^4+120\theta^3+69\theta^2+30\theta+5)
\\[-2pt] &\;
+2^{12}z^4(2\theta+1)^4
\\[-8pt]
\endaligned &\cr
& &\span\hrulefill&\cr
& &&
A_n=\binom{2n}n\sum_k(-1)^k\binom{2n}k^4
&\cr
\noalign{\hrule}
& \eqnnol{211}{271} &&
\aligned
\\[-12pt]
&
D=\theta^4
+2^4z(704\theta^4+928\theta^3+612\theta^2+148\theta+13)
\\[-2pt] &\;
+2^{12}z^2(5856\theta^4+4704\theta^3-1632\theta^2-972\theta-121)
\\[-2pt] &\;
+2^{20}\cdot5z^3(2752\theta^4+96\theta^3-60\theta^2+24\theta+7)
\\[-2pt] &\;
+2^{28}\cdot5^2z^4(2\theta+1)^4
\\[-8pt]
\endaligned &\cr
& &\span\hrulefill&\cr
& &&
\gathered
\text{the reflection of \#210 at infinity}
\\[-8pt]
\endgathered &\cr
& &\span\hrulefill&\cr
& &&
\aligned
A_n
&=(-1)^n\binom{2n}n^4\biggl(\sum_{k=0}^n(-1)^k
\binom nk^2\binom{2k}k\binom{4n-2k}{2n-k}\binom{n+k}n^{-2}\binom{2n}k^{-1}
\\[-2pt] &\quad
+\sum_{k=1}^n\binom nk^2\binom{2n+k}{2n}\binom{4n+2k}{2n+k}
\binom{n+k}n^{-2}\binom{2k}k^{-1}\biggr)
\\[2pt]
\endaligned &\cr
\noalign{\hrule}
}}\hss}

\newpage

\hbox to\hsize{\hss\vbox{\offinterlineskip
\halign to120mm{\strut\tabskip=100pt minus 100pt
\strut\vrule#&\hbox to6.5mm{\hss$#$\hss}&%
\vrule#&\hbox to113mm{\hfil$\dsize#$\hfil}&%
\vrule#\tabskip=0pt\cr\noalign{\hrule}
& \# &%
& \text{differential operator $D$ and coefficients $A_n$, $n=0,1,2,\dots$} &\cr
\noalign{\hrule\vskip1pt\hrule}
& \eqnnol{212}{264} &&
\aligned
\\[-12pt]
&
D=7^2\theta^4
-14z(134\theta^4+286\theta^3+234\theta^2+91\theta+14)
\\[-2pt] &\;
-2^2z^2(3183\theta^4+10266\theta^3+13501\theta^2+8225\theta+1918)
\\[-2pt] &\;
-2^3z^3(2588\theta^4+8400\theta^3+10256\theta^2+5649\theta+1190)
\\[-2pt] &\;
-2^4\cdot3z^4(256\theta^4+848\theta^3+1141\theta^2+717\theta+174)
\\[-2pt] &\;
-2^8\cdot3^2z^5(\theta+1)^4
\\[-8pt]
\endaligned &\cr
& &\span\hrulefill&\cr
& &&
\gathered
\text{the reflection of \#117 at infinity}
\\[-8pt]
\endgathered &\cr
& &\span\hrulefill&\cr
& &&
A_n=\sum_{i,j}\binom ni^2\binom nj^2\binom{i+j}j\binom{2n-i-j}n
&\cr
\noalign{\hrule}
& \eqnnol{213}{265} &&
\aligned
\\[-12pt]
&
D=17^2\theta^4
-34z(647\theta^4+1240\theta^3+977\theta^2+357\theta+51)
\\[-2pt] &\;
-2^2z^2(14437\theta^4+89752\theta^3+147734\theta^2+92123\theta+20400)
\\[-2pt] &\;
+2^2z^3(64614\theta^4+77040\theta^3-125937\theta^2-168453\theta-52326)
\\[-2pt] &\;
+2^3z^4(51920\theta^4+166384\theta^3-83149\theta^2-217017\theta-79362)
\\[-2pt] &\;
+2^4\cdot3z^5(-9360\theta^4+26784\theta^3+43813\theta^2+21965\theta+3496)
\\[-2pt] &\;
+2^5\cdot3z^6(-10160\theta^4+96\theta^3+10535\theta^2+5385\theta+438)
\\[-2pt] &\;
-2^8\cdot3^2z^7(288\theta^4+864\theta^3+1082\theta^2+641\theta+147)
\\[-2pt] &\;
-2^{11}\cdot3^2z^8(\theta+1)^2(4\theta+3)(4\theta+5)
\\[-8pt]
\endaligned &\cr
& &\span\hrulefill&\cr
& &&
A_n=\sum_{i,j}\binom ni^2\binom nj^2\binom{i+j}j\binom{2i}n
&\cr
\noalign{\hrule}
& \eqnnol{214}{266} &&
\aligned
\\[-12pt]
&
D=\theta^4
-2z(90\theta^4+188\theta^3+141\theta^2+47\theta+6)
\\[-2pt] &\;
-2^2z^2(564\theta^4+1520\theta^3+1705\theta^2+934\theta+192)
\\[-2pt] &\;
-2^4z^3(2\theta+1)(286\theta^3+813\theta^2+851\theta+294)
\\[-2pt] &\;
-2^6\cdot3z^4(2\theta+1)(2\theta+3)(4\theta+3)(4\theta+5)
\\[-8pt]
\endaligned &\cr
& &\span\hrulefill&\cr
& &&
A_n=\binom{2n}n\sum_k\binom nk^2\binom{n+k}n\binom{3k}n
&\cr
\noalign{\hrule}
& \eqnnol{215}{267} &&
\aligned
\\[-12pt]
&
D=3^2\theta^4
-12z(268\theta^4+632\theta^3+463\theta^2+147\theta+18)
\\[-2pt] &\;
+2^7z^2(-448\theta^4+1616\theta^3+4280\theta^2+2418\theta+441)
\\[-2pt] &\;
+2^{12}z^3(416\theta^4+2016\theta^3+756\theta^2-288\theta-135)
\\[-2pt] &\;
+2^{19}z^4(2\theta+1)^2(8\theta^2-28\theta-33)
-2^{24}z^5(2\theta+1)^2(2\theta+3)^2
\\[-8pt]
\endaligned &\cr
& &\span\hrulefill&\cr
& &&
A_n=\binom{2n}n^2\sum_k\binom nk^2\binom{4k}{2n}
&\cr
\noalign{\hrule}
}}\hss}

\newpage

\hbox to\hsize{\hss\vbox{\offinterlineskip
\halign to120mm{\strut\tabskip=100pt minus 100pt
\strut\vrule#&\hbox to6.5mm{\hss$#$\hss}&%
\vrule#&\hbox to113mm{\hfil$\dsize#$\hfil}&%
\vrule#\tabskip=0pt\cr\noalign{\hrule}
& \# &%
& \text{differential operator $D$ and coefficients $A_n$, $n=0,1,2,\dots$} &\cr
\noalign{\hrule\vskip1pt\hrule}
& \eqnnol{216}{268} &&
\aligned
\\[-12pt]
&
D=\theta^4
-z\theta(27\theta^3+18\theta^2+11\theta+2)
\\[-2pt] &\;
-2\cdot3^3z^2(72\theta^4+414\theta^3+603\theta^2+330\theta+64)
\\[-2pt] &\;
+2^2\cdot3^5z^3(93\theta^4-720\theta^3-708\theta^2-184)
\\[-2pt] &\;
+2^3\cdot3^7z^4(2\theta+1)(54\theta^3+405\theta^2+544\theta+200)
\\[-2pt] &\;
-2^4\cdot3^{10}z^5(2\theta+1)(2\theta+3)(3\theta+2)(3\theta+4)
\\[-8pt]
\endaligned &\cr
& &\span\hrulefill&\cr
& &&
A_n=\binom{2n}n\sum_k\binom nk^2\binom{3k}n\binom{3n-3k}n
&\cr
\noalign{\hrule}
& \eqnnol{217}{269} &&
\aligned
\\[-12pt]
&
D=7^2\theta^4
+7z\theta(13\theta^3-118\theta^2-73\theta-14)
\\[-2pt] &\;
-2^3\cdot3z^2(3378\theta^4+13446\theta^3+18869\theta^2+11158\theta+2352)
\\[-2pt] &\;
-2^4\cdot3^3z^3(3628\theta^4+17920\theta^3+31668\theta^2+22596\theta+5383)
\\[-2pt] &\;
-2^8\cdot3^3z^4(2\theta+1)(572\theta^3+2370\theta^2+2896\theta+1095)
\\[-2pt] &\;
-2^{10}\cdot3^4z^5(2\theta+1)(2\theta+3)(6\theta+5)(6\theta+7)
\\[-8pt]
\endaligned &\cr
& &\span\hrulefill&\cr
& &&
A_n=\binom{2n}n\sum_k\binom nk^2\binom{3k}n\binom{2n-2k}n
&\cr
\noalign{\hrule}
& \eqnnol{218}{270} &&
\aligned
\\[-12pt]
&
D=7^2\theta^4
-2\cdot3\cdot7z(192\theta^4+396\theta^3+303\theta^2+105\theta+14)
\\[-2pt] &\;
+2^2\cdot3z^2(1188\theta^4+11736\theta^3+20431\theta^2+12152\theta+2436)
\\[-2pt] &\;
+2^2\cdot3^3z^3(532\theta^4+504\theta^3-3455\theta^2-3829\theta-1036)
\\[-2pt] &\;
-2^4\cdot3^4z^4(2\theta+1)(36\theta^3+306\theta^2+421\theta+156)
\\[-2pt] &\;
-2^6\cdot3^4z^5(2\theta+1)(2\theta+3)(3\theta+2)(3\theta+4)
\\[-8pt]
\endaligned &\cr
& &\span\hrulefill&\cr
& &&
A_n=\binom{2n}n\sum_k\binom nk^2\binom{2k}n\binom{3k}n
&\cr
\noalign{\hrule}
& \eqnnol{219}{274} &&
\aligned
\\[-12pt]
&
D=5^2\theta^4
-2\cdot5z(464\theta^4+1036\theta^3+763\theta^2+245\theta+30)
\\[-2pt] &\;
-2^2\cdot3^2z^2(7064\theta^4+22472\theta^3+26699\theta^2+13200\theta+2340)
\\[-2pt] &\;
-2^4\cdot3^4z^3(3440\theta^4+13320\theta^3+18784\theta^2+10665\theta+2070)
\\[-2pt] &\;
-2^6\cdot3^8z^4(2\theta+1)^2(19\theta^2+59\theta+45)
-2^8\cdot3^9z^5(2\theta+1)^2(2\theta+3)^2
\\[-8pt]
\endaligned &\cr
& &\span\hrulefill&\cr
& &&
A_n=\binom{2n}n^2\sum_k\binom nk^2\binom{3k}n
&\cr
\noalign{\hrule}
& \eqnnol{220}{291} &&
\aligned
\\[-12pt]
&
D=\theta^4
-2^4z(20\theta^4+56\theta^3+38\theta^2+10\theta+1)
\\[-2pt] &\;
-2^{10}z^2(84\theta^4+240\theta^3+261\theta^2+134\theta+25)
\\[-2pt] &\;
-2^{16}z^3(2\theta+1)^2(23\theta^2+55\theta+39)
\\[-2pt] &\;
-2^{23}z^4(2\theta+1)^2(2\theta+3)^2
\\[-8pt]
\endaligned &\cr
& &\span\hrulefill&\cr
& &&
\aligned
A_n
&=\binom{2n}n^3\biggl((-1)^n\sum_{k=0}^{[n/2]}\binom nk^2
\binom n{2k}\binom{2n}{4k}^{-1}
\\[-2pt] &\quad
+\sum_{k=[n/2]+1}^n\binom nk^2
\binom{4k}{2n}\binom{2k}n^{-1}\biggr)
\\[2pt]
\endaligned &\cr
\noalign{\hrule}
}}\hss}

\newpage

\hbox to\hsize{\hss\vbox{\offinterlineskip
\halign to120mm{\strut\tabskip=100pt minus 100pt
\strut\vrule#&\hbox to6.5mm{\hss$#$\hss}&%
\vrule#&\hbox to113mm{\hfil$\dsize#$\hfil}&%
\vrule#\tabskip=0pt\cr\noalign{\hrule}
& \# &%
& \text{differential operator $D$ and coefficients $A_n$, $n=0,1,2,\dots$} &\cr
\noalign{\hrule\vskip1pt\hrule}
& \eqnnol{221}{292} &&
\aligned
\\[-12pt]
&
D=5^2\theta^4
-5\cdot 2^2z(404\theta^4+1096\theta^3+773\theta^2+225\theta+25)
\\[-2pt] &\;
-2^4z^2(66896\theta^4+137408\theta^3+1010960\theta^2+52800\theta+11625)
\\[-2pt] &\;
-2^8\cdot3\cdot5z^3(2\theta+1)(5672\theta^3+9500\theta^2+8422\theta+2689)
\\[-2pt] &\;
-2^{15}\cdot3^2z^4(2\theta+1)(1208\theta^3+2892\theta^2+2842\theta+969)
\\[-2pt] &\;
-2^{20}\cdot3^3z^5(2\theta+1)(2\theta+3)(6\theta+5)(6\theta+7)
\\[-8pt]
\endaligned &\cr
& &\span\hrulefill&\cr
& &&
\aligned
A_n
&=\binom{2n}n^2\biggl((-1)^n\sum_{k=0}^{[n/2]}\binom nk^2
\binom{n+k}n\binom n{2k}\binom{2n}{4k}^{-1}
\\[-2pt] &\quad
+\sum_{k=[n/2]+1}^n\binom nk^2\binom{n+k}n
\binom{4k}{2n}\binom{2k}n^{-1}\biggr)
\\[2pt]
\endaligned &\cr
\noalign{\hrule}
& \eqnnol{222}{279} &&
\aligned
\\[-12pt]
&
D=5^2\theta^4
-5z(407\theta^4+1198\theta^3+909\theta^2+310\theta+40)
\\[-2pt] &\;
-2^7z^2(2103\theta^4+6999\theta^3+8358\theta^2+4050\theta+680)
\\[-2pt] &\;
-2^{12}z^3(1387\theta^4+3840\theta^3+3081\theta^2+960\theta+100)
-2^{21}z^4(2\theta+1)^4
\\[-8pt]
\endaligned &\cr
& &\span\hrulefill&\cr
& &&
A_n=\binom{2n}n^2\sum_k\binom nk\binom{2k}n\binom{2n}{n-k}
&\cr
\noalign{\hrule}
& \eqnnol{223}{280} &&
\aligned
\\[-12pt]
&
D=\theta^4+6z\theta(48\theta^3-12\theta^2-7\theta-1)
\\[-2pt] &\;
+2^23^2z^2(392\theta^4+488\theta^3+775\theta^2+376\theta+64)
\\[-2pt] &\;
+2^43^5z^3(1184\theta^4+3288\theta^3+3512\theta^2+1635\theta+278)
\\[-2pt] &\;
+2^63^8z^4(2\theta+1)^2(169\theta^2+361\theta+238)
+2^{11}3^{11}z^5(2\theta+1)^2(2\theta+3)^2
\\[-8pt]
\endaligned &\cr
& &\span\hrulefill&\cr
& &&
A_n=\binom{2n}n\sum_k(-1)^k3^{2n-3k}\binom{2n}{3k}\binom{2k}n
\frac{(3k)!}{k!^3}
&\cr
\noalign{\hrule}
& \eqnnol{224}{281} &&
\aligned
\\[-12pt]
&
D=5^2\theta^4
-5z(1057\theta^4+1058\theta^3+819\theta^2+290\theta+40)
\\[-2pt] &\;
+2^5z^2(10123\theta^4+11419\theta^3+5838\theta^2+1510\theta+180)
\\[-2pt] &\;
-2^8z^3(3098\theta^4+46560\theta^3+48211\theta^2+25500\theta+5100)
\\[-2pt] &\;
+2^{14}\cdot11z^4(2\theta+1)(234\theta^3+591\theta^2+581\theta+202)
\\[-2pt] &\;
-2^{20}\cdot11^2z^5(\theta+1)^2(2\theta+1)(2\theta+3)
\\[-8pt]
\endaligned &\cr
& &\span\hrulefill&\cr
& &&
A_n=2^{-n}\binom{2n}n^2\sum_k(-1)^{n+k}\binom nk\binom{2k}k
\binom{2n-2k}{n-k}\binom{2k}n\binom{2n}{2k}^{-1}
&\cr
\noalign{\hrule}
& \eqnnol{225}{282} &&
\aligned
\\[-12pt]
&
D=\theta^4
+2^4z(22192\theta^4-17056\theta^3-9576\theta^2-1048\theta-49)
\\[-2pt] &\;
+2^{20}z^2(33648\theta^4-44688\theta^3+16224\theta^2+1764\theta+17)
\\[-2pt] &\;
+5\cdot2^{34}z^3(6512\theta^4-6144\theta^3-4440\theta^2-1536\theta-193)
\\[-2pt] &\;
-5^2\cdot2^{55}z^4(2\theta+1)^4
\\[-8pt]
\endaligned &\cr
& &\span\hrulefill&\cr
& &&
\gathered
\text{the reflection of \#222 at infinity}
\\[-8pt]
\endgathered &\cr
& &\span\hrulefill&\cr
& &&
\aligned
A_n
&=\binom{2n}n^2\sum_k\binom nk\binom{n+2k}k
\binom{3n-2k}{n-k}\binom{2n+4k}{n+2k}\binom{6n-4k}{3n-2k}
\\[-2pt] &\;\times
\bigl(1+k(-2H_k+2H_{n-k}-H_{n+k}+H_{2n-k}-2H_{n+2k}
\\[-2pt] &\;\quad
+2H_{3n-2k}+4H_{2n+4k}-4H_{6n-4k})\bigr)
\\[2pt]
\endaligned &\cr
\noalign{\hrule}
}}\hss}

\newpage

\hbox to\hsize{\hss\vbox{\offinterlineskip
\halign to120mm{\strut\tabskip=100pt minus 100pt
\strut\vrule#&\hbox to6.5mm{\hss$#$\hss}&%
\vrule#&\hbox to113mm{\hfil$\dsize#$\hfil}&%
\vrule#\tabskip=0pt\cr\noalign{\hrule}
& \# &%
& \text{differential operator $D$ and coefficients $A_n$, $n=0,1,2,\dots$} &\cr
\noalign{\hrule\vskip1pt\hrule}
& \eqnnol{226}{283} &&
\aligned
\\[-12pt]
&
D=5^2\theta^4
-2\cdot5z(328\theta^4+692\theta^3+551\theta^2+205\theta+30)
\\[-2pt] &\;
+2^2\cdot3z^2(5352\theta^4+25416\theta^3+38387\theta^2+23020\theta+4860)
\\[-2pt] &\;
-2^4\cdot3^3z^3(352\theta^4+4520\theta^3+12108\theta^2+10205\theta+2630)
\\[-2pt] &\;
-2^6\cdot3^3z^4(2\theta+1)(586\theta^4+3039\theta^2+3947\theta+1527)
\\[-2pt] &\;
-2^8\cdot3^4z^5(2\theta+1)(2\theta+3)(6\theta+5)(6\theta+7)
\\[-8pt]
\endaligned &\cr
& &\span\hrulefill&\cr
& &&
A_n=\binom{2n}n\sum_k(-1)^{n+k}\binom nk\binom{2k}k
\binom{2n-2k}{n-k}\binom{3k}n
&\cr
\noalign{\hrule}
& \eqnnol{227}{284} &&
\aligned
\\[-12pt]
&
D=\theta^4
-2^2\cdot3^2z(132\theta^4+264\theta^3+201\theta^2+69\theta+10)
\\[-2pt] &\;
+2^9\cdot3^6z^2(20\theta^4+80\theta^3+107\theta^2+54\theta+10)
\\[-2pt] &\;
+2^{12}\cdot3^{10}z^3(2\theta+1)^2(2\theta+5)^2
\\[-8pt]
\endaligned &\cr
& &\span\hrulefill&\cr
& &&
A_n=432^n\binom{2n}n\sum_k(-1)^k\binom nk\binom{3k}n\binom{-1/6}k\binom{-5/6}k
&\cr
\noalign{\hrule}
& \eqnnol{228}{285} &&
\aligned
\\[-12pt]
&
D=\theta^4
-2^2z(176\theta^4+352\theta^3+289\theta^2+113\theta+18)
\\[-2pt] &\;
+2^{11}z^2(80\theta^4+320\theta^3+449\theta^2+258\theta+54)
\\[-2pt] &\;
-3\cdot2^{16}z^3(2\theta+1)(2\theta+5)(4\theta+3)(4\theta+9)
\\[-8pt]
\endaligned &\cr
& &\span\hrulefill&\cr
& &&
A_n=64^n\binom{2n}n\sum_k(-1)^{n+k}\binom nk\binom{3k}n
\binom{-1/4}k\binom{-3/4}k
&\cr
\noalign{\hrule}
& \eqnnol{229}{293} &&
\aligned
\\[-12pt]
&
D=\theta^4
-2^2z(256\theta^4+728\theta^3+506\theta^2+142\theta+15)
\\[-2pt] &\;
+2^4\cdot3^2z^2(-2336\theta^4-2336\theta^3+1768\theta^2+1176\theta+189)
\\[-2pt] &\;
+2^9\cdot3^4z^3(-512\theta^4+432\theta^3+404\theta^2+108\theta+9)
\\[-2pt] &\;
+2^{12}\cdot3^8z^4(2\theta+1)^4
\\[-8pt]
\endaligned &\cr
& &\span\hrulefill&\cr
& &&
A_n=\binom{2n}n^2
\sum_{i+j+k=2n}\biggl(\frac{(2n)!}{i!j!k!}\biggr)^2
&\cr
\noalign{\hrule}
& \eqnnol{230}{294} &&
\aligned
\\[-12pt]
&
D=\theta^4
+3z(945\theta^4-162\theta^3-49\theta^2+32\theta+8)
\\[-2pt] &\;
+2\cdot3^2z^2(17928\theta^4+2970\theta^3+10187\theta^2+3376\theta+408)
\\[-2pt] &\;
+2^2\cdot3^7z^3(156285\theta^4+200016\theta^3+19630\theta^2+84378\theta+13964)
\\[-2pt] &\;
+2^4\cdot3^{10}\cdot19z^4(2\theta+1)^2(4743\theta^2+8199\theta+4922)
\\[-2pt] &\;
+2^9\cdot3^{15}\cdot19^2z^5(2\theta+1)^2(2\theta+3)^2
\\[-8pt]
\endaligned &\cr
& &\span\hrulefill&\cr
& &&
A_n=27^n\binom{2n}n^2\sum_k(-1)^k\binom nk\binom{2k}n\binom{n+k}n^{-1}
\binom{-1/3}k\binom{-2/3}k
&\cr
\noalign{\hrule}
& \eqnnol{231}{295} &&
\aligned
\\[-12pt]
&
D=9\theta^4
-2^2z(84\theta^4+3048\theta^3+2217\theta^2+693\theta+90)
\\[-2pt] &\;
+2^9z^2(-1168\theta^4+968\theta^3+9518\theta^2+5325\theta+1005)
\\[-2pt] &\;
+2^{16}z^3(988\theta^4+8208\theta^3-743\theta^2-4230\theta-1245)
\\[-2pt] &\;
+2^{24}\cdot5z^4(2\theta+1)^2(9\theta^2-279\theta-277)
\\[-2pt] &\;
-2^{33}\cdot5^2z^5(2\theta+1)^2(2\theta+3)^2
\\[-8pt]
\endaligned &\cr
& &\span\hrulefill&\cr
& &&
\aligned
&
A_n=\binom{2n}n^3\biggl(\sum_{k=0}^{[n/2]}(-1)^k\binom nk
\binom{2k}k\binom{2n-2k}{n-k}\binom n{2k}
\binom{2n}{4k}^{-1}\binom{n+k}n^{-1}
\\[-2pt] &\quad
+\sum_{k=[n/2]+1}^n(-1)^{n+k}\binom nk\binom{2k}k\binom{2n-2k}{n-k}
\binom{4k}{2n}\binom{2k}n^{-1}\binom{n+k}n^{-1}\biggr)
\\[2pt]
\endaligned &\cr
\noalign{\hrule}
}}\hss}

\newpage

\hbox to\hsize{\hss\vbox{\offinterlineskip
\halign to120mm{\strut\tabskip=100pt minus 100pt
\strut\vrule#&\hbox to6.5mm{\hss$#$\hss}&%
\vrule#&\hbox to113mm{\hfil$\dsize#$\hfil}&%
\vrule#\tabskip=0pt\cr\noalign{\hrule}
& \# &%
& \text{differential operator $D$ and coefficients $A_n$, $n=0,1,2,\dots$} &\cr
\noalign{\hrule\vskip1pt\hrule}
& \eqnnol{232}{296} &&
\aligned
\\[-12pt]
&
D=5^2\theta^4
-5z(2617\theta^4+4658\theta^3+3379\theta^2+1050\theta+120)
\\[-2pt] &\;
-2^6\cdot3z^2(-673\theta^4+4871\theta^3+10282\theta^2+5410\theta+860)
\\[-2pt] &\;
+2^{10}\cdot3^2z^3(955\theta^4+4320\theta^3+3477\theta^2+1020\theta+100)
\\[-2pt] &\;
-2^{17}\cdot3^3z^4(2\theta+1)^2(3\theta+1)(3\theta+2)
\\[-8pt]
\endaligned &\cr
& &\span\hrulefill&\cr
& &&
A_n=\binom{2n}n^2\sum_k\binom nk^2\binom{3n}{n+k}
&\cr
\noalign{\hrule}
& \eqnnol{233}{297} &&
\aligned
\\[-12pt]
&
D=\theta^4
-2^4z(83\theta^4+94\theta^3+71\theta^2+24\theta+3)
\\[-2pt] &\;
+3\cdot 2^{11}z^2(101\theta^4+191\theta^3+174\theta^2+71\theta+10)
\\[-2pt] &\;
-3^22^{16}z^3(203\theta^4+432\theta^3+333\theta^2+102\theta+11)
\\[-2pt] &\;
+3^32^{23}z^4(2\theta+1)^2(3\theta+1)(3\theta+2)
\\[-8pt]
\endaligned &\cr
& &\span\hrulefill&\cr
& &&
A_n=\binom{2n}n^3\sum_k\binom nk^2\binom{3n}{n+k}\binom{2n}{2k}^{-1}
&\cr
\noalign{\hrule}
& \eqnnol{234}{298} &&
\aligned
\\[-12pt]
&
D=7^2\theta^4
-7\cdot 2z\theta(192\theta^3+60\theta^2+37\theta+7)
\\[-2pt] &\;
-2^2z^2(17608\theta^4+115144\theta^3+166715\theta^2+94556\theta+18816)
\\[-2pt] &\;
+2^4\cdot3^2z^3(20288\theta^4+57288\theta^3+27524\theta^2-7455\theta-5026)
\\[-2pt] &\;
-2^6\cdot3^5z^4(2\theta+1)(458\theta^3-657\theta^2-1799\theta-846)
\\[-2pt] &\;
-2^{12}\cdot3^8z^5(\theta+1)^2(2\theta+1)(2\theta+3)
\\[-8pt]
\endaligned &\cr
& &\span\hrulefill&\cr
& &&
A_n=\binom{2n}n\sum_k\binom nk^2\binom{2k}{n-k}\binom{2n-2k}k
&\cr
\noalign{\hrule}
& \eqnnol{235}{299} &&
\aligned
\\[-12pt]
&
D=7^2\theta^4
-2\cdot7z\theta(46\theta^3+52\theta^2+33\theta+7)
\\[-2pt] &\;
-2^2z^2(7332\theta^4+28848\theta^3+42633\theta^2+26670\theta+6272)
\\[-2pt] &\;
-2^4z^3(2860\theta^4+44760\theta^3+120483\theta^2+111279\theta+35098)
\\[-2pt] &\;
+2^9z^4(2230\theta^4+5920\theta^3-741\theta^2-6509\theta-3049)
\\[-2pt] &\;
+2^{14}z^5(174\theta^4+1320\theta^3+1971\theta^2+1095\theta+190)
\\[-2pt] &\;
+2^{19}z^6(-22\theta^4-24\theta^3+9\theta^2+21\theta+7)
\\[-2pt] &\;
-2^{25}z^7(\theta+1)^4
\\[-8pt]
\endaligned &\cr
& &\span\hrulefill&\cr
& &&
A_n=\sum_k\binom nk\binom{2k}k\binom{2n-2k}{n-k}
\binom{2k}{n-k}\binom{2n-2k}k
&\cr
\noalign{\hrule}
& \eqnnol{236}{300} &&
\aligned
\\[-12pt]
&
D=\theta^4
+2^4z(22\theta^4+64\theta^3+51\theta^2+19\theta+3)
\\[-2pt] &\;
+2^9z^2(-174\theta^4+624\theta^3+945\theta^2+417\theta+80)
\\[-2pt] &\;
+2^{14}z^3(-2230\theta^4-3000\theta^3+5121\theta^2+3813\theta+971)
\\[-2pt] &\;
+2^{19}z^4(2860\theta^4-33320\theta^3+3363\theta^2+6847\theta+2402)
\\[-2pt] &\;
+2^{27}z^5(7332\theta^4+480\theta^3+81\theta^2+1380\theta+719)
\\[-2pt] &\;
+2^{36}\cdot7z^6(\theta+1)(46\theta^3+86\theta^2+67\theta+20)
\\[-2pt] &\;
+2^{45}\cdot7^2z^7(\theta+1)^4
\\[-8pt]
\endaligned &\cr
& &\span\hrulefill&\cr
& &&
\gathered
\text{the reflection of \#235 at infinity}
\\[-8pt]
\endgathered &\cr
& &\span\hrulefill&\cr
& &&
\text{a formula for $A_n$ is not known}
&\cr
\noalign{\hrule}
}}\hss}

\newpage

\hbox to\hsize{\hss\vbox{\offinterlineskip
\halign to120mm{\strut\tabskip=100pt minus 100pt
\strut\vrule#&\hbox to6.5mm{\hss$#$\hss}&%
\vrule#&\hbox to113mm{\hfil$\dsize#$\hfil}&%
\vrule#\tabskip=0pt\cr\noalign{\hrule}
& \# &%
& \text{differential operator $D$ and coefficients $A_n$, $n=0,1,2,\dots$} &\cr
\noalign{\hrule\vskip1pt\hrule}
& \eqnnol{237}{301} &&
\aligned
\\[-12pt]
&
D=\theta^4
-2^4z(46\theta^4+128\theta^3+91\theta^2+27\theta+3)
\\[-2pt] &\;
+2^9\cdot3z^2(-74\theta^4+16\theta^3+231\theta^2+127\theta+20)
\\[-2pt] &\;
+2^{14}\cdot3^2z^3(14\theta^4+216\theta^3+175\theta^2+51\theta+5)
\\[-2pt] &\;
+2^{19}\cdot3^3z^4(2\theta+1)^2(3\theta+1)(3\theta+2)
\\[-8pt]
\endaligned &\cr
& &\span\hrulefill&\cr
& &&
A_n=\binom{2n}n^2\sum_k\binom nk^3\binom{2n+2k}{n+k}
\binom{4n-2k}{2n-k}\binom{2n}k^{-1}\binom{2n}{n-k}^{-1}
&\cr
\noalign{\hrule}
& \eqnnol{238}{302} &&
\aligned
\\[-12pt]
&
D=\theta^4
+2^2z(500\theta^4+976\theta^3+677\theta^2+189\theta+19)
\\[-2pt] &\;
+2^4z^2(3968\theta^4+3968\theta^3-1336\theta^2-1164\theta-177)
\\[-2pt] &\;
+2^{10}z^3(500\theta^4+24\theta^3-37\theta^2+6\theta+3)
+2^{12}z^4(2\theta+1)^4
\\[-8pt]
\endaligned &\cr
& &\span\hrulefill&\cr
& &&
A_n=\binom{2n}n\sum_k\binom nk\binom{n+k}n
\binom{2n+2k}{n+k}\binom{2n+k}{2n-k}
&\cr
\noalign{\hrule}
& \eqnnoltmp{239}{303} &&
\aligned
\\[-12pt]
&
D=\theta^4
-2^4\cdot3z(-9\theta^4+198\theta^3+131\theta^2+32\theta+39)
\\[-2pt] &\;
-2^{11}\cdot3^2z^2(486\theta^4+1215\theta^3+81\theta^2-27\theta-5)
\\[-2pt] &\;
-2^{16}\cdot3^5z^3(891\theta^4+972\theta^3+675\theta^2+216\theta+25)
\\[-2pt] &\;
-2^{23}\cdot3^8z^4(3\theta+1)^2(3\theta+2)^2
\\[-8pt]
\endaligned &\cr
& &\span\hrulefill&\cr
& &&
A_n=\frac{(3n)!}{n!^3}\sum_k\binom nk\binom{2n+2k}{n+k}\binom{4n-2k}{2n-k}
&\cr
\noalign{\hrule}
& \eqnnoltmp{240}{304} &&
\aligned
\\[-12pt]
&
D=13^2\theta^4
-13z(1449\theta^4+4050\theta^3+3143\theta^2+1118\theta+156)
\\[-2pt] &\;
+2^4z^2(-22760\theta^4+27112\theta^3+121046\theta^2+82316\theta+17589)
\\[-2pt] &\;
+2^8z^3(3824\theta^4+39936\theta^3-34292\theta^2-63492\theta-19539)
\\[-2pt] &\;
-2^{16}\cdot3z^4(2\theta+1)(40\theta^3+684\theta^2+1013\theta+399)
\\[-2pt] &\;
-2^{20}\cdot3^2z^5(2\theta+1)(2\theta+3)(4\theta+3)(4\theta+5)
\\[-8pt]
\endaligned &\cr
& &\span\hrulefill&\cr
& &&
A_n=\sum_k\binom nk\binom{2k}k\binom{2n-2k}{n-k}\binom{n+2k}n\binom{3n-2k}n
&\cr
\noalign{\hrule}
& \eqnnoltmp{241}{305} &&
\aligned
\\[-12pt]
&
D=\theta^4
-2^4z(152\theta^4+160\theta^3+110\theta^2+30\theta+3)
\\[-2pt] &\;
+2^{10}\cdot3z^2(428\theta^4+176\theta^3-299\theta^2-170\theta-25)
\\[-2pt] &\;
+2^{17}\cdot3^2z^3(-136\theta^4+216\theta^3+180\theta^2+51\theta+5)
\\[-2pt] &\;
-2^{24}\cdot3^3z^4(2\theta+1)^2(3\theta+1)(3\theta+2)
\\[-8pt]
\endaligned &\cr
& &\span\hrulefill&\cr
& &&
A_n=\binom{2n}n\sum_k\binom{2n+2k}{n+k}\binom{4n-2k}{2n-k}
\binom{n+k}{n-k}\binom{2n-k}k
&\cr
\noalign{\hrule}
}}\hss}

\newpage

\hbox to\hsize{\hss\vbox{\offinterlineskip
\halign to120mm{\strut\tabskip=100pt minus 100pt
\strut\vrule#&\hbox to6.5mm{\hss$#$\hss}&%
\vrule#&\hbox to113mm{\hfil$\dsize#$\hfil}&%
\vrule#\tabskip=0pt\cr\noalign{\hrule}
& \# &%
& \text{differential operator $D$ and coefficients $A_n$, $n=0,1,2,\dots$} &\cr
\noalign{\hrule\vskip1pt\hrule}
& \eqnnoltmp{242}{306} &&
\aligned
\\[-12pt]
&
D=\theta^4
+2\cdot3z(72\theta^4+108\theta^3+91\theta^2+37\theta+6)
\\[-2pt] &\;
+2^2\cdot3^3z^2(648\theta^2+1800\theta^3+2211\theta^2+1248\theta+260)
\\[-2pt] &\;
+2^4\cdot3^5z^3(1344\theta^4+4968\theta^3+7320\theta^2+4749\theta+1072)
\\[-2pt] &\;
+2^6\cdot3^7z^4(2\theta+1)(630\theta^3+2241\theta^2+2617\theta+971)
\\[-2pt] &\;
+2^8\cdot3^{10}z^5(2\theta+1)(2\theta+3)(6\theta+5)(6\theta+7)
\\[-8pt]
\endaligned &\cr
& &\span\hrulefill&\cr
& &&
\aligned
A_n
&=\binom{2n}n\sum_{k=[n/3]}^{[2n/3]}\binom nk^2
\binom{2k}k\binom{2n-2k}{n-k}\binom{3k}n\binom{3n-3k}n
\\[-2pt] &\;\quad\times
\bigl(1+k(-4H_k+4H_{n-k}+2H_{2k}-2H_{2n-2k}
\\[-2pt] &\;\qquad
+3H_{3k}-3H_{3n-3k}+3H_{2n-3k}-3H_{3k-n})\bigr)
\\[-2pt] &\;
+3\binom{2n}n\sum_{k=0}^{[(n-1)/3]}(-1)^{n+k}\frac{n-2k}{n-3k}\binom nk^2
\\[-2pt] &\;\quad\times
\binom{2k}k\binom{2n-2k}{n-k}\binom{3n-3k}n\binom n{3k}^{-1}
\\[2pt]
\endaligned &\cr
\noalign{\hrule}
& \eqnnol{243}{222} &&
\aligned
\\[-12pt]
&
D=\theta^4
+z(295\theta^4+572\theta^3+424\theta^2+138\theta+17)
\\[-2pt] &\;
+2z^2(843\theta^4+744\theta^2-473\theta^2-481\theta-101)
\\[-2pt] &\;
+2z^3(1129\theta^4-516\theta^3-725\theta^2-159\theta+4)
\\[-2pt] &\;
-3z^4(173\theta^4+352\theta^3+290\theta^2+114\theta+18)
-3^2z^5(\theta+1)^4
\\[-8pt]
\endaligned &\cr
& &\span\hrulefill&\cr
& &&
\gathered
\text{the reflection of \#27 at infinity}
\\[-8pt]
\endgathered &\cr
& &\span\hrulefill&\cr
& &&
A_n=(-1)^n\sum_{k,l}\binom nk\binom nl\binom{n+k}n\binom{n+l}n
\binom{n+k+l}n\binom n{l-k}
&\cr
\noalign{\hrule}
& \eqnnoltmp{244}{307} &&
\aligned
\\[-12pt]
&
D=\theta^4
+z(416\theta^4+104\theta^3+36\theta^2-16\theta-11)
\\[-2pt] &\;
+z^2(63168\theta^4+31584\theta^3+12158\theta^2+4828\theta+6973)
\\[-2pt] &\;
+2^2z^3(990080\theta^4+742560\theta^3+311828\theta^2+192658\theta-33281)
\\[-2pt] &\;
+z^4(62035456\theta^4+62035456\theta^3+26395808\theta^2
+5661536\theta+7849233)
\\[-2pt] &\;
-2^4\cdot3z^5(35642880\theta^4+44553600\theta^3+24468128\theta^2
\\[-2pt] &\;\quad
+1916112\theta+4807463)
\\[-2pt] &\;
+2^5\cdot3^2z^6(40932864\theta^4+61399296\theta^3+38293776\theta^2
\\[-2pt] &\;\quad
+9496512\theta+4807463)
\\[-2pt] &\;
-2^8\cdot3^5z^7(539136\theta^4+943488\theta^3+647568\theta^2+229320\theta+23377)
\\[-2pt] &\;
+2^8\cdot3^8z^8(12\theta+5)^2(12\theta+7)^2
\\[-8pt]
\endaligned &\cr
& &\span\hrulefill&\cr
& &&
\gathered
\text{the pullback of the 5th-order differential equation $D'y=0$, where}
\\[-8pt]
\endgathered &\cr
& &\span\hrulefill&\cr
& &&
\aligned
\\[-12pt]
&
D'=\theta^5
+2z(2\theta+1)(26\theta^4+52\theta^3+44\theta^2+18\theta+3)
\\[-2pt] &\;
-12z^2(\theta+1)^3(6\theta+5)(6\theta+7)
\\[-8pt]
\endaligned &\cr
& &\span\hrulefill&\cr
& &&
\aligned
A_n'
&=\sum_k\binom nk^6\binom{2k}k\binom{2n-2k}{n-k}
\\[-2pt] &\quad\times
\bigl(1+k(-8H_k+8H_{n-k}+2H_{2k}-2H_{2n-2k})\bigr)
\\[2pt]
\endaligned &\cr
\noalign{\hrule}
}}\hss}

\newpage

\hbox to\hsize{\hss\vbox{\offinterlineskip
\halign to120mm{\strut\tabskip=100pt minus 100pt
\strut\vrule#&\hbox to6.5mm{\hss$#$\hss}&%
\vrule#&\hbox to113mm{\hfil$\dsize#$\hfil}&%
\vrule#\tabskip=0pt\cr\noalign{\hrule}
& \# &%
& \text{differential operator $D$ and coefficients $A_n$, $n=0,1,2,\dots$} &\cr
\noalign{\hrule\vskip1pt\hrule}
& \eqnnoltmp{245}{308} &&
\aligned
\\[-12pt]
&
D=\theta^4
-3z(144\theta^4+72\theta^3+32\theta^2-4\theta-5)
\\[-2pt] &\;
+3^2z^2(7776\theta^4+7776\theta^3+4086\theta^2+828\theta+589)
\\[-2pt] &\;
-2^3\cdot3^7z^3(288\theta^4+432\theta^3+262\theta^2+81\theta+6)
\\[-2pt] &\;
+3^{12}z^4(4\theta+1)^2(4\theta+3)^2
\\[-8pt]
\endaligned &\cr
& &\span\hrulefill&\cr
& &&
\gathered
\text{the pullback of the 5th-order differential equation $D'y=0$, where}
\\[-8pt]
\endgathered &\cr
& &\span\hrulefill&\cr
& &&
\aligned
\\[-12pt]
&
D'=\theta^5
-6z(2\theta+1)(18\theta^4+36\theta^3+34\theta^2+16\theta+3)
\\[-2pt] &\;
+2916z^2(\theta+1)^3(2\theta+1)(2\theta+3)
\\[-8pt]
\endaligned &\cr
& &\span\hrulefill&\cr
& &&
A_n'=3\binom{2n}n^2\sum_{k=0}^{[n/3]}(-1)^k\frac{n-2k}{2n-3k}
\binom nk^4\binom{3n-3k}{2n}\binom{2n}{3k}^{-1}
&\cr
\noalign{\hrule}
& \eqnnoltmp{246}{309} &&
\aligned
\\[-12pt]
&
D=5^2\theta^4
-2^2\cdot5z(12\theta^4+48\theta^3+49\theta^2+25\theta+5)
\\[-2pt] &\;
-2^4z^2(544\theta^4+1792\theta^3+2444\theta^2+1580\theta+405)
\\[-2pt] &\;
+2^9z^3(112\theta^4+960\theta^3+2306\theta^2+2130\theta+685)
\\[-2pt] &\;
+2^{12}z^4(144\theta^4+768\theta^3+1308\theta^2+924\theta+235)
\\[-2pt] &\;
+2^{20}z^5(\theta+1)^4
\\[-8pt]
\endaligned &\cr
& &\span\hrulefill&\cr
& &&
\aligned
A_n
&=\sum_k\binom nk^3\binom{2k}k^2\binom{2n-2k}{n-k}^2
\\[-2pt] &\quad\times
\bigl(1+k(-7H_k+7H_{n-k}+4H_{2k}-4H_{2n-2k})\bigr)
\\[2pt]
\endaligned &\cr
\noalign{\hrule}
& \eqnnoltmp{247}{310} &&
\aligned
\\[-12pt]
&
D=\theta^4
-2^4z(-144\theta^4+192\theta^3+132\theta^2+36\theta+5)
\\[-2pt] &\;
+2^{13}z^2(112\theta^4-512\theta^3+98\theta^2+50\theta+13)
\\[-2pt] &\;
-2^{20}z^3(544\theta^4+384\theta^3+332\theta^2+108\theta+21)
\\[-2pt] &\;
+2^{30}\cdot5z^4(-12\theta^4+23\theta^3+23\theta+7)
\\[-2pt] &\;
+2^{40}\cdot5^2z^5(\theta+1)^4
\\[-8pt]
\endaligned &\cr
& &\span\hrulefill&\cr
& &&
\gathered
\text{the reflection of \#246 at infinity}
\\[-8pt]
\endgathered &\cr
& &\span\hrulefill&\cr
& &&
\aligned
A_n
&=4^{-n}\sum_k\binom nk^{-3}\binom{2k}k^5\binom{2n-2k}{n-k}^5
\\[-2pt] &\quad\times
\bigl(1+k(-7H_k+7H_{n-k}+10H_{2k}-10H_{2n-2k})\bigr)
\\[2pt]
\endaligned &\cr
\noalign{\hrule}
& \eqnnoltmp{248}{311} &&
\aligned
\\[-12pt]
&
D=3^2\theta^4
-3z(106\theta^4+146\theta^3+115\theta^2+42\theta+6)
\\[-2pt] &\;
-z^2(4511^4+24314\theta^3+37829\theta^2+23598\theta+5286)
\\[-2pt] &\;
+2^2z^3(10457\theta^4+32184\theta^3+24449\theta^2+3627\theta-1317)
\\[-2pt] &\;
-2^2\cdot11z^4(1596\theta^4+2040\theta^3-101\theta^2-1085\theta-386)
\\[-2pt] &\;
-2^4\cdot11^2z^5(\theta+1)^2(4\theta+3)(4\theta+5)
\\[-8pt]
\endaligned &\cr
& &\span\hrulefill&\cr
& &&
\aligned
A_n
&=(-1)^n\sum_k\binom nk^3\binom{n+k}n
\binom{2n-k}n\binom{2k}k\binom{2n-2k}{n-k}
\\[-2pt] &\quad\times
\bigl(1+k(-6H_k+6H_{n-k}+H_{n+k}-H_{2n-k}+2H_{2k}-2H_{2n-2k})\bigr)
\\[2pt]
\endaligned &\cr
\noalign{\hrule}
}}\hss}

\newpage

\hbox to\hsize{\hss\vbox{\offinterlineskip
\halign to120mm{\strut\tabskip=100pt minus 100pt
\strut\vrule#&\hbox to6.5mm{\hss$#$\hss}&%
\vrule#&\hbox to113mm{\hfil$\dsize#$\hfil}&%
\vrule#\tabskip=0pt\cr\noalign{\hrule}
& \# &%
& \text{differential operator $D$ and coefficients $A_n$, $n=0,1,2,\dots$} &\cr
\noalign{\hrule\vskip1pt\hrule}
& \eqnnoltmp{249}{312} &&
\aligned
\\[-12pt]
&
D=5^2\theta^4
+2^2\cdot5z(148\theta^4+392\theta^3+341\theta^2+145\theta+25)
\\[-2pt] &\;
+2^4z^2(4096\theta^4+32128\theta^3+57016\theta^2+37920\theta+9175)
\\[-2pt] &\;
+2^8z^3(-6656\theta^4-7680\theta^3+36960\theta^2+49920\theta+16985)
\\[-2pt] &\;
-2^{15}z^4(512\theta^4+4864\theta^3+9136^2+6464\theta+1587)
\\[-2pt] &\;
+2^{20}z^5(4\theta+3)^2(4\theta+5)^2
\\[-8pt]
\endaligned &\cr
& &\span\hrulefill&\cr
& &&
\aligned
A_n
&=(-1)^n\sum_k\binom nk\binom{n+k}n\binom{2n-k}n\binom{2k}k^2
\binom{2n-2k}{n-k}^2
\\[-2pt] &\quad\times
\bigl(1+k(-6H_k+6H_{n-k}+H_{n+k}-H_{2n-k}+4H_{2k}-4H_{2n-k})\bigr)
\\[2pt]
\endaligned &\cr
\noalign{\hrule}
& \eqnnoltmp{250}{313} &&
\aligned
\\[-12pt]
&
D=23^2\theta^4
-23z(3271\theta^4+5078\theta^3+3896\theta^2+1357\theta+184)
\\[-2pt] &\;
+z^2(1357863\theta^4+999924\theta^3-787393\theta^2-850862\theta-205712)
\\[-2pt] &\;
+2^3z^3(-775799\theta^4+272481\theta^3+218821\theta^2-176709\theta-100234)
\\[-2pt] &\;
+2^4\cdot61z^4(-1005\theta^4+15654\theta^3+36317\theta^2+27938\theta+7304)
\\[-2pt] &\;
-2^9\cdot61^2z^5(\theta+1)^2(4\theta+3)(4\theta+5)
\\[-8pt]
\endaligned &\cr
& &\span\hrulefill&\cr
& &&
\aligned
A_n
&=(-1)^n\sum_k\binom nk^3\binom{2k}k\binom{2n-2k}{n-k}
\binom{n+2k}n\binom{3n-2k}n
\\[-2pt] &\quad\times
\bigl(1+k(-5H_k+5H_{n-k}+2H_{n+2k}-2H_{3n-2k})\bigr)
\\[2pt]
\endaligned &\cr
\noalign{\hrule}
& \eqnnoltmp{251}{314} &&
\aligned
\\[-12pt]
&
D=\theta^4
-3z\theta(27\theta^3+18\theta^2+11\theta+2)
\\[-2pt] &\;
-2\cdot3^2z^2(39\theta^4+480\theta^3+474\theta^2+276\theta+64)
\\[-2pt] &\;
+2^3\cdot3^4z^3(348\theta^4+1152\theta^3+1759\theta^2+1110\theta+260)
\\[-2pt] &\;
-2^3\cdot3^5z^4(3420\theta^4+15912\theta^3+28437\theta^2+20544\theta+5296)
\\[-2pt] &\;
+2^4\cdot3^7z^5(1125\theta^4+12546\theta^3+31089\theta^2+26448\theta+7480)
\\[-2pt] &\;
+2^5\cdot3^9z^6(1395\theta^4+3240\theta^3-3378\theta^2-7146\theta-2696)
\\[-2pt] &\;
-2^7\cdot3^{11}z^7(351\theta^4+2646\theta^3+4767\theta^2+3309\theta+800)
\\[-2pt] &\;
-2^7\cdot3^{13}z^8(3\theta+2)(3\theta+4)(6\theta+5)(6\theta+7)
\\[-8pt]
\endaligned &\cr
& &\span\hrulefill&\cr
& &&
A_n=\sum_k\binom nk\binom{2k}k\binom{2n-2k}{n-k}\binom{3k}n\binom{3n-3k}n
&\cr
\noalign{\hrule}
& \eqnnoltmp{252}{315} &&
\aligned
\\[-12pt]
&
D=5^2\theta^4
-5z(-36\theta^4+636\theta^3+488\theta^2+170\theta+25)
\\[-2pt] &\;
+2^4z^2(-21301\theta^4-27148\theta^3+86889\theta^2+63110\theta+14975)
\\[-2pt] &\;
+2^8z^3(-19535\theta^4+294315\theta^3+126425\theta^2-54390\theta-35755)
\\[-2pt] &\;
-2^{10}\cdot59z^4(-10981\theta^4+29878\theta^3+89811\theta^2+70372\theta+17759)
\\[-2pt] &\;
-2^{15}\cdot3\cdot59^2z^5(3\theta+2)(3\theta+4)(4\theta+3)(4\theta+5)
\\[-8pt]
\endaligned &\cr
& &\span\hrulefill&\cr
& &&
\aligned
&
A_n=\sum_k\binom nk\binom{2k}k^2\binom{2n-2k}{n-k}^2
\binom{n+2k}n\binom{3n-2k}n
\\[-2pt] &\quad\times
\bigl(1+k(-5H_k+5H_{n-k}+2H_{2k}-2H_{2n-2k}+2H_{n+2k}-2H_{3n-2k})\bigr)
\\[2pt]
\endaligned &\cr
\noalign{\hrule}
}}\hss}

\newpage

\hbox to\hsize{\hss\vbox{\offinterlineskip
\halign to120mm{\strut\tabskip=100pt minus 100pt
\strut\vrule#&\hbox to6.5mm{\hss$#$\hss}&%
\vrule#&\hbox to113mm{\hfil$\dsize#$\hfil}&%
\vrule#\tabskip=0pt\cr\noalign{\hrule}
& \# &%
& \text{differential operator $D$ and coefficients $A_n$, $n=0,1,2,\dots$} &\cr
\noalign{\hrule\vskip1pt\hrule}
& \eqnnoltmp{253}{316} &&
\aligned
\\[-12pt]
&
D=\theta^4
-2^2z(64\theta^4+32\theta^3+15\theta^2-\theta-2)
\\[-2pt] &\;
+2^8z^2(96\theta^4+96\theta^3+53\theta^2+13\theta+8)
\\[-2pt] &\;
-2^{12}z^3(256\theta^4+384\theta^3+244\theta^2+84\theta+7)
\\[-2pt] &\;
+2^{18}z^4(2\theta+1)^2(4\theta+1)(4\theta+3)
\\[-8pt]
\endaligned &\cr
& &\span\hrulefill&\cr
& &&
\gathered
\text{the pullback of the 5th-order differential equation $D'y=0$, where}
\\[-8pt]
\endgathered &\cr
& &\span\hrulefill&\cr
& &&
\aligned
\\[-12pt]
&
D'=\theta^5
-4z(2\theta+1)(16\theta^4+32\theta^3+31\theta^2+15\theta+3)
\\[-2pt] &\;
+16z^2(\theta+1)(4\theta+3)^2(4\theta+5)^2
\\[-8pt]
\endaligned &\cr
& &\span\hrulefill&\cr
& &&
\aligned
A_n'
&=\sum_k\binom nk^2\binom{n+k}n\binom{2n-k}n
\binom{2k}k^2\binom{2n-2k}{n-k}^2
\\[-2pt] &\quad\times
\bigl(1+k(-7H_k+7H_{n-k}+H_{n+k}-H_{2n-k}+4H_{2k}-4H_{2n-2k})\bigr)
\\[2pt]
\endaligned &\cr
\noalign{\hrule}
& \eqnnoltmp{254}{317} &&
\aligned
\\[-12pt]
&
D=\theta^4
+2^4z(-2608\theta^4+544\theta^3+200\theta^2-72\theta-15)
\\[-2pt] &\;
+2^{15}\cdot3z^2(6128\theta^4-208\theta^3+2328\theta^2+452\theta+25)
\\[-2pt] &\;
-2^{24}\cdot3^2\cdot5z^3(4592\theta^4+3456\theta^3+2632\theta^2+816\theta+95)
\\[-2pt] &\;
+2^{38}\cdot3^3\cdot5^2z^4(2\theta+1)^2(3\theta+1)(3\theta+2)
\\[-8pt]
\endaligned &\cr
& &\span\hrulefill&\cr
& &&
\aligned
&
A_n=\sum_k\binom nk\binom{n+k}n\binom{2n-k}n
\binom{2n+2k}{n+k}^2\binom{4n-2k}{2n-k}^2
\\[-2pt] &\quad\times
\bigl(1+k(-2H_k+2H_{n-k}-3H_{n+k}+3H_{2n-k}+4H_{2n+2k}-4H_{4n-2k})\bigr)
\\[2pt]
\endaligned &\cr
\noalign{\hrule}
& \eqnnoltmp{255}{318} &&
\aligned
\\[-12pt]
&
D=\theta^4
+2^2z(256\theta^4+128\theta^3+77\theta^2+13\theta-2)
\\[-2pt] &\;
+2^7z^2(3072\theta^4+3072\theta^3+1960\theta^2+536\theta+141)
\\[-2pt] &\;
+2^{12}z^3(16384\theta^4+24576\theta^3+16576\theta^2+5184\theta+491)
\\[-2pt] &\;
+2^{22}z^4(4\theta+1)(4\theta+3)(8\theta+3)(8\theta+5)
\\[-8pt]
\endaligned &\cr
& &\span\hrulefill&\cr
& &&
\gathered
\text{the pullback of the 5th-order differential equation $D'y=0$, where}
\\[-8pt]
\endgathered &\cr
& &\span\hrulefill&\cr
& &&
\aligned
\\[-12pt]
&
D'=\theta^5
+4z(2\theta+1)(64\theta^4+128\theta^3+141\theta^2+77\theta+17)
\\[-2pt] &\;
+16z^2(\theta+1)(8\theta+5)(8\theta+7)(8\theta+9)(8\theta+11)
\\[-8pt]
\endaligned &\cr
& &\span\hrulefill&\cr
& &&
\aligned
&
A_n'=\sum_k\binom nk^2\binom{n+k}n\binom{2n-k}n
\frac{(4k)!}{k!^2(2k)!}\frac{(4n-4k)!}{(n-k)!^2(2n-2k)!}
\\[-2pt] &\quad\times
\bigl(1+k(-5H_k+5H_{n-k}+H_{n+k}-H_{2n-k}-2H_{2k}+2H_{2n-2k}+4H_{4k}-4H_{4n-4k})\bigr)
\\[2pt]
\endaligned &\cr
\noalign{\hrule}
& \eqnnoltmp{256}{319} &&
\aligned
\\[-12pt]
&
D=\theta^4
+2^5z(24\theta^4+42\theta^3+30\theta^2+9\theta+1)
\\[-2pt] &\;
+2^8z^2(164\theta^4+104\theta^3-144\theta^2-100\theta-17)
\\[-2pt] &\;
+2^{14}z^3(28\theta^4-48\theta^3-44\theta^2-12\theta-1)
-2^{18}z^4(2\theta+1)^4
\\[-8pt]
\endaligned &\cr
& &\span\hrulefill&\cr
& &&
\aligned
&
A_n=\sum_k\binom nk^3\binom{n+k}n\binom{2n-k}n
\binom{2n+2k}{n+k}\binom{4n-2k}{2n-k}
\\[-2pt] &\quad\times
\bigl(1+k(-4H_k+4H_{n-k}-H_{n+k}+H_{2n-k}+2H_{2n+2k}-2H_{4n-2k})\bigr)
\\[2pt]
\endaligned &\cr
\noalign{\hrule}
}}\hss}

\newpage

\hbox to\hsize{\hss\vbox{\offinterlineskip
\halign to120mm{\strut\tabskip=100pt minus 100pt
\strut\vrule#&\hbox to6.5mm{\hss$#$\hss}&%
\vrule#&\hbox to113mm{\hfil$\dsize#$\hfil}&%
\vrule#\tabskip=0pt\cr\noalign{\hrule}
& \# &%
& \text{differential operator $D$ and coefficients $A_n$, $n=0,1,2,\dots$} &\cr
\noalign{\hrule\vskip1pt\hrule}
& \eqnnoltmp{257}{320} &&
\aligned
\\[-12pt]
&
D=\theta^4
-2^4z(112\theta^4+416\theta^3+280\theta^2+72\theta+7)
\\[-2pt] &\;
+2^{12}z^2(-656\theta^4-896\theta^3+216\theta^2+160\theta+23)
\\[-2pt] &\;
-2^{23}z^3(96\theta^4+24\theta^3+12\theta^2+6\theta+1)
-2^{30}z^4(2\theta+1)^4
\\[-8pt]
\endaligned &\cr
& &\span\hrulefill&\cr
& &&
\gathered
\text{the reflection of \#256 at infinity}
\\[-8pt]
\endgathered &\cr
& &\span\hrulefill&\cr
& &&
\aligned
A_n
&=(-1)^n\binom{2n}n^2\sum_k\binom nk\binom{2n+2k}{n+k}
\binom{4n-2k}{2n-k}\binom{n+k}{n-k}\binom{2n-k}k
\\[-2pt] &\quad\times
\bigl(1+k(-2H_k+2H_{n-k}-H_{n+k}+H_{2n-k}-2H_{2k}
\\[-2pt] &\quad\qquad
+2H_{2n-2k}+2H_{2n+2k}-2H_{4n-2k})\bigr)
\\[2pt]
\endaligned &\cr
\noalign{\hrule}
& \eqnnoltmp{258}{321} &&
\aligned
\\[-12pt]
&
D=\theta^4
-2^4z(16\theta^4+224\theta^3+156\theta^2+44\theta+5)
\\[-2pt] &\;
+2^{14}z^2(-48\theta^4-48\theta^3+120\theta^2+66\theta+11)
\\[-2pt] &\;
+2^{22}z^3(-16\theta^4+192\theta^3+156\theta^2+48\theta+5)
+2^{32}z^4(2\theta+1)^4
\\[-8pt]
\endaligned &\cr
& &\span\hrulefill&\cr
& &&
\gathered
\text{the case is self-dual at infinity}
\\[-8pt]
\endgathered &\cr
& &\span\hrulefill&\cr
& &&
\aligned
A_n
&=\sum_k\binom nk\binom{n+k}n\binom{2n-k}n\binom{2k}k
\binom{2n-2k}{n-k}\binom{2n+2k}{n+k}\binom{4n-2k}{2n-k}
\\[-2pt] &\quad\times
\bigl(1+k(-4H_k+4H_{n-k}-H_{n+k}+H_{2n-k}+2H_{2k}
\\[-2pt] &\quad\qquad
-2H_{2n-2k}+2H_{2n+2k}-2H_{4n-2k})\bigr)
\\[2pt]
\endaligned &\cr
\noalign{\hrule}
& \eqnnoltmp{259}{322} &&
\aligned
\\[-12pt]
&
D=\theta^4+10z(40000\theta^4-17500\theta^3-8125\theta^2+625\theta+238)
\\[-2pt] &\;
+2^2\cdot5^6z^2(835000\theta^4-365000\theta^3+371125\theta^2+58500\theta+2116)
\\[-2pt] &\;
+2^4\cdot5^{11}z^3(3130000\theta^4+1815000\theta^3+1662000\theta^2+625875\theta+96914)
\\[-2pt] &\;
+2^6\cdot5^{19}\cdot13z^4(2\theta+1)^2(625\theta^2+745\theta+351)
\\[-2pt] &\;
+2^8\cdot5^{25}\cdot13^2z^5(2\theta+1)^2(2\theta+3)^2
\\[-8pt]
\endaligned &\cr
& &\span\hrulefill&\cr
& &&
\aligned
A_n
&=(-1)^n\sum_k\binom{2n}k\binom{2n}{n-k}\frac{(5k)!\,(5n-5k)!}{n!^3k!^2(n-k)!^2}
\\[-2pt] &\quad\times
\bigl(1+k(-3H_k+3H_{n-k}-H_{n+k}+H_{2n-k}+5H_{5k}-5H_{5n-5k})\bigr)
\\[2pt]
\endaligned &\cr
\noalign{\hrule}
}}\hss}

\newpage

\hbox to\hsize{\hss\vbox{\offinterlineskip
\halign to120mm{\strut\tabskip=100pt minus 100pt
\strut\vrule#&\hbox to6.5mm{\hss$#$\hss}&%
\vrule#&\hbox to113mm{\hfil$\dsize#$\hfil}&%
\vrule#\tabskip=0pt\cr\noalign{\hrule}
& \# &%
& \text{differential operator $D$ and coefficients $A_n$, $n=0,1,2,\dots$} &\cr
\noalign{\hrule\vskip1pt\hrule}
& \eqnnoltmp{260}{323} &&
\aligned
\\[-12pt]
&
D=5^2\theta^4
+2^2\cdot5z(596\theta^4+544\theta^3+397\theta^2+125\theta+15)
\\[-2pt] &\;
+2^4\cdot3z^2(30048\theta^4+14784\theta^3-13312\theta^2-10940\theta-2115)
\\[-2pt] &\;
+2^8\cdot3^3z^3(6368\theta^4-6720\theta^3-9052\theta^2-4080\theta-655)
\\[-2pt] &\;
-2^{12}\cdot3^6z^4(2\theta+1)^2(76\theta^2+196\theta+139)
\\[-2pt] &\;
-2^{16}\cdot3^9z^5(2\theta+1)^2(2\theta+3)^2
\\[-8pt]
\endaligned &\cr
& &\span\hrulefill&\cr
& &&
\aligned
&
A_n=\binom{2n}n^2\sum_{k=[n/3]}^{[2n/3]}\binom nk^3
\binom{2k}{n-k}\binom{2n-2k}k
\\[-2pt] &\quad\times
\bigl(1+k(-4H_k+4H_{n-k}+2H_{2k}-2H_{2n-2k}+3H_{2n-3k}-3H_{3k-n})\bigr)
\\[-2pt] &\;
+3\binom{2n}n^2\sum_{k=0}^{[(n-1)/3]}(-1)^{n+k}\frac{n-2k}{n-3k}
\binom nk^3\binom{2n-2k}k\binom{n-k}{2k}^{-1}
\\[2pt]
\endaligned &\cr
\noalign{\hrule}
& \eqnnoltmp{261}{324} &&
\aligned
\\[-12pt]
&
D=5^2\theta^4
+2^2\cdot5z(292\theta^4+368\theta^3+289\theta^2+105\theta+15)
\\[-2pt] &\;
+2^4z^2(24736\theta^4+43648\theta^3+38936\theta^2+18980\theta+3735)
\\[-2pt] &\;
+2^9\cdot3^2z^3(2512\theta^4+5760\theta^3+6328\theta^2+3330\theta+655)
\\[-2pt] &\;
+2^{12}\cdot3^4z^4(2\theta+1)(232\theta^3+588\theta^2+590\theta+207)
\\[-2pt] &\;
+2^{18}\cdot3^6z^5(\theta+1)^2(2\theta+1)(2\theta+3)
\\[-8pt]
\endaligned &\cr
& &\span\hrulefill&\cr
& &&
\aligned
&
A_n=\binom{2n}n\sum_{k=[n/3]}^{[2n/3]}\binom nk^2\binom{2k}k
\binom{2n-2k}{n-k}\binom{2k}{n-k}\binom{2n-2k}k
\\[-2pt] &\;\quad\times
\bigl(1+k(-5H_k+5H_{n-k}+4H_{2k}-4H_{2n-2k}+3H_{2n-3k}-3H_{3k-n})\bigr)
\\[-2pt] &\;
+3\binom{2n}n\sum_{k=0}^{[(n-1)/3]}(-1)^{n+k}\frac{n-2k}{n-3k}\binom nk^2
\\[-2pt] &\;\quad\times
\binom{2k}k\binom{2n-2k}{n-k}\binom{2n-2k}k\binom{n-k}{2k}^{-1}
\\[2pt]
\endaligned &\cr
\noalign{\hrule}
& \eqnnoltmp{262}{325} &&
\aligned
\\[-12pt]
&
D=5^2\theta^4
+2^2\cdot5z(136\theta^4+224\theta^3+197\theta^2+85\theta+15)
\\[-2pt] &\;
+2^4z^2(5584\theta^4+16192\theta^3+21924\theta^2+14800\theta+3955)
\\[-2pt] &\;
+2^{11}z^3(608\theta^4+2280\theta^3+3642\theta^2+2745\theta+780)
\\[-2pt] &\;
+2^{14}z^4(464\theta^4+1888\theta^3+2956\theta^2+2012\theta+501)
+2^{24}z^5(\theta+1)^4
\\[-8pt]
\endaligned &\cr
& &\span\hrulefill&\cr
& &&
\aligned
&
A_n=\sum_{k=[n/3]}^{[2n/3]}\binom nk\binom{2k}k^2\binom{2n-2k}{n-k}^2
\binom{2k}{n-k}\binom{2n-2k}k
\\[-2pt] &\quad\times
\bigl(1+k(-6H_k+6H_{n-k}+6H_{2k}-6H_{2n-2k}+3H_{2n-3k}-3H_{3k-n})\bigr)
\\[-2pt] &\;
+3\sum_{k=0}^{[(n-1)/3]}(-1)^{n+k}\frac{n-2k}{n-3k}\binom nk\binom{2k}k^2
\binom{2n-2k}{n-k}^2\binom{2n-2k}k\binom{n-k}{2k}^{-1}
\\[2pt]
\endaligned &\cr
\noalign{\hrule}
& \eqnnoltmp{263}{326} &&
\aligned
\\[-12pt]
&
D=\theta^4
+2^4z(464\theta^4-32\theta^3+76\theta^2+92\theta+21)
\\[-2pt] &\;
+2^{15}z^2(608\theta^4+152\theta^3+450\theta^2+131\theta+5)
\\[-2pt] &\;
+2^{22}z^3(5584\theta^4+6144\theta^3+6852\theta^2+2808\theta+471)
\\[-2pt] &\;
+2^{34}\cdot5z^4(136\theta^4+320\theta^3+341\theta^2+181\theta+39)
+2^{46}z^5(\theta+1)^4
\\[-8pt]
\endaligned &\cr
& &\span\hrulefill&\cr
& &&
\gathered
\text{the reflection of \#262 at infinity}
\\[-8pt]
\endgathered &\cr
& &\span\hrulefill&\cr
& &&
\text{a formula for $A_n$ is not known}
&\cr
\noalign{\hrule}
}}\hss}

\newpage

\hbox to\hsize{\hss\vbox{\offinterlineskip
\halign to120mm{\strut\tabskip=100pt minus 100pt
\strut\vrule#&\hbox to6.5mm{\hss$#$\hss}&%
\vrule#&\hbox to113mm{\hfil$\dsize#$\hfil}&%
\vrule#\tabskip=0pt\cr\noalign{\hrule}
& \# &%
& \text{differential operator $D$ and coefficients $A_n$, $n=0,1,2,\dots$} &\cr
\noalign{\hrule\vskip1pt\hrule}
& \eqnnoltmp{264}{327} &&
\aligned
\\[-12pt]
&
D=\theta^4
-2^4z(-3392\theta^4+9344\theta^3+5764\theta^2+1092\theta+93)
\\[-2pt] &\;
+2^{17}\cdot3z^2(-1952\theta^4-15200\theta^3+7758\theta^2+2593\theta+323)
\\[-2pt] &\;
+2^{26}\cdot3^2\cdot7z^3(-11584\theta^4+6912\theta^3+5364\theta^2+1632\theta+167)
\\[-2pt] &\;
+2^{42}\cdot3^3\cdot7^2z^4(2\theta+1)^2(3\theta+1)(3\theta+2)
\\[-8pt]
\endaligned &\cr
& &\span\hrulefill&\cr
& &&
\aligned
A_n
&=16^{-n}\binom{2n}n^2\biggl(\sum_{k=0}^n\binom nk
\binom{2k}k\binom{2n-2k}{n-k}\binom{2n+2k}{n+k}^2\binom{4n-2k}{2n-k}^2
\\[-2pt] &\quad\times
\binom{2n}k^{-1}\binom{2n}{n-k}^{-1}
\bigl(1+k(-2H_k+2H_{n-k}-3H_{n+k}+3H_{2n-k}
\\[-2pt] &\quad\qquad
+2H_{2k}-2H_{2n-2k}+4H_{2n+2k}-4H_{4n-2k})\bigr)
\\[-2pt] &\;
+\sum_{k=1}^n\frac{n+2k}k\binom{2n+k}{2n}\binom{2n+2k}{n+k}
\binom{2n-2k}{n-k}^2\binom{4n+2k}{2n+k}^2
\\[-2pt] &\quad\times
\binom{2k}k^{-1}\binom{n+k}n^{-1}\binom{2n}{n+k}^{-1}\biggr)
\\[2pt]
\endaligned &\cr
\noalign{\hrule}
& \eqnnoltmp{265}{328} &&
\aligned
\\[-12pt]
&
D=\theta^4
-2^4\cdot3z(-96\theta^4+96\theta^3+6\theta^2+12\theta+1)
\\[-2pt] &\;
+2^{13}\cdot3z^2(288\theta^4-144\theta^3+526\theta^2+206\theta+27)
\\[-2pt] &\;
+2^{20}\cdot3^3z^3(288\theta^4+864\theta^3+652\theta^2+204\theta+23)
\\[-2pt] &\;
+2^{30}\cdot3^5z^4(2\theta+1)^2(3\theta+1)(3\theta+2)
\\[-8pt]
\endaligned &\cr
& &\span\hrulefill&\cr
& &&
\aligned
A_n
&=\sum_k\binom nk^3\binom{n+k}n\binom{2n-k}n\binom{2n+2k}{n+k}^2
\binom{4n-2k}{2n-k}^2
\\[-2pt] &\quad\times
\binom{2n}k^{-1}\binom{2n}{n-k}^{-1}
\cdot\bigl(1+k(-3H_k+3H_{n-k}-2H_{n+k}
\\[-2pt] &\quad\qquad
+2H_{2n-k}+4H_{2n+2k}-4H_{4n-2k)}\bigr)
\\[2pt]
\endaligned &\cr
\noalign{\hrule}
& \eqnnoltmp{266}{329} &&
\aligned
\\[-12pt]
&
D=5^2\theta^4
-3\cdot5z(27\theta^4+108\theta^3+124\theta^2+70\theta+15)
\\[-2pt] &\;
-3^2z^2(2754\theta^4+9072\theta^2+13014\theta^2+8910\theta+2440)
\\[-2pt] &\;
+3^5z^3(1134\theta^4+9720\theta^3+23166\theta^2+21330\theta+6890)
\\[-2pt] &\;
+3^8z^4(729\theta^4+3888\theta^3+6606\theta^2+4662\theta+1184)
+3^{15}z^5(\theta+1)^4
\\[-8pt]
\endaligned &\cr
& &\span\hrulefill&\cr
& &&
\aligned
A_n
&=\sum_k\binom nk^3\frac{(3k)!}{k!^3}\frac{(3n-3k)!}{(n-k)!^3}
\\[-2pt] &\quad\times
\bigl(1+k(-6H_k+6H_{n-k}+3H_{3k}-3H_{3n-3k})\bigr)
\\[2pt]
\endaligned &\cr
\noalign{\hrule}
& \eqnnoltmp{267}{330} &&
\aligned
\\[-12pt]
&
D=\theta^4
+3^2z(729\theta^4-972\theta^3-684\theta^2-198\theta-31)
\\[-2pt] &\;
+2\cdot3^8z^2(567\theta^4-2592\theta^3+405\theta^2+189\theta+70)
\\[-2pt] &\;
-2\cdot3^{14}z^3(1377\theta^4+972\theta^3+1161\theta^2+459\theta+113)
\\[-2pt] &\;
+3^{22}\cdot5z^4(-27\theta^4+38\theta^2+38\theta+12)
+3^{30}\cdot5^2z^5(\theta+1)^4
\\[-8pt]
\endaligned &\cr
& &\span\hrulefill&\cr
& &&
\gathered
\text{the reflection of \#266 at infinity}
\\[-8pt]
\endgathered &\cr
& &\span\hrulefill&\cr
& &&
\text{a formula for $A_n$ is not known}
&\cr
\noalign{\hrule}
}}\hss}

\newpage

\hbox to\hsize{\hss\vbox{\offinterlineskip
\halign to120mm{\strut\tabskip=100pt minus 100pt
\strut\vrule#&\hbox to6.5mm{\hss$#$\hss}&%
\vrule#&\hbox to113mm{\hfil$\dsize#$\hfil}&%
\vrule#\tabskip=0pt\cr\noalign{\hrule}
& \# &%
& \text{differential operator $D$ and coefficients $A_n$, $n=0,1,2,\dots$} &\cr
\noalign{\hrule\vskip1pt\hrule}
& \eqnnoltmp{268}{331} &&
\aligned
\\[-12pt]
&
D=5^2\theta^4
-2^2\cdot3\cdot5z(108\theta^4+432\theta^3+661\theta^2+445\theta+105)
\\[-2pt] &\;
-2^4\cdot3^2z^2(44064\theta^4+14515\theta^3+239004\theta^2+186300\theta+58045)
\\[-2pt] &\;
+2^9\cdot3^5z^3(9072\theta^4+7760\theta^3+180954\theta^2+164970\theta+53965)
\\[-2pt] &\;
+2^{12}\cdot3^8z^4(11664\theta^4+62208\theta^3+104940\theta^2+73836\theta+18659)
\\[-2pt] &\;
+2^{20}\cdot3^{15}z^5(\theta+1)^4
\\[-8pt]
\endaligned &\cr
& &\span\hrulefill&\cr
& &&
\aligned
A_n
&=\sum_k\binom nk^3\frac{(6k)!}{k!(2k)!(3k)!}
\frac{(6n-6k)!}{(n-k)!(2n-2k)!(3n-3k)!}
\\[-2pt] &\quad\times
\bigl(1+k(-4H_k+4H_{n-k}+6H_{6k}-6H_{6n-6k}-3H_{3k}
\\[-2pt] &\quad\qquad
+3H_{3n-3k}-2H_{2k}+2H_{2n-2k})\bigr)
\\[2pt]
\endaligned &\cr
\noalign{\hrule}
& \eqnnoltmp{269}{332} &&
\aligned
\\[-12pt]
&
D=\theta^4
+2^4\cdot3^2z(11664\theta^4-15552\theta^3-11700\theta^2-3924\theta-781)
\\[-2pt] &\;
+2^{13}\cdot3^8z^2(9072\theta^4-41472\theta^3+2106^2-54\theta+1261)
\\[-2pt] &\;
-2^{20}\cdot3^{14}z^3(44064\theta^4+31104\theta^3+67932\theta^2+32508\theta+9661)
\\[-2pt] &\;
-2^{30}\cdot3^{22}\cdot5z^4(108\theta^4+13\theta^2+13\theta-3)
\\[-2pt] &\;
+2^{40}\cdot3^{30}\cdot5^2z^5(\theta+1)^4
\\[-8pt]
\endaligned &\cr
& &\span\hrulefill&\cr
& &&
\gathered
\text{the reflection of \#268 at infinity}
\\[-8pt]
\endgathered &\cr
& &\span\hrulefill&\cr
& &&
\text{a formula for $A_n$ is not known}
&\cr
\noalign{\hrule}
& \eqnnoltmp{270}{333} &&
\aligned
\\[-12pt]
&
D=5^2\theta^4
-2^2\cdot5z(48\theta^4+192\theta^3+251\theta^2+155\theta+35)
\\[-2pt] &\;
-2^4z^2(8704\theta^4+28672\theta^3+43664\theta^2+31760\theta+9265)
\\[-2pt] &\;
+2^{11}z^3(1792\theta^4+15360\theta^3+36248\theta^2+33240\theta+10795)
\\[-2pt] &\;
+2^{16}z^4(2304\theta^4+12288\theta^3+20816\theta^2+14672\theta+3719)
+2^{30}z^5(\theta+1)^4
\\[-8pt]
\endaligned &\cr
& &\span\hrulefill&\cr
& &&
\aligned
&
A_n=\sum_k\binom nk^3\frac{(4k)!}{k!^2(2k)!}
\frac{(4n-4k)!}{(n-k)!^2(2n-2k)!}
\\[-2pt] &\quad\times
\bigl(1+k(-5H_k+5H_{n-k}+4H_{4k}-4H_{4n-4k}+2H_{2n-2k}-2H_{2k})\bigr)
\\[2pt]
\endaligned &\cr
\noalign{\hrule}
& \eqnnoltmp{271}{334} &&
\aligned
\\[-12pt]
&
D=\theta^4
-2^4z(-2304\theta^4+3072\theta^3+2224\theta^2+668\theta+121)
\\[-2pt] &\;
+2^{17}z^2(1792\theta^4-8192\theta^3+920\theta^2+344\theta+235)
\\[-2pt] &\;
-2^{28}z^3(8704\theta^4+6144\theta^3+9872\theta^2+4368\theta+1201)
\\[-2pt] &\;
+2^{44}\cdot5z^4(-48\theta^4+37\theta^2+37\theta+13)
+2^{60}\cdot5^2z^5(\theta+1)^4
\\[-8pt]
\endaligned &\cr
& &\span\hrulefill&\cr
& &&
\gathered
\text{the reflection of \#270 at infinity}
\\[-8pt]
\endgathered &\cr
& &\span\hrulefill&\cr
& &&
\text{a formula for $A_n$ is not known}
&\cr
\noalign{\hrule}
& \eqnnoltmp{272}{335} &&
\aligned
\\[-12pt]
&
D=5^2\theta^4
-2^2\cdot3\cdot5z(1332\theta^4+3528\theta^3+3289\theta^2+1525\theta+285)
\\[-2pt] &\;
+2^4\cdot3^2z^2(331776\theta^4+1602368\theta^3+453333\theta^2+2996640\theta+724415)
\\[-2pt] &\;
+2^8\cdot3^5z^3(539136\theta^4+622080\theta^3-3024864\theta^2-4008960\theta-1315985)
\\[-2pt] &\;
-2^{15}\cdot3^8z^4(41472\theta^4+393984\theta^3+735984\theta^2+510912\theta+120811)
\\[-2pt] &\;
-2^{20}\cdot3^{11}z^5(12\theta+7)(12\theta+11)(12\theta+13)(12\theta+17)
\\[-8pt]
\endaligned &\cr
& &\span\hrulefill&\cr
& &&
\aligned
&
A_n=\sum_k\binom nk\binom{n+k}n\binom{2n-k}n\frac{(6k)!}{k!(2k)!(3k)!}
\\[-2pt] &\;\times
\frac{(6n-6k)!}{(n-k)!(2n-2k)!(3n-3k)!}
\cdot\bigl(1+k(-3H_k+3H_{n-k}+H_{n+k}
\\[-2pt] &\;\quad
-H_{2n-k}+6H_{6k}-6H_{6n-6k}
-3H_{3k}+3H_{3n-3k}-2H_{2k}+2H_{2n-2k})\bigr)
\\[2pt]
\endaligned &\cr
\noalign{\hrule}
}}\hss}

\newpage

\hbox to\hsize{\hss\vbox{\offinterlineskip
\halign to120mm{\strut\tabskip=100pt minus 100pt
\strut\vrule#&\hbox to6.5mm{\hss$#$\hss}&%
\vrule#&\hbox to113mm{\hfil$\dsize#$\hfil}&%
\vrule#\tabskip=0pt\cr\noalign{\hrule}
& \# &%
& \text{differential operator $D$ and coefficients $A_n$, $n=0,1,2,\dots$} &\cr
\noalign{\hrule\vskip1pt\hrule}
& \eqnnoltmp{273}{336} &&
\aligned
\\[-12pt]
&
D=5^2\theta^4
-3\cdot5z(333\theta^4+882\theta^3+781\theta^2+340\theta+60)
\\[-2pt] &\;
+2^2\cdot3^2z^2(5184\theta^4+40662\theta^3+71829\theta^2+47700\theta+11540)
\\[-2pt] &\;
+2^2\cdot3^5z^3(8424\theta^4+9720\theta^3-46899\theta^2-63045\theta-21260)
\\[-2pt] &\;
-2^4\cdot3^8z^4(1296\theta^4+12312\theta^3+23094\theta^2+16263\theta+3956)
\\[-2pt] &\;
-2^6\cdot3^{11}z^5(3\theta+2)(3\theta+4)(6\theta+5)(6\theta+7)
\\[-8pt]
\endaligned &\cr
& &\span\hrulefill&\cr
& &&
\aligned
&
A_n=\sum_k\binom nk\binom{n+k}n\binom{2n-k}n
\frac{(3k)!}{k!^3}\frac{(3n-3k)!}{(n-k)!^3}
\\[-2pt] &\quad\times
\bigl(1+k(-5H_k+5H_{n-k}+H_{n+k}-H_{2n-k}+3H_{3k}-3H_{3n-3k})\bigr)
\\[2pt]
\endaligned &\cr
\noalign{\hrule}
& \eqnnoltmp{274}{337} &&
\aligned
\\[-12pt]
&
D=5^2\theta^4
-5z(757\theta^4+1298\theta^3+1049\theta^2+400\theta+60)
\\[-2pt] &\;
+2^2\cdot3^2z^2(5456\theta^4+17498\theta^3+22121\theta^2+11940\theta+2340)
\\[-2pt] &\;
-2^2\cdot3^4z^3(15128\theta^4+68040\theta^3+112171\theta^2+73845\theta+16380)
\\[-2pt] &\;
+2^4\cdot3^8z^4(2\theta+1)(216\theta^3+864\theta^2+1015\theta+356)
\\[-2pt] &\;
-2^6\cdot3^{10}z^5(2\theta+1)(2\theta+3)(3\theta+2)(3\theta+4)
\\[-8pt]
\endaligned &\cr
& &\span\hrulefill&\cr
& &&
\aligned
&
A_n=(-1)^n9^{-n}\binom{2n}n\sum_{k=[n/3]}^{[2n/3]}\binom{3k}n\binom{3n-3k}n
\frac{(3k)!}{k!^3}\frac{(3n-3k)!}{(n-k)!^3}
\\[-2pt] &\quad\times
\bigl(1+k(-3H_k+3H_{n-k}+6H_{3k}-6H_{3n-3k}+3H_{2n-3k}-3H_{3k-n})\bigr)
\\[-2pt] &\;
+3\cdot9^{-n}\binom{2n}n\sum_{k=0}^{[(n-1)/3]}(-1)^k\frac{n-2k}{n-3k}
\frac{(3k)!}{k!^3}\frac{(3n-3k)!}{(n-k)!^3}\binom{3n-3k}n\binom n{3k}^{-1}
\\[2pt]
\endaligned &\cr
\noalign{\hrule}
& \eqnnoltmp{275}{338} &&
\aligned
\\[-12pt]
&
D=5^2\theta^4
-2^2\cdot5z(592\theta^4+1568\theta^3+1419\theta^2+635\theta+115)
\\[-2pt] &\;
+2^4z^2(6553\theta^4+514048\theta^3+902816\theta^2+598400\theta+144735)
\\[-2pt] &\;
+2^{10}z^3(106496\theta^4+122880\theta^3-594816\theta^2-794880\theta-265065)
\\[-2pt] &\;
-2^{19}z^4(8192\theta^4+77824\theta^3+145728\theta^2+102016\theta+24527)
\\[-2pt] &\;
-2^{26}z^5(8\theta+5)(8\theta+7)(8\theta+9)(8\theta+11)
\\[-8pt]
\endaligned &\cr
& &\span\hrulefill&\cr
& &&
\aligned
A_n
&=\sum_k\binom nk\binom{n+k}n\binom{2n-k}n\frac{(4k)!}{k!^2(2k)!}
\frac{(4n-4k)!}{(n-k)!^2(2n-2k)!}
\\[-2pt] &\quad\times
\bigl(1+k(-4H_k+4H_{n-k}+H_{n+k}-H_{2n-k}
\\[-2pt] &\quad\qquad
+4H_{4k}-4H_{4n-4k}+2H_{2n-2k}-2H_{2k})\bigr)
\\[2pt]
\endaligned &\cr
\noalign{\hrule}
& \eqnnoltmp{276}{339} &&
\aligned
\\[-12pt]
&
D=\theta^4
+2^4\cdot3z(-18432\theta^4+4608\theta^3+1024\theta^2-1280\theta-221)
\\[-2pt] &\;
+2^{17}\cdot3^4z^2(25344\theta^4-2304\theta^3+11680\theta^2+1472\theta-33)
\\[-2pt] &\;
-2^{28}\cdot3^8z^3(18432\theta^4+13824\theta^3+11392\theta^2+3264\theta+359)
\\[-2pt] &\;
+2^{46}\cdot3^{12}z^4(2\theta+1)^2(3\theta+1)(3\theta+2)
\\[-8pt]
\endaligned &\cr
& &\span\hrulefill&\cr
& &&
\aligned
&
A_n=\binom{2n}n\sum_k\binom{2n+2k}{n+k}\binom{4n-2k}{2n-k}
\frac{(6k)!}{k!(2k)!(3k)!}
\\[-2pt] &\;\times
\frac{(6n-6k)!}{(n-k)!(2n-2k)!(3n-3k)!}
\cdot\bigl(1+k(-H_k+H_{n-k}-2H_{2k}+2H_{2n-2k}
\\[-2pt] &\;\quad
-3H_{3k}+3H_{3n-3k}
+6H_{6k}-6H_{6n-6k}+2H_{2n+2k}-2H_{4n-2k})\bigr)
\\[2pt]
\endaligned &\cr
\noalign{\hrule}
}}\hss}

\newpage

\hbox to\hsize{\hss\vbox{\offinterlineskip
\halign to120mm{\strut\tabskip=100pt minus 100pt
\strut\vrule#&\hbox to6.5mm{\hss$#$\hss}&%
\vrule#&\hbox to113mm{\hfil$\dsize#$\hfil}&%
\vrule#\tabskip=0pt\cr\noalign{\hrule}
& \# &%
& \text{differential operator $D$ and coefficients $A_n$, $n=0,1,2,\dots$} &\cr
\noalign{\hrule\vskip1pt\hrule}
& \eqnnol{277}{223} &&
\aligned
\\[-12pt]
&
D=\theta^4
-2^4z(-576\theta^4+1152\theta^3+724\theta^2+148\theta+13)
\\[-2pt] &\;
+2^{17}z^2(-32\theta^4-992\theta^2+166\theta^2+57\theta+6)
\\[-2pt] &\;
-2^{26}\cdot3z^3(832\theta^4+768\theta^3+556\theta^2+192\theta+25)
-2^{40}\cdot3^2z^4(2\theta+1)^4
\\[-8pt]
\endaligned &\cr
& &\span\hrulefill&\cr
& &&
\gathered
\text{the reflection of \#55 at infinity}
\\[-8pt]
\endgathered &\cr
& &\span\hrulefill&\cr
& &&
\aligned
A_n
&=\binom{2n}n^2\sum_k\binom nk\binom{2k}k\binom{2n-2k}{n-k}
\binom{2n+2k}{n+k}\binom{4n-2k}{2n-k}
\\[-2pt] &\quad\times
\bigl(1+k(-3H_k+3H_{n-k}+2H_{2k}-2H_{2n-2k}+2H_{2n-k}
\\[-2pt] &\quad\quad
+2H_{n+k}+2H_{2n+2k}-2H_{4n-2k})\bigr)
\\[2pt]
\endaligned &\cr
\noalign{\hrule}
& \eqnnoltmp{278}{340} &&
\aligned
\\[-12pt]
&
D=\theta^4
-3z(279\theta^4+882\theta^3+641\theta^2+200\theta+24)
\\[-2pt] &\;
+2\cdot3^5z^2(-72\theta^4+1710\theta^3+3665\theta^2+1864\theta+296)
\\[-2pt] &\;
+2^2\cdot3^9z^3(909\theta^4+3888\theta^3+3082\theta^2+918\theta+92)
\\[-2pt] &\;
+2^4\cdot3^{15}z^4(2\theta+1)^2(3\theta+1)(3\theta+2)
\\[-8pt]
\endaligned &\cr
& &\span\hrulefill&\cr
& &&
A_n=\binom{2n}n^2\binom{3n}n\sum_k\frac{(3k)!}{k!^3}
\frac{(3n-3k)!}{(n-k)!^3}\binom{n+k}n^{-1}\binom{2n-k}n^{-1}
&\cr
\noalign{\hrule}
& \eqnnoltmp{279}{341} &&
\aligned
\\[-12pt]
&
D=17^2\theta^4
+17z(286\theta^4+734\theta^3+656\theta^2+289\theta+51)
\\[-2pt] &\;
+3^2z^2(4110\theta^4+22074\theta^3+37209\theta^2+26265\theta+6800)
\\[-2pt] &\;
-3^5z^3(1521\theta^4+7344\theta^3+12936\theta^2+9945\theta+2822)
\\[-2pt] &\;
+3^8z^4(123\theta^4+552\theta^3+879\theta^2+603\theta+152)
-3^{12}z^5(\theta+1)^4
\\[-8pt]
\endaligned &\cr
& &\span\hrulefill&\cr
& &&
A_n=3\sum_{k=[2(n+1)/3]}^n(-1)^k\binom nk\frac{n-2k}{n-3k}\frac{(3k)!}{k!^3}
\frac{(3n-3k)!}{(n-k)!^3}\binom n{3n-3k}\binom{3k}n^{-1}
&\cr
\noalign{\hrule}
& \eqnnoltmp{280}{342} &&
\aligned
\\[-12pt]
&
D=\theta^4
+3^2z(123\theta^4-60\theta^3-39\theta^2-9\theta-1)
\\[-2pt] &\;
+3^5z^2(1521\theta^2-1260\theta^3+30\theta^2-21\theta-10)
\\[-2pt] &\;
+3^8z^3(4110\theta^4-5634\theta^3-4353\theta^2-1629\theta-220)
\\[-2pt] &\;
-3^{12}\cdot17z^4(286\theta^4+410\theta^3+170\theta^2-35\theta-30)
+3^{18}\cdot17^2z^5(\theta+1)^4
\\[-8pt]
\endaligned &\cr
& &\span\hrulefill&\cr
& &&
\gathered
\text{the reflection of \#279 at infinity}
\\[-8pt]
\endgathered &\cr
& &\span\hrulefill&\cr
& &&
\text{a formula for $A_n$ is not known}
&\cr
\noalign{\hrule}
}}\hss}

\newpage

\hbox to\hsize{\hss\vbox{\offinterlineskip
\halign to120mm{\strut\tabskip=100pt minus 100pt
\strut\vrule#&\hbox to6.5mm{\hss$#$\hss}&%
\vrule#&\hbox to113mm{\hfil$\dsize#$\hfil}&%
\vrule#\tabskip=0pt\cr\noalign{\hrule}
& \# &%
& \text{differential operator $D$ and coefficients $A_n$, $n=0,1,2,\dots$} &\cr
\noalign{\hrule\vskip1pt\hrule}
& \eqnnoltmp{281}{343} &&
\aligned
\\[-12pt]
&
D=\theta^4
+z(328\theta^4+82\theta^3+33\theta^2-8\theta-7)
\\[-2pt] &\;
+z^2(52844\theta^4+26422\theta^3+11771\theta^2-3472\theta-3676)
\\[-2pt] &\;
+z^3(5280472\theta^4+3960354\theta^3+1932893\theta^2-37864\theta+197467)
\\[-2pt] &\;
+z^4(355955926\theta^4+355955926\theta^3+188833541\theta^2
\\[-2pt] &\;\quad
+39149936\theta+53236566)
\\[-2pt] &\;
+5z^5(3300295000\theta^4+4125368750\theta^3+2363161375\theta^2
\\[-2pt] &\;\quad
+789584250\theta+300065632)
\\[-2pt] &\;
+5^6z^6(33027500\theta^4+49541250\theta^3+30476125\theta^2+10997750\theta
+305303)
\\[-2pt] &\;
+5^{11}z^7(205000\theta^4+358750\theta^3+235875\theta^2+71750\theta+6104)
\\[-2pt] &\;
+5^{16}z^8(5\theta+1)(5\theta+2)(5\theta+3)(5\theta+4)
\\[-8pt]
\endaligned &\cr
& &\span\hrulefill&\cr
& &&
\gathered
\text{the pullback of the 5th-order differential equation $D'y=0$, where}
\\[-8pt]
\endgathered &\cr
& &\span\hrulefill&\cr
& &&
\aligned
\\[-12pt]
&
D'=\theta^5
+z(2\theta+1)(41\theta^4+82\theta^3+74\theta^2+33\theta+6)
\\[-2pt] &\;
+5z^2(\theta+1)(5\theta+3)(5\theta+4)(5\theta+6)(5\theta+7)
\\[-8pt]
\endaligned &\cr
& &\span\hrulefill&\cr
& &&
A_n'
=3\sum_{k=0}^{[n/3]}(-1)^k\binom nk^3\frac{n-2k}{2n-3k}
\binom{n+k}n\binom{2n-k}n\binom{2n-3k}n
&\cr
\noalign{\hrule}
& \eqnnoltmp{282}{344} &&
\aligned
\\[-12pt]
&
D=5^2\theta^4
-2^2\cdot5z(1348\theta^4+752\theta^3+521\theta^2+145\theta+15)
\\[-2pt] &\;
+2^4\cdot3^4z^2(5696\theta^4-1792\theta^3-7304\theta^2-3740\theta-585)
\\[-2pt] &\;
+2^{10}\cdot3^8z^3(-20\theta^4+360\theta^3+289\theta^2+90\theta+10)
\\[-2pt] &\;
-2^{12}\cdot3^{13}z^4(2\theta+1)^4
\\[-8pt]
\endaligned &\cr
& &\span\hrulefill&\cr
& &&
\aligned
&
A_n=(-1)^n\binom{2n}n\sum_{k=[n/3]}^{[2n/3]}\binom{2n}{3k}
\binom{2n}{3n-3k}\frac{(3k)!}{k!^3}\frac{(3n-3k)!}{(n-k)!^3}
\\[-2pt] &\quad\times
\bigl(1+k(-3H_k+3H_{n-k}+3H_{2n-3k}-3H_{3k-n})\bigr)
\\[-2pt] &\;
+3\binom{2n}n\sum_{k=0}^{[(n-1)/3]}(-1)^{n+k}\frac{n-2k}{2n-3k}
\frac{(3k)!}{k!^3}\frac{(3n-3k)!}{(n-k)!^3}
\binom{2n}{3n-3k}\binom{3k}{2n}^{-1}
\\[2pt]
\endaligned &\cr
\noalign{\hrule}
& \eqnnoltmp{283}{345} &&
\aligned
\\[-12pt]
&
D=\theta^4
+2^2z(20\theta^4+400\theta^3+281\theta^2+81\theta+9)
\\[-2pt] &\;
-2^4\cdot3z^2(5696\theta^4+13184\theta^3+3928\theta^2+628\theta+39)
\\[-2pt] &\;
+2^{10}\cdot3^2\cdot5z^3(1348\theta^4+1944\theta^4+1415\theta^2+486\theta+63)
\\[-2pt] &\;
-2^{12}\cdot3^7\cdot5^2z^4(2\theta+1)^4
\\[-8pt]
\endaligned &\cr
& &\span\hrulefill&\cr
& &&
\gathered
\text{the reflection of \#282 at infinity}
\\[-8pt]
\endgathered &\cr
& &\span\hrulefill&\cr
& &&
\text{a formula for $A_n$ is not known}
&\cr
\noalign{\hrule}
}}\hss}

\newpage

\hbox to\hsize{\hss\vbox{\offinterlineskip
\halign to120mm{\strut\tabskip=100pt minus 100pt
\strut\vrule#&\hbox to6.5mm{\hss$#$\hss}&%
\vrule#&\hbox to113mm{\hfil$\dsize#$\hfil}&%
\vrule#\tabskip=0pt\cr\noalign{\hrule}
& \# &%
& \text{differential operator $D$ and coefficients $A_n$, $n=0,1,2,\dots$} &\cr
\noalign{\hrule\vskip1pt\hrule}
& \eqnnoltmp{284}{346} &&
\aligned
\\[-12pt]
&
D=38^2\theta^4
-38z(3014\theta^4+5878\theta^3+4725\theta^2+1786\theta+266)
\\[-2pt] &\;
+z^2(689717\theta^4+2305502\theta^3+3057079\theta^2+1810054\theta+402002)
\\[-2pt] &\;
-z^3(1438808\theta^4+5812350\theta^3+9142457\theta^2+6295992\theta+1576582)
\\[-2pt] &\;
+z^4(1395491\theta^4+5075392\theta^3+6297445\theta^2+3375833\theta+663471)
\\[-2pt] &\;
-z^5(834163\theta^4+2657224\theta^3+2407768\theta^2+604005\theta-52928)
\\[-2pt] &\;
+z^6(277543\theta^4+692484\theta^3+572576\theta^2+148359\theta-4832)
\\[-2pt] &\;
+11z^7(4625\theta^4+9100\theta^3+6395\theta^2+1845\theta+178)
-11^2z^8(\theta+1)^4
\\[-8pt]
\endaligned &\cr
& &\span\hrulefill&\cr
& &&
A_n=\sum_{k,l}\binom nk\binom nl\binom{n+k}n\binom{n+l}n
\binom n{k+l}\binom{n+k-l}{n-l}
&\cr
\noalign{\hrule}
& \eqnnoltmp{285}{347} &&
\aligned
\\[-12pt]
&
D=11^2\theta^4
+11z(4625\theta^4+9400\theta^3+6845\theta^2+2145\theta+253)
\\[-2pt] &\;
-z^2(834163\theta^4+417688\theta^3+160382\theta^2+29513\theta+4444)
\\[-2pt] &\;
+z^3(834163\theta^4+679428\theta^3-558926\theta^2-423489\theta-72226)
\\[-2pt] &\;
-z^4(1395491\theta^4+506572\theta^3-555785\theta^2-425155\theta-94818)
\\[-2pt] &\;
+z^5(1438808\theta^4+57118\theta^3+338255\theta^2+307104\theta+49505)
\\[-2pt] &\;
-z^6(689717\theta^4+453366\theta^3+278875\theta^2+146466\theta+33242)
\\[-2pt] &\;
+38z^7(3014\theta^4+6178\theta^3+5175\theta^2+2086\theta+341)
-38^2z^8(\theta+1)^4
\\[-8pt]
\endaligned &\cr
& &\span\hrulefill&\cr
& &&
\gathered
\text{the reflection of \#284 at infinity}
\\[-8pt]
\endgathered &\cr
& &\span\hrulefill&\cr
& &&
\text{a formula for $A_n$ is not known}
&\cr
\noalign{\hrule}
& \eqnnoltmp{286}{348} &&
\aligned
\\[-12pt]
&
D=3^2\theta^4
-3^2z(38\theta^4+82\theta^3+67\theta^2+26\theta+4)
\\[-2pt] &\;
+3z^2(2045\theta^4+5702\theta^3+7535\theta^2+4170\theta+852)
\\[-2pt] &\;
+2^3\cdot3z^3(2208\theta^4+5925\theta^3+7925\theta^2+5607\theta+1512)
\\[-2pt] &\;
+2^3z^4(60287\theta^4+56374\theta^3-215983\theta^2-268986\theta-85452)
\\[-2pt] &\;
-2^4z^5(205651\theta^4+605608\theta^3+603579\theta^2+204622\theta+8104)
\\[-2pt] &\;
+2^7z^6(-51414\theta^4+273267\theta^3+502700\theta^2+305649\theta+63398)
\\[-2pt] &\;
+2^8\cdot37z^7(7909\theta^4+18122\theta^3+17595\theta^2+8462\theta+1672)
\\[-2pt] &\;
-2^{13}\cdot37^2z^8(\theta+1)^2(4\theta+3)(4\theta+5)
\\[-8pt]
\endaligned &\cr
& &\span\hrulefill&\cr
& &&
A_n=\sum_{k,l}\binom nk\binom nl\binom{2n-k}n
\binom n{k-l}\binom{2k}n\binom{2l}n
&\cr
\noalign{\hrule}
& \eqnnoltmp{287}{349} &&
\aligned
\\[-12pt]
&
D=21^2\theta^4
-21z(3289\theta^4+6098\theta^3+4645\theta^2+1596\theta+210)
\\[-2pt] &\;
+2^2z^2(38560\theta^4-230840\theta^3-534425\theta^2-337050\theta-68145)
\\[-2pt] &\;
+2^4z^3(106636\theta^4+493416\theta^3+420211\theta^2+116361\theta+6090)
\\[-2pt] &\;
-5\cdot2^8z^4(2\theta+1)(1916\theta^3+2622\theta^2+1077\theta+91)
\\[-2pt] &\;
-5^2\cdot2^{12}z^5(\theta+1)^2(2\theta+1)(2\theta+3)
\\[-8pt]
\endaligned &\cr
& &\span\hrulefill&\cr
& &&
A_n=\binom{2n}n\sum_{i,j}\binom ni\binom nj\binom{n+j}n
\binom{i+j}n\binom n{i-j}
&\cr
\noalign{\hrule}
}}\hss}

\newpage

\hbox to\hsize{\hss\vbox{\offinterlineskip
\halign to120mm{\strut\tabskip=100pt minus 100pt
\strut\vrule#&\hbox to6.5mm{\hss$#$\hss}&%
\vrule#&\hbox to113mm{\hfil$\dsize#$\hfil}&%
\vrule#\tabskip=0pt\cr\noalign{\hrule}
& \# &%
& \text{differential operator $D$ and coefficients $A_n$, $n=0,1,2,\dots$} &\cr
\noalign{\hrule\vskip1pt\hrule}
& \eqnnoltmp{288}{350} &&
\aligned
\\[-12pt]
&
D=\theta^4
-2^4z(496\theta^4+1568\theta^3+1060\theta^2+276\theta+27)
\\[-2pt] &\;
+3\cdot 2^{15}z^2(-32\theta^4+760\theta^3+1570\theta^2+651\theta+81)
\\[-2pt] &\;
+3^22^{22}z^3(1616\theta^4+6912\theta^3+5092\theta^2+1416\theta+135)
\\[-2pt] &\;
+3^32^{34}z^4(3\theta+1)(3\theta+2)(4\theta+1)(4\theta+3)
\\[-8pt]
\endaligned &\cr
& &\span\hrulefill&\cr
& &&
\aligned
A_n
&==\binom{2n}n\binom{3n}n\binom{4n}{2n}\sum_k\frac{(4k)!}{k!^2(2k)!}\frac{(4n-4k)!}{(n-k)!^2(2n-2k)!}
\\[-2pt] &\qquad\times
\binom{n+k}n^{-1}\binom{2n-k}n^{-1}
\\[2pt]
\endaligned &\cr
\noalign{\hrule}
& \eqnnoltmp{289}{351} &&
\aligned
\\[-12pt]
&
D=\theta^4
-2^4z(400\theta^4+2720\theta^3+1752\theta^2+392\theta+33)
\\[-2pt] &\;
+2^{15}z^2(-4272\theta^4-6288\theta^3+3184\theta^2+1484\theta+177)
\\[-2pt] &\;
+2^{24}\cdot5z^3(-4688\theta^4+1536\theta^3+1384\theta^2+336\theta+27)
\\[-2pt] &\;
+2^{36}\cdot5^2z^4(2\theta+1)^2(4\theta+1)(4\theta+3)
\\[-8pt]
\endaligned &\cr
& &\span\hrulefill&\cr
& &&
\aligned
&
A_n=(-1)^n\binom{2n}n\binom{4n}{2n}\sum_k\binom nk\binom{n+k}n
\binom{2n-k}n\binom{2n}k\binom{2n}{n-k}
\\[-2pt] &\quad\times
\binom{2n+2k}{n+k}\binom{4n-2k}{2n-k}\binom{2k}k^{-1}\binom{2n-2k}{n-k}^{-1}
\bigl(1+k(-H_k+H_{n-k}
\\[-2pt] &\quad\qquad
-2H_{n+k}+2H_{2n-k}-2H_{2k}+2H_{2n-2k}
+2H_{2n+2k}-2H_{4n-2k})\bigr)
\\[2pt]
\endaligned &\cr
\noalign{\hrule}
& \eqnnol{290}{221} &&
\aligned
\\[-12pt]
&
D=\theta^4
+3z(-279\theta^4+252\theta^3+160\theta^2+34\theta+3)
\\[-2pt] &\;
+2\cdot3^5z^2(423\theta^4-468\theta^3+457\theta^2+215\theta+37)
\\[-2pt] &\;
-2\cdot3^9z^3(531\theta^4+1296\theta^3+1243\theta^2+567\theta+104)
\\[-2pt] &\;
+3^{15}\cdot5z^4(51\theta^4+120\theta^3+126\theta^2+66\theta+14)
\\[-2pt] &\;
-3^{20}\cdot5^2z^5(\theta+1)^4
\\[-8pt]
\endaligned &\cr
& &\span\hrulefill&\cr
& &&
\gathered
\text{the reflection of \#17 at infinity}
\\[-8pt]
\endgathered &\cr
& &\span\hrulefill&\cr
& &&
\text{a formula for $A_n$ is not known}
&\cr
\noalign{\hrule}
& \eqnnol{291}{226} &&
\aligned
\\[-12pt]
&
D=\theta^4
-z(566\theta^4+34\theta^3+62\theta^2+45\theta+9)
\\[-2pt] &\;
+3z^2(39370\theta^4+17302\theta^3+22493\theta^2+8369\theta+1140)
\\[-2pt] &\;
-3^2z^3(1215215\theta^4+1432728\theta^3+1274122\theta^2+538245\theta+93222)
\\[-2pt] &\;
+3^7\cdot61z^4(3029\theta^4+6544\theta^3+6135\theta^2+2863\theta+548)
\\[-2pt] &\;
-3^{12}\cdot61^2z^5(\theta+1)^4
\\[-8pt]
\endaligned &\cr
& &\span\hrulefill&\cr
& &&
\gathered
\text{the reflection of \#124 at infinity}
\\[-8pt]
\endgathered &\cr
& &\span\hrulefill&\cr
& &&
\text{a formula for $A_n$ is not known}
&\cr
\noalign{\hrule}
& \eqnnoltmp{292}{352} &&
\aligned
\\[-12pt]
&
D=3^2\theta^4
-2^2\cdot3z(4636\theta^2+7928\theta^3+5347\theta^2+1383\theta+126)
\\[-2pt] &\;
+2^9z^2(59048\theta^4+50888\theta^3-26248\theta^2-16827\theta-2205)
\\[-2pt] &\;
+2^{16}\cdot7z^3(-9004\theta^4+2304\theta^3+2511\theta^2+504\theta+27)
\\[-2pt] &\;
-2^{24}\cdot7^2z^4(2\theta+1)^2(4\theta+1)(4\theta+3)
\\[-8pt]
\endaligned &\cr
& &\span\hrulefill&\cr
& &&
A_n=\binom{2n}n\sum_k(-1)^{n+k}4^{n-k}\binom nk
\binom{2n+k}n^2\binom{4n+2k}{2n+k}
&\cr
\noalign{\hrule}
}}\hss}

\newpage

\hbox to\hsize{\hss\vbox{\offinterlineskip
\halign to120mm{\strut\tabskip=100pt minus 100pt
\strut\vrule#&\hbox to6.5mm{\hss$#$\hss}&%
\vrule#&\hbox to113mm{\hfil$\dsize#$\hfil}&%
\vrule#\tabskip=0pt\cr\noalign{\hrule}
& \# &%
& \text{differential operator $D$ and coefficients $A_n$, $n=0,1,2,\dots$} &\cr
\noalign{\hrule\vskip1pt\hrule}
& \eqnnoltmp{293}{353} &&
\aligned
\\[-12pt]
&
D=\theta^4
-2^2z(54\theta^4+66\theta^3+49\theta^2+16\theta+2)
\\[-2pt] &\;
+2^4z^2(417\theta^4-306\theta^3-1219\theta^2-776\theta-154)
\\[-2pt] &\;
+2^8z^3(166\theta^4+1920\theta^3+1589\theta^2+432\theta+23)
\\[-2pt] &\;
-2^{12}\cdot7z^4(2\theta+1)(38\theta^3+45\theta^2+12\theta-2)
\\[-2pt] &\;
-2^{14}\cdot7^2z^5(\theta+1)^2(2\theta+1)(2\theta+3)
\\[-8pt]
\endaligned &\cr
& &\span\hrulefill&\cr
& &&
A_n=\sum_{k,l}\binom nk\binom nl\binom{k+l}n
\binom{2k}n\binom{2l}n\binom{2n}{k+l}
&\cr
\noalign{\hrule}
& \eqnnoltmp{294}{354} &&
\aligned
\\[-12pt]
&
D=\theta^4
+2^4z(-18800\theta^4+14624\theta^3+8184\theta^2+872\theta+33)
\\[-2pt] &\;
+2^{18}z^2(101744\theta^4-107920\theta^3+74968\theta^2+15100\theta+1191)
\\[-2pt] &\;
-2^{30}\cdot17z^3(40048\theta^4+49152\theta^3+35848\theta^2+10752\theta+1143)
\\[-2pt] &\;
+2^{50}\cdot17^2z^4(2\theta+1)^2(4\theta+1)(4\theta+3)
\\[-8pt]
\endaligned &\cr
& &\span\hrulefill&\cr
& &&
\aligned
&
A_n=\frac{(4n)!}{n!^2(2n)!}\sum_k\binom{2n}{2k}\binom{n+2k}k
\binom{3n-2k}{n-k}\binom{2n+4k}{n+2k}\binom{6n-4k}{3n-3k}
\\[-2pt] &\quad\times
\bigl(1+k(-H_k+H_{n-k}-H_{n+k}+H_{2n-k}-2H_{2k}+2H_{2n-2k}
\\[-2pt] &\quad\;
+2H_{3n-2k}-2H_{n+2k}+4H_{2n+4k}-4H_{6n-4k})\bigr)
\\[2pt]
\endaligned &\cr
\noalign{\hrule}
& \eqnnoltmp{295}{355} &&
\aligned
\\[-12pt]
&
D=\theta^4
+2^4z(816\theta^4-1440\theta^3-904\theta^2-184\theta-17)
\\[-2pt] &\;
+2^{18}z^2(80\theta^4-592\theta^3+432\theta^2+164\theta+23)
\\[-2pt] &\;
+2^{30}z^3(-80\theta^4+384\theta^3+296\theta^2+96\theta+11)
+2^{45}z^4(2\theta+1)^4
\\[-8pt]
\endaligned &\cr
& &\span\hrulefill&\cr
& &&
\aligned
&
A_n=\binom{3n}n\sum_k\binom{2k}k\binom{2n-2k}{n-k}
\binom{2n}{2k}\binom{n+2k}k
\\[-2pt] &\;\times
\binom{3n-2k}{n-k}\binom{2n+4k}{n+2k}
\binom{6n-4k}{3n-2k}\binom{3n}{n+k}^{-1}
\\[-2pt] &\;\times
\bigl(1+k(-3H_k+3H_{n-k}-2H_{n+2k}+2H_{3n-2k}
+4H_{2n+4k}-4H_{6n-4k})\bigr)
\\[2pt]
\endaligned &\cr
\noalign{\hrule}
& \eqnnoltmp{296}{356} &&
\aligned
\\[-12pt]
&
D=\theta^4
-2^4z(5\theta^4+34\theta^3+25\theta^2+8\theta+1)
\\[-2pt] &\;
+2^{11}z^2(5\theta^4+47\theta^3+90\theta^2+47\theta+8)
\\[-2pt] &\;
+2^{16}z^3(51\theta^4+192\theta^3+155\theta^2+48\theta+5)
+2^{23}z^4(2\theta+1)^4
\\[-8pt]
\endaligned &\cr
& &\span\hrulefill&\cr
& &&
\gathered
\text{the reflection of \#295 at infinity}
\\[-8pt]
\endgathered &\cr
& &\span\hrulefill&\cr
& &&
A_n=\binom{2n}n^3\sum_k(-1)^{n+k}4^{n-k}\binom nk^2\binom{2k}n
\binom{2n}{n-k}\binom{2n}{2k}^{-1}\binom{2n-2k}{n-k}^{-1}
&\cr
\noalign{\hrule}
& \eqnnoltmp{297}{357} &&
\aligned
\\[-12pt]
&
D=7^2\theta^4
-2\cdot7z\theta(520\theta^3+68\theta^2+41\theta+7)
\\[-2pt] &\;
-2^2\cdot3z^2(9480\theta^4+153912\theta^3+212893\theta^2+108080\theta+18816)
\\[-2pt] &\;
+2^4\cdot3^3z^3(93968\theta^4+341544\theta^3+319592\theta^2+125853\theta+18242)
\\[-2pt] &\;
-2^6\cdot3^7z^4(2\theta+1)^2(2257\theta^2+3601\theta+1942)
\\[-2pt] &\;
+2^{11}\cdot3^{11}z^5(2\theta+1)^2(2\theta+3)^2
\\[-8pt]
\endaligned &\cr
& &\span\hrulefill&\cr
& &&
A_n=\binom{2n}n\sum_k\binom nk\binom{2n}{2k}\binom{2k}{n-k}\binom{2n-2k}k
&\cr
\noalign{\hrule}
}}\hss}

\newpage

\hbox to\hsize{\hss\vbox{\offinterlineskip
\halign to120mm{\strut\tabskip=100pt minus 100pt
\strut\vrule#&\hbox to6.5mm{\hss$#$\hss}&%
\vrule#&\hbox to113mm{\hfil$\dsize#$\hfil}&%
\vrule#\tabskip=0pt\cr\noalign{\hrule}
& \# &%
& \text{differential operator $D$ and coefficients $A_n$, $n=0,1,2,\dots$} &\cr
\noalign{\hrule\vskip1pt\hrule}
& \eqnnoltmp{298}{358} &&
\aligned
\\[-12pt]
&
D=3^4\theta^4
-3^2z(1993\theta^4+3218\theta^3+2437\theta^2+828\theta+108)
\\[-2pt] &\;
+2^5z^2(17486\theta^4+25184\theta^3+12239\theta^2+2790\theta+297)
\\[-2pt] &\;
-2^8z^3(23620\theta^4+34776\theta^3+28905\theta^2+12447\theta+2106)
\\[-2pt] &\;
+2^{15}z^4(2\theta+1)(340\theta^3+618\theta^2+455\theta0129)
\\[-2pt] &\;
-2^{22}z^5(\theta+1)^2(2\theta+1)(2\theta+3)
\\[-8pt]
\endaligned &\cr
& &\span\hrulefill&\cr
& &&
A_n=\sum_{k,l}\binom nk\binom nl\binom{k+l}k\binom{k+l}n^2\binom{2n}{k+l}
&\cr
\noalign{\hrule}
& \eqnnoltmp{299}{359} &&
\aligned
\\[-12pt]
&
D=\theta^4
-6z(144\theta^4+36\theta^3+47\theta^2+29\theta+6)
\\[-2pt] &\;
+6^2z^2(8376\theta^4+6648\theta^3+8157\theta^2+3900\theta+724)
\\[-2pt] &\;
-6^4z^3(42672\theta^4+68616\theta^3+81056\theta^2+44841\theta+9964)
\\[-2pt] &\;
+2^6\cdot3^5z^4(374028\theta^4+962040\theta^3+1262091\theta^2+794463\theta+195335)
\\[-2pt] &\;
-2^8\cdot3^7z^5(633840\theta^4+2243328\theta^3
+3405968\theta^2+2385208\theta+6529129)
\\[-2pt] &\;
+2^{12}\cdot3^8z^6(438960\theta^4+1884384\theta^3
+3176664\theta^2+2380392\theta+652943)
\\[-2pt] &\;
-2^{19}\cdot3^{10}z^7(5760\theta^4+25128\theta^3+39548\theta^2+26606\theta+6517)
\\[-2pt] &\;
+12^{11}z^8(6\theta+5)^2(6\theta+7)^2
\\[-8pt]
\endaligned &\cr
& &\span\hrulefill&\cr
& &&
\aligned
&
A_n=(-1)^n\binom{2n}n\sum_{k=[n/3]}^{[2n/3]}\binom nk^3
\binom{2k}k\binom{2n-2k}{n-k}\binom{3k}n\binom{3n-3k}n
\binom{2n}{2k}^{-1}
\\[-2pt] &\quad\times
\bigl(1+k(-5H_k+5H_{n-k}+4H_{2k}-4H_{2n-2k}
\\[-2pt] &\quad\quad
+3H_{3k}-3H_{3k-n}+3H_{2n-3k}-3H_{3n-3k})\bigr)
\\[-2pt] &\;
+3\binom{2n}n\sum_{k=0}^{[(n-1)/3]}(-1)^k\frac{n-2k}{n-3k}
\binom nk^3\binom{2k}k\binom{2n-2k}{n-k}
\\[-2pt] &\;\quad\times
\binom{3n-3k}n\binom{2n}{2k}^{-1}\binom n{3k}^{-1}
\\[2pt]
\endaligned &\cr
\noalign{\hrule}
& \eqnnoltmp{300}{360} &&
\aligned
\\[-12pt]
&
D=\theta^4
+2^4z(371\theta^4+862\theta^3+591\theta^2+160\theta+15)
\\[-2pt] &\;
+5\cdot2^{11}z^2(224\theta^4+2069\theta^3+3277\theta^2+1363\theta+159)
\\[-2pt] &\;
-5^2\cdot2^{16}z^3(2089\theta^4+7500\theta^3+5533\theta^2+1500\theta+135)
\\[-2pt] &\;
+5^3\cdot2^{23}z^4(5\theta+1)(5\theta+2)(5\theta+3)(5\theta+4)
\\[-8pt]
\endaligned &\cr
& &\span\hrulefill&\cr
& &&
A_n=\binom{2n}n\sum_k(-1)^k4^{n-k}\binom nk\binom{2n+k}n
\binom{2k}n\binom{4n+2k}{2n+k}
&\cr
\noalign{\hrule}
& \eqnnoltmp{301}{361} &&
\aligned
\\[-12pt]
&
D=11^2\theta^4
-11z(1517\theta^4+3136\theta^3+2393\theta^2+825\theta+110)
\\[-2pt] &\;
-z^2(90362\theta^4+207620\theta^3+202166\theta^2+106953\theta+24266)
\\[-2pt] &\;
-z^3(245714\theta^4+507996\theta^3+415082\theta^2+217437\theta+53130)
\\[-2pt] &\;
-z^4(407863\theta^4+785972\theta^3+564786\theta^2+183269\theta+15226)
\\[-2pt] &\;
-z^5(434831\theta^4+790148\theta^3+728323\theta^2+279826\theta+25160)
\\[-2pt] &\;
-2^3z^6(36361\theta^4+70281\theta^3+73343\theta^2+37947\theta+7644)
\\[-2pt] &\;
-2^4\cdot5z^7(1307\theta^4+3430\theta^3+3877\theta^2+2162\theta+488)
\\[-2pt] &\;
-2^9\cdot5^2z^8(\theta+1)^4
\\[-8pt]
\endaligned &\cr
& &\span\hrulefill&\cr
& &&
A_n=\sum_{k,l}\binom nk\binom nl\binom{n+k}n\binom{k+l}k
\binom{2l}n\binom n{l-k}
&\cr
\noalign{\hrule}
}}\hss}

\newpage

\hbox to\hsize{\hss\vbox{\offinterlineskip
\halign to120mm{\strut\tabskip=100pt minus 100pt
\strut\vrule#&\hbox to6.5mm{\hss$#$\hss}&%
\vrule#&\hbox to113mm{\hfil$\dsize#$\hfil}&%
\vrule#\tabskip=0pt\cr\noalign{\hrule}
& \# &%
& \text{differential operator $D$ and coefficients $A_n$, $n=0,1,2,\dots$} &\cr
\noalign{\hrule\vskip1pt\hrule}
& \eqnnoltmp{302}{362} &&
\aligned
\\[-12pt]
&
D=5^2\theta^4
+5z(1307\theta^4+1798\theta^3+1429\theta^2+530\theta+80)
\\[-2pt] &\;
+2^4z^2(36361\theta^4+75163\theta^3+80666\theta^2+43340\theta+9120)
\\[-2pt] &\;
+2^6z^3(434831\theta^4+949176\theta^3+966865\theta^2+545700\theta+118340)
\\[-2pt] &\;
+2^{11}z^4(407863\theta^4+845480\theta^3+654048\theta^2+219839\theta+18634)
\\[-2pt] &\;
+2^{16}z^5(245714\theta^4+474860\theta^3+365378\theta^2+71595\theta-11507)
\\[-2pt] &\;
+2^{21}z^6(90362\theta^4+153828\theta^3+121478\theta^2+35967\theta+2221)
\\[-2pt] &\;
+2^{26}\cdot11z^7(1517\theta^4+2932\theta^3+2087\theta^2+621\theta+59)
\\[-2pt] &\;
-2^{31}\cdot11^2z^8(\theta+1)^4
\\[-8pt]
\endaligned &\cr
& &\span\hrulefill&\cr
& &&
\gathered
\text{the reflection of \#301 at infinity}
\\[-8pt]
\endgathered &\cr
& &\span\hrulefill&\cr
& &&
\text{a formula for $A_n$ is not known}
&\cr
\noalign{\hrule}
& \eqnnoltmp{303}{363} &&
\aligned
\\[-12pt]
&
D=13^2\theta^4
-13z(1505\theta^4+2746\theta^3+2127\theta^2+754\theta+104)
\\[-2pt] &\;
+2^2z^2(22961\theta^4-2086\theta^3-55741\theta^2-41574\theta-9256)
\\[-2pt] &\;
+2^5z^3(7524\theta^4+28098\theta^3+16131\theta^2+2691\theta-52)
\\[-2pt] &\;
-2^7z^4(7241\theta^4+6214\theta^3+17522\theta^2+15423\theta+4146)
\\[-2pt] &\;
+2^8z^5(-6087\theta^4-1806\theta^3+3905\theta^2+3796\theta+1036)
\\[-2pt] &\;
+2^{10}z^6(553\theta^4+4062\theta^3+4405\theta^2+1752\theta+220)
\\[-2pt] &\;
+2^{14}z^7(82\theta^4+230\theta^3+275\theta^2+160\theta+37)
\\[-2pt] &\;
+2^{18}z^8(\theta+1)^4
\\[-8pt]
\endaligned &\cr
& &\span\hrulefill&\cr
& &&
A_n=\sum_{k,l}\binom nk\binom nl\binom{k+l}n
\binom{2k}n\binom{2l}n\binom n{l-k}
&\cr
\noalign{\hrule}
& \eqnnoltmp{304}{364} &&
\aligned
\\[-12pt]
&
D=\theta^4
+z(82\theta^4+98\theta^3+77\theta^2+28\theta+4)
\\[-2pt] &\;
-z^2(-553\theta^4+1850\theta^3+4463\theta^2+2916\theta+636)
\\[-2pt] &\;
-2^2z^3(6087\theta^4+22542\theta^3+27199\theta^2+14916\theta+3136)
\\[-2pt] &\;
-2^5z^4(7241\theta^4+22750\theta^3+42326\theta^2+29943\theta+7272)
\\[-2pt] &\;
+2^7z^5(7524\theta^4+1998\theta^3-23019\theta^2-24627\theta-7186)
\\[-2pt] &\;
+2^8z^6(22961\theta^4+93930\theta^3+88283\theta^2+28194\theta+1624)
\\[-2pt] &\;
-2^{10}\cdot13z^7(1505\theta^4+3274\theta^3+2919\theta^2+1282\theta+236)
\\[-2pt] &\;
+2^{14}\cdot13^2z^8(\theta+1)^4
\\[-8pt]
\endaligned &\cr
& &\span\hrulefill&\cr
& &&
\gathered
\text{the reflection of \#303 at infinity}
\\[-8pt]
\endgathered &\cr
& &\span\hrulefill&\cr
& &&
A_n=(-1)^n\sum_{k,l}\binom nk^2\binom nl\binom{n+l}n\binom{2n-2k}n\binom{2k}{n-l}
&\cr
\noalign{\hrule}
}}\hss}

\newpage

\hbox to\hsize{\hss\vbox{\offinterlineskip
\halign to120mm{\strut\tabskip=100pt minus 100pt
\strut\vrule#&\hbox to6.5mm{\hss$#$\hss}&%
\vrule#&\hbox to113mm{\hfil$\dsize#$\hfil}&%
\vrule#\tabskip=0pt\cr\noalign{\hrule}
& \# &%
& \text{differential operator $D$ and coefficients $A_n$, $n=0,1,2,\dots$} &\cr
\noalign{\hrule\vskip1pt\hrule}
& \eqnnoltmp{305}{365} &&
\aligned
\\[-12pt]
&
D=\theta^4
+2^4\cdot3z(81552\theta^4-94944\theta^3-53688\theta^2-6216\theta-379)
\\[-2pt] &\;
+2^{20}\cdot3z^2(1091952\theta^4-2917008\theta^3+1388032\theta^2+225284\theta+19545)
\\[-2pt] &\;
+2^{34}\cdot3^3\cdot7z^3(-207504\theta^4+221184\theta^3+157480\theta^2+52224\theta+5855)
\\[-2pt] &\;
+2^{59}\cdot3^5\cdot7z^4(2\theta+1)^2(3\theta+1)(3\theta+2)
\\[-8pt]
\endaligned &\cr
& &\span\hrulefill&\cr
& &&
\aligned
&
A_n=\binom{2n}n^2\sum_{k=0}^n\binom{n+2k}k\binom{3n-2k}{n-k}
\binom{2n+4k}{n+2k}\binom{6n-4k}{3n-2k}\binom{3n}{n+k}
\\[-2pt] &\quad\times
\bigl(1+k(-H_k+H_{n-k}-2H_{n+k}+2H_{2n-k}+2H_{3n-2k}
\\[-2pt] &\quad\quad
-2H_{n+2k}+4H_{2n+4k}-4H_{6n-4k})\bigr)
\\[-2pt] &\;
+\binom{2n}n^2\sum_{k=1}^{[(n-1)/2]}(-1)^k\frac{n+2k}k\binom{3n+2k}{n+k}
\binom{2n-4k}{n-2k}
\\[-2pt] &\;\quad\times
\binom{6n+4k}{3n+2k}\binom{3n}{n-k}\binom{n-k}{n-2k}^{-1}
\\[-2pt] &\;
+\binom{2n}n^2\sum_{k=[(n+1)/2]}^n(-1)^{n+k}\frac{n+2k}k\binom{3n+2k}{n+k}
\binom k{n-k}
\\[-2pt] &\;\quad\times
\binom{6n+4k}{3n+2k}\binom{3n}{n-k}\binom{4k-2n}{2k-n}^{-1}
\\[2pt]
\endaligned &\cr
\noalign{\hrule}
& \eqnnoltmp{306}{366} &&
\aligned
\\[-12pt]
&
D=3^2\theta^4
-3z(592\theta^4+1100\theta^3+829\theta^2+279\theta+36)
\\[-2pt] &\;
+z^2(13801\theta^4+6652\theta^3-18041\theta^2-14904\theta-3312)
\\[-2pt] &\;
-2z^3\theta(8461\theta^3-29160\theta^2-28365\theta-7236)
\\[-2pt] &\;
-2^2\cdot3\cdot7z^4(513\theta^4+864\theta^3+487\theta^2+64\theta-16)
\\[-2pt] &\;
+2^3\cdot3\cdot7^2z^5(\theta+1)^2(3\theta+2)(3\theta+4)
\\[-8pt]
\endaligned &\cr
& &\span\hrulefill&\cr
& &&
A_n=\sum_{k,l}\binom nk\binom nl\binom{n+k}n\binom{n+l}n
\binom{2l}n\binom n{l-k}
&\cr
\noalign{\hrule}
& \eqnnoltmp{307}{367} &&
\aligned
\\[-12pt]
&
D=11^2\theta^4
-11z(1083\theta^4+1590\theta^3+1257\theta^2+462\theta+66)
\\[-2pt] &\;
+2^2z^2(47008\theta^4+45904\theta^3-3251\theta^2-17094\theta-4851)
\\[-2pt] &\;
-2^4\cdot3z^3(31436\theta^4+86856\theta^3+160363\theta^2+122133\theta+30294)
\\[-2pt] &\;
+2^9\cdot3^2z^4(2\theta+1)(1252\theta^3+5442\theta^2+6767\theta+2625)
\\[-2pt] &\;
-2^{14}\cdot3^6z^5(\theta+1)^2(2\theta+1)(2\theta+3)
\\[-8pt]
\endaligned &\cr
& &\span\hrulefill&\cr
& &&
A_n=\binom{2n}n\sum_{k,l}\binom nk\binom nl
\binom n{k+l}^2\binom{k+l}k
&\cr
\noalign{\hrule}
& \eqnnoltmp{308}{368} &&
\aligned
\\[-12pt]
&
D=29^2\theta^4
-2\cdot29z(1318\theta^4+2336\theta^3+1806\theta^2+638\theta+87)
\\[-2pt] &\;
-2^2z^2(90996\theta^4+744384\theta^3+1267526\theta^2+791584\theta-168345)
\\[-2pt] &\;
+2^2\cdot5^2z^3(34172\theta^4+77256\theta^3-46701\theta^2-110403\theta-36540)
\\[-2pt] &\;
+2^4\cdot5^4z^4(2\theta+1)(68\theta^3+1842\theta^2+2899\theta+1215)
\\[-2pt] &\;
-2^6\cdot5^7z^5(\theta+1)^2(2\theta+1)(2\theta+3)
\\[-8pt]
\endaligned &\cr
& &\span\hrulefill&\cr
& &&
A_n=\binom{2n}n\sum_{k,l}\binom nk\binom kl
\binom{n+k-l}n\binom{2l}l\binom{2l}{k-l}
&\cr
\noalign{\hrule}
}}\hss}

\newpage

\hbox to\hsize{\hss\vbox{\offinterlineskip
\halign to120mm{\strut\tabskip=100pt minus 100pt
\strut\vrule#&\hbox to6.5mm{\hss$#$\hss}&%
\vrule#&\hbox to113mm{\hfil$\dsize#$\hfil}&%
\vrule#\tabskip=0pt\cr\noalign{\hrule}
& \# &%
& \text{differential operator $D$ and coefficients $A_n$, $n=0,1,2,\dots$} &\cr
\noalign{\hrule\vskip1pt\hrule}
& \eqnnoltmp{309}{369} &&
\aligned
\\[-12pt]
&
D=3^4\theta^4
-3^2z(1993\theta^4+3218\theta^3+2437\theta^2+828\theta+108)
\\[-2pt] &\;
+2^5z^2(17486\theta^4+25184\theta^3+12239\theta^2+2790\theta+297)
\\[-2pt] &\;
-2^8z^3(23620\theta^4+34776\theta^3+28905\theta^2+12447\theta+2106)
\\[-2pt] &\;
+2^{15}z^4(2\theta+1)(340\theta^3+618\theta^2+455\theta+129)
\\[-2pt] &\;
-2^{22}z^5(\theta+1)^2(2\theta+1)(2\theta+3)
\\[-8pt]
\endaligned &\cr
& &\span\hrulefill&\cr
& &&
A_n=\binom{2n}n\sum_{k,l}\binom nk^2\binom nl\binom{n+l}n\binom{2k}{n-l}
&\cr
\noalign{\hrule}
& \eqnnoltmp{310}{370} &&
\aligned
\\[-12pt]
&
D=11^2\theta^4
-11z(1673\theta^4+3046\theta^3+2337\theta^2+814\theta+110)
\\[-2pt] &\;
+2\cdot5z^2(19247\theta^4+28298\theta^3+13285\theta^2+3454\theta+660)
\\[-2pt] &\;
-2^2z^3(167497\theta^4+245982\theta^3+227451\theta^2+115434\theta+22968)
\\[-2pt] &\;
+2^3\cdot5^2z^4(4079\theta^4+10270\theta^3+11427\theta^2+6226\theta+1340)
\\[-2pt] &\;
-2^5\cdot5^4z^5(\theta+1)^2(4\theta+3)(4\theta+5)
\\[-8pt]
\endaligned &\cr
& &\span\hrulefill&\cr
& &&
A_n=\sum_{k,l}\binom nk^2\binom nl\binom{n+k}n
\binom{n+l}n\binom{2k}{n-l}
&\cr
\noalign{\hrule}
& \eqnnoltmp{311}{371} &&
\aligned
\\[-12pt]
&
D=13^2\theta^4
-13z(327\theta^4+1038\theta^3+857\theta^2+338\theta+52)
\\[-2pt] &\;
-2^4z^2(12848\theta^4+42008\theta^3+52082\theta^2+28548\theta+5707)
\\[-2pt] &\;
+2^{11}z^3(-122\theta^4+1872\theta^3+6341\theta^2+5772\theta+1547)
\\[-2pt] &\;
+2^{16}z^4(2\theta+1)(76\theta^3+426\theta^2+570\theta+227)
\\[-2pt] &\;
+2^{23}z^5(\theta+1)^2(2\theta+1)(2\theta+3)
\\[-8pt]
\endaligned &\cr
& &\span\hrulefill&\cr
& &&
A_n=\binom{2n}n\sum_{k,l}\binom nk^2\binom nl\binom{2l}n\binom{2k}{n+l}
&\cr
\noalign{\hrule}
& \eqnnoltmp{312}{372} &&
\aligned
\\[-12pt]
&
D=7^2\theta^4
-7z(39\theta^4+234\theta^3+201\theta^2+84\theta+14)
\\[-2pt] &\;
-2z^2(12073\theta^4+43222\theta^3+57461\theta^2+34328\theta+7756)
\\[-2pt] &\;
+2^2z^3(-28923\theta^4-48421\theta^3+33393\theta^2+80976\theta+32032)
\\[-2pt] &\;
+2^3\cdot13z^4(359\theta^4+9790\theta^3+20805\theta^2+15784\theta+4124)
\\[-2pt] &\;
+2^5\cdot3\cdot13^2z^5(\theta+1)^2(6\theta+5)(6\theta+7)
\\[-8pt]
\endaligned &\cr
& &\span\hrulefill&\cr
& &&
A_n=\sum_{k,l}\binom nk^2\binom nl\binom{n+k-l}{n-l}
\binom{2l}n\binom{2k}{n+l}
&\cr
\noalign{\hrule}
& \eqnnoltmp{313}{373} &&
\aligned
\\[-12pt]
&
D=\theta^4
-z(\theta+1)(285\theta^3+321\theta^2+128\theta+18)
\\[-2pt] &\;
+2z^2(-1640\theta^4-1322\theta^3+1337\theta^2+1178\theta+240)
\\[-2pt] &\;
+2^2\cdot3^2z^3(-213\theta^4+256\theta^3+286\theta^2+80\theta+5)
\\[-2pt] &\;
+2^3\cdot3^3z^4(2\theta+1)(22\theta^3+37\theta^2+24\theta+6)
\\[-2pt] &\;
+2^4\cdot3^3z^5(\theta+1)^2(2\theta+1)(2\theta+3)
\\[-8pt]
\endaligned &\cr
& &\span\hrulefill&\cr
& &&
A_n=\binom{2n}n\sum_{k,l}\binom nk^2\binom nl
\binom{n+k-l}{n-l}\binom{3k}{n+l}
&\cr
\noalign{\hrule}
}}\hss}

\newpage

\hbox to\hsize{\hss\vbox{\offinterlineskip
\halign to120mm{\strut\tabskip=100pt minus 100pt
\strut\vrule#&\hbox to6.5mm{\hss$#$\hss}&%
\vrule#&\hbox to113mm{\hfil$\dsize#$\hfil}&%
\vrule#\tabskip=0pt\cr\noalign{\hrule}
& \# &%
& \text{differential operator $D$ and coefficients $A_n$, $n=0,1,2,\dots$} &\cr
\noalign{\hrule\vskip1pt\hrule}
& \eqnnoltmp{314}{374} &&
\aligned
\\[-12pt]
&
D=2^4\theta^4
-2^2z(1282\theta^4+2618\theta^3+1909\theta^2+600\theta+72)
\\[-2pt] &\;
-3^2z^2(9503\theta^4+26810\theta^3+31755\theta^2+15944\theta+2936)
\\[-2pt] &\;
-3^4z^3(-15627\theta^4+18288\theta^3+91412\theta^2+53256\theta+9688)
\\[-2pt] &\;
+2\cdot3^6z^4(15106\theta^4+20300\theta^3-20421\theta^2-23443\theta-5907)
\\[-2pt] &\;
+2^2\cdot3^8z^5(-2072\theta^4+18256\theta^3+2563\theta^2-4626\theta-1495)
\\[-2pt] &\;
-2^2\cdot3^{10}z^6(6204\theta^4+360\theta^3-281\theta^2+1017\theta+434)
\\[-2pt] &\;
-2^5\cdot3^{12}z^7(2\theta+1)(100\theta^3+162\theta^2+95\theta+21)
\\[-2pt] &\;
+2^8\cdot3^{14}z^8(\theta+1)^2(2\theta+1)(2\theta+3)
\\[-8pt]
\endaligned &\cr
& &\span\hrulefill&\cr
& &&
A_n=\binom{2n}n\sum_{k,l}\binom nk^2\binom nl\binom{2n-l}n\binom{3k}{n+l}
&\cr
\noalign{\hrule}
& \eqnnoltmp{315}{375} &&
\aligned
\\[-12pt]
&
D=5^2\theta^4
-5^2z(239\theta^4+496\theta^3+368\theta^2+120\theta+15)
\\[-2pt] &\;
-3\cdot5z^2(3454\theta^4+6412\theta^3+4682\theta^2+2180\theta+490)
\\[-2pt] &\;
-3^2\cdot5z^3(1519\theta^4+7338\theta^3+14271\theta^2+8340\theta+1690)
\\[-2pt] &\;
-3^3z^4(-10358\theta^4+16622\theta^3+49763\theta^2+37900\theta+10210)
\\[-2pt] &\;
+3^4z^5(4610\theta^4+17630\theta^3-6785\theta^2-15140\theta-5155)
\\[-2pt] &\;
+3^5z^6(-1219\theta^4+6030\theta^3+6441\theta^2+1740\theta-160)
\\[-2pt] &\;
-2^2\cdot3^6z^7(162\theta^4+234\theta^3+65\theta^2-52\theta-25)
\\[-2pt] &\;
-2^4\cdot3^8z^8(\theta+1)^4
\\[-8pt]
\endaligned &\cr
& &\span\hrulefill&\cr
& &&
A_n=\sum_{k,l}\binom nk^2\binom nl\binom{2n-l}n
\binom{n+k-l}{n-l}\binom{3k}{n+l}
&\cr
\noalign{\hrule}
& \eqnnoltmp{316}{376} &&
\aligned
\\[-12pt]
&
D=11^2\theta^4
-2^2\cdot11z(1092\theta^4+2472\theta^3+1792\theta^3+561\theta+66)
\\[-2pt] &\;
-2^4z^2(124328\theta^4+168992\theta^3+24998\theta^2-12804\theta-3168)
\\[-2pt] &\;
-2^4\cdot3z^3(484016\theta^4+474144\theta^3+366952\theta^2+161832\theta+27027)
\\[-2pt] &\;
-2^{11}\cdot3^2z^4(2\theta+1)^2(964\theta^2+1360\theta+669)
\\[-2pt] &\;
-2^{16}\cdot3^4z^5(2\theta+1)^2(2\theta+3)^2
\\[-8pt]
\endaligned &\cr
& &\span\hrulefill&\cr
& &&
A_n=\binom{2n}n^2\sum_{k,l}\binom nk^2\binom nl\binom{3k}{n+l}
&\cr
\noalign{\hrule}
& \eqnnoltmp{317}{377} &&
\aligned
\\[-12pt]
&
D=2^4\theta^4
-2^2\cdot3z(162\theta^4+414\theta^3+335\theta^2+128\theta+20)
\\[-2pt] &\;
+3^3z^2(1219\theta^4+10906\theta^3+18963\theta^2+11824\theta+2708)
\\[-2pt] &\;
+3^5\cdot5z^3(-922\theta^4-162\theta^3+6403\theta^2+6576\theta+1964)
\\[-2pt] &\;
-3^7z^4(10358\theta^4+58054\theta^3+62251\theta^2+29672\theta+4907)
\\[-2pt] &\;
+3^9\cdot5z^5(1519\theta^4-1262\theta^3+1371\theta^2+4264\theta+1802)
\\[-2pt] &\;
+2\cdot3^{11}\cdot5z^6(1727\theta^4+3702\theta^3+3085\theta^2+882\theta+17)
\\[-2pt] &\;
+3^{13}\cdot5^2z^7(239\theta^4+460^3+314\theta^2+84\theta+6)
\\[-2pt] &\;
-3^{16}\cdot5^2z^8(\theta+1)^4
\\[-8pt]
\endaligned &\cr
& &\span\hrulefill&\cr
& &&
\gathered
\text{the reflection of \#315 at infinity}
\\[-8pt]
\endgathered &\cr
& &\span\hrulefill&\cr
& &&
\text{a formula for $A_n$ is not known}
&\cr
\noalign{\hrule}
}}\hss}

\newpage

\hbox to\hsize{\hss\vbox{\offinterlineskip
\halign to120mm{\strut\tabskip=100pt minus 100pt
\strut\vrule#&\hbox to6.5mm{\hss$#$\hss}&%
\vrule#&\hbox to113mm{\hfil$\dsize#$\hfil}&%
\vrule#\tabskip=0pt\cr\noalign{\hrule}
& \# &%
& \text{differential operator $D$ and coefficients $A_n$, $n=0,1,2,\dots$} &\cr
\noalign{\hrule\vskip1pt\hrule}
& \eqnnoltmp{318}{378} &&
\aligned
\\[-12pt]
&
D=5^2\theta^4
-5z(473\theta^4+892\theta^3+696\theta^2+250\theta+35)
\\[-2pt] &\;
-2z^2(-1973\theta^4+4636\theta^3+14417\theta^2+10895\theta+2745)
\\[-2pt] &\;
+2\cdot3^2z^3(343\theta^4+1920\theta^3+1147\theta^2-345\theta-320)
\\[-2pt] &\;
+3^4z^4(-83\theta^4+104\theta^3+458\theta^2+406\theta+114)
\\[-2pt] &\;
-3^8z^5(\theta+1)^4
\\[-8pt]
\endaligned &\cr
& &\span\hrulefill&\cr
& &&
A_n=\sum_{k,l}\binom nk^2\binom nl^2
\binom{n+k-l}{n-l}\binom{2k+l}{n+l}
&\cr
\noalign{\hrule}
& \eqnnoltmp{319}{379} &&
\aligned
\\[-12pt]
&
D=\theta^4
+z(83\theta^4+436\theta^3+352\theta^2+134\theta+21)
\\[-2pt] &\;
+2\cdot3^2z^2(-343\theta^4+548\theta^3+2555\theta^2+1749\theta+405)
\\[-2pt] &\;
-3^4z^3(3946\theta^4+25056\theta^3+22658\theta^2+7722\theta+684)
\\[-2pt] &\;
+3^8\cdot5z^4(473\theta^4+1000\theta^3+858\theta^2+358\theta+62)
\\[-2pt] &\;
-3^{12}\cdot5^2z^5(\theta+1)^4
\\[-8pt]
\endaligned &\cr
& &\span\hrulefill&\cr
& &&
\gathered
\text{the reflection of \#318 at infinity}
\\[-8pt]
\endgathered &\cr
& &\span\hrulefill&\cr
& &&
\text{a formula for $A_n$ is not known}
&\cr
\noalign{\hrule}
& \eqnnoltmp{320}{380} &&
\aligned
\\[-12pt]
&
D=11^2\theta^4
-11z(4843\theta^4+8918\theta^3+6505\theta^2+2046\theta+242)
\\[-2pt] &\;
+2^2z^2(312184\theta^4+343456\theta^3-23371\theta^2-73942\theta-14883)
\\[-2pt] &\;
-2^4z^3(511972\theta^4+256344\theta^3+144969\theta^2+78639\theta+15642)
\\[-2pt] &\;
+2^{11}z^4(2\theta+1)(1964\theta^3+3078\theta^2+1853\theta+419)
\\[-2pt] &\;
-2^{18}z^5(\theta+1)^2(2\theta+1)(2\theta+3)
\\[-8pt]
\endaligned &\cr
& &\span\hrulefill&\cr
& &&
A_n=\binom{2n}n\sum_{k,l}\binom nk^2\binom nl\binom{2n-l}n\binom{n+2k}{n+l}
&\cr
\noalign{\hrule}
& \eqnnoltmp{321}{381} &&
\aligned
\\[-12pt]
&
D=3^4\theta^4
-3^2z(191\theta^4+862\theta^3+683\theta^2+252\theta+36)
\\[-2pt] &\;
-2^5z^2(7225\theta^4+24835\theta^3+30634\theta^2+16173\theta+3069)
\\[-2pt] &\;
-2^8z^3(13251\theta^4+35856\theta^3+27641\theta^2+6966\theta+180)
\\[-2pt] &\;
-2^{12}\cdot5z^4(2\theta+1)(314\theta^3+363\theta^2+68\theta-31)
\\[-2pt] &\;
+2^{16}\cdot5^2z^5(\theta+1)^2(2\theta+1)(2\theta+3)
\\[-8pt]
\endaligned &\cr
& &\span\hrulefill&\cr
& &&
A_n=\binom{2n}n\sum_{k,l}\binom nk^2\binom nl
\binom{2l}n\binom{n+k}{n+l}
&\cr
\noalign{\hrule}
& \eqnnoltmp{322}{382} &&
\aligned
\\[-12pt]
&
D=3^2\theta^4
-3z(-5\theta^4+122\theta^3+100\theta^2+39\theta+6)
\\[-2pt] &\;
-z^2(8603\theta^4+32600\theta^3+41729\theta^2+23736\theta+5052)
\\[-2pt] &\;
-2^2z^3(33304\theta^4+108297\theta^3+122347\theta^2+61470\theta+11712)
\\[-2pt] &\;
-2^2z^4(180401\theta^4+547606\theta^3+638125\theta^2+339248\theta+69036)
\\[-2pt] &\;
-2^4z^5(94934\theta^4+298745\theta^3+355667\theta^2+189660\theta+38224)
\\[-2pt] &\;
-2^4z^6(73291\theta^4+204216\theta^3+190453\theta^2+68916\theta+6964)
\\[-2pt] &\;
-2^7\cdot3z^7(811\theta^4+1886\theta^3+1804\theta^2+861\theta+174)
\\[-2pt] &\;
-2^{10}\cdot3^2z^8(\theta+1)^4
\\[-8pt]
\endaligned &\cr
& &\span\hrulefill&\cr
& &&
A_n=\sum_{k,l}\binom nk^2\binom nl\binom{n+k-l}{n-l}
\binom{2l}n\binom{n+k}{n+l}
&\cr
\noalign{\hrule}
}}\hss}

\newpage

\hbox to\hsize{\hss\vbox{\offinterlineskip
\halign to120mm{\strut\tabskip=100pt minus 100pt
\strut\vrule#&\hbox to6.5mm{\hss$#$\hss}&%
\vrule#&\hbox to113mm{\hfil$\dsize#$\hfil}&%
\vrule#\tabskip=0pt\cr\noalign{\hrule}
& \# &%
& \text{differential operator $D$ and coefficients $A_n$, $n=0,1,2,\dots$} &\cr
\noalign{\hrule\vskip1pt\hrule}
& \eqnnoltmp{323}{383} &&
\aligned
\\[-12pt]
&
D=3^2\theta^4
+3z(811\theta^4+1358\theta^3+1012\theta^2+333\theta+42)
\\[-2pt] &\;
+z^2(73291\theta^4+88948\theta^3+17551\theta^2-7494\theta-2424)
\\[-2pt] &\;
+2^3z^3(94934\theta^4+80991\theta^3+29036\theta^2+5175\theta+420)
\\[-2pt] &\;
+2^4z^4(180401\theta^4+173998\theta^3+77713\theta^2+15788\theta+708)
\\[-2pt] &\;
+2^7z^5(33304\theta^4+24919\theta^3-2720\theta^2-8451\theta-2404)
\\[-2pt] &\;
+2^8z^6(8603\theta^4+1812\theta^3-4453\theta^2-3666\theta-952)
\\[-2pt] &\;
-2^{11}\cdot3z^7(5\theta^4+142\theta^3+296\theta^2+225\theta+60)
\\[-2pt] &\;
-2^{14}\cdot3^2z^8(\theta+1)^4
\\[-8pt]
\endaligned &\cr
& &\span\hrulefill&\cr
& &&
\gathered
\text{the reflection of \#322 at infinity}
\\[-8pt]
\endgathered &\cr
& &\span\hrulefill&\cr
& &&
A_n=(-1)^n\sum_{k,l}(-4)^{n-k}\binom nk\binom nl\binom kl\binom{2k}k\binom{n+2k}n\binom{n+k-l}{n-l}
&\cr
\noalign{\hrule}
& \eqnnoltmp{324}{384} &&
\aligned
\\[-12pt]
&
D=11^2\theta^4
-2^2\cdot11z(432\theta^4+624\theta^3+477\theta^2+165\theta+22)
\\[-2pt] &\;
+2^5z^2(-12944\theta^4+4736\theta^3-15491\theta^2-12914\theta-2860)
\\[-2pt] &\;
+2^45z^3(-10688\theta^4+114048\theta^3+159132\theta^2+83028\theta+15455)
\\[-2pt] &\;
-2^{11}\cdot5^2z^4(2\theta+1)(4\theta+3)(76\theta^2+189\theta+125)
\\[-2pt] &\;
+2^{14}\cdot5^2z^5(2\theta+1)(2\theta+3)(4\theta+3)(4\theta+5)
\\[-8pt]
\endaligned &\cr
& &\span\hrulefill&\cr
& &&
\gathered
\text{$A_n$ is the constant term of $S^{2n}$
(Batyrev \#11.7658), where}
\\[-2pt]
S=x+\frac1x+\frac yx+\frac zx+\frac x{yz}+\frac tx(1+y+z+yz)
+\frac xt\biggl(\frac1y+\frac1z+\frac1{yz}\biggr)
\endgathered &\cr
\noalign{\hrule}
& \eqnnoltmp{325}{385} &&
\aligned
\\[-12pt]
&
D=19^2\theta^4
-19z(4333\theta^4+6212\theta^3+4778\theta^2+1672\theta+228)
\\[-2pt] &\;
+z^2(4307495\theta^4+7600484\theta^3+6216406\theta^2+2802424\theta+530556)
\\[-2pt] &\;
-z^3(93739369\theta^4+213316800\theta^3+236037196\theta^2
\\[-2pt] &\;\quad
+125748612\theta+25260804)
\\[-2pt] &\;
+z^4(240813800\theta^4+778529200\theta^3+1041447759\theta^2
\\[-2pt] &\;\quad
+631802809\theta+138510993)
\\[-2pt] &\;
-2^2\cdot409z^5(\theta+1)(2851324\theta^3+100355\theta^2+11221241\theta+3481470)
\\[-2pt] &\;
+2^2\cdot3^2\cdot19^2\cdot409^2z^6(\theta+1)(\theta+2)(2\theta+1)(2\theta+5)
\\[-8pt]
\endaligned &\cr
& &\span\hrulefill&\cr
& &&
\gathered
\text{$A_n$ is the constant term of $S^{2n}$
(Batyrev \#11.7661), where}
\\[-2pt]
S=\frac1x+\frac yx+\frac xy+\frac zx+\frac xz+\frac{yz}x
+\frac x{yz}+\frac tx(1+y+yz)
+\frac xt\biggl(1+\frac1y+\frac1{yz}\biggr)
\endgathered &\cr
\noalign{\hrule}
& \eqnnoltmp{326}{386} &&
\aligned
\\[-12pt]
&
D=13^2\theta^4
-13z\theta(56\theta^3+178\theta^2+115\theta+26)
\\[-2pt] &\;
-z^2(28466\theta^4+109442\theta^3+165603\theta^2+117338\theta+32448)
\\[-2pt] &\;
-z^3(233114\theta^4+1257906\theta^3+2622815\theta^2+2467842\theta+872352)
\\[-2pt] &\;
-z^4(989585\theta^4+6852298\theta^3+17737939\theta^2+19969754\theta+8108448)
\\[-2pt] &\;
-z^5(\theta+1)(2458967\theta^3+18007287\theta^2+44047582\theta+35386584)
\\[-2pt] &\;
-9z^6(\theta+1)(\theta+2)(393163\theta^2+2539029\theta+4164444)
\\[-2pt] &\;
-297z^7(\theta+1)(\theta+2)(\theta+3)(8683\theta+34604)
\\[-2pt] &\;
+3^3\cdot11^2\cdot13\cdot17z^8(\theta+1)(\theta+2)(\theta+3)(\theta+4)
\\[-8pt]
\endaligned &\cr
& &\span\hrulefill&\cr
& &&
\gathered
\text{$A_n$ is the constant term of $S^n$
(Batyrev \#20.1454; Kreuzer $X_{117,114}^{37}$),}
\\[-2pt]
\aligned
\\[-13pt]
\text{where}\;
S&=x+\frac1x+y+\frac1y+z+\frac1z+t+\frac 1t
+\frac xz+\frac tx+\frac zy+\frac xy
\\[-2pt] &\quad
+\frac1{xy}+\frac1{yt}+\frac1{zt}+\frac z{xy}
+\frac{zt}x+\frac x{zt}+\frac{yt}x+\frac z{yt}
\endaligned
\endgathered &\cr
\noalign{\hrule}
}}\hss}

\newpage

\hbox to\hsize{\hss\vbox{\offinterlineskip
\halign to120mm{\strut\tabskip=100pt minus 100pt
\strut\vrule#&\hbox to6.5mm{\hss$#$\hss}&%
\vrule#&\hbox to113mm{\hfil$\dsize#$\hfil}&%
\vrule#\tabskip=0pt\cr\noalign{\hrule}
& \# &%
& \text{differential operator $D$ and coefficients $A_n$, $n=0,1,2,\dots$} &\cr
\noalign{\hrule\vskip1pt\hrule}
& \eqnnoltmp{327}{387} &&
\aligned
\\[-12pt]
&
D=29^2\theta^4
+2\cdot29z\theta(24\theta^3-198\theta^2-128\theta-29)
\\[-2pt] &\;
-2^2z^2(44284\theta^4+172954\theta^3+248589\theta^2+172057\theta+47096)
\\[-2pt] &\;
-2^2z^3(525708\theta^4+2414772\theta^3+4447643\theta^2+3839049\theta+1275594)
\\[-2pt] &\;
-2^3z^4(1415624\theta^4+7911004\theta^3+17395449\theta^2+17396359\theta+6496262)
\\[-2pt] &\;
-2^4z^5(\theta+1)(2152040\theta^3+12186636\theta^2+24179373\theta+16560506)
\\[-2pt] &\;
-2^5z^6(\theta+1)(\theta+2)(1912256\theta^2+9108540\theta+11349571)
\\[-2pt] &\;
-2^8\cdot41z^7(\theta+1)(\theta+2)(\theta+39(5671\theta+16301)
\\[-2pt] &\;
-2^8\cdot3\cdot19\cdot41^2z^8(\theta+1)(\theta+2)(\theta+3)(\theta+4)
\\[-8pt]
\endaligned &\cr
& &\span\hrulefill&\cr
& &&
\gathered
\text{$A_n$ is the constant term of $S^n$
(Batyrev \#20.1268; Kreuzer $X_{116,116}^{41}$),}
\\[-2pt]
\aligned
\\[-13pt]
\text{where}\;
S&=x+\frac1x+y+\frac1y+z+\frac1z+t+\frac1t
+\frac xt+\frac ty+\frac tx+zt
\\[-2pt] &\quad
+\frac1{xy}+\frac1{xz}+\frac{zt}x+\frac{zt}y+\frac y{zt}
+\frac t{xy} +\frac{zt}{xy}+\frac{xy}{zt}
\endaligned
\endgathered &\cr
\noalign{\hrule}
& \eqnnoltmp{328}{388} &&
\aligned
\\[-12pt]
&
D=\theta^4-2^3z(33\theta^4+58\theta^3+44\theta^2+15\theta+2)
\\[-2pt] &\;
+2^6z^2(174\theta^4+448\theta^3+493\theta^2+262\theta+52)
\\[-2pt] &\;
-2^9z^3(2\theta+1)(166\theta^3+465\theta^2+477\theta+158)
\\[-2pt] &\;
+2^{13}\cdot3z^4(2\theta+1)(2\theta+3)(3\theta+2)(3\theta+4)
\\[-8pt]
\endaligned &\cr
& &\span\hrulefill&\cr
& &&
A_n=\binom{2n}n\sum_k(-1)^{n+k}4^{n-k}\binom nk\binom{2k}k
\binom{n+k}n\binom{2k}n
&\cr
\noalign{\hrule}
& \eqnnoltmp{329}{389} &&
\aligned
\\[-12pt]
&
D=\theta^4-2^4z(8\theta^4+34\theta^3+25\theta^2+8\theta+1)
\\[-2pt] &\;
-2^8z^2(87\theta^4+150\theta^3+32\theta^2+2\theta+1)
\\[-2pt] &\;
-2^{12}z^3(202\theta^4+240\theta^3+211\theta^2+102\theta+19)
\\[-2pt] &\;
-2^{16}\cdot3z^4(2\theta+1)(22\theta^3+45\theta^2+38\theta+12)
\\[-2pt] &\;
-2^{20}\cdot3^2z^5(\theta+1)^2(2\theta+1)(2\theta+3)
\\[-8pt]
\endaligned &\cr
& &\span\hrulefill&\cr
& &&
A_n=\binom{2n}n\sum_k(-1)^{n+k}4^{n-k}\binom nk\binom{2k}k\binom{2k}n^2
&\cr
\noalign{\hrule}
& \eqnnoltmp{330}{390} &&
\aligned
\\[-12pt]
&
D=\theta^4+2^4z(112\theta^4-64\theta^3-32\theta^2+1)
\\[-2pt] &\;
+2^{14}z^2(56\theta^4-64\theta^3+3\theta^2-10\theta-4)
\\[-2pt] &\;
+2^{20}z^3(32\theta^4-384\theta^3-436\theta^2-264\theta-55)
\\[-2pt] &\;
-2^{29}\cdot3z^4(2\theta+1)(10\theta+7)(2\theta^2+4\theta+3)
\\[-2pt] &\;
-2^{38}\cdot3^2z^5(\theta+1)^2(2\theta+1)(2\theta+3)
\\[-8pt]
\endaligned &\cr
& &\span\hrulefill&\cr
& &&
A_n=\binom{2n}n\sum_k(-1)^{n+k}4^{n-k}\binom nk^{-1}
\binom{2k}k^3\binom{2n-2k}{n-k}^2
&\cr
\noalign{\hrule}
& \eqnnoltmp{331}{391} &&
\aligned
\\[-12pt]
&
D=\theta^4+2^4z(-18\theta^4+48\theta^3+33\theta^2+9\theta+1)
\\[-2pt] &\;
-2^9z^2(86\theta^4+512\theta^3+125\theta^2+45\theta+10)
\\[-2pt] &\;
+2^{14}z^3(1138\theta^4+2040\theta^3+1883\theta^2+879\theta+157)
\\[-2pt] &\;
-2^{19}\cdot7z^4(2\theta+1)(186\theta^3+375\theta^2+317\theta+100)
\\[-2pt] &\;
+2^{27}\cdot7^2z^5(\theta+1)^2(2\theta+1)(2\theta+3)
\\[-8pt]
\endaligned &\cr
& &\span\hrulefill&\cr
& &&
A_n=\binom{2n}n\sum_k4^{n-k}\binom{2k}k^2\binom{2n-2k}{n-k}\binom{2k}n
&\cr
\noalign{\hrule}
}}\hss}

\newpage

\hbox to\hsize{\hss\vbox{\offinterlineskip
\halign to120mm{\strut\tabskip=100pt minus 100pt
\strut\vrule#&\hbox to6.5mm{\hss$#$\hss}&%
\vrule#&\hbox to113mm{\hfil$\dsize#$\hfil}&%
\vrule#\tabskip=0pt\cr\noalign{\hrule}
& \# &%
& \text{differential operator $D$ and coefficients $A_n$, $n=0,1,2,\dots$} &\cr
\noalign{\hrule\vskip1pt\hrule}
& \eqnnoltmp{332}{392} &&
\aligned
\\[-12pt]
&
D=3^2\theta^4
+2^2\cdot3z(67\theta^4+122\theta^3+100\theta^2+39\theta+6)
\\[-2pt] &\;
+2^5z^2(1172\theta^4+4298\theta^3+5831\theta^2+3315\theta+678)
\\[-2pt] &\;
+2^8z^3(302\theta^4+15912\theta^3+29314\theta^2+20925\theta+4926)
\\[-2pt] &\;
+2^{11}z^4(2\theta+1)(826\theta^3+3543\theta^2+4321\theta+1594)
\\[-2pt] &\;
+2^{16}z^5(2\theta+1)^2(4\theta+3)(4\theta+5)
\\[-8pt]
\endaligned &\cr
& &\span\hrulefill&\cr
& &&
A_n=2^n\binom{2n}n\sum_k(-1)^{n+k}4^{-k}\binom nk\binom{2k}k
\binom{2n-2k}n\binom{n+k}n
&\cr
\noalign{\hrule}
& \eqnnoltmp{333}{393} &&
\aligned
\\[-12pt]
&
D=\theta^4
+z\theta^2(71\theta^2-2\theta-1)
\\[-2pt] &\;
+2^3\cdot3z^2(154\theta^4+334\theta^3+461\theta^2+248\theta+48)
\\[-2pt] &\;
+2^6\cdot3^2z^3(5\theta+3)(31\theta^3+39\theta^2-25\theta-21)
\\[-2pt] &\;
+2^9\cdot3^6z^4(2\theta+1)(2\theta^3-33\theta^2-56\theta-24)
\\[-2pt] &\;
-12^6z^5(\theta+1)^2(2\theta+1)(2\theta+3)
\\[-8pt]
\endaligned &\cr
& &\span\hrulefill&\cr
& &&
A_n=\sum_k(-1)^k3^{2n-3k}\binom{2n}{3k}\binom{2k}n^2\frac{(3k)!}{k!^3}
&\cr
\noalign{\hrule}
& \eqnnoltmp{334}{394} &&
\aligned
\\[-12pt]
&
D=3^2\theta^4-3z(166\theta^4+320\theta^3+271\theta^2+111\theta+18)
\\[-2pt] &\;
+z^2(11155\theta^4+42652\theta^3+60463\theta^2+36876\theta+8172)
\\[-2pt] &\;
-9z^3(4705\theta^4+23418\theta^3+42217\theta^2+31152\theta+7932)
\\[-2pt] &\;
+12z^4(3514\theta^4+16403\theta^3+25581\theta^2+16442\theta+3744)
\\[-2pt] &\;
-20z^5(5\theta+3)(5\theta+4)(5\theta+6)(5\theta+7)
\\[-8pt]
\endaligned &\cr
& &\span\hrulefill&\cr
& &&
A_n=\sum_{k,l}(-1)^l3^{n-3l}\binom nk^2\binom n{3l}\binom{k+l}k\frac{(3l)!}{l!^3}
&\cr
\noalign{\hrule}
& \eqnnoltmp{335}{395} &&
\aligned
\\[-12pt]
&
D=\theta^4
-z(61\theta^4+122\theta^3+125\theta^2+64\theta+12)
\\[-2pt] &\;
-2^3z^2(193\theta^4+772\theta^3+1033\theta^2+522\theta+72)
\\[-2pt] &\;
+2^9\cdot3z^3(146\theta^4+876\theta^3+1838\theta^2+1572\theta+405)
\\[-2pt] &\;
-2^{12}\cdot3^2z^4(204\theta^4+1632\theta^3+4449\theta^2+4740\theta+1400)
\\[-2pt] &\;
+2^{16}\cdot3^3z^5(2\theta+5)^2(16\theta^2+80\theta+35)
\\[-2pt] &\;
-2^{19}\cdot3^4z^6(2\theta+1)(2\theta+5)(2\theta+7)(2\theta+11)
\\[-8pt]
\endaligned &\cr
& &\span\hrulefill&\cr
& &&
A_n=\binom{2n}n\sum_{k,l}(-8)^{n-k}\binom nk\binom kl^3\binom{3k}n
&\cr
\noalign{\hrule}
}}\hss}

\newpage

\hbox to\hsize{\hss\vbox{\offinterlineskip
\halign to120mm{\strut\tabskip=100pt minus 100pt
\strut\vrule#&\hbox to6.5mm{\hss$#$\hss}&%
\vrule#&\hbox to113mm{\hfil$\dsize#$\hfil}&%
\vrule#\tabskip=0pt\cr\noalign{\hrule}
& \# &%
& \text{differential operator $D$ and coefficients $A_n$, $n=0,1,2,\dots$} &\cr
\noalign{\hrule\vskip1pt\hrule}
& \eqnnoltmp{336}{396} &&
\aligned
\\[-12pt]
&
D=\theta^4
-z(205\theta^4+410\theta^3+305\theta^2+100\theta+12)
\\[-2pt] &\;
-2^5z^2(127\theta^4+508\theta^3+742^2+468\theta+99)
\\[-2pt] &\;
-2^2\cdot3z^3(2588\theta^4+15528\theta^3+32639\theta^2+28041\theta+7290)
\\[-2pt] &\;
-2^6\cdot3^2z^4(204\theta^4+1632\theta^3+4449\theta^2+4740\theta+1400)
\\[-2pt] &\;
-2^7\cdot3^3z^5(2\theta+5)^2(16\theta^2+80\theta+35)
\\[-2pt] &\;
-2^7\cdot3^4z^6(2\theta+1)(2\theta+5)(2\theta+7)(2\theta+11)
\\[-8pt]
\endaligned &\cr
& &\span\hrulefill&\cr
& &&
A_n=\binom{2n}n\sum_{k,l}\binom nk\binom kl^3\binom{3k}n
&\cr
\noalign{\hrule}
& \eqnnoltmp{337}{397} &&
\aligned
\\[-12pt]
&
D=5^2\theta^4
-5\cdot3z(3483\theta^4+6102\theta^3+4241\theta^2+1190\theta+120)
\\[-2pt] &\;
+2^5\cdot3^2z^2(31428\theta^4+35559\theta^3+243\theta^2-4320\theta-740)
\\[-2pt] &\;
-2^8\cdot3^5z^3(7371\theta^4+4860\theta^3+2997\theta^2+1080\theta+140)
\\[-2pt] &\;
+2^{13}\cdot3^8z^4(3\theta+1)^2(3\theta+2)^2
\\[-8pt]
\endaligned &\cr
& &\span\hrulefill&\cr
& &&
A_n=\binom{2n}n^2\sum_k(-1)^{n+k}4^{n-k}\binom nk\binom{2n+k}{2n}\binom{3n+2k}{n+k}
&\cr
\noalign{\hrule}
& \eqnnoltmp{338}{398} &&
\aligned
\\[-12pt]
&
D=3^2\theta^4
+2^2\cdot3z(278\theta^4+424\theta^3+311\theta^2+99\theta+12)
\\[-2pt] &\;
+2^3z^2(20840\theta^4+15776\theta^3-10540\theta^2-9732\theta-1968)
\\[-2pt] &\;
+2^8z^3(8190\theta^4-3528\theta^3-3991\theta^2-585\theta+114)
\\[-2pt] &\;
-2^{11}\cdot11z^4(2\theta+1)(86\theta^3+57\theta^2-39\theta-329)
\\[-2pt] &\;
+2^{15}\cdot11^2z^5(\theta+1)^2(2\theta+1)(2\theta+3)
\\[-8pt]
\endaligned &\cr
& &\span\hrulefill&\cr
& &&
A_n=\binom{2n}n\sum_k(-1)^k4^{n-k}\binom nk\binom{2k}k\binom{n+k}n\binom{n+2k}n
&\cr
\noalign{\hrule}
& \eqnnoltmp{339}{399} &&
\aligned
\\[-12pt]
&
D=\theta^4
-2^2z(10\theta^4+50\theta^3+39\theta^2+14\theta+2)
\\[-2pt] &\;
+2^4z^2(177\theta^4+1158\theta^3+2007\theta^2+1158\theta+230)
\\[-2pt] &\;
+2^8z^3(539\theta^4+1344\theta^3-300\theta^2-1068\theta-340)
\\[-2pt] &\;
+2^{10}\cdot5z^4(2\theta+1)(4\theta^3-642\theta^2-1002\theta-385)
\\[-2pt] &\;
-2^{13}\cdot3\cdot5^2z^5(2\theta+1)(2\theta+3)(3\theta+2)(3\theta+4)
\\[-8pt]
\endaligned &\cr
& &\span\hrulefill&\cr
& &&
A_n=2^{-n}\binom{2n}n\sum_k(-1)^k4^{n-k}\binom nk\binom{2k}k\binom{2n-2k}n\binom{n+2k}n
&\cr
\noalign{\hrule}
& \eqnnoltmp{340}{400} &&
\aligned
\\[-12pt]
&
D=3^2\theta^4
-2^2\cdot3z(124\theta^4+1064\theta^3+769\theta^2+237\theta+30)
\\[-2pt] &\;
+2^7z^2(-8092\theta^4-5848\theta^3+22175\theta^2+13869\theta+2751)
\\[-2pt] &\;
+2^{12}z^3(-5412\theta^4+92376\theta^3+67609\theta^2+15615\theta+96)
\\[-2pt] &\;
+2^{17}\cdot17z^4(2\theta+1)(2242\theta^3+1419\theta^2-1047\theta-733)
\\[-2pt] &\;
-2^{23}\cdot3\cdot17^2z^5(2\theta+1)(2\theta+3)(3\theta+2)(3\theta+4)
\\[-8pt]
\endaligned &\cr
& &\span\hrulefill&\cr
& &&
A_n=\binom{2n}n^2\sum_k4^{n-k}\binom nk^2\binom{n+k}n\binom{n+2k}k\binom{2n}{2k}^{-1}
&\cr
\noalign{\hrule}
}}\hss}

\newpage

\hbox to\hsize{\hss\vbox{\offinterlineskip
\halign to120mm{\strut\tabskip=100pt minus 100pt
\strut\vrule#&\hbox to6.5mm{\hss$#$\hss}&%
\vrule#&\hbox to113mm{\hfil$\dsize#$\hfil}&%
\vrule#\tabskip=0pt\cr\noalign{\hrule}
& \# &%
& \text{differential operator $D$ and coefficients $A_n$, $n=0,1,2,\dots$} &\cr
\noalign{\hrule\vskip1pt\hrule}
& \eqnnoltmp{341}{401} &&
\aligned
\\[-12pt]
&
D=13^2\theta^4
-13z(1217\theta^4+1474\theta^3+1127\theta^2+390\theta+52)
\\[-2pt] &\;
-2^4z^2(5134\theta^4+83956\theta^3+142024\theta^2+83616\theta+16575)
\\[-2pt] &\;
+2^6z^3(14292\theta^4+565032\theta^3+604615\theta^2+269841\theta+44070)
\\[-2pt] &\;
-2^{11}\cdot5z^4(2\theta+1)(4324\theta^3+10698\theta^2+9903\theta+3110)
\\[-2pt] &\;
+2^{16}\cdot3\cdot5^2z^5(2\theta+1)(2\theta+3)(3\theta+2)(3\theta+4)
\\[-8pt]
\endaligned &\cr
& &\span\hrulefill&\cr
& &&
A_n=\binom{2n}n\sum_k\binom nk^2\binom{2k}n\binom{3n-2k}n
&\cr
\noalign{\hrule}
& \eqnnoltmp{342}{402} &&
\aligned
\\[-12pt]
&
D=\theta^4
+2z(50\theta^4+64\theta^3+52\theta^2+20\theta+3)
\\[-2pt] &\;
+2^2\cdot3z^2(380\theta^4+992\theta^3+1166\theta^2+612\theta+117)
\\[-2pt] &\;
+2^2\cdot3^2z^3(2140\theta^4+5832\theta^3+5651\theta^22349\theta+360)
\\[-2pt] &\;
+2^4\cdot3^6z^4(2\theta+1)(20\theta^3+42\theta^2+35\theta+11)
\\[-2pt] &\;
+2^6\cdot3^7z^5(\theta+1)^2(2\theta+1)(2\theta+3)
\\[-8pt]
\endaligned &\cr
& &\span\hrulefill&\cr
& &&
A_n=\binom{2n}n\sum_{k,l}(-1)^{n+k}3^{n-3k}\binom n{3k}\binom{2k}l\frac{(3k)!}{k!^3}\binom nl
&\cr
\noalign{\hrule}
& \eqnnoltmp{343}{403} &&
\aligned
\\[-12pt]
&
D=\theta^4
+z(91\theta^4+116\theta^3+96\theta^2+38\theta+6)
\\[-2pt] &\;
+z^2(3649\theta^4+9388\theta^3+11076\theta^2+5950\theta+1218)
\\[-2pt] &\;
+z^3(63585\theta^4+203832\theta^3+258070\theta^2+148542\theta+32814)
\\[-2pt] &\;
+2z^4(244543\theta^4+938432\theta^3+1417427\theta^2+933049\theta+226317)
\\[-2pt] &\;
+2^2z^5(374407\theta^4+1908784\theta^3+3407293\theta^2+2501538\theta+653454)
\\[-2pt] &\;
+2^2\cdot3z^6(130530\theta^4+686256\theta^3+1382165\theta^2+1159645\theta+333030)
\\[-2pt] &\;
-2^3z^7(-276464\theta^4+92912\theta^3+3194335\theta^2+3755703\theta+1224450)
\\[-2pt] &\;
+2^4z^8(341712\theta^4+1614816\theta^3+1576879\theta^2+219863\theta-145632)
\\[-2pt] &\;
-2^5z^9(29968\theta^4+412128\theta^3+489227\theta^2+156573\theta-3258)
\\[-2pt] &\;
+2^8\cdot3z^{10}(6368\theta^4+13600\theta^3+11014\theta^2+4187\theta+681)
\\[-2pt] &\;
-2^{11}\cdot3^3z^{11}(\theta+1)^2(4\theta+3)(4\theta+5)
\\[-8pt]
\endaligned &\cr
& &\span\hrulefill&\cr
& &&
A_n=\sum_{k,l}(-1)^{n+k}3^{n-3k}\binom n{3k}\binom{2k}l\frac{(3k)!}{k!^3}\binom nl\binom{2n-l}n
&\cr
\noalign{\hrule}
& \eqnnoltmp{344}{404} &&
\aligned
\\[-12pt]
&
D=7^2\theta^4
-7z\theta(29\theta^3-50\theta^2-32\theta-7)
\\[-2pt] &\;
+3z^2\theta(1235\theta^3+512\theta^2+1165\theta+532)
\\[-2pt] &\;
-2\cdot3^2z^3(5373\theta^4+29040\theta^3+61493\theta^2+51786\theta+15876)
\\[-2pt] &\;
+2^2\cdot3^3z^4(10813\theta^4+68120\theta^3+160529\theta^2+154570\theta+53396)
\\[-2pt] &\;
-2^3\cdot3^4z^5(13929\theta^4+84348\theta^3+181015\theta^2+171080\theta+59172)
\\[-2pt] &\;
+2^5\cdot3^5z^6(6160\theta^4+35964\theta^3+69935\theta^2+58677\theta+18110)
\\[-2pt] &\;
-2^8\cdot3^6z^7(944\theta^4+5308\theta^3+10916\theta^2+9657\theta+3109)
\\[-2pt] &\;
+2^{11}\cdot3^7z^8(\theta+1)^2(96\theta^2+300\theta+265)
\\[-2pt] &\;
-2^{15}\cdot3^9z^9(\theta+1)^2(\theta+2)^2
\\[-8pt]
\endaligned &\cr
& &\span\hrulefill&\cr
& &&
A_n=\sum_{k,l}(-1)^{n+k}3^{n-3k}\binom n{3k}\binom{2k}l
\frac{(3k)!}{k!^3}\binom nl\binom{2l}n
&\cr
\noalign{\hrule}
}}\hss}

\newpage

\hbox to\hsize{\hss\vbox{\offinterlineskip
\halign to120mm{\strut\tabskip=100pt minus 100pt
\strut\vrule#&\hbox to6.5mm{\hss$#$\hss}&%
\vrule#&\hbox to113mm{\hfil$\dsize#$\hfil}&%
\vrule#\tabskip=0pt\cr\noalign{\hrule}
& \# &%
& \text{differential operator $D$ and coefficients $A_n$, $n=0,1,2,\dots$} &\cr
\noalign{\hrule\vskip1pt\hrule}
& \eqnnoltmp{345}{405} &&
\aligned
\\[-12pt]
&
D=11^2\theta^4
-11\cdot3z(113\theta^4+184\theta^3+158\theta^2+66\theta+11)
\\[-2pt] &\;
+2z^2(28397^4+95138\theta^3+128420\theta^2+77715\theta+17622)
\\[-2pt] &\;
-3z^3(3165\theta^4+180822\theta^3+560611\theta^2+539022\theta+167508)
\\[-2pt] &\;
-3z^4(233330\theta^4+1052614\theta^3+1424797\theta^2+774518\theta+145896)
\\[-2pt] &\;
-3^2z^5(12866\theta^4-98902\theta^3-52127\theta^2+102028\theta+63723)
\\[-2pt] &\;
+3^2z^6(183763\theta^4+473778\theta^3+427847\theta^2+147060\theta+11268)
\\[-2pt] &\;
-2^3\cdot3^3z^7(5006\theta^4+13414\theta^3+14935\theta^2+8228\theta+1869)
\\[-2pt] &\;
+2^6\cdot3^7z^8(\theta+1)^4
\\[-8pt]
\endaligned &\cr
& &\span\hrulefill&\cr
& &&
A_n=\sum_{k,l}(-1)^{n+k}3^{n-3k}\binom n{3k}\frac{(3k)!}{k!^3}\binom{k+l}k\binom nl\binom{2k}{n-l}
&\cr
\noalign{\hrule}
& \eqnnoltmp{346}{406} &&
\aligned
\\[-12pt]
&
D=8^2\theta^4
-8z(5006\theta^4+6610\theta^3+4729\theta^2+1424\theta+168)
\\[-2pt] &\;
+3^3z^2(183763\theta^4+261274\theta^3+109091\theta^2+22352\theta+2040)
\\[-2pt] &\;
-3^7z^3(12866\theta^4+150366\theta^3+321775\theta^2+141888\theta+21336)
\\[-2pt] &\;
+3^{10}z^4(-233330\theta^4+119294\theta^3+333065\theta^2+149446\theta+23109)
\\[-2pt] &\;
-3^{14}z^5(3165\theta^4-168162\theta^3+37135\theta^2+52394\theta+11440)
\\[-2pt] &\;
+2\cdot3^{17}z^6(28397\theta^4+18450\theta^3+13388\theta^2+7299\theta+1586)
\\[-2pt] &\;
+3^{22}\cdot11z^7(113\theta^4+268\theta^3+284\theta^2+150\theta+32)
\\[-2pt] &\;
+3^{25}\cdot11^2z^8(\theta+1)^4
\\[-8pt]
\endaligned &\cr
& &\span\hrulefill&\cr
& &&
\gathered
\text{the reflection of \#345 at infinity}
\\[-8pt]
\endgathered &\cr
& &\span\hrulefill&\cr
& &&
\text{a formula for $A_n$ is not known}
&\cr
\noalign{\hrule}
& \eqnnoltmp{347}{407} &&
\aligned
\\[-12pt]
&
D=\theta^4-3z(213\theta^4+186\theta^3+149\theta^2+56\theta+8)
\\[-2pt] &\;
+2^3\cdot3^3z^2(702\theta^4+1078\theta^3+949\theta^2+392\theta+60)
\\[-2pt] &\;
-2^6\cdot3^3z^3(9277\theta^4+18432\theta^3+16008\theta^2+6000\theta+840)
\\[-2pt] &\;
+2^{13}\cdot3^4\cdot5z^4(2\theta+1)^2(51\theta^2+69\theta+32)
\\[-2pt] &\;
-2^{14}\cdot3^6\cdot5^2z^5(2\theta+1)^2(2\theta+3)^2
\\[-8pt]
\endaligned &\cr
& &\span\hrulefill&\cr
& &&
A_0=1, \quad
A_n=6\binom{2n}n^2\sum_{k=0}^{[n/6]}\frac{n-2k}{5n-6k}\binom nk^2\binom{5n-6k}{4n}
&\cr
\noalign{\hrule}
& \eqnnoltmp{348}{408} &&
\aligned
\\[-12pt]
&
D=\theta^4
+2^2z(70\theta^4+194\theta^3+145\theta^2+48\theta+6)
\\[-2pt] &\;
+2^4\cdot3z^2(-141\theta^4+858\theta^3+2111\theta^2+1192\theta+206)
\\[-2pt] &\;
+2^8\cdot3^2z^3(-18\theta^4+324\theta^3+2364\theta^2+1953\theta+403)
\\[-2pt] &\;
-2^{10}\cdot3^4z^4(3\theta+1)(3\theta+2)(42\theta^2+258\theta+223)
\\[-2pt] &\;
+2^{14}\cdot3^6z^5(3\theta+1)(3\theta+2)(3\theta+4)(3\theta+5)
\\[-8pt]
\endaligned &\cr
& &\span\hrulefill&\cr
& &&
A_n=\binom{3n}n\sum_{k,l}(-1)^{k+l}\binom nk\binom nl\binom{2n}{k+l}\binom{2n-2k}n\binom{2l}n
&\cr
\noalign{\hrule}
}}\hss}

\newpage

\hbox to\hsize{\hss\vbox{\offinterlineskip
\halign to120mm{\strut\tabskip=100pt minus 100pt
\strut\vrule#&\hbox to6.5mm{\hss$#$\hss}&%
\vrule#&\hbox to113mm{\hfil$\dsize#$\hfil}&%
\vrule#\tabskip=0pt\cr\noalign{\hrule}
& \# &%
& \text{differential operator $D$ and coefficients $A_n$, $n=0,1,2,\dots$} &\cr
\noalign{\hrule\vskip1pt\hrule}
& \eqnnoltmp{349}{409} &&
\aligned
\\[-12pt]
&
D=7^2\theta^4
+2\cdot3\cdot7z(111\theta^4+120\theta^3+102\theta^2+42\theta+7)
\\[-2pt] &\;
+3z^2(67308\theta^4+136032\theta^3+153856\theta^2+83384\theta+17976)
\\[-2pt] &\;
+3^3z^3(178553\theta^4+439878\theta^3+528099\theta^2+313502\theta+74536)
\\[-2pt] &\;
+2\cdot3^3z^4(1355053\theta^4+3438698\theta^3+3854711\theta^2+2221354\theta+519896)
\\[-2pt] &\;
+2^2\cdot3^4z^5(2406561\theta^4+5708802\theta^3+5082043\theta^2+2161754\theta+336752)
\\[-2pt] &\;
+2^3\cdot3^5z^6(3133411\theta^4+6625998\theta^3+4266961\theta^2+238710\theta-485736)
\\[-2pt] &\;
+2^6\cdot3^6z^7(746186\theta^4+1366021\theta^3+743388\theta^2-203279\theta-212552)
\\[-2pt] &\;
+2^7\cdot3^7z^8(506499\theta^4+760668\theta^3+404459\theta^2-112958\theta-117216)
\\[-2pt] &\;
+2^{11}\cdot3^8z^9(27992\theta^4+34962\theta^3+7197\theta^2-14685\theta-7604)
\\[-2pt] &\;
+2^{14}\cdot3^9z^{10}(1381\theta^4+1244\theta^3-2460\theta^2-4030\theta-1571)
\\[-2pt] &\;
-2^{18}\cdot3^{10}z^{11}(\theta+1)^2(22\theta^2+98\theta+105)
\\[-2pt] &\;
-2^{22}\cdot3^{11}z^{12}(\theta+1)^2(\theta+2)^2
\\[-8pt]
\endaligned &\cr
& &\span\hrulefill&\cr
& &&
A_n=\sum_{k,l}(-1)^{n+k}3^{n-3k}\binom n{3k}\frac{(3k)!}{k!^3}\binom{2k}{n-l}\binom nl\binom{2l}n
&\cr
\noalign{\hrule}
& \eqnnoltmp{350}{410} &&
\aligned
\\[-12pt]
&
D=\theta^4
-z(289\theta^4+722\theta^3+545\theta^2+184\theta+24)
\\[-2pt] &\;
+2^3\cdot3z^2(214\theta^4+2734\theta^3+4861\theta^2+2640\theta+468)
\\[-2pt] &\;
+2^6\cdot3^2z^3(1391\theta^4+5184\theta^3+4252\theta^2+1296\theta+126)
\\[-2pt] &\;
+2^{10}\cdot3^6z^4(2\theta+1)^4
\\[-8pt]
\endaligned &\cr
& &\span\hrulefill&\cr
& &&
A_n=3\binom{2n}n^3\sum_{k=0}^{[n/3]}(-1)^k\frac{n-2k}{2n-3k}
\binom nk^2\binom{2n-2k}{n+k}\binom{2n-k}{2k}^{-1}
&\cr
\noalign{\hrule}
& \eqnnoltmp{351}{411} &&
\aligned
\\[-12pt]
&
D=\theta^4
+2^4z(22256\theta^4-38432\theta^3-23000\theta^2-3784\theta-321)
\\[-2pt] &\;
+2^{18}\cdot3^3z^2(1712\theta^4-18448\theta^3+8648\theta^2+2220\theta+279)
\\[-2pt] &\;
+2^{30}\cdot3^6z^3(-4624\theta^4+2304\theta^3+1672\theta^2+576\theta+63)
\\[-2pt] &\;
+2^{46}\cdot3^{10}z^4(2\theta+1)^4
\\[-8pt]
\endaligned &\cr
& &\span\hrulefill&\cr
& &&
\gathered
\text{the reflection of \#350 at infinity}
\\[-8pt]
\endgathered &\cr
& &\span\hrulefill&\cr
& &&
\text{a formula for $A_n$ is not known}
&\cr
\noalign{\hrule}
& \eqnnoltmp{352}{412} &&
\aligned
\\[-12pt]
&
D=\theta^4
-z(70\theta^4+86\theta^3+77\theta^2+34\theta+6)
\\[-2pt] &\;
+3z^2(675\theta^4+1602\theta^3+1933\theta^2+1130\theta+258)
\\[-2pt] &\;
-2^2\cdot3^3z^3(271\theta^4+888\theta^3+1259\theta^2+831\theta+207)
\\[-2pt] &\;
+2^2\cdot3^5z^4(212\theta^4+808\theta^3+1189\theta^2+773\theta+186)
\\[-2pt] &\;
-2^4\cdot3^7z^5(\theta+1)^2(4\theta+3)(4\theta+5)
\\[-8pt]
\endaligned &\cr
& &\span\hrulefill&\cr
& &&
A_n=3\sum_{k=0}^{[n/3]}(-1)^k\frac{n-2k}{2n-3k}\binom nk^2\binom{2k}k
\binom{2n-2k}{n-k}\binom{2n-3k}n
&\cr
\noalign{\hrule}
}}\hss}

\newpage

\hbox to\hsize{\hss\vbox{\offinterlineskip
\halign to120mm{\strut\tabskip=100pt minus 100pt
\strut\vrule#&\hbox to6.5mm{\hss$#$\hss}&%
\vrule#&\hbox to113mm{\hfil$\dsize#$\hfil}&%
\vrule#\tabskip=0pt\cr\noalign{\hrule}
& \# &%
& \text{differential operator $D$ and coefficients $A_n$, $n=0,1,2,\dots$} &\cr
\noalign{\hrule\vskip1pt\hrule}
& \eqnnoltmp{353}{413} &&
\aligned
\\[-12pt]
&
D=\theta^4
-2^2z(52\theta^4+40\theta^3+37\theta^2+17\theta+3)
\\[-2pt] &\;
+2^4z^2(960\theta^4+1536\theta^3+1512\theta^2+688\theta+123)
\\[-2pt] &\;
-2^8z^3(1792\theta^4+4608\theta^3+5184\theta^2+2816\theta-597)
\\[-2pt] &\;
-2^{14}z^4(4\theta+3)^2(4\theta+5)^2
\\[-8pt]
\endaligned &\cr
& &\span\hrulefill&\cr
& &&
\aligned
A_n&=3\binom{2n}n\sum_{k=0}^{[n/3]}(-1)^k\frac{n-2k}{2n-3k}\binom nk^3\binom{n+k}n\binom{2n-k}n
\\[-2pt] &\;\times
\binom{2n-2k}{n+k}\binom{2n-k}{2k}^{-1}\binom{2n}{2k}^{-1}
\\[2pt]
\endaligned &\cr
\noalign{\hrule}
& \eqnnoltmp{354}{414} &&
\aligned
\\[-12pt]
&
D=\theta^4-5z(170\theta^4+160\theta^3+125\theta^2+45\theta+6)
\\[-2pt] &\;
+3\cdot5^3z^2(725\theta^4+1220\theta^3+1105\theta^2+460\theta+68)
\\[-2pt] &\;
-3^2\cdot5^5z^3(1421\theta^4+3186\theta^3+3053\theta^2+1272\theta+188)
\\[-2pt] &\;
+2^2\cdot3^3\cdot5^7z^4(3\theta+1)(3\theta+2)(34\theta^2+61\theta+36)
\\[-2pt] &\;
-2^2\cdot3^4\cdot5^9z^5(3\theta+1)(3\theta+2)(3\theta+4)(3\theta+5)
\\[-8pt]
\endaligned &\cr
& &\span\hrulefill&\cr
& &&
A_0=1, \;
A_n=5\binom{2n}n\binom{3n}n\sum_{k=0}^{[n/5]}(-1)^k\frac{n-2k}{4n-5k}
\binom nk^2\binom{4n-5k}{3n}
&\cr
\noalign{\hrule}
& \eqnnoltmp{355}{415} &&
\aligned
\\[-12pt]
&
D=\theta^4
-z\bigl(344(\theta+\tfrac12)^{4}+326(\theta+\tfrac12)^2+11\bigr)
\\[-2pt] &\;
+z^2\bigl(43408(\theta+1)^4+51724(\theta+1)^2+4816\bigr)
\\[-2pt] &\;
-48z^3\bigl(49536(\theta+\tfrac32)^4+43504(\theta+\tfrac32)^2+861\bigr)
\\[-2pt] &\;
+2^{14}\cdot3^2z^4(3\theta+5)(3\theta+7)(6\theta+11)(6\theta+13)
\\[-8pt]
\endaligned &\cr
& &\span\hrulefill&\cr
& &&
\gathered
\text{the YY-pullback of the 5th-order differential equation $D'y=0$, where}
\\[-8pt]
\endgathered &\cr
& &\span\hrulefill&\cr
& &&
\aligned
\\[-12pt]
&
D'=\theta^5
-2z(2\theta+1)(43\theta^4+86\theta^3+77\theta^2+34\theta+6)
\\[-2pt] &\;
+48z^2(\theta+1)(2\theta+1)(2\theta+3)(6\theta+5)(6\theta+7)
\\[-8pt]
\endaligned &\cr
& &\span\hrulefill&\cr
& &&
\text{a formula for $A_n'$ is not known}
&\cr
\noalign{\hrule}
& \eqnnoltmp{356}{416} &&
\aligned
\\[-12pt]
&
D=\theta^4
-z\bigl(472(\theta+\tfrac12)^4+446(\theta+\tfrac12)^2+15\bigr)
\\[-2pt] &\;
+z^2\bigl(83344(\theta+1)^4+102060(\theta+1)^2+9252\bigr)
\\[-2pt] &\;
-z^3\bigl(6524928(\theta+\tfrac32)^4+5576448(\theta+\tfrac32)^2+69888\bigr)
\\[-2pt] &\;
+9216z^4(12\theta+19)(12\theta+23)(12\theta+25)(12\theta+29)
\\[-8pt]
\endaligned &\cr
& &\span\hrulefill&\cr
& &&
\gathered
\text{the YY-pullback of the 5th-order differential equation $D'y=0$, where}
\\[-8pt]
\endgathered &\cr
& &\span\hrulefill&\cr
& &&
\aligned
\\[-12pt]
&
D'=\theta^5
-2z(2\theta+1)(59\theta^4+118\theta^3+105\theta^2+46\theta+8)
\\[-2pt] &\;
+384z^2(\theta+1)(2\theta+1)(2\theta+3)(3\theta+2)(3\theta+4)
\\[-8pt]
\endaligned &\cr
& &\span\hrulefill&\cr
& &&
A_n'=4\binom{2n}{n}\sum_{k=0}^{[n/4]}\frac{n-2k}{3n-4k}\binom n{4k}
\binom{4n-4k}n^{-1}\binom{2n}{2k}^{-1}\frac{(4k)!}{k!^4}\,\frac{(4n-4k)!}{(n-k)!^4}
&\cr
\noalign{\hrule}
}}\hss}

\newpage

\hbox to\hsize{\hss\vbox{\offinterlineskip
\halign to120mm{\strut\tabskip=100pt minus 100pt
\strut\vrule#&\hbox to6.5mm{\hss$#$\hss}&%
\vrule#&\hbox to113mm{\hfil$\dsize#$\hfil}&%
\vrule#\tabskip=0pt\cr\noalign{\hrule}
& \# &%
& \text{differential operator $D$ and coefficients $A_n$, $n=0,1,2,\dots$} &\cr
\noalign{\hrule\vskip1pt\hrule}
& \eqnnoltmp{357}{417} &&
\aligned
\\[-12pt]
&
D=13^2\theta^4
-13z(441\theta^4+690\theta^3+631\theta^2+286\theta+52)
\\[-2pt] &\;
+2^4z^2(5121\theta^4+15576\theta^3+21215\theta^2+13702\theta+3445)
\\[-2pt] &\;
-2^{10}z^3(640\theta^4+2847\theta^3+5078\theta^2+4056\theta+1196)
\\[-2pt] &\;
+2^{14}z^4(125\theta^4+562\theta^3+905\theta^2+624\theta+157)
\\[-2pt] &\;
-2^{21}z^5(\theta+1)^4
\\[-8pt]
\endaligned &\cr
& &\span\hrulefill&\cr
& &&
A_0=1, \quad
A_n=4\sum_{k=0}^{[n/4]}\frac{n-2k}{3n-4k}\binom nk^4\binom{3n-4k}{2n}
&\cr
\noalign{\hrule}
& \eqnnoltmp{358}{418} &&
\aligned
\\[-12pt]
&
D=\theta^4
-2^4z(125\theta^4-62\theta^3-31\theta^2+1)
\\[-2pt] &\;
+2^{11}z^2(640\theta^4-287\theta^3+377\theta^2+119\theta+11)
\\[-2pt] &\;
-2^{16}z^3(5121\theta^4+4908\theta^3+5213\theta^2+2484\theta+503)
\\[-2pt] &\;
+2^{23}\cdot13z^4(441\theta^4+1074\theta^3+1207\theta^2+670\theta+148)
\\[-2pt] &\;
-2^{34}\cdot13^2z^5(\theta+1)^4
\\[-8pt]
\endaligned &\cr
& &\span\hrulefill&\cr
& &&
\gathered
\text{the reflection of \#357 at infinity}
\\[-8pt]
\endgathered &\cr
& &\span\hrulefill&\cr
& &&
\text{a formula for $A_n$ is not known}
&\cr
\noalign{\hrule}
& \eqnnoltmp{359}{419} &&
\aligned
\\[-12pt]
&
D=5^2\theta^4
+30z(51\theta^4+84\theta^3+72\theta^2+30\theta+5)
\\[-2pt] &\;
+2^2\cdot3z^2(3297\theta^4+10236\theta^3+13562\theta^2+8110\theta+1830)
\\[-2pt] &\;
+2^2\cdot3^3z^3(3866\theta^4+14088\theta^3+21137\theta^2+14355\theta+3600)
\\[-2pt] &\;
+2^3\cdot3^3z^4(11680\theta^4+38792\theta^3+45641\theta^2+24205\theta+4854)
\\[-2pt] &\;
+2^4\cdot3^5z^5(2624\theta^4+8240\theta^3+8275\theta^2+2971\theta+216)
\\[-2pt] &\;
+2^5\cdot3^5z^6(3248\theta^4+8832\theta^3+9739\theta^2+4803\theta+882)
\\[-2pt] &\;
+2^7\cdot3^7z^7(144\theta^4+384\theta^3+428\theta^2+233\theta+51)
\\[-2pt] &\;
+2^9\cdot3^7z^8(\theta+1)^2(4\theta+3)(4\theta+5)
\\[-8pt]
\endaligned &\cr
& &\span\hrulefill&\cr
& &&
A_{n}=\sum_{k,l}(-1)^{n+k}3^{n-3k}\binom{n}{3k}\binom{n}{l}
\binom{k}{n-l}\binom{2l}{n}\frac{(3k)!}{k!^{3}}
&\cr
\noalign{\hrule}
& \eqnnoltmp{360}{420} &&
\aligned
\\[-12pt]
&
D=17^2\theta^4-17z(-10622\theta^4+19904\theta^3+13913\theta^2+3961\theta+510)
\\[-2pt] &\;
+3^2z^2(1596891\theta^4-10821444\theta^3+10580847\theta^2+6358884\theta+1355036)
\\[-2pt] &\;
+3^5z^3(-5472387\theta^4+81131922\theta^3+52565469\theta^2+9898488\theta-1434596)
\\[-2pt] &\;
+2^2\cdot3^8\cdot127z^4(318018\theta^4+157911\theta^3-445563\theta^2-476706\theta-130792)
\\[-2pt] &\;
-2^2\cdot3^{12}\cdot5\cdot127^2z^5(5\theta+3)(5\theta+4)(5\theta+6)(5\theta+7)
\\[-8pt]
\endaligned &\cr
& &\span\hrulefill&\cr
& &&
\aligned
A_n
&=\sum_k\binom nk\binom{n+3k}n\binom{4n-3k}n\frac{(3k)!}{k!^3}\,\frac{(3n-3k)!}{(n-k)!^3}
\\[-2pt] &\quad\times
\bigl(1+k(-4H_{k}+4H_{n-k}+3H_{n+3k}-3H_{4n-3k})\bigr)
\\[2pt]
\endaligned &\cr
\noalign{\hrule}
& \eqnnoltmp{361}{421} &&
\aligned
\\[-12pt]
&
D=\theta^4-2^2z(80\theta^4+88\theta^3+67\theta^2+23\theta+3)
\\[-2pt] &\;
+2^4\cdot3z^2(928\theta^4+2080\theta^3+2176\theta^2+972\theta+153)
\\[-2pt] &\;
-2^{10}\cdot3^2z^3(272\theta^4+648\theta^3+511\theta^2+162\theta+18)
\\[-2pt] &\;
+2^{12}\cdot3^6z^4(2\theta+1)^4
\\[-8pt]
\endaligned &\cr
& &\span\hrulefill&\cr
& &&
A_0=1, \quad
A_n=3\binom{2n}n^2\sum_{k=0}^{[n/3]}(-1)^k\frac{n-2k}{2n-3k}
\binom nk\binom{2n}{2k}\binom{2n-3k}n
&\cr
\noalign{\hrule}
}}\hss}

\newpage

\hbox to\hsize{\hss\vbox{\offinterlineskip
\halign to120mm{\strut\tabskip=100pt minus 100pt
\strut\vrule#&\hbox to6.5mm{\hss$#$\hss}&%
\vrule#&\hbox to113mm{\hfil$\dsize#$\hfil}&%
\vrule#\tabskip=0pt\cr\noalign{\hrule}
& \# &%
& \text{differential operator $D$ and coefficients $A_n$, $n=0,1,2,\dots$} &\cr
\noalign{\hrule\vskip1pt\hrule}
& \eqnnoltmp{362}{422} &&
\aligned
\\[-12pt]
&
D=\theta^4+2^4z(-1088\theta^4+416\theta^3+212\theta^2+4\theta-3)
\\[-2pt] &\;
+2^{12}\cdot3^3z^2(928\theta^4-224\theta^3+448\theta^2+108\theta+9)
\\[-2pt] &\;
-2^{20}\cdot3^6z^3(320\theta^4+288\theta^3+220\theta^2+72\theta+9)
\\[-2pt] &\;
+2^{28}\cdot3^{10}z^4(2\theta+1)^4
\\[-8pt]
\endaligned &\cr
& &\span\hrulefill&\cr
& &&
\gathered
\text{the reflection of \#361 at infinity}
\\[-8pt]
\endgathered &\cr
& &\span\hrulefill&\cr
& &&
\text{a formula for $A_n$ is not known}
&\cr
\noalign{\hrule}
& \eqnnoltmp{363}{423} &&
\aligned
\\[-12pt]
&
D=\theta^4+3^2z(231\theta^4+318\theta^3+231\theta^2+72\theta+8)
\\[-2pt] &\;
+2^3\cdot3^5z^2(774\theta^4+1854\theta^3+1869\theta^2+768\theta+100)
\\[-2pt] &\;
+2^6\cdot3^8z^3(951\theta^4+2304\theta^3+1740\theta^2+504\theta+50)
\\[-2pt] &\;
+2^{10}\cdot3^{12}z^4(2\theta+1)^2(4\theta+1)(4\theta+3)
\\[-8pt]
\endaligned &\cr
& &\span\hrulefill&\cr
& &&
A_0=1, \;
A_n=3\frac{(4n)!}{n!^2(2n)!}\sum_{k=0}^{[n/3]}(-1)^{n+k}\frac{n-2k}{2n-3k}
\binom{2n-3k}n\binom{2n}k\binom{2n}{n-k}
&\cr
\noalign{\hrule}
& \eqnnoltmp{364}{424} &&
\aligned
\\[-12pt]
&
D=5^2\theta^4-5z(553\theta^4+722\theta^3+611\theta^2+250\theta+40)
\\[-2pt] &\;
+2^6z^2(1914\theta^4+4722\theta^3+5519\theta^2+3010\theta+610)
\\[-2pt] &\;
-2^{12}z^3(685\theta^4+2400\theta^3+3466\theta^2+2220\theta+500)
\\[-2pt] &\;
+2^{19}z^4(2\theta+1)(30\theta^3+105\theta^2+122\theta+46)
\\[-2pt] &\;
-2^{25}z^5(\theta+1)^2(2\theta+1)(2\theta+3)
\\[-8pt]
\endaligned &\cr
& &\span\hrulefill&\cr
& &&
A_{n}=4\binom{2n}n\sum_{k=0}^{[n/4]}\frac{n-2k}{3n-4k}\binom nk^3\binom{3n-4k}{2n}
&\cr
\noalign{\hrule}
& \eqnnoltmp{365}{425} &&
\aligned
\\[-12pt]
&
D=\theta^4-2^2z(99\theta^4+78\theta^3+65\theta^2+26\theta+4)
\\[-2pt] &\;
+2^6z^2(938\theta^4+1382\theta^3+1269\theta^2+554\theta+92)
\\[-2pt] &\;
-2^{10}z^3(4171\theta^4+8736\theta^3+8690\theta^2+3948\theta+680)
\\[-2pt] &\;
+2^{15}\cdot5z^4(2\theta+1)(418\theta^3+951\theta^2+846\theta+260)
\\[-2pt] &\;
-2^{19}\cdot3\cdot 5^2z^5(2\theta+1)(2\theta+3)(3\theta+2)(3\theta+4)
\\[-8pt]
\endaligned &\cr
& &\span\hrulefill&\cr
& &&
A_0=1, \;
A_n=4\binom{2n}n\sum_{k=0}^{[n/4]}\frac{n-2k}{3n-4k}\binom nk
\binom{2k}k\binom{2n-2k}{n-k}\binom{3n-4k}{2n}
&\cr
\noalign{\hrule}
& \eqnnoltmp{366}{426} &&
\aligned
\\[-12pt]
&
D=\theta^4+z\theta(39\theta^3-30\theta^2-19\theta-4)
\\[-2pt] &\;
+2z^2(16\theta^4-1070\theta^3-1057\theta^2-676\theta-192)
\\[-2pt] &\;
-2^23^2z^3(3\theta+2)(171\theta^3+566\theta^2+600\theta+316)
\\[-2pt] &\;
-2^53^3z^4(384\theta^4+1542\theta^3+2635\theta^2+2173\theta+702)
\\[-2pt] &\;
-2^63^3z^5(\theta+1)(1393\theta^3+5571\theta^2+8378\theta+4584)
\\[-2pt] &\;
-2^{10}3^5z^6(\theta+1)(\theta+2)(31\theta^2+105\theta+98)
\\[-2pt] &\;
-2^{12}3^7z^7(\theta+1)(\theta+2)^2(\theta+3)
\\[-8pt]
\endaligned &\cr
& &\span\hrulefill&\cr
& &&
\aligned
A_n
&=\sum_{i,j,k,l,m}\binom{2i}{i}\binom{2j}{j}\binom{2k}{k}\binom{l+m}{m}
\binom{2(n-i-j-k)}{n-i-j-k}\binom{n}{2(n-i-j-k)}
\\[-2pt] &\quad\times
\binom{2(n-i-j-k)}{l+m}\binom{2i+2j+2k-n}{n-2i-l-m}\binom{4i+2j+2k+l+m-2n}{2i+2j+m-n}
\\[2pt]
\endaligned &\cr
\noalign{\hrule}
}}\hss}

\newpage

\hbox to\hsize{\hss\vbox{\offinterlineskip
\halign to120mm{\strut\tabskip=100pt minus 100pt
\strut\vrule#&\hbox to6.5mm{\hss$#$\hss}&%
\vrule#&\hbox to113mm{\hfil$\dsize#$\hfil}&%
\vrule#\tabskip=0pt\cr\noalign{\hrule}
& \# &%
& \text{differential operator $D$ and coefficients $A_n$, $n=0,1,2,\dots$} &\cr
\noalign{\hrule\vskip1pt\hrule}
& \eqnnoltmp{367}{427} &&
\aligned
\\[-12pt]
&
D=3^2\theta^4-2^23z(760\theta^4+2048\theta^3+1423\theta^2+399\theta+42)
\\[-2pt] &\;
+2^7z^2(-20440\theta^4-25216\theta^3+4415\theta^2+4845\theta+795)
\\[-2pt] &\;
+2^{12}z^3(39928\theta^4+16512\theta^3+23719\theta^2+11637\theta+1830)
\\[-2pt] &\;
+2^{17}z^4(2928\theta^4-41856\theta^3-42871\theta^2-16873\theta-2425)
\\[-2pt] &\;
+2^{23}z^5(608\theta^4+3968\theta^3+10676\theta^2+6177\theta+1089)
\\[-2pt] &\;
+2^{29}z^6(272\theta^4+1056\theta^3+861\theta^2+264\theta+27)
\\[-2pt] &\;
+2^{35}z^7(2\theta+1)^4
\\[-8pt]
\endaligned &\cr
& &\span\hrulefill&\cr
& &&
A_{n}=\binom{2n}{n}\sum_{k}(-1)^{k}4^{n-k}\binom{n}{k}\binom{n+k}{n}\binom{n+2k}{n}\binom{2n+2k}{n+k}
&\cr
\noalign{\hrule}
& \eqnnoltmp{368}{428} &&
\aligned
\\[-12pt]
&
D=\theta^4+2^4z(1088\theta^4-2048\theta^3-1260\theta^2-236\theta-19)
\\[-2pt] &\;
+2^{15}z^2(1216\theta^4-5504\theta^3+11272\theta^2+3654\theta+423)
\\[-2pt] &\;
+2^{24}z^3(11712\theta^4+190848\theta^3+97220\theta^2+27432\theta+2835)
\\[-2pt] &\;
+2^{35}z^4(159712\theta^4+253376\theta^3+235372^2+78648\theta+9491)
\\[-2pt] &\;
+2^{46}z^5(-81760\theta^4-62656\theta^3+46316\theta^2+33048\theta+5403)
\\[-2pt] &\;
+2^{57}3z^6(3040\theta^4-2112\theta^3-2036\theta^2-528\theta-41)
\\[-2pt] &\;
+2^{69}3^2z^7(2\theta+1)^4
\\[-8pt]
\endaligned &\cr
& &\span\hrulefill&\cr
& &&
\gathered
\text{the reflection of \#367 at infinity}
\\[-8pt]
\endgathered &\cr
& &\span\hrulefill&\cr
& &&
\text{a formula for $A_n$ is not known}
&\cr
\noalign{\hrule}
& \eqnnoltmp{369}{429} &&
\aligned
\\[-12pt]
&
D=3^2\theta^4-3z(112\theta^4+140\theta^3+133\theta^2+63\theta+12)
\\[-2pt] &\;
+z^2(4393\theta^4+9340\theta^3+10903\theta^2+6360\theta+1488)
\\[-2pt] &\;
-2z^3(11669\theta^4+27720\theta^3+27019\theta^2+8460\theta-912)
\\[-2pt] &\;
+2^2z^4(6799\theta^4-10288\theta^3-82183\theta^2-119168\theta-52672)
\\[-2pt] &\;
-2^37z^5(\theta+1)(2611\theta^3+15537\theta^2+26998\theta+14360)
\\[-2pt] &\;
-2^67^2z^6(\theta+1)(\theta+2)(83\theta^2+105\theta-66)
\\[-2pt] &\;
-2^{10}7^3z^7(\theta+1)(\theta+2)^2(\theta+3)
\\[-8pt]
\endaligned &\cr
& &\span\hrulefill&\cr
& &&
\gathered
\text{the Hurwitz product (a)$\circ$(a)}
\\[-8pt]
\endgathered &\cr
& &\span\hrulefill&\cr
& &&
A_{n}=\sum_{k=0}^{n}\sum_{j=0}^{k}\sum_{l=0}^{n-k}\binom{n}{k}\binom{k}{j}^3\binom{n-k}{l}^3
&\cr
\noalign{\hrule}
& \eqnnoltmp{370}{430} &&
\aligned
\\[-12pt]
&
D=3^2\theta^4-3z(176\theta^4+220\theta^3+206\theta^2+96\theta+18)
\\[-2pt] &\;
+z^2(11692\theta^4+26440\theta^3+32164\theta^2+19632\theta+4824)
\\[-2pt] &\;
-z^3(123365\theta^4+374814\theta^3+519741\theta^2+346176\theta+89676)
\\[-2pt] &\;
+2z^4(309657\theta^4+1102938\theta^3+1591157\theta^2+1032920\theta+249740)
\\[-2pt] &\;
-2^3\cdot11z^5(\theta+1)(12897\theta^3+35469\theta^2+31181\theta+8042)
\\[-2pt] &\;
-2^3\cdot11^2z^6(\theta+1)(\theta+2)(355\theta^2+1047\theta+806)
\\[-2pt] &\;
-2^4\cdot11^3z^7(\theta+1)(\theta+2)^2(\theta+3)
\\[-8pt]
\endaligned &\cr
& &\span\hrulefill&\cr
& &&
\gathered
\text{the Hurwitz product (b)$\circ$(b)}
\\[-8pt]
\endgathered &\cr
& &\span\hrulefill&\cr
& &&
A_{n}=\sum_{k=0}^{n}\sum_{j=0}^{k}\sum_{l=0}^{n-k}\binom{n}{k}\binom{k}{j}^2\binom{k+j}{k}
\binom{n-k}{l}^2\binom{n-k+l}{n-k}
&\cr
\noalign{\hrule}
}}\hss}

\newpage

\hbox to\hsize{\hss\vbox{\offinterlineskip
\halign to120mm{\strut\tabskip=100pt minus 100pt
\strut\vrule#&\hbox to6.5mm{\hss$#$\hss}&%
\vrule#&\hbox to113mm{\hfil$\dsize#$\hfil}&%
\vrule#\tabskip=0pt\cr\noalign{\hrule}
& \# &%
& \text{differential operator $D$ and coefficients $A_n$, $n=0,1,2,\dots$} &\cr
\noalign{\hrule\vskip1pt\hrule}
& \eqnnoltmp{371}{431} &&
\aligned
\\[-12pt]
&
D=3^2\theta^4-3z(272\theta^4+340\theta^3+347\theta^2+177\theta+36)
\\[-2pt] &\;
+z^2(31273\theta^4+76540\theta^3+103783\theta^2+71112\theta+19728)
\\[-2pt] &\;
-2z^3(328219\theta^4+1181160\theta^3+1977957\theta^2+1620036\theta+522288)
\\[-2pt] &\;
+2^2z^4(2036999\theta^4+9602752\theta^3+19022113\theta^2+17726192\theta+6309408)
\\[-2pt] &\;
-2^317z^5(\theta+1)(439669\theta^3+2114103^2+3708554\theta+2306280)
\\[-2pt] &\;
+2^63^317^2z^6(\theta+1)(\theta+2)(481\theta^2+1875\theta+1962)
\\[-2pt] &\;
-2^{10}3^417^3z^7(\theta+1)(\theta+2)^2(\theta+3)
\\[-8pt]
\endaligned &\cr
& &\span\hrulefill&\cr
& &&
\gathered
\text{the Hurwitz product (g)$\circ$(g)}
\\[-8pt]
\endgathered &\cr
& &\span\hrulefill&\cr
& &&
A_{n}=\sum_{k=0}^{n}\sum_{i,j,l,m}(-1)^{i+l}8^{n-i-l}\binom{n}{k}\binom{k}{i}
\binom{i}{j}^3\binom{n-k}{l}\binom{l}{m}^3
&\cr
\noalign{\hrule}
& \eqnnoltmp{372}{432} &&
\aligned
\\[-12pt]
&
D=\theta^4-3z(96\theta^4+120^3+127\theta^2+67\theta+14)
\\[-2pt] &\;
+3^2z^2(3897\theta^4+9540\theta^3+13209^2+9246\theta+2608)
\\[-2pt] &\;
-2\cdot3^4z^3(14445\theta^4+52002\theta^3+88179\theta^2+73278\theta+23920)
\\[-2pt] &\;
+2^23^6z^4(31671\theta^4+149364\theta^3+298089\theta^2+280512\theta+100780)
\\[-2pt] &\;
-2^33^{12}z^5(\theta+1)(507\theta^3+2439\theta^2+4306\theta+2704)
\\[-2pt] &\;
+2^63^{14}z^6(\theta+1)(\theta+2)(90\theta^2+351\theta+370)
\\[-2pt] &\;
-2^73^{19}z^7(\theta+1)(\theta+2)^2(\theta+3)
\\[-8pt]
\endaligned &\cr
& &\span\hrulefill&\cr
& &&
\gathered
\text{the Hurwitz product (h)$\circ$(h)}
\\[-8pt]
\endgathered &\cr
& &\span\hrulefill&\cr
& &&
A_{n}=27^{n}\sum_{k=0}^{n}\sum_{j=0}^{k}\sum_{l=0}^{n-k}(-1)^{j+l}\binom{n}{k}
\binom{-2/3}{j}\binom{-2/3}{l}\binom{-1/3}{k-j}^2\binom{-1/3}{n-k-l}^2
&\cr
\noalign{\hrule}
& \eqnnoltmp{373}{433} &&
\aligned
\\[-12pt]
&
D=\theta^4-2z(190\theta^4+308\theta^3+227\theta^2+73\theta+9)
\\[-2pt] &\;
+2^2z^2(4780\theta^4+6304\theta^3+2395\theta^2+642\theta+135)
\\[-2pt] &\;
-2^43z^3(6700\theta^4+8472\theta^3+7607\theta^2+3615\theta+648)
\\[-2pt] &\;
+2^73^2z^4(2\theta+1)(760\theta^3+1464\theta^2+1211\theta+375)
\\[-2pt] &\;
-2^{10}3^6z^5(\theta+1)^2(2\theta+1)(2\theta+3)
\\[-8pt]
\endaligned &\cr
& &\span\hrulefill&\cr
& &&
A_{n}=\binom{2n}{n}\sum_{k,l}\binom{n}{k}^2\binom{n}{l}^2\binom{2k+2l}{2k}
&\cr
\noalign{\hrule}
}}\hss}

\newpage

\hbox to\hsize{\hss\vbox{\offinterlineskip
\halign to120mm{\strut\tabskip=100pt minus 100pt
\strut\vrule#&\hbox to6.5mm{\hss$#$\hss}&%
\vrule#&\hbox to113mm{\hfil$\dsize#$\hfil}&%
\vrule#\tabskip=0pt\cr\noalign{\hrule}
& \# &%
& \text{differential operator $D$ and coefficients $A_n$, $n=0,1,2,\dots$} &\cr
\noalign{\hrule\vskip1pt\hrule}
& \eqnnoltmp{374}{434} &&
\aligned
\\[-12pt]
&
D=97^2\theta^4-97z\theta(-1727\theta^3+2018\theta^2+1300\theta+291)
\\[-2pt] &\;
-z^2(1652135\theta^4+13428812\theta^3+16174393\theta^2+10216234\theta+2709792)
\\[-2pt] &\;
-3z^3(27251145\theta^4+121375398\theta^3+189546499\theta^2
\\[-2pt] &\;\quad
+147705198\theta+46000116)
\\[-2pt] &\;
-2z^4(587751431\theta^4+2711697232\theta^3+5003189285\theta^2
\\[-2pt] &\;\quad
+4434707760\theta+1524637512)
\\[-2pt] &\;
-z^5(9726250397\theta^4+50507429234\theta^3+106108023451\theta^2
\\[-2pt] &\;\quad
+103964102350\theta+38537290992)
\\[-2pt] &\;
-2\cdot3z^6(8793822649\theta^4+52062405804\theta^3+122175610025\theta^2
\\[-2pt] &\;\quad
+130254629814\theta+51340027968)
\\[-2pt] &\;
-2^2\cdot3^2z^7(5429262053\theta^4+36477756530\theta^3+94431307279\theta^2
\\[-2pt] &\;\quad
+108363704338\theta+44982230808)
\\[-2pt] &\;
-2^4\cdot3^2z^8(\theta+1)(3432647479\theta^3+22487363787\theta^2
\\[-2pt] &\;\quad
+50808614711\theta+38959393614)
\\[-2pt] &\;
-2^4\cdot3^3z^9(\theta+1)(\theta+2)(1903493629\theta^2+10262864555\theta+14314039440)
\\[-2pt] &\;
-2^5\cdot3^4\cdot13^2z^{10}(\theta+1)(\theta+2)(\theta+3)(1862987\theta+5992902)
\\[-2pt] &\;
-2^6\cdot3^3\cdot13^4\cdot7457z^{11}(\theta+1)(\theta+2)(\theta+3)(\theta+4)
\\[-8pt]
\endaligned &\cr
& &\span\hrulefill&\cr
& &&
\gathered
\text{$A_n$ is the constant term of $S^n$, where}
\\[-2pt]
\aligned
\\[-13pt]
S&=x+y+z+\frac{1}{x}+\frac{1}{y}+\frac{1}{z}+\frac{x}{z}+\frac{y}{x}+\frac{z}{y}
\\[-2pt] &\quad
+t\biggl(1+\frac{1}{x}+\frac{1}{y}+\frac{1}{xy}+\frac{z}{x}+\frac{z}{y}+\frac{z}{xy}\biggr)
\\[-2pt] &\quad
+\frac{1}{t}\biggl(1+x+y+\frac{1}{z}+\frac{x}{z}+\frac{y}{z}+\frac{xy}{z}\biggr)
\endaligned
\endgathered &\cr
\noalign{\hrule}
}}\hss}

\newpage

\hbox to\hsize{\hss\vbox{\offinterlineskip
\halign to120mm{\strut\tabskip=100pt minus 100pt
\strut\vrule#&\hbox to6.5mm{\hss$#$\hss}&%
\vrule#&\hbox to113mm{\hfil$\dsize#$\hfil}&%
\vrule#\tabskip=0pt\cr\noalign{\hrule}
& \# &%
& \text{differential operator $D$ and coefficients $A_n$, $n=0,1,2,\dots$} &\cr
\noalign{\hrule\vskip1pt\hrule}
& \eqnnoltmp{375}{435} &&
\aligned
\\[-12pt]
&
D=7^2\cdot13^2\theta^4-7\cdot13z\theta(-782\theta^3+1874\theta^2+1210\theta+273)
\\[-2pt] &\;
-z^2(2515785\theta^4+11622522\theta^3+15227939\theta^2+9962953\theta+2649920)
\\[-2pt] &\;
-z^3(59827597\theta^4+258678126\theta^3+432607868\theta^2
\\[-2pt] &\;\quad
+348819198\theta+110445426)
\\[-2pt] &\;
-2z^4(306021521\theta^4+1499440609\theta^32950997910\theta^2
\\[-2pt] &\;\quad
+2719866190\theta+957861945)
\\[-2pt] &\;
-3z^5(1254280114\theta^4+7075609686\theta^3+15834414271\theta^2
\\[-2pt] &\;\quad
+16174233521\theta+6159865002)
\\[-2pt] &\;
-z^6(15265487382\theta^4+98210309094\theta^3+244753624741\theta^2
\\[-2pt] &\;\quad
+271941545379\theta+110147546634)
\\[-2pt] &\;
-2z^7(21051636001\theta^4+152243816141\theta^3+415982528557\theta^2
\\[-2pt] &\;\quad
+495914741301\theta+211134581226)
\\[-2pt] &\;
-z^8(\theta+1)(39253400626\theta^3+275108963001\theta^2
\\[-2pt] &\;\quad
+654332416678\theta+521254338620)
\\[-2pt] &\;
-z^9(\theta+1)(\theta+2)(9498735517\theta^2+545340710193\theta+799002779040)
\\[-2pt] &\;
-2^2\cdot5\cdot7\cdot11z^{10}(\theta+1)(\theta+2)(\theta+3)(43765159\theta+149264765)
\\[-2pt] &\;
-2^23\cdot5^2\cdot7^2\cdot11^2\cdot11971z^{11}(\theta+1)(\theta+2)(\theta+3)(\theta+4)
\\[-8pt]
\endaligned &\cr
& &\span\hrulefill&\cr
& &&
\gathered
\text{$A_n$ is the constant term of $S^n$, where}
\\[-2pt]
\aligned
\\[-13pt]
S&=x+y+z+\frac{1}{x}+\frac{1}{y}+\frac{1}{z}+yz+\frac{1}{xy}+\frac{1}{xz}+\frac{1}{yz}+\frac{1}{xyz}
\\[-2pt] &\quad
+t(1+x+z+yz+xz+xyz)
\\[-2pt] &\quad
+\frac{1}{t}\biggl(1+\frac{1}{x}+\frac{1}{y}+\frac{1}{z}+\frac{1}{xy}+\frac{1}{xz}+\frac{1}{yz}+\frac{1}{xyz}\biggr)
\endaligned
\endgathered &\cr
\noalign{\hrule}
& \eqnnoltmp{376}{436} &&
\aligned
\\[-12pt]
&
D=2^4\theta^4-2^2z\theta(2\theta^3+82\theta^2+53\theta+12)
\\[-2pt] &\;
-z^2(4895\theta^4+18410\theta^3+26199\theta^2+18308\theta+5120)
\\[-2pt] &\;
-z^3(60679\theta^4+272424\theta^3+497452^2+430092\theta+143808)
\\[-2pt] &\;
-z^4(344527^4+1870838\theta^3+4034628\theta^2+3987101\theta+1478544)
\\[-2pt] &\;
-z^5(\theta+1)(1076509\theta^3+5847783\theta^3+11226106\theta+7492832)
\\[-2pt] &\;
-2z^6(\theta+1)(\theta+2)(944887\theta^2+4249317\theta+5045304)
\\[-2pt] &\;
-2^8\cdot13z^7(\theta+1)(\theta+2)(\theta+3)(518\theta+1381)
\\[-2pt] &\;
-2^5\cdot5\cdot13^2\cdot23z^8(\theta+1)(\theta+2)(\theta+3)(\theta+4)
\\[-8pt]
\endaligned &\cr
& &\span\hrulefill&\cr
& &&
\gathered
\text{$A_n$ is the constant term of $S^n$, where}
\\[-2pt]
\aligned
\\[-13pt]
S&=x+y+z+\frac{1}{x}+\frac{1}{y}+\frac{1}{z}+\frac{x}{y}+\frac{y}{x}+\frac{x}{z}+\frac{z}{x}+\frac{x}{yz}
\\[-2pt] &\quad
+t\biggl(1+z+\frac{1}{y}+\frac{z}{x}+\frac{z}{y}+\frac{z}{xy}\biggr)
+\frac{1}{t}\biggl(1+y+\frac{x}{z}+\frac{y}{z}+\frac{y}{x}+\frac{xy}{z}\biggr)
\endaligned
\endgathered &\cr
\noalign{\hrule}
& \eqnnoltmp{377}{437} &&
\aligned
\\[-12pt]
&
D=3^2\theta^4
-2^33z(61\theta^4+74\theta^3+58\theta^2+21\theta+3)
\\[-2pt] &\;
+2^4z^2(3883\theta^4+5356\theta^3+3451\theta^2+1278\theta+228)
\\[-2pt] &\;
-2^7z^3(8067\theta^4+13410\theta^3+12875\theta^2+6336\theta+1236)
\\[-2pt] &\;
+2^{14}z^4(413\theta^4+1069\theta^3+1206\theta^2+658\theta+140)
\\[-2pt] &\;
-2^{19}3z^5(\theta+1)^2(3\theta+2)(3\theta+4)
\\[-8pt]
\endaligned &\cr
& &\span\hrulefill&\cr
& &&
A_n=\sum_{k,l}(-1)^{n+k}4^{n-k}\binom{n}{k}\binom{n}{l}\binom{2k}{k}
\binom{k}{l}\binom{n+k}{n}\binom{n+k}{n}
&\cr
\noalign{\hrule}
}}\hss}

\newpage

\hbox to\hsize{\hss\vbox{\offinterlineskip
\halign to120mm{\strut\tabskip=100pt minus 100pt
\strut\vrule#&\hbox to6.5mm{\hss$#$\hss}&%
\vrule#&\hbox to113mm{\hfil$\dsize#$\hfil}&%
\vrule#\tabskip=0pt\cr\noalign{\hrule}
& \# &%
& \text{differential operator $D$ and coefficients $A_n$, $n=0,1,2,\dots$} &\cr
\noalign{\hrule\vskip1pt\hrule}
& \eqnnoltmp{378}{438} &&
\aligned
\\[-12pt]
&
D=3^2\theta^4
-2^23^2z(23\theta^4+58\theta^3+44\theta^2+15\theta+2)
\\[-2pt] &\;
-2^53z^2(254\theta^4+662\theta^3+623\theta^2+309\theta+66)
\\[-2pt] &\;
-2^8z^3(1707\theta^4+3276\theta^3+1806\theta^2+855\theta+234)
\\[-2pt] &\;
-2^{11}z^4(2266\theta^4+4076\theta^3+2167\theta^2+537\theta+18)
\\[-2pt] &\;
-2^{16}z^5(519\theta^4+798\theta^3+821\theta^2+391\theta+62)
\\[-2pt] &\;
-2^{19}z^6(305\theta^4+558\theta^3+625\theta^2+360\theta+82)
\\[-2pt] &\;
-2^{26}z^7(26\theta^4+70\theta^3+83\theta^2+48\theta+11)
\\[-2pt] &\;
-2^{29}z^8(\theta+1)^4
\\[-8pt]
\endaligned &\cr
& &\span\hrulefill&\cr
& &&
A_n=\sum_{k,l}(-1)^{n+k}4^{n-k}\binom{n}{k}\binom{n}{l}\binom{2k}{k}
\binom{k}{l}\binom{2k}{n}\binom{2l}{n}
&\cr
\noalign{\hrule}
& \eqnnoltmp{379}{439} &&
\aligned
\\[-12pt]
&
D=7^2\theta^4
-2\cdot 7z(452\theta^4+640\theta^3+509\theta^2+189\theta+28)
\\[-2pt] &\;
+2^2z^2(47156\theta^4+78224\theta^3+63963\theta^2+31010\theta+7000)
\\[-2pt] &\;
-2^5z^3(77224\theta^4+1509366\theta^3+155876\theta^2+86751\theta+19838)
\\[-2pt] &\;
+2^8z^4(65988\theta^4+160584\theta^3+193653\theta^2+117501\theta+28198)
\\[-2pt] &\;
-2^{12}z^5(15712\theta^4+46888\theta^3+63382\theta^2+41163\theta+10338)
\\[-2pt] &\;
+2^{16}z^6(2088\theta^4+7272\theta^3+10589\theta^2+7140\theta+1828)
\\[-2pt] &\;
-2^{22}z^7(36\theta^4+138\theta^3+206\theta^2+137\theta+34)
\\[-2pt] &\;
+2^{27}z^8(\theta+1)^4
\\[-8pt]
\endaligned &\cr
& &\span\hrulefill&\cr
& &&
A_n=\sum_{k,l}(-1)^{n+k}4^{n-k}\binom{n}{k}\binom{n}{l}\binom{2k}{k}
\binom{k}{l}\binom{k+l}{l}\binom{2l}{n}
&\cr
\noalign{\hrule}
& \eqnnoltmp{380}{440} &&
\aligned
\\[-12pt]
&
D=\theta^4
-2z(60\theta^4+90\theta^3+68\theta^2+23\theta+3)
\\[-2pt] &\;
+2^2z^2(313\theta^4-398\theta^3-1417\theta^2-1033\theta-252)
\\[-2pt] &\;
+2^3z^3(654\theta^4+5064\theta^3+3574\theta^2+129\theta-405)
\\[-2pt] &\;
+2^45z^4(-628\theta^4+40\theta^3+1699\theta^2+1661\theta+480)
\\[-2pt] &\;
-2^63\cdot 5^2z^5(\theta+1)^2(6\theta+5)(6\theta+7)
\\[-8pt]
\endaligned &\cr
& &\span\hrulefill&\cr
& &&
A_n=\sum_{k,l}(-1)^{n+k}4^{n-k}\binom{n}{k}\binom{n}{l}\binom{2k}{k}
\binom{k}{l}\binom{k+l}{l}\binom{n+l}{n}
&\cr
\noalign{\hrule}
& \eqnnoltmp{381}{441} &&
\aligned
\\[-12pt]
&
D=5^2\theta^4
+2^25z(19\theta^4+86\theta^3+73\theta^2+30\theta+5)
\\[-2pt] &\;
+2^4z^2(709\theta^4+4252\theta^3+7339\theta^2+4830\theta+1165)
\\[-2pt] &\;
+2^8z^3(-420\theta^4-114\theta^3+3294\theta^2+3960\theta+1325)
\\[-2pt] &\;
-2^{10}z^4(949\theta^4+6782\theta^3+11350\theta^2+7719\theta+1889)
\\[-2pt] &\;
+2^{12}z^5(1315\theta^4+4282\theta^3+7199\theta^2+5744\theta+1691)
\\[-2pt] &\;
+2^{14}z^6(613\theta^4+1560\theta^3+973\theta^2-216\theta-249)
\\[-2pt] &\;
+2^{18}z^7(11\theta^4-2\theta^3-40\theta^2-39\theta-11)
\\[-2pt] &\;
-2^{20}z^8(\theta+1)^4
\\[-8pt]
\endaligned &\cr
& &\span\hrulefill&\cr
& &&
A_n=\sum_{k,l}(-1)^{n+k}4^{n-k}\binom{n}{k}\binom{n}{l}\binom{2k}{k}
\binom{k}{l}\binom{2n-2k}{n}\binom{n+l}{n}
&\cr
\noalign{\hrule}
}}\hss}

\newpage

\hbox to\hsize{\hss\vbox{\offinterlineskip
\halign to120mm{\strut\tabskip=100pt minus 100pt
\strut\vrule#&\hbox to6.5mm{\hss$#$\hss}&%
\vrule#&\hbox to113mm{\hfil$\dsize#$\hfil}&%
\vrule#\tabskip=0pt\cr\noalign{\hrule}
& \# &%
& \text{differential operator $D$ and coefficients $A_n$, $n=0,1,2,\dots$} &\cr
\noalign{\hrule\vskip1pt\hrule}
& \eqnnoltmp{382}{442} &&
\aligned
\\[-12pt]
&
D=\theta^4
+2^2z(26\theta^4+34\theta^3+29\theta^2+12\theta+2)
\\[-2pt] &\;
+2^4z^2(305\theta^4+662\theta^3+781\theta^2+436\theta+94)
\\[-2pt] &\;
+2^8z^3(519\theta^4+1278\theta^3+1541\theta^2+933\theta+213)
\\[-2pt] &\;
-2^{10}z^4(2266\theta^4+4988\theta^3+3535\theta^2+633\theta-162)
\\[-2pt] &\;
+2^{14}3z^5(569\theta^4+1184\theta^3+740\theta^2-81\theta-128)
\\[-2pt] &\;
+2^{18}3z^6(254\theta^4+354\theta^3+161\theta^2-33\theta-28)
\\[-2pt] &\;
+2^{22}3^2z^7(23\theta^4+34\theta^3+8\theta^2-9\theta-4)
\\[-2pt] &\;
-2^{27}3^2z^8(\theta+1)^4
\\[-8pt]
\endaligned &\cr
& &\span\hrulefill&\cr
& &&
\gathered
\text{the reflection of \#378 at infinity}
\\[-8pt]
\endgathered &\cr
& &\span\hrulefill&\cr
& &&
\text{a formula for $A_n$ is not known}
&\cr
\noalign{\hrule}
& \eqnnoltmp{383}{443} &&
\aligned
\\[-12pt]
&
D=\theta^4
-2^5z(36\theta^4+6\theta^3+8\theta^2+5\theta+1)
\\[-2pt] &\;
+2^8z^2(2088\theta^4+1080\theta^3+1301\theta^2+574\theta+93)
\\[-2pt] &\;
-2^{13}z^3(15712\theta^4+15960\theta^3+16990\theta^2+7785\theta+1381)
\\[-2pt] &\;
+2^{18}z^4(65988\theta^4+103368\theta^3+107829\theta^2+52005\theta+9754)
\\[-2pt] &\;
-2^{24}z^5(77224\theta^4+157960\theta^3+166412\theta^2+81089\theta+15251)
\\[-2pt] &\;
+2^{30}z^6(47156\theta^4+110400\theta^3+112227\theta^2+50868\theta+8885)
\\[-2pt] &\;
-2^{38}7z^7(452\theta^4+1168\theta^3+1301\theta^2+717\theta+160)
\\[-2pt] &\;
+2^{46}7^2z^8(\theta+1)^4
\\[-8pt]
\endaligned &\cr
& &\span\hrulefill&\cr
& &&
\gathered
\text{the reflection of \#379 at infinity}
\\[-8pt]
\endaligned &\cr
& &\span\hrulefill&\cr
& &&
\text{a formula for $A_n$ is not known}
&\cr
\noalign{\hrule}
& \eqnnoltmp{384}{444} &&
\aligned
\\[-12pt]
&
D=\theta^4
-2^5z(11\theta^4+46\theta^3+32\theta^2+9\theta+1)
\\[-2pt] &\;
+2^8z^2(-613\theta^4-892\theta^3+29\theta^2+66\theta+7)
\\[-2pt] &\;
-2^{13}z^3(1315\theta^4+978\theta^3+2243\theta^2+1068\theta+179)
\\[-2pt] &\;
+2^{18}z^4(949\theta^4-2986\theta^3-3302\theta^2-1569\theta-313)
\\[-2pt] &\;
+2^{23}z^5(420\theta^4+1566\theta^3-1116\theta^2-1290\theta-353)
\\[-2pt] &\;
+2^{26}z^6(-709\theta^4+1416\theta^3+1163\theta^2+72\theta-131)
\\[-2pt] &\;
+2^{31}5z^7(-19\theta^4+10\theta^3+71\theta^2+66\theta+19)
\\[-2pt] &\;
-2^{36}5^2z^8(\theta+1)^4
\\[-8pt]
\endaligned &\cr
& &\span\hrulefill&\cr
& &&
\gathered
\text{the reflection of \#381 at infinity}
\\[-8pt]
\endaligned &\cr
& &\span\hrulefill&\cr
& &&
\text{a formula for $A_n$ is not known}
&\cr
\noalign{\hrule}
& \eqnnoltmp{385}{445} &&
\aligned
\\[-12pt]
&
D=\theta^4
-3z(42\theta^4+84\theta^3+77\theta^2+35\theta+6)
\\[-2pt] &\;
+3^3z^2(291\theta^4+1164\theta^3+1747\theta^2+1166\theta+264)
\\[-2pt] &\;
-2^2\cdot3^5z^3(360\theta^4+2160\theta^3+4553\theta^2+3939\theta+1035)
\\[-2pt] &\;
+2^3\cdot3^8z^4(204\theta^4+1632\theta^3+4449\theta^2+4740\theta+1400)
\\[-2pt] &\;
-2^4\cdot3^{11}z^5(2\theta+5)^2(16\theta^2+80\theta+35)
\\[-2pt] &\;
+2^4\cdot3^{14}z^6(2\theta+1)(2\theta+5)(2\theta+7)(2\theta+11)
\\[-8pt]
\endaligned &\cr
& &\span\hrulefill&\cr
& &&
A_n=\binom{2n}{n}\sum_{k,l}(-9)^{n-k}\binom nk{\binom kl}^{2}\binom{2l}l\binom{3k}n
&\cr
\noalign{\hrule}
}}\hss}

\newpage

\hbox to\hsize{\hss\vbox{\offinterlineskip
\halign to120mm{\strut\tabskip=100pt minus 100pt
\strut\vrule#&\hbox to6.5mm{\hss$#$\hss}&%
\vrule#&\hbox to113mm{\hfil$\dsize#$\hfil}&%
\vrule#\tabskip=0pt\cr\noalign{\hrule}
& \# &%
& \text{differential operator $D$ and coefficients $A_n$, $n=0,1,2,\dots$} &\cr
\noalign{\hrule\vskip1pt\hrule}
& \eqnnoltmp{386}{446} &&
\aligned
\\[-12pt]
&
D=\theta^4
-2z(422\theta^4+844\theta^3+751\theta^2+327\theta+57)
\\[-2pt] &\;
+2^2\cdot3^4z^2(\theta+1)^2(716\theta^2+1432\theta+579)
\\[-2pt] &\;
-2^4\cdot3^{8}\cdot7^2z^3(\theta+1)(\theta+2)(2\theta+1)(2\theta+5)
\\[-8pt]
\endaligned &\cr
& &\span\hrulefill&\cr
& &&
\text{a formula for $A_n$ is not known; \cite[\#\,$G_{7/9}$]{BR}}
&\cr
\noalign{\hrule}
& \eqnnoltmp{387}{447} &&
\aligned
\\[-12pt]
&
D=\theta^4
-2z(456\theta^4+912\theta^3+770\theta^2+314\theta+57)
\\[-2pt] &\;
+2^{11}z^2(\theta+1)^2(132\theta^2+264\theta+109)
\\[-2pt] &\;
-2^{18}\cdot5^2z^3(\theta+1)(\theta+2)(2\theta+1)(2\theta+5)
\\[-8pt]
\endaligned &\cr
& &\span\hrulefill&\cr
& &&
\text{a formula for $A_n$ is not known; \cite[\#\,$G_{5/4}$]{BR}}
&\cr
\noalign{\hrule}
& \eqnnoltmp{388}{448} &&
\aligned
\\[-12pt]
&
D=\theta^4
-2z(582\theta^4+1164\theta^3+815\theta^2+233\theta+25)
\\[-2pt] &\;
+4z^2(\theta+1)^2(2316\theta^2+4632\theta+1907)
\\[-2pt] &\;
-2^4\cdot17^2z^3(\theta+1)(\theta+2)(2\theta+1)(2\theta+5)
\\[-8pt]
\endaligned &\cr
& &\span\hrulefill&\cr
& &&
\text{a formula for $A_n$ is not known; \cite[\#\,$G_{17}$]{BR}}
&\cr
\noalign{\hrule}
& \eqnnoltmp{389}{449} &&
\aligned
\\[-12pt]
&
D=\theta^4
-2z(742\theta^4+1484\theta^3+1295\theta^2+553\theta+95)
\\[-2pt] &\;
+500z^2(\theta+1)^2(1468\theta^2+2936\theta+1211)
\\[-2pt] &\;
-2^4\cdot5^{6}\cdot11^2z^3(\theta+1)(\theta+2)(2\theta+1)(2\theta+5)
\\[-8pt]
\endaligned &\cr
& &\span\hrulefill&\cr
& &&
\text{a formula for $A_n$ is not known; \cite[\#\,$G_{11\sqrt{5}/25}$]{BR}}
&\cr
\noalign{\hrule}
& \eqnnoltmp{390}{450} &&
\aligned
\\[-12pt]
&
D=\theta^4
-z(561\theta^4+1122\theta^3+975\theta^2+414\theta+70)
\\[-2pt] &\;
+196z^2(\theta+1)^2(534\theta^2+10646\theta+433)
\\[-2pt] &\;
-2^2\cdot7^4\cdot13^2z^3(\theta+1)(\theta+2)(2\theta+1)(2\theta+5)
\\[-8pt]
\endaligned &\cr
& &\span\hrulefill&\cr
& &&
\text{a formula for $A_n$ is not known; \cite[\#\,$J_{13i\sqrt3/9}$]{BR}}
&\cr
\noalign{\hrule}
& \eqnnoltmp{391}{451} &&
\aligned
\\[-12pt]
&
D=\theta^4
-2z(6720\theta^4+11536\theta^3+8770\theta^2+3002\theta+372)
\\[-2pt] &\;
+2^{10}\cdot3^2z^2(4\theta+3)(1732\theta^3+4475\theta^2+3531\theta+645)
\\[-2pt] &\;
-2^{14}\cdot3^4\cdot17^2z^3(4\theta+1)(4\theta+3)(4\theta+7)(4\theta+9)
\\[-8pt]
\endaligned &\cr
& &\span\hrulefill&\cr
& &&
\text{a formula for $A_n$ is not known; \cite[\#\,$H_{17}$]{BR}}
&\cr
\noalign{\hrule}
& \eqnnoltmp{392}{452} &&
\aligned
\\[-12pt]
&
D=\theta^4
-2z(230\theta^4+446\theta^3+323\theta^2+75\theta+6)
\\[-2pt] &\;
-12z(6\theta+5)(1866\theta^3+5341\theta^2+4760\theta+1084)
\\[-2pt] &\;
-2^4\cdot3^2\cdot13^2z^3(3\theta+1)(3\theta+7)(6\theta+5)(6\theta+11)
\\[-8pt]
\endaligned &\cr
& &\span\hrulefill&\cr
& &&
\text{a formula for $A_n$ is not known; \cite[\#\,$G_{13i\sqrt3/9}$]{BR}}
&\cr
\noalign{\hrule}
& \eqnnoltmp{393}{453} &&
\aligned
\\[-12pt]
&
D=\theta^4
-2z(1264\theta^4+2240\theta^3+1792\theta^2+672\theta+96)
\\[-2pt] &\;
+768z(6\theta+5)(462\theta^3+1255\theta^2+1052\theta+235)
\\[-2pt] &\;
+2^{13}\cdot3^2\cdot5^2z^3(3\theta+1)(3\theta+7)(6\theta+5)(6\theta+11)
\\[-8pt]
\endaligned &\cr
& &\span\hrulefill&\cr
& &&
\text{a formula for $A_n$ is not known; \cite[\#\,$G_{5\sqrt3/9}$]{BR}}
&\cr
\noalign{\hrule}
}}\hss}

\newpage

\hbox to\hsize{\hss\vbox{\offinterlineskip
\halign to120mm{\strut\tabskip=100pt minus 100pt
\strut\vrule#&\hbox to6.5mm{\hss$#$\hss}&%
\vrule#&\hbox to113mm{\hfil$\dsize#$\hfil}&%
\vrule#\tabskip=0pt\cr\noalign{\hrule}
& \# &%
& \text{differential operator $D$ and coefficients $A_n$, $n=0,1,2,\dots$} &\cr
\noalign{\hrule\vskip1pt\hrule}
& \eqnnoltmp{394}{454} &&
\aligned
\\[-12pt]
&
D=3^4\theta^4
-3^3z(367\theta^4+398\theta^3+295\theta^2+96\theta+12)
\\[-2pt] &\;
-2^4\cdot3^3z^2(200\theta^4+2081^3+3614\theta^2+2009\theta+392)
\\[-2pt] &\;
+2^6\cdot3z^3(72449\theta^4+102684\theta^3-48579\theta^2-77922\theta-22536)
\\[-2pt] &\;
+2^{10}z^4(109873\theta^4+619970\theta^3+56260\theta^2-219027\theta-78216)
\\[-2pt] &\;
+2^{14}\cdot7z^5(-40669\theta^4+18266\theta^3+36570\theta^2+16190\theta+1955)
\\[-2pt] &\;
-2^{17}\cdot7z^6(80805\theta^4+76590\theta^3+51265^2+23076\theta+4780)
\\[-2pt] &\;
-2^{24}\cdot7^2z^7(437\theta^4+1117\theta^3+1236\theta^2+664\theta+140)
\\[-2pt] &\;
-2^{29}\cdot3\cdot7^2z^8(\theta+1)^2(3\theta+2)(3\theta+4)
\\[-8pt]
\endaligned &\cr
& &\span\hrulefill&\cr
& &&
A_n=\sum_{k,l}(-4)^{n-k}\binom nk\binom nl\binom kl\binom{2k}k\binom{n+2k}n\binom{2n-2l}n
&\cr
\noalign{\hrule}
& \eqnnoltmp{395}{455} &&
\aligned
\\[-12pt]
&
D=\theta^4
-2^2z\theta(22\theta^3+8\theta^2+5\theta+1)
\\[-2pt] &\;
+2^5z^2(34\theta^4-152\theta^3-265\theta^2-163\theta-36)
\\[-2pt] &\;
+2^8z^3(142\theta^4+600\theta^3+335\theta^2-39\theta-54)
\\[-2pt] &\;
+2^{11}\cdot3z^4(-68^4+56\theta^3+295\theta^2+261\theta+72)
\\[-2pt] &\;
-2^{15}3^2z^5(\theta+1)^2(4\theta+3)(4\theta+5)
\\[-8pt]
\endaligned &\cr
& &\span\hrulefill&\cr
& &&
A_n=\sum_{k,l}(-4)^{n-k}\binom nk\binom nl\binom kl\binom{2k}k\binom{n+k}n\binom{2n-2l}n
&\cr
\noalign{\hrule}
& \eqnnoltmp{396}{456} &&
\aligned
\\[-12pt]
&
D=5^2\theta^4
-2^2\cdot5z(197\theta^4+418\theta^3+319\theta^2+110\theta+15)
\\[-2pt] &\;
+2^4z^2(181\theta^4+5068^3+10291^2+6750\theta+1585)
\\[-2pt] &\;
+2^6z^3(-1727\theta^4+4758\theta^3+11365\theta^2+4560\theta+345)
\\[-2pt] &\;
+2^9z^4(2351^4+4552^3-11125\theta^2-12552\theta-3833)
\\[-2pt] &\;
-2^{12}z^5(527^4+1448\theta^3+16\theta^2-1811\theta-887)
\\[-2pt] &\;
+2^{15}z^6(493^4-1527\theta^3-789\theta^2-363\theta-116)
\\[-2pt] &\;
-2^{17}z^7(780\theta^4-282\theta^3+865\theta^2+1459\theta+563)
\\[-2pt] &\;
+2^{20}z^8(151\theta^4-104^3-291\theta^2-239\theta-65)
\\[-2pt] &\;
-2^{22}z^9(23\theta^4+24\theta^3+85\theta^2+132\theta+55)
\\[-2pt] &\;
+2^{25}z^{10}(\theta+1)(7\theta^3+31\theta^2+35\theta+12)
\\[-2pt] &\;
-2^{28}z^{11}(\theta+1)^4
\\[-8pt]
\endaligned &\cr
& &\span\hrulefill&\cr
& &&
A_n=\sum_{k,l}(-4)^{n-k}\binom nk\binom nl\binom kl\binom{2k}k\binom{n+k-l}{n-l}\binom{2k}n
&\cr
\noalign{\hrule}
}}\hss}

\newpage

\hbox to\hsize{\hss\vbox{\offinterlineskip
\halign to120mm{\strut\tabskip=100pt minus 100pt
\strut\vrule#&\hbox to6.5mm{\hss$#$\hss}&%
\vrule#&\hbox to113mm{\hfil$\dsize#$\hfil}&%
\vrule#\tabskip=0pt\cr\noalign{\hrule}
& \# &%
& \text{differential operator $D$ and coefficients $A_n$, $n=0,1,2,\dots$} &\cr
\noalign{\hrule\vskip1pt\hrule}
& \eqnnoltmp{397}{457} &&
\aligned
\\[-12pt]
&
D=\theta^4
-2^4z\theta(7\theta^3-10\theta^2-6\theta-1)
\\[-2pt] &\;
+2^8z^2(23\theta^4+68\theta^3+151\theta^2+58\theta+7)
\\[-2pt] &\;
-2^{13}z^3(151\theta^4+708\theta^3+927\theta^2+573\theta+138)
\\[-2pt] &\;
+2^{17}z^4(780\theta^4+3402\theta^3+6391\theta^2+4237\theta+1031)
\\[-2pt] &\;
-2^{22}z^5(493^4+3499\theta^3+6750\theta^2+5338\theta+1478)
\\[-2pt] &\;
+2^{26}z^6(527\theta^4+660\theta^3-1166\theta^2-393\theta+19)
\\[-2pt] &\;
+2^{30}z^7(-2351\theta^4-4852\theta^3+10675\theta^2+13950\theta+4607)
\\[-2pt] &\;
+2^{34}z^8(1727\theta^4+11666^3+13271\theta^2+3012\theta-665)
\\[-2pt] &\;
+2^{39}z^9(-181\theta^4+4344\theta^3+3827\theta^2+648\theta-239)
\\[-2pt] &\;
+2^{44}\cdot5z^{10}(197\theta^4+370\theta^3+247\theta^2+62\theta+3)
\\[-2pt] &\;
-2^{49}\cdot5^2z^{11}(\theta+1)^4
\\[-8pt]
\endaligned &\cr
& &\span\hrulefill&\cr
& &&
\gathered
\text{the reflection of \#396 at infinity}
\\[-8pt]
\endgathered &\cr
& &\span\hrulefill&\cr
& &&
\text{a formula for $A_n$ is not known}
&\cr
\noalign{\hrule}
& \eqnnoltmp{398}{458} &&
\aligned
\\[-12pt]
&
D=13^2\theta^4
-13z(5249\theta^4+4930\theta^3+3687\theta^2+1222\theta+156)
\\[-2pt] &\;
+2^4z^2(526803\theta^4+564192\theta^3+270729\theta^2+58266\theta+4641)
\\[-2pt] &\;
-2^7z^3(3336915\theta^4+3777024\theta^3+2377229\theta^2+746148\theta+94185)
\\[-2pt] &\;
+2^{10}z^4(8591694\theta^4+11872968\theta^3+7381951\theta^2+2132674\theta+236280)
\\[-2pt] &\;
-2^{12}z^5(15421829\theta^4+18326342\theta^3+7032841\theta^2+833608\theta-2718)
\\[-2pt] &\;
+2^{16}\cdot3^2z^6(334895\theta^4+615600\theta^3+867965\theta^2+590850\theta+138536)
\\[-2pt] &\;
-2^{19}\cdot3^4\cdot7z^7(\theta+1)(\theta+2)(646\theta^2+1715\theta+1044)
\\[-2pt] &\;
+2^{22}\cdot3^6\cdot7^2z^8(\theta+1)^2(2\theta+1)(2\theta+3)
\\[-8pt]
\endaligned &\cr
& &\span\hrulefill&\cr
& &&
A_n=\binom{2n}n\sum_{k,l}(-4)^{n-k}\binom nk\binom nl\binom kl\binom{2k}k\binom{n+k+l}n
&\cr
\noalign{\hrule}
& \eqnnoltmp{399}{459} &&
\aligned
\\[-12pt]
&
D=13^2\theta^4
-13z(5041\theta^4+7634\theta^3+5767\theta^2+1950\theta+260)
\\[-2pt] &\;
+2^3z^2(744635\theta^4+1560842\theta^3+1510101\theta^2+768170\theta+156078)
\\[-2pt] &\;
-2^6z^3(3698955\theta^4+10227906\theta^3+12569064\theta^2+7257627\theta+1555242)
\\[-2pt] &\;
+2^9z^4(9225025\theta^4+33675338\theta^3+49289090\theta^2+318449807\theta+7296732)
\\[-2pt] &\;
-2^{12}\cdot3\cdot17z^5(\theta+1)(222704\theta^3+833160\theta^2+989659\theta+317310)
\\[-2pt] &\;
+2^{15}\cdot3^3\cdot17^2\cdot23^2z^6(\theta+1)(\theta+2)(2\theta+1)(2\theta+5)
\\[-8pt]
\endaligned &\cr
& &\span\hrulefill&\cr
& &&
A_n=\binom{2n}n\sum_{k,l}(-4)^{n-k}\binom nk\binom nl\binom kl\binom{2k}k\binom{2k+l}n
&\cr
\noalign{\hrule}
& \eqnnoltmp{400}{460} &&
\aligned
\\[-12pt]
&
D=3^2\theta^4
-2^2\cdot3z(29\theta^4+178\theta^3+134\theta^2+45\theta+6)
\\[-2pt] &\;
-2^5z^2(2233\theta^4+2536\theta^3+607\theta^2+132\theta+12)
\\[-2pt] &\;
-2^{10}z^3(1274\theta^4+7425\theta^3+20002\theta^2+12717\theta+2670)
\\[-2pt] &\;
+2^{13}z^4(2539\theta^4-36538\theta^3-52775\theta^2-31122\theta-6192)
\\[-2pt] &\;
+2^{20}z^5(1617\theta^4+9771\theta^3+4484\theta^2-674\theta-556)
\\[-2pt] &\;
+2^{25}z^6(1135\theta^4+4272\theta^3+3439\theta^2+858\theta+16)
\\[-2pt] &\;
-2^{31}\cdot3z^7(2\theta+1)(110\theta^3+225\theta^2+184\theta+57)
\\[-2pt] &\;
+2^{37}\cdot3^2z^8(\theta+1)^2(2\theta+1)(2\theta+3)
\\[-8pt]
\endaligned &\cr
& &\span\hrulefill&\cr
& &&
A_n=\binom{2n}n\sum_{k,l}(-4)^{n-k}\binom nk\binom nl\binom kl\binom{2k}k\binom{2k-2l}n
&\cr
\noalign{\hrule}
}}\hss}

\newpage

\hbox to\hsize{\hss\vbox{\offinterlineskip
\halign to120mm{\strut\tabskip=100pt minus 100pt
\strut\vrule#&\hbox to6.5mm{\hss$#$\hss}&%
\vrule#&\hbox to113mm{\hfil$\dsize#$\hfil}&%
\vrule#\tabskip=0pt\cr\noalign{\hrule}
& \# &%
& \text{differential operator $D$ and coefficients $A_n$, $n=0,1,2,\dots$} &\cr
\noalign{\hrule\vskip1pt\hrule}
& \eqnnoltmp{401}{461} &&
\aligned
\\[-12pt]
&
D=7^2\theta^4-2\cdot7z(1488\theta^4+1452\theta^3+1125\theta^2+399\theta+56)
\\[-2pt] &\;
+2^2z^2(766392\theta^4+1184952\theta^3+1010797\theta^2+454076\theta+83776)
\\[-2pt] &\;
-2^4z^3(12943616\theta^4+28354200\theta^3+30710572\theta^2+16054731\theta+3215254)
\\[-2pt] &\;
+2^6z^4(105973188\theta^4+333359304\theta^3+436182381\theta^2
\\[-2pt] &\;\qquad
+261265857\theta+57189166)
\\[-2pt] &\;
-2^{11}\cdot127z^5(\theta+1)(390972\theta^3+1350660\theta^2+1486781\theta+460439)
\\[-2pt] &\;
+2^{14}\cdot23^2\cdot127^2z^6(\theta+1)(\theta+2)(2\theta+1)(2\theta+5)
\\[-8pt]
\endaligned &\cr
& &\span\hrulefill&\cr
& &&
A_n=\binom{2n}{n}\sum_{k,l}(-4)^{n-k}\binom{n}{k}\binom{n}{l}\binom{k}{l}\binom{2k}{k}\binom{k+2l}{n}
&\cr
\noalign{\hrule}
& \eqnnoltmp{402}{462} &&
\aligned
\\[-12pt]
&
D=\theta^4+2z(2\theta+1)^2(3\theta^2+3\theta+1)
\\[-2pt] &\;
-4z^2(2\theta+1)(2\theta+3)(47\theta^2+94\theta+51)
\\[-2pt] &\;
+112z^3(2\theta+1)(2\theta+3)^2(2\theta+5)
\\[-8pt]
\endaligned &\cr
& &\span\hrulefill&\cr
& &&
\text{a formula for $A_n$ is not known; \cite[\#\,$\widetilde{C}_9$]{Co}}
&\cr
\noalign{\hrule}
& \eqnnoltmp{403}{463} &&
\aligned
\\[-12pt]
&
D=\theta^4+2z(2\theta+1)^2(7\theta^2+7\theta+3)
\\[-2pt] &\;
+4z^2(2\theta+1)(2\theta+3)(29\theta^2+58\theta+33)
\\[-2pt] &\;
+240z^3(2\theta+1)(2\theta+3)^2(2\theta+5)
\\[-8pt]
\endaligned &\cr
& &\span\hrulefill&\cr
& &&
\text{a formula for $A_n$ is not known; \cite[\#\,$\widetilde{C}_{17}$]{Co}}
&\cr
\noalign{\hrule}
& \eqnnoltmp{403}{463} &&
\aligned
\\[-12pt]
&
D=5^2\theta^4+5z(487\theta^4+878\theta^3+709\theta^2+270\theta+40)
\\[-2pt] &\;
+2^5z^2(1013\theta^4+2639\theta^3+2943\theta^2+1520\theta+280)
\\[-2pt] &\;
-2^8z^3(2169\theta^4+144880\theta^3+30789\theta^2+22440\theta+5240)
\\[-2pt] &\;
-2^{12}z^4(2\theta+1)(518\theta^3+2397\theta^2+2940\theta+1048)
\\[-2pt] &\;
-2^{16}\cdot3z^5(2\theta+1)(2\theta+3)^2(2\theta+5)
\\[-8pt]
\endaligned &\cr
& &\span\hrulefill&\cr
& &&
\text{a formula for $A_n$ is not known; \cite[\#\,$\widetilde{C}_{25}$]{Co}}
&\cr
\noalign{\hrule}
}}\hss}

\begin{remarka}
In the table there are 32 pairs of differential equations dual
at~0 and $\infty$, respectively. For the pairs
(\#22,\#118), (\#21,\#71), (\#23,\#56)
we know a simple formula for the coefficients.
Then we have the sporadic pairs
(\#193,\#198), (\#210,\#211), (\#117,\#212), (\#222,\#225),
(\#246,\#247), (\#55,\allowbreak\#277), (\#295,\#296),
where a formula for $A_n$ at~$\infty$ was found by chance. It would be
desirable to find formulas for $A_n$ at~$\infty$ also for the remaining
cases. An even more interesting question is:
What is the geometric meaning of the instanton numbers at~$\infty$, e.g.,
for the cases \#17, \#21, \#22, \#23, \#27,
where we know the manifold?
\end{remarka}

\begin{remarka}
Most of the differential equations \#242--290 were found by using Maple's
\texttt{Zeilberger} on
\[
\text{``}A_n\text{''}
=\sum_k(n-2k)C(n,k)
\]
(cf.~Subsection~\ref{subsec:howtosum}),
where a binomial expression $C(n,k)$ sarisfies
$C(n,n-k)=C(n,k)$.
Using this symmetry one sees that $\text{``}A_n\text{''}=0$,
but Maple finds several bona fide
differential equations and \emph{usually\/}
the correct coefficients are found by differentiation
\[
A_n=\frac{\d\text{``}A_n\text{''}}{\d k},
\]
which is the case in \#244, \#246, \#247, etc.,
leading to harmonic sums containing
$H_n=\sum_{j=1}^nj^{-1}$.
However in some cases
(like \#242, \#245, \#259, \#260, \#261, \#262, \#264,
\#274, \#279, \#281, \#282, \#299)
this is not enough. This will be explained in a forthcoming
paper~\cite{AK}.
\end{remarka}

\newpage

\cleardoublepage
\section{Table of powers}

We found that in most of the cases $q(z)/z$ was a high power
of a power series in $\Z[[z]]$ (and $z(q)/q$ the same power
of a power series in $\Z[[q]]$). Similarly, $y_0(z)$ was
a (not so high) power of an integral power series. In the following
table we present (conjectured) powers $r,s$ in
$$
q=z(1+C_1z+\dotsb)^r, \qquad
y_0=(1+D_1z+\dotsb)^s;
$$
the first column is reserved for numeration of cases.

\medskip
\hbox to\hsize{\hss\vbox{\offinterlineskip
\halign to29mm{\strut\tabskip=100pt minus 100pt
\strut\vrule\vphantom{\vrule height9.8pt}#&\hbox to5.5mm{\hfil$#$\kern-.2pt}&%
\vrule#&\hbox to10.2mm{\hfil$#$\hfil}&%
\vrule#&\hbox to10.2mm{\hfil$#$\hfil}&%
\vrule#\tabskip=0pt\cr\noalign{\hrule}
& \#\, && r && s &\cr
\noalign{\hrule\vskip1pt\hrule}
&  \eqnno{1} &&   10 &&   4 &\cr
&  \eqnno{2} &&  960 &&  24 &\cr
&  \eqnno{3} &&   64 &&   8 &\cr
&  \eqnno{4} &&  180 &&   6 &\cr
&  \eqnno{5} &&  108 &&   4 &\cr
&  \eqnno{6} &&   32 &&   8 &\cr
&  \eqnno{7} &&   64 &&   8 &\cr
&  \eqnno{8} &&   64 &&  12 &\cr
&  \eqnno{9} &&  576 &&  24 &\cr
& \eqnno{10} &&  960 &&  24 &\cr
& \eqnno{11} &&   60 &&  12 &\cr
& \eqnno{12} &&   96 &&  24 &\cr
& \eqnno{13} && 2880 && 120 &\cr
& \eqnno{14} && 1728 &&   4 &\cr
& \eqnno{15} &&    6 &&   2 &\cr
& \eqnno{16} &&    4 &&   4 &\cr
& \eqnno{17} &&    6 &&   1 &\cr
& \eqnno{18} &&   12 &&   2 &\cr
& \eqnno{19} &&    1 &&   2 &\cr
& \eqnno{20} &&    3 &&   2 &\cr
& \eqnno{21} &&    2 &&   2 &\cr
& \eqnno{22} &&    5 &&   1 &\cr
& \eqnno{23} &&    4 &&   4 &\cr
& \eqnno{24} &&    3 &&   3 &\cr
& \eqnno{25} &&    4 &&   2 &\cr
& \eqnno{26} &&    2 &&   4 &\cr
& \eqnno{27} &&   14 &&   1 &\cr
& \eqnno{28} &&    8 &&   1 &\cr
& \eqnno{29} &&    2 &&   1 &\cr
& 2^* && 160 &&  16 &\cr
& \eqnnol{\hbox{$3^*$}}{35} &&   16 &&   8 &\cr
& \eqnnol{\hbox{$4^*$}}{186} &&   54 &&   1 &\cr
& \eqnnol{\hbox{$4^{**}$}}{46} && 18 &&   4 &\cr
& \eqnnol{\hbox{$6^*$}}{220} &&   64 &&   8 &\cr
\noalign{\hrule}
}}\hss
\vbox{\offinterlineskip
\halign to29mm{\strut\tabskip=100pt minus 100pt
\strut\vrule\vphantom{\vrule height9.8pt}#&\hbox to5.5mm{\hfil$#$\kern-.2pt}&%
\vrule#&\hbox to10.2mm{\hfil$#$\hfil}&%
\vrule#&\hbox to10.2mm{\hfil$#$\hfil}&%
\vrule#\tabskip=0pt\cr\noalign{\hrule}
& \#\, && r && s &\cr
\noalign{\hrule\vskip1pt\hrule}
& \eqnnol{\hbox{$7^*$}}{192} &&   32 &&   8 &\cr
& \eqnnol{\hbox{$7^{**}$}}{189} && 32 &&   8 &\cr
& \eqnnol{\hbox{$8^*$}}{193} &&   72 &&   2 &\cr
& \eqnnol{\hbox{$8^{**}$}}{141} &&216 &&   2 &\cr
& \eqnnol{\hbox{$9^*$}}{191} &&  288 &&   8 &\cr
& \eqnnol{\hbox{$9^{**}$}}{190} &&864 &&   8 &\cr
& \eqnnol{\hbox{$10^*$}}{185} &&   32 &&   8 &\cr
& \eqnnol{\hbox{$10^{**}$}}{111} \hss&& 32 &&   8 &\cr
& \eqnnol{\hbox{$13^*$}}{188} &&  432 &&   8 &\cr
& \eqnnol{\hbox{$13^{**}$}}{47} \hss&&144 &&   8 &\cr
& 14^* &&   48 &&  24 &\cr
& \eqnnol{\hbox{$\wh1$}}{206} &&  10 &&   - &\cr
& \eqnnol{\hbox{$\wh2$}}{207} && 320 &&   8 &\cr
& \eqnnol{\hbox{$\wh3$}}{208} &&  64 &&   8 &\cr
& \eqnnol{\hbox{$\wh4$}}{209} &&  18 &&   - &\cr
& \eqnnol{\hbox{$\wh5$}}{210} &&  24 &&   6 &\cr
& \eqnnol{\hbox{$\wh6$}}{211} &&  64 &&   8 &\cr
& \eqnnol{\hbox{$\wh7$}}{212} &&  64 &&   8 &\cr
& \eqnnol{\hbox{$\wh8$}}{213} &&  72 &&   6 &\cr
& \eqnnol{\hbox{$\wh9$}}{214} && 576 &&   8 &\cr
& \eqnnol{\hbox{$\wh{10}$}}{215} && 64 &&   8 &\cr
& \eqnnol{\hbox{$\wh{11}$}}{216} && 24 &&   2 &\cr
& \eqnnol{\hbox{$\wh{12}$}}{217} &&192 &&   8 &\cr
& \eqnnol{\hbox{$\wh{13}$}}{218} &&576 &&  24 &\cr
& \eqnnol{\hbox{$\wh{14}$}}{219} &&192 &&  24 &\cr
& \eqnno{30} &&   64 &&   8 &\cr
& \eqnno{31} &&   64 &&   8 &\cr
& \eqnno{32} &&   78 &&   1 &\cr
& \eqnno{33} &&    4 &&   4 &\cr
& \eqnno{34} &&    8 &&   1 &\cr
& \eqnno{35} &&   12 &&   4 &\cr
& \eqnno{36} &&   24 &&   8 &\cr
& \eqnno{37} &&   24 &&   2 &\cr
& \eqnno{38} &&    8 &&   8 &\cr
\noalign{\hrule}
}}\hss
\vbox{\offinterlineskip
\halign to29mm{\strut\tabskip=100pt minus 100pt
\strut\vrule\vphantom{\vrule height9.8pt}#&\hbox to5.5mm{\hfil$#$\kern-.2pt}&%
\vrule#&\hbox to10.2mm{\hfil$#$\hfil}&%
\vrule#&\hbox to10.2mm{\hfil$#$\hfil}&%
\vrule#\tabskip=0pt\cr\noalign{\hrule}
& \#\, && r && s &\cr
\noalign{\hrule\vskip1pt\hrule}
& \eqnno{39} &&    8 &&   2 &\cr
& \eqnno{40} &&   32 &&   8 &\cr
& \eqnno{41} &&    2 &&   1 &\cr
& \eqnno{42} &&    4 &&   4 &\cr
& \eqnno{43} &&   32 &&   8 &\cr
& \eqnno{44} &&    4 &&  12 &\cr
& \eqnno{45} &&    4 &&   4 &\cr
& \eqnno{46} &&   18 &&   4 &\cr
& \eqnno{47} &&  144 &&   8 &\cr
& \eqnno{48} &&   12 &&   4 &\cr
& \eqnno{49} &&   72 &&   2 &\cr
& \eqnno{50} &&    6 &&   3 &\cr
& \eqnno{51} &&    4 &&   2 &\cr
& \eqnno{52} &&   12 &&   4 &\cr
& \eqnno{53} &&    6 &&   3 &\cr
& \eqnno{54} &&    1 &&   1 &\cr
& \eqnno{55} &&   12 &&   4 &\cr
& \eqnno{56} &&    8 &&   8 &\cr
& \eqnno{57} &&    2 &&   2 &\cr
& \eqnno{58} &&    8 &&   2 &\cr
& \eqnno{59} &&    6 &&   4 &\cr
& \eqnno{60} &&    2 &&   4 &\cr
& \eqnno{61} &&  288 &&  24 &\cr
& \eqnno{62} &&   12 &&   4 &\cr
& \eqnno{63} &&   12 &&   6 &\cr
& \eqnno{64} &&   24 &&   6 &\cr
& \eqnno{65} &&   24 &&   8 &\cr
& \eqnno{66} &&   24 &&   6 &\cr
& \eqnno{67} &&  288 &&   8 &\cr
& \eqnno{68} &&    4 &&   4 &\cr
& \eqnno{69} &&   24 &&   6 &\cr
& \eqnno{70} &&    3 &&   3 &\cr
& \eqnno{71} &&   32 &&   8 &\cr
& \eqnno{72} &&   64 &&   1 &\cr
\noalign{\hrule}
}}\hss
\vbox{\offinterlineskip
\halign to29mm{\strut\tabskip=100pt minus 100pt
\strut\vrule\vphantom{\vrule height9.8pt}#&\hbox to5.5mm{\hfil$#$\kern-.2pt}&%
\vrule#&\hbox to10.2mm{\hfil$#$\hfil}&%
\vrule#&\hbox to10.2mm{\hfil$#$\hfil}&%
\vrule#\tabskip=0pt\cr\noalign{\hrule}
& \#\, && r && s &\cr
\noalign{\hrule\vskip1pt\hrule}
& \eqnno{73} &&   18 &&   6 &\cr
& \eqnno{74} &&    6 &&   4 &\cr
& \eqnno{75} &&    2 &&   1 &\cr
& \eqnno{76} &&    1 &&   1 &\cr
& \eqnno{77} &&    2 &&   1 &\cr
& \eqnno{78} &&    1 &&   8 &\cr
& \eqnno{79} &&    1 &&   1 &\cr
& \eqnno{80} &&    2 &&   1 &\cr
& \eqnno{81} &&    1 &&   1 &\cr
& \eqnno{82} &&    1 &&   6 &\cr
& \eqnno{83} &&    8 &&   8 &\cr
& \eqnno{84} &&    4 &&   4 &\cr
& \eqnno{85} &&    2 &&   1 &\cr
& \eqnno{86} &&    1 &&   1 &\cr
& \eqnno{87} &&    2 &&   1 &\cr
& \eqnno{88} &&    2 &&   1 &\cr
& \eqnno{89} &&    1 &&  24 &\cr
& \eqnno{90} &&    2 &&   1 &\cr
& \eqnno{91} &&    1 &&   1 &\cr
& \eqnno{92} &&    2 &&   2 &\cr
& \eqnno{93} &&  108 &&   1 &\cr
& \eqnno{94} &&    1 &&   2 &\cr
& \eqnno{95} &&    1 &&   1 &\cr
& \eqnno{96} &&   32 &&   1 &\cr
& \eqnno{97} &&  192 &&   1 &\cr
& \eqnno{98} &&   24 &&   1 &\cr
& \eqnno{99} &&    2 &&   2 &\cr
&\eqnno{100} &&   12 &&   2 &\cr
&\eqnno{101} &&   30 &&   1 &\cr
&\eqnno{102} &&    1 &&   1 &\cr
&\eqnno{103} &&   24 &&   3 &\cr
&\eqnno{104} &&    1 &&   1 &\cr
&\eqnno{105} &&    1 &&   1 &\cr
&\eqnno{106} &&    2 &&   2 &\cr
\noalign{\hrule}
}}\hss}

\newpage

\hbox to\hsize{\hss\vbox{\offinterlineskip
\halign to29mm{\strut\tabskip=100pt minus 100pt
\strut\vrule\vphantom{\vrule height9.8pt}#&\hbox to5.5mm{\hfil$#$\kern-.2pt}&%
\vrule#&\hbox to10.2mm{\hfil$#$\hfil}&%
\vrule#&\hbox to10.2mm{\hfil$#$\hfil}&%
\vrule#\tabskip=0pt\cr\noalign{\hrule}
& \#\, && r && s &\cr
\noalign{\hrule\vskip1pt\hrule}
&\eqnno{107} &&   16 &&   8 &\cr
&\eqnno{108} &&    1 &&   1 &\cr
&\eqnno{109} &&    6 &&   8 &\cr
&\eqnno{110} &&   12 &&  12 &\cr
&\eqnno{111} &&   32 &&   8 &\cr
&\eqnno{112} &&   96 &&  24 &\cr
&\eqnno{113} &&    1 &&   3 &\cr
&\eqnno{114} &&    4 &&   4 &\cr
&\eqnno{115} &&   64 &&   8 &\cr
&\eqnno{116} &&   48 &&   8 &\cr
&\eqnno{117} &&    4 &&   4 &\cr
&\eqnno{118} &&   10 &&   4 &\cr
&\eqnno{119} &&    4 &&   4 &\cr
&\eqnno{120} &&    8 &&   6 &\cr
&\eqnno{121} &&    4 &&   2 &\cr
&\eqnno{122} &&    8 &&   8 &\cr
&\eqnno{123} &&    2 &&   2 &\cr
&\eqnno{124} &&    1 &&   1 &\cr
&\eqnno{125} &&    1 &&   1 &\cr
&\eqnno{126} &&    2 &&   1 &\cr
&\eqnno{127} &&   72 &&   1 &\cr
&\eqnno{128} &&    2 &&   1 &\cr
&\eqnno{129} &&    1 &&   6 &\cr
&\eqnno{130} &&   12 &&   1 &\cr
&\eqnno{131} &&    1 &&   1 &\cr
&\eqnno{132} &&   10 &&   1 &\cr
&\eqnno{133} &&   12 &&   2 &\cr
&\eqnno{134} &&    9 &&   3 &\cr
&\eqnno{135} &&   12 &&   6 &\cr
&\eqnno{136} &&   36 &&   6 &\cr
&\eqnno{137} &&    4 &&   4 &\cr
&\eqnno{138} &&    6 &&   6 &\cr
&\eqnno{139} &&   12 &&  12 &\cr
&\eqnno{140} &&   12 &&  12 &\cr
&\eqnno{141} &&  216 &&   2 &\cr
&\eqnno{142} &&    9 &&   3 &\cr
&\eqnno{143} &&   72 &&   6 &\cr
&\eqnno{144} &&   12 &&   6 &\cr
&\eqnno{145} &&   72 &&   3 &\cr
&\eqnno{146} &&    2 &&   2 &\cr
&\eqnno{147} &&    8 &&   4 &\cr
&\eqnno{148} &&    5 &&   1 &\cr
&\eqnno{149} &&   12 &&  12 &\cr
\noalign{\hrule}
}}\hss
\vbox{\offinterlineskip
\halign to29mm{\strut\tabskip=100pt minus 100pt
\strut\vrule\vphantom{\vrule height9.8pt}#&\hbox to5.5mm{\hfil$#$\kern-.2pt}&%
\vrule#&\hbox to10.2mm{\hfil$#$\hfil}&%
\vrule#&\hbox to10.2mm{\hfil$#$\hfil}&%
\vrule#\tabskip=0pt\cr\noalign{\hrule}
& \#\, && r && s &\cr
\noalign{\hrule\vskip1pt\hrule}
&\eqnno{150} &&   12 &&   4 &\cr
&\eqnno{151} &&    6 &&   1 &\cr
&\eqnno{152} &&    4 &&   4 &\cr
&\eqnno{153} &&   12 &&   4 &\cr
&\eqnno{154} &&  180 &&   6 &\cr
&\eqnno{155} &&  192 &&   8 &\cr
&\eqnno{156} &&    2 &&   1 &\cr
&\eqnno{157} &&    1 &&   1 &\cr
&\eqnno{158} &&    2 &&   1 &\cr
&\eqnno{159} &&    2 &&   1 &\cr
&\eqnno{160} &&    3 &&   1 &\cr
&\eqnno{161} &&    6 &&   3 &\cr
&\eqnno{162} &&    3 &&   3 &\cr
&\eqnno{163} &&    6 &&   2 &\cr
&\eqnno{164} &&   12 &&   6 &\cr
&\eqnno{165} &&    6 &&   3 &\cr
&\eqnno{166} && 2880 &&  24 &\cr
&\eqnno{167} &&    3 &&   1 &\cr
&\eqnno{168} &&    3 &&   3 &\cr
&\eqnno{169} &&   12 &&   3 &\cr
&\eqnno{170} &&    6 &&   2 &\cr
&\eqnno{171} &&   24 &&   6 &\cr
&\eqnno{172} &&    9 &&   3 &\cr
&\eqnno{173} &&    4 &&   2 &\cr
&\eqnno{174} &&    1 &&   1 &\cr
&\eqnno{175} &&    3 &&   3 &\cr
&\eqnno{176} &&    4 &&   4 &\cr
&\eqnno{177} &&    4 &&   4 &\cr
&\eqnno{178} &&    3 &&   3 &\cr
&\eqnno{179} &&    3 &&   3 &\cr
&\eqnno{180} &&   24 &&   8 &\cr
&\eqnno{181} &&    6 &&   3 &\cr
&\eqnno{182} &&    2 &&   1 &\cr
&\eqnno{183} &&    4 &&   2 &\cr
&\eqnnol{184}{187} &&    2 &&   1 &\cr
&\eqnnol{185}{184} &&    6 &&   1 &\cr
&\eqnnol{186}{247} &&    2 &&   1 &\cr
&\eqnnol{187}{248} &&    6 &&   9 &\cr
&\eqnnol{188}{245} &&    2 &&   2 &\cr
&\eqnnol{189}{246} &&    2 &&   1 &\cr
&\eqnnol{190}{286} &&  320 &&   1 &\cr
&\eqnnol{191}{287} &&   24 &&   1 &\cr
&\eqnnol{192}{288} &&  192 &&   1 &\cr
\noalign{\hrule}
}}\hss
\vbox{\offinterlineskip
\halign to29mm{\strut\tabskip=100pt minus 100pt
\strut\vrule\vphantom{\vrule height9.8pt}#&\hbox to5.5mm{\hfil$#$\kern-.2pt}&%
\vrule#&\hbox to10.2mm{\hfil$#$\hfil}&%
\vrule#&\hbox to10.2mm{\hfil$#$\hfil}&%
\vrule#\tabskip=0pt\cr\noalign{\hrule}
& \#\, && r && s &\cr
\noalign{\hrule\vskip1pt\hrule}
&\eqnnol{193}{254} &&    3 &&   1 &\cr
&\eqnnol{194}{255} &&    4 &&   1 &\cr
&\eqnnol{195}{256} &&    1 &&   4 &\cr
&\eqnnol{196}{257} &&    1 &&   1 &\cr
&\eqnnol{197}{258} &&    2 &&   2 &\cr
&\eqnnol{198}{259} &&    3 &&   1 &\cr
&\eqnnol{199}{260} &&    4 &&   1 &\cr
&\eqnnol{200}{261} &&    1 &&   1 &\cr
&\eqnnol{201}{262} &&    8 &&   8 &\cr
&\eqnnol{202}{272} &&    1 &&   1 &\cr
&\eqnnol{203}{273} &&    2 &&   4 &\cr
&\eqnnol{204}{200} &&  320 &&   8 &\cr
&\eqnno{205} &&    2 &&   4 &\cr
&\eqnnoltmp{206}{244} &&    4 &&   4 &\cr
&\eqnnol{207}{225} &&   64 &&   8 &\cr
&\eqnnol{208}{289} &&    6 &&   8 &\cr
&\eqnnol{209}{290} &&    2 &&   1 &\cr
&\eqnnol{210}{263} &&   12 &&   2 &\cr
&\eqnnol{211}{271} &&   96 &&   8 &\cr
&\eqnnol{212}{264} &&    2 &&   2 &\cr
&\eqnnol{213}{265} &&    2 &&   1 &\cr
&\eqnnol{214}{266} &&    2 &&   2 &\cr
&\eqnnol{215}{267} &&    4 &&   4 &\cr
&\eqnnol{216}{268} &&    6 &&  12 &\cr
&\eqnnol{217}{269} &&    2 &&  12 &\cr
&\eqnnol{218}{270} &&    6 &&   2 &\cr
&\eqnnol{219}{274} &&    2 &&   2 &\cr
&\eqnnol{220}{291} &&   48 &&   8 &\cr
&\eqnnol{221}{292} &&    4 &&   2 &\cr
&\eqnnol{222}{279} &&    6 &&   4 &\cr
&\eqnnol{223}{280} &&    6 &&  24 &\cr
&\eqnnol{224}{281} &&    2 &&   4 &\cr
&\eqnnol{225}{282} &&  192 &&   8 &\cr
&\eqnnol{226}{283} &&    2 &&   2 &\cr
&\eqnnol{227}{284} &&   36 &&  12 &\cr
&\eqnnol{228}{285} &&    4 &&  12 &\cr
&\eqnnol{229}{293} &&    8 &&   2 &\cr
&\eqnnol{230}{294} &&  360 &&   4 &\cr
&\eqnnol{231}{295} &&    4 &&   4 &\cr
&\eqnnol{232}{296} &&    6 &&   4 &\cr
&\eqnnol{233}{297} &&   96 &&   8 &\cr
&\eqnnol{234}{298} &&    2 &&  16 &\cr
&\eqnnol{235}{299} &&    2 &&  16 &\cr
\noalign{\hrule}
}}\hss
\vbox{\offinterlineskip
\halign to29mm{\strut\tabskip=100pt minus 100pt
\strut\vrule\vphantom{\vrule height9.8pt}#&\hbox to5.5mm{\hfil$#$\kern-.2pt}&%
\vrule#&\hbox to10.2mm{\hfil$#$\hfil}&%
\vrule#&\hbox to10.2mm{\hfil$#$\hfil}&%
\vrule#\tabskip=0pt\cr\noalign{\hrule}
& \#\, && r && s &\cr
\noalign{\hrule\vskip1pt\hrule}
&\eqnnol{236}{300} &&    8 &&   8 &\cr
&\eqnnol{237}{301} &&   24 &&   8 &\cr
&\eqnnol{238}{302} &&    4 &&   2 &\cr
&\eqnnoltmp{239}{303} &&  480 &&  24 &\cr
&\eqnnoltmp{240}{304} &&    2 &&   2 &\cr
&\eqnnoltmp{241}{305} &&   48 &&   8 &\cr
&\eqnnoltmp{242}{306} &&    6 &&   6 &\cr
&\eqnnol{243}{222} &&   14 &&   1 &\cr
&\eqnnoltmp{244}{307} &&   28 &&   1 &\cr
&\eqnnoltmp{245}{308} &&   24 &&   1 &\cr
&\eqnnoltmp{246}{309} &&    4 &&   2 &\cr
&\eqnnoltmp{247}{310} &&   32 &&   8 &\cr
&\eqnnoltmp{248}{311} &&    2 &&   1 &\cr
&\eqnnoltmp{249}{312} &&   12 &&   2 &\cr
&\eqnnoltmp{250}{313} &&    1 &&   2 &\cr
&\eqnnoltmp{251}{314} &&    6 &&  12 &\cr
&\eqnnoltmp{252}{315} &&    2 &&   2 &\cr
&\eqnnoltmp{253}{316} &&    4 &&   4 &\cr
&\eqnnoltmp{254}{317} &&  192 &&   8 &\cr
&\eqnnoltmp{255}{318} &&   12 &&   4 &\cr
&\eqnnoltmp{256}{319} &&   16 &&   8 &\cr
&\eqnnoltmp{257}{320} &&   64 &&   8 &\cr
&\eqnnoltmp{258}{321} &&   96 &&   8 &\cr
&\eqnnoltmp{259}{322} &&   30 &&   2 &\cr
&\eqnnoltmp{260}{323} &&    4 &&   2 &\cr
&\eqnnoltmp{261}{324} &&    4 &&   2 &\cr
&\eqnnoltmp{262}{325} &&    4 &&   2 &\cr
&\eqnnoltmp{263}{326} &&   32 &&   8 &\cr
&\eqnnoltmp{264}{327} &&   96 &&   8 &\cr
&\eqnnoltmp{265}{328} &&   96 &&   8 &\cr
&\eqnnoltmp{266}{329} &&    6 &&   3 &\cr
&\eqnnoltmp{267}{330} &&   18 &&   3 &\cr
&\eqnnoltmp{268}{331} &&   60 &&   6 &\cr
&\eqnnoltmp{269}{332} && 1440 &&  24 &\cr
&\eqnnoltmp{270}{333} &&   12 &&   2 &\cr
&\eqnnoltmp{271}{334} &&  192 &&   8 &\cr
&\eqnnoltmp{272}{335} &&   12 &&   6 &\cr
&\eqnnoltmp{273}{336} &&    6 &&   6 &\cr
&\eqnnoltmp{274}{337} &&    4 &&   2 &\cr
&\eqnnoltmp{275}{338} &&    4 &&   2 &\cr
&\eqnnoltmp{276}{339} &&  576 &&   8 &\cr
&\eqnnol{277}{223} &&   96 &&   8 &\cr
&\eqnnoltmp{278}{340} &&   12 &&  12 &\cr
\noalign{\hrule}
}}\hss}

\newpage

\hbox to\hsize{\hss\vbox{\offinterlineskip
\halign to29mm{\strut\tabskip=100pt minus 100pt
\strut\vrule\vphantom{\vrule height9.8pt}#&\hbox to5.5mm{\hfil$#$\kern-.2pt}&%
\vrule#&\hbox to10.2mm{\hfil$#$\hfil}&%
\vrule#&\hbox to10.2mm{\hfil$#$\hfil}&%
\vrule#\tabskip=0pt\cr\noalign{\hrule}
& \#\, && r && s &\cr
\noalign{\hrule\vskip1pt\hrule}
&\eqnnoltmp{279}{341} &&    1 &&   1 &\cr
&\eqnnoltmp{280}{342} &&    9 &&   3 &\cr
&\eqnnoltmp{281}{343} &&   10 &&   1 &\cr
&\eqnnoltmp{282}{344} &&    4 &&   2 &\cr
&\eqnnoltmp{283}{345} &&    4 &&   6 &\cr
&\eqnnoltmp{284}{346} &&    1 &&   1 &\cr
&\eqnnoltmp{285}{347} &&    1 &&   1 &\cr
&\eqnnoltmp{286}{348} &&    1 &&   2 &\cr
&\eqnnoltmp{287}{349} &&    2 &&   1 &\cr
&\eqnnoltmp{288}{350} &&   96 &&  72 &\cr
&\eqnnoltmp{289}{351} &&   64 &&   8 &\cr
&\eqnnol{290}{221} &&    9 &&   3 &\cr
&\eqnnol{291}{226} &&    1 &&   3 &\cr
&\eqnnoltmp{292}{352} &&    4 &&   4 &\cr
&\eqnnoltmp{293}{353} &&   16 &&   4 &\cr
&\eqnnoltmp{294}{354} &&   64 &&   8 &\cr
&\eqnnoltmp{295}{355} &&   64 &&   8 &\cr
&\eqnnoltmp{296}{356} &&   32 &&   8 &\cr
&\eqnnoltmp{297}{357} &&    2 &&  48 &\cr
&\eqnnoltmp{298}{358} &&    4 &&   2 &\cr
&\eqnnoltmp{299}{359} &&    6 &&   6 &\cr
&\eqnnoltmp{300}{360} &&  160 &&   8 &\cr
&\eqnnoltmp{301}{361} &&    1 &&   1 &\cr
&\eqnnoltmp{302}{362} &&    2 &&   4 &\cr
&\eqnnoltmp{303}{363} &&    2 &&   2 &\cr
&\eqnnoltmp{304}{364} &&    2 &&   2 &\cr
&\eqnnoltmp{305}{365} &&  192 &&   8 &\cr
&\eqnnoltmp{306}{366} &&    9 &&   2 &\cr
&\eqnnoltmp{307}{367} &&    6 &&   1 &\cr
&\eqnnoltmp{308}{368} &&    2 &&   1 &\cr
&\eqnnoltmp{309}{369} &&    4 &&   2 &\cr
&\eqnnoltmp{310}{370} &&    2 &&   1 &\cr
\noalign{\hrule}
}}\hss
\vbox{\offinterlineskip
\halign to29mm{\strut\tabskip=100pt minus 100pt
\strut\vrule\vphantom{\vrule height9.8pt}#&\hbox to5.5mm{\hfil$#$\kern-.2pt}&%
\vrule#&\hbox to10.2mm{\hfil$#$\hfil}&%
\vrule#&\hbox to10.2mm{\hfil$#$\hfil}&%
\vrule#\tabskip=0pt\cr\noalign{\hrule}
& \#\, && r && s &\cr
\noalign{\hrule\vskip1pt\hrule}
&\eqnnoltmp{311}{371} &&    2 &&   2 &\cr
&\eqnnoltmp{312}{372} &&    2 &&   1 &\cr
&\eqnnoltmp{313}{373} &&    2 &&   1 &\cr
&\eqnnoltmp{314}{374} &&    6 &&   3 &\cr
&\eqnnoltmp{315}{375} &&    6 &&   1 &\cr
&\eqnnoltmp{316}{376} &&   12 &&   4 &\cr
&\eqnnoltmp{317}{377} &&    6 &&   1 &\cr
&\eqnnoltmp{318}{378} &&    2 &&   1 &\cr
&\eqnnoltmp{319}{379} &&    2 &&   1 &\cr
&\eqnnoltmp{320}{380} &&    2 &&   1 &\cr
&\eqnnoltmp{321}{381} &&    4 &&   2 &\cr
&\eqnnoltmp{322}{382} &&    1 &&   1 &\cr
&\eqnnoltmp{323}{383} &&    1 &&   1 &\cr
&\eqnnoltmp{324}{384} &&    4 &&   4 &\cr
&\eqnnoltmp{325}{385} &&    4 &&   2 &\cr
&\eqnnoltmp{326}{386} &&    1 &&   1 &\cr
&\eqnnoltmp{327}{387} &&    1 &&   1 &\cr
&\eqnnoltmp{328}{388} &&    4 &&   8 &\cr
&\eqnnoltmp{329}{389} &&   16 &&   8 &\cr
&\eqnnoltmp{330}{390} &&   64 &&   8 &\cr
&\eqnnoltmp{331}{391} &&    8 &&   8 &\cr
&\eqnnoltmp{332}{392} &&    4 &&   4 &\cr
&\eqnnoltmp{333}{393} &&    4 &&  12 &\cr
&\eqnnoltmp{334}{394} &&    1 &&   1 &\cr
&\eqnnoltmp{335}{395} &&    2 &&   2 &\cr
&\eqnnoltmp{336}{396} &&    2 &&   2 &\cr
&\eqnnoltmp{337}{397} &&    6 &&  12 &\cr
&\eqnnoltmp{338}{398} &&    4 &&   4 &\cr
&\eqnnoltmp{339}{399} &&   12 &&   4 &\cr
&\eqnnoltmp{340}{400} &&   12 &&   4 &\cr
&\eqnnoltmp{341}{401} &&    2 &&   2 &\cr
&\eqnnoltmp{342}{402} &&    2 &&   1 &\cr
\noalign{\hrule}
}}\hss
\vbox{\offinterlineskip
\halign to29mm{\strut\tabskip=100pt minus 100pt
\strut\vrule\vphantom{\vrule height9.8pt}#&\hbox to5.5mm{\hfil$#$\kern-.2pt}&%
\vrule#&\hbox to10.2mm{\hfil$#$\hfil}&%
\vrule#&\hbox to10.2mm{\hfil$#$\hfil}&%
\vrule#\tabskip=0pt\cr\noalign{\hrule}
& \#\, && r && s &\cr
\noalign{\hrule\vskip1pt\hrule}
&\eqnnoltmp{343}{403} &&    2 &&   1 &\cr
&\eqnnoltmp{344}{404} &&    1 &&  12 &\cr
&\eqnnoltmp{345}{405} &&    2 &&   1 &\cr
&\eqnnoltmp{346}{406} &&    2 &&   1 &\cr
&\eqnnoltmp{347}{407} &&   24 &&   4 &\cr
&\eqnnoltmp{348}{408} &&   24 &&   4 &\cr
&\eqnnoltmp{349}{409} &&    6 &&   1 &\cr
&\eqnnoltmp{350}{410} &&    8 &&   4 &\cr
&\eqnnoltmp{351}{411} &&   64 &&   8 &\cr
&\eqnnoltmp{352}{412} &&    1 &&   1 &\cr
&\eqnnoltmp{353}{413} &&    4 &&   2 &\cr
&\eqnnoltmp{354}{414} &&   15 &&   1 &\cr
&\eqnnoltmp{355}{415} &&    8 &&   1 &\cr
&\eqnnoltmp{356}{416} &&    2 &&   1 &\cr
&\eqnnoltmp{357}{417} &&    2 &&   2 &\cr
&\eqnnoltmp{358}{418} &&   32 &&   8 &\cr
&\eqnnoltmp{359}{419} &&   12 &&   1 &\cr
&\eqnnoltmp{360}{420} &&    1 &&   1 &\cr
&\eqnnoltmp{361}{421} &&    4 &&   2 &\cr
&\eqnnoltmp{362}{422} &&   32 &&   8 &\cr
&\eqnnoltmp{363}{423} &&   72 &&  12 &\cr
&\eqnnoltmp{364}{424} &&    2 &&   4 &\cr
&\eqnnoltmp{365}{425} &&    4 &&   8 &\cr
&\eqnnoltmp{366}{426} &&    2 &&   4 &\cr
&\eqnnoltmp{367}{427} &&    4 &&   4 &\cr
&\eqnnoltmp{368}{428} &&   32 &&   8 &\cr
&\eqnnoltmp{369}{429} &&    1 &&   2 &\cr
&\eqnnoltmp{370}{430} &&    1 &&   1 &\cr
&\eqnnoltmp{371}{431} &&    1 &&   2 &\cr
&\eqnnoltmp{372}{432} &&    3 &&   1 &\cr
&\eqnnoltmp{373}{433} &&    2 &&   3 &\cr
&\eqnnoltmp{374}{434} &&    1 &&   1 &\cr
\noalign{\hrule}
}}\hss
\vbox{\offinterlineskip
\halign to29mm{\strut\tabskip=100pt minus 100pt
\strut\vrule\vphantom{\vrule height9.8pt}#&\hbox to5.5mm{\hfil$#$\kern-.2pt}&%
\vrule#&\hbox to10.2mm{\hfil$#$\hfil}&%
\vrule#&\hbox to10.2mm{\hfil$#$\hfil}&%
\vrule#\tabskip=0pt\cr\noalign{\hrule}
& \#\, && r && s &\cr
\noalign{\hrule\vskip1pt\hrule}
&\eqnnoltmp{375}{435} &&    1 &&   1 &\cr
&\eqnnoltmp{376}{436} &&    1 &&   2 &\cr
&\eqnnoltmp{377}{437} &&    6 &&   2 &\cr
&\eqnnoltmp{378}{438} &&    4 &&   4 &\cr
&\eqnnoltmp{379}{439} &&    2 &&   4 &\cr
&\eqnnoltmp{380}{440} &&    2 &&   1 &\cr
&\eqnnoltmp{381}{441} &&    2 &&   2 &\cr
&\eqnnoltmp{382}{442} &&    4 &&   4 &\cr
&\eqnnoltmp{383}{443} &&   16 &&   8 &\cr
&\eqnnoltmp{384}{444} &&   16 &&   8 &\cr
&\eqnnoltmp{385}{445} &&    3 &&   3 &\cr
&\eqnnoltmp{386}{446} &&    2 &&   1 &\cr
&\eqnnoltmp{387}{447} &&    2 &&   1 &\cr
&\eqnnoltmp{388}{448} &&    2 &&   1 &\cr
&\eqnnoltmp{389}{449} &&    2 &&   1 &\cr
&\eqnnoltmp{390}{450} &&    2 &&   1 &\cr
&\eqnnoltmp{391}{451} &&    6 &&   2 &\cr
&\eqnnoltmp{392}{452} &&    6 &&   2 &\cr
&\eqnnoltmp{393}{453} &&   48 &&   8 &\cr
&\eqnnoltmp{394}{454} &&    2 &&   2 &\cr
&\eqnnoltmp{395}{455} &&    4 &&  12 &\cr
&\eqnnoltmp{396}{456} &&    2 &&   2 &\cr
&\eqnnoltmp{397}{457} &&    8 &&   8 &\cr
&\eqnnoltmp{398}{458} &&    2 &&   2 &\cr
&\eqnnoltmp{399}{459} &&    2 &&   2 &\cr
&\eqnnoltmp{400}{460} &&    4 &&   4 &\cr
&\eqnnoltmp{401}{461} &&    2 &&   8 &\cr
&\eqnnoltmp{402}{462} &&    2 &&   1 &\cr
&\eqnnoltmp{403}{463} &&    2 &&   1 &\cr
&\eqnnoltmp{404}{464} &&    2 &&   4 &\cr
&    &&      &&     &\cr
&    &&      &&     &\cr
\noalign{\hrule}
}}\hss}

\newpage

\cleardoublepage
\section{Superseeker of Calabi--Yau differential equations}

\newcount\superscount
\def\supers{\global\advance\superscount1\relax\number\superscount}

\hbox to\hsize{\hss\vbox{\offinterlineskip
\halign to127.5mm{\strut\tabskip=100pt minus 100pt
\strut\vrule\vphantom{\vrule height9.8pt}#&\hbox to4.8mm{\hss$#$\kern-1.5pt}&%
\vrule#&\hbox to4.2mm{\hss$#$\kern-1.2pt\hss}&%
\vrule#&\hbox to12.2mm{\hss$#$\kern-1.2pt\hss}&%
\vrule#&\hbox to32.2mm{\hss$#$\kern-1.2pt\hss}&%
\vrule#&\hbox to29.2mm{\hss#\kern-1.2pt\hss}&%
\vrule#&\hbox to25.2mm{\hss#\kern-1.2pt\hss}&%
\vrule#\tabskip=0pt\cr\noalign{\hrule}
& && N_0 && |N_1| && |N_3| && {\small Hadamard product}&%
&\kern-3pt{\small\# in Table~A}\kern-5pt&\cr
\noalign{\hrule\vskip1pt\hrule}
& \supers &&  1 &&       2 &&                  8 &&  && \eqnnol{184}{187} &\cr
& \supers &&  6 &&       2 &&                 10 && \hbox{(a)$\circ$(a)} && \eqnnoltmp{369}{429} &\cr
& \supers &&  3 &&       2 &&                 13 && \hbox{(b)$\circ$(b)} && \eqnnoltmp{370}{430} &\cr
& \supers &&  1 &&       2 &&                104 &&  && \eqnno{41} &\cr
& \supers &&  3 &&       3 &&                 28 &&  && \eqnno{34} &\cr
& \supers &&  1 &&       3 &&                 64 &&  && \eqnnoltmp{366}{426} &\cr
& \supers &&  3 &&       3 &&               2668 &&  && \eqnnoltmp{333}{393} &\cr
& \supers && 12 &&       4 &&                 20 && \hbox{(g)$\circ$(g)} && \eqnnoltmp{371}{431} &\cr
& \supers &&  1 &&       4 &&                 44 &&  && \eqnnoltmp{253}{316} &\cr
& \supers &&  3 &&       4 &&                 44 &&  && \eqnno{23} &\cr
& \supers &&  1 &&       4 &&                 84 &&  && \eqnno{84} &\cr
& \supers &&  5 &&       4 &&                108 &&  && \eqnnoltmp{246}{309} &\cr
& \supers &&  1 &&       4 &&                940 &&  && \eqnnoltmp{395}{455} &\cr
& \supers &&  1 &&       4 &&               3252 && (p)$*$(p) && \eqnnol{\hbox{$3^*$}}{35} &\cr
& \supers &&  1 &&       5 &&                454 && \eqnnol{(a)$*$(a)}{100} && \eqnno{100} &\cr
& \supers &&  1 &&       6 &&                104 && \eqnnol{(k)$*$(k)}{275} && \eqnnol{\hbox{$4^*$}}{186} &\cr
& \supers &&  1 &&       6 &&                170 &&  && \eqnnoltmp{245}{308} &\cr
& \supers &&  1 &&       6 &&                325 && \eqnnol{(a)$*$(f)}{160} && \eqnno{160} &\cr
& \supers && 13 &&       7 &&                 21 &&  && \eqnnoltmp{357}{417} &\cr
& \supers &&  8 &&       8 &&                 96 &&  && \eqnnoltmp{352}{412} &\cr
& \supers &&  2 &&       8 &&                280 &&  && \eqnnoltmp{382}{442} &\cr
& \supers &&  1 &&       9 &&                748 && \eqnnol{(f)$*$(f)}{165} && \eqnno{165} &\cr
& \supers && 17 &&      10 &&                170 &&  && \eqnnoltmp{279}{341} &\cr
& \supers &&  7 &&      10 &&                295 &&  && \eqnno{22} &\cr
& \supers &&  7 &&      10 &&                508 &&  && \eqnnol{235}{299} &\cr
& \supers &&  1 &&      10 &&                664 && \eqnnol{(c)$*$(c)}{103} && \eqnno{103} &\cr
& \supers &&  1 &&      10 &&                870 &&  && \eqnno{60} &\cr
& \supers &&  1 &&      10 &&              18328 &&  && \eqnnoltmp{386}{446} &\cr
& \supers &&  5 &&      11 &&                 71 &&  && \eqnnoltmp{364}{424} &\cr
& \supers &&  6 &&      12 &&                140 &&  && \eqnno{130} &\cr
& \supers &&  2 &&      12 &&                208 &&  && \eqnno{46} &\cr
& \supers && 10 &&      12 &&                236 &&  && \eqnno{17} &\cr
& \supers &&  3 &&      12 &&                644 &&  && \eqnno{16} &\cr
& \supers &&  1 &&      12 &&               3204 && \eqnnol{(A)$*$(a)}{45} && \eqnno{45} &\cr
& \supers &&  1 &&      13 &&               2650 && \eqnnol{(b)$*$(b)}{101} && \eqnno{101} &\cr
& \supers &&  7 &&      14 &&                756 &&  && \eqnnoltmp{402}{462} &\cr
& \supers &&  2 &&      16 &&                208 && \eqnnol{(a)$*$(d)}{105} && \eqnno{105} &\cr
& \supers && 10 &&      16 &&                304 &&  && \eqnno{21} &\cr
& \supers && 12 &&      16 &&                380 &&  && \eqnnoltmp{322}{382} &\cr
& \supers &&  1 &&      16 &&               1232 && \eqnnol{(A)$*$(d)}{36} && \eqnno{36} &\cr
\noalign{\hrule}
}}\hss}

\newpage

\hbox to\hsize{\hss\vbox{\offinterlineskip
\halign to127.5mm{\strut\tabskip=100pt minus 100pt
\strut\vrule\vphantom{\vrule height9.8pt}#&\hbox to4.8mm{\hss$#$\kern-1.5pt}&%
\vrule#&\hbox to4.2mm{\hss$#$\kern-1.2pt\hss}&%
\vrule#&\hbox to12.2mm{\hss$#$\kern-1.2pt\hss}&%
\vrule#&\hbox to32.2mm{\hss$#$\kern-1.2pt\hss}&%
\vrule#&\hbox to29.2mm{\hss#\kern-1.2pt\hss}&%
\vrule#&\hbox to25.2mm{\hss#\kern-1.2pt\hss}&%
\vrule#\tabskip=0pt\cr\noalign{\hrule}
& && N_0 && |N_1| && |N_3| && {\small Hadamard product}&%
&\kern-3pt{\small\# in Table~A}\kern-5pt&\cr
\noalign{\hrule\vskip1pt\hrule}
& \supers &&  1 &&      16 &&               1744 && \eqnnol{(d)$*$(d)}{107} && \eqnno{107} &\cr
& \supers &&  2 &&      16 &&               2000 &&  && \eqnno{42} &\cr
& \supers &&  2 &&      16 &&               2106 && \eqnnol{(a)$*$(b)}{102} && \eqnno{102} &\cr
& \supers &&  1 &&      16 &&               3280 &&  && \eqnno{56} &\cr
& \supers &&  4 &&      16 &&               5072 &&  && \eqnnoltmp{365}{425} &\cr
& \supers &&  6 &&      18 &&                490 &&  && \eqnno{20} &\cr
& \supers &&  9 &&      18 &&               3820 &&  && \eqnnol{199}{260} &\cr
& \supers &&  7 &&      18 &&               5676 &&  && \eqnnol{234}{298} &\cr
& \supers &&  1 &&      19 &&               4455 &&  && \eqnnoltmp{390}{450} &\cr
& \supers && 20 &&      20 &&                100 &&  && \eqnno{205} &\cr
& \supers && 28 &&      20 &&                192 &&  && \eqnnoltmp{312}{372} &\cr
& \supers &&  4 &&      20 &&               1680 &&  && \eqnnoltmp{281}{343} &\cr
& \supers &&  5 &&      20 &&               1820 &&  && \eqnno{18} &\cr
& \supers &&  2 &&      20 &&               2036 &&  && \eqnnoltmp{244}{307} &\cr
& \supers &&  1 &&      20 &&               5924 && \eqnnol{(a)$*$(e)}{114} && \eqnno{114}, \eqnno{150} &\cr
& \supers &&  1 &&      20 &&               8220 &&  && \eqnno{25} &\cr
& \supers &&  1 &&      21 &&              15894 &&  && \eqnno{15} &\cr
& \supers && 29 &&      24 &&                284 &&  && \eqnnoltmp{327}{387} &\cr
& \supers &&  3 &&      24 &&               1552 &&  && \eqnnol{188}{245} &\cr
& \supers &&  7 &&      26 &&              55644 &&  && \eqnnoltmp{297}{357} &\cr
& \supers &&  3 &&      27 &&                217 &&  && \eqnnoltmp{385}{445} &\cr
& \supers &&  3 &&      27 &&              14201 &&  && \eqnnol{216}{268} &\cr
& \supers &&  1 &&      27 &&              18089 && (B)$*(\iota)$ && \eqnno{70} &\cr
& \supers &&  1 &&      28 &&               1036 && \eqnnol{(b)$*$(e)}{121} && \eqnno{121} &\cr
& \supers &&  5 &&      28 &&               1268 &&  && \eqnnoltmp{262}{325} &\cr
& \supers &&  6 &&      28 &&               1820 &&  && \eqnno{27} &\cr
& \supers &&  3 &&      28 &&               3892 &&  && \eqnno{119} &\cr
& \supers &&  5 &&      29 &&               1481 &&  && \eqnnoltmp{404}{464} &\cr
& \supers && 15 &&      30 &&               1540 &&  && \eqnnoltmp{403}{463} &\cr
& \supers &&  6 &&      32 &&                416 &&  && \eqnnoltmp{332}{392} &\cr
& \supers &&  1 &&      32 &&                608 &&  && \eqnnoltmp{397}{457} &\cr
& \supers &&  1 &&      32 &&               1440 && \eqnnol{(l)$*$(l)}{276}, (A)$*(\beta)$ && \eqnnol{\hbox{$10^{**}$}}{185}, \eqnno{40} &\cr
& \supers &&  1 &&      32 &&               7584 &&  && \eqnnol{201}{262} &\cr
& \supers &&  1 &&      32 &&              26016 && \eqnnol{(C)$*$(e)}{30}, (A)$*(\theta)$ && \eqnno{3}, \eqnno{30}, \eqnno{31}, \eqnno{72} &\cr
& \supers &&  3 &&      32 &&              38880 && \eqnnol{(A)$*$(e)}{111} && \eqnno{111} &\cr
& \supers && 39 &&      33 &&                385 &&  && \eqnnoltmp{326}{386} &\cr
& \supers &&  3 &&      33 &&               3422 &&  && \eqnnoltmp{335}{395} &\cr
& \supers &&  3 &&      33 &&               3600 && \eqnnol{(b)$*$(c)}{113} && \eqnno{113} &\cr
& \supers &&  7 &&      35 &&               2184 &&  && \eqnno{28} &\cr
& \supers &&  9 &&      36 &&                556 &&  && \eqnno{183} &\cr
& \supers && 63 &&      36 &&                955 &&  && \eqnnoltmp{344}{404} &\cr
& \supers && 12 &&      36 &&                980 && \hbox{(h)$\circ$(h)} && \eqnnoltmp{372}{432} &\cr
& \supers && 10 &&      36 &&               1284 &&  && \eqnnoltmp{266}{329} &\cr
& \supers &&  9 &&      36 &&               1580 &&  && \eqnnoltmp{353}{413} &\cr
\noalign{\hrule}
}}\hss}

\newpage

\hbox to\hsize{\hss\vbox{\offinterlineskip
\halign to127.5mm{\strut\tabskip=100pt minus 100pt
\strut\vrule\vphantom{\vrule height9.8pt}#&\hbox to4.8mm{\hss$#$\kern-1.5pt}&%
\vrule#&\hbox to4.2mm{\hss$#$\kern-1.2pt\hss}&%
\vrule#&\hbox to12.2mm{\hss$#$\kern-1.2pt\hss}&%
\vrule#&\hbox to32.2mm{\hss$#$\kern-1.2pt\hss}&%
\vrule#&\hbox to29.2mm{\hss#\kern-1.2pt\hss}&%
\vrule#&\hbox to25.2mm{\hss#\kern-1.2pt\hss}&%
\vrule#\tabskip=0pt\cr\noalign{\hrule}
& && N_0 && |N_1| && |N_3| && {\small Hadamard product}&%
&\kern-3pt{\small\# in Table~A}\kern-5pt&\cr
\noalign{\hrule\vskip1pt\hrule}
& \supers &&  3 &&      36 &&               3020 && \eqnnol{(d)$*$(f)}{163} && \eqnno{163} &\cr
& \supers &&  3 &&      36 &&               3284 && \eqnnol{(A)$*$(f)}{133} && \eqnno{133} &\cr
& \supers &&  6 &&      36 &&               3648 &&  && \eqnnol{185}{184} &\cr
& \supers &&  9 &&      36 &&               3856 &&  && \eqnnoltmp{342}{402} &\cr
& \supers &&  1 &&      36 &&               8076 && \eqnnol{(B)$*$(e)}{110} && \eqnno{110} &\cr
& \supers &&  1 &&      36 &&              41421 && \eqnnol{(B)$*$(b)}{24} && \eqnno{24} &\cr
& \supers &&  1 &&      36 &&          128217204 && \eqnnol{(B)$*$(j)}{195} &&  &\cr
& \supers &&  1 &&      37 &&              15270 && \eqnnol{(g)$*$(g)}{144} && \eqnno{144} &\cr
& \supers && 13 &&      39 &&               1621 &&  && \eqnnol{197}{258} &\cr
& \supers &&  8 &&      40 &&               5128 &&  && \eqnnoltmp{304}{364} &\cr
& \supers &&  2 &&      40 &&              26376 &&  && \eqnnoltmp{293}{353} &\cr
& \supers &&  6 &&      42 &&               2542 && \eqnnol{(a)$*$(c)}{104} && \eqnno{104} &\cr
& \supers && 44 &&      44 &&                308 &&  && \eqnno{182} &\cr
& \supers &&  5 &&      44 &&               2980 &&  && \eqnnoltmp{249}{312} &\cr
& \supers && 11 &&      44 &&               3124 &&  && \eqnnoltmp{206}{244} &\cr
& \supers &&  1 &&      44 &&              22500 && \eqnnol{(C)$*$(g)}{139} && \eqnno{139} &\cr
& \supers &&  1 &&      45 &&              43531 &&  && \eqnnoltmp{313}{373} &\cr
& \supers &&  4 &&      48 &&                112 &&  && \eqnnoltmp{339}{399} &\cr
& \supers && 20 &&      48 &&                400 &&  && \eqnnoltmp{381}{441} &\cr
& \supers &&  4 &&      48 &&               1424 && \eqnnol{(b)$*$(d)}{106} && \eqnno{106} &\cr
& \supers &&  1 &&      48 &&               2864 && {\small(B)$*(\beta)$, \eqnnol{(A)$*$(h)}{141}, (k)$*$(m)} && \eqnnol{\hbox{$8^{**}$}}{141}, \eqnno{49}, \eqnno{141} &\cr
& \supers &&  1 &&      48 &&               9104 &&  && \eqnnoltmp{296}{356} &\cr
& \supers &&  4 &&      48 &&               9280 &&  && \eqnnoltmp{380}{440} &\cr
& \supers &&  3 &&      48 &&              11056 && \eqnnol{(A)$*$(c)}{58} && \eqnno{58} &\cr
& \supers &&  2 &&      48 &&              11664 && \eqnnol{(B)$*$(d)}{48} && \eqnno{48} &\cr
& \supers &&  1 &&      48 &&              25200 &&  && \eqnnoltmp{329}{389} &\cr
& \supers &&  1 &&      48 &&              32368 && \eqnnol{(d)$*$(e)}{122}, (A)$*(\varepsilon)$ && \eqnno{122} &\cr
& \supers &&  1 &&      48 &&              73328 && \eqnnol{(C)$*$(d)}{38} && \eqnno{38} &\cr
& \supers && 19 &&      49 &&               1761 &&  && \eqnnol{186}{247} &\cr
& \supers &&  1 &&      50 &&              68472 &&  && \eqnnoltmp{373}{433} &\cr
& \supers && 21 &&      51 &&               5095 &&  && \eqnnol{217}{269} &\cr
& \supers && 20 &&      52 &&               1356 &&  && \eqnnol{203}{273} &\cr
& \supers &&  3 &&      52 &&              52284 &&  && \eqnno{117} &\cr
& \supers &&  1 &&      52 &&             220220 &&  && \eqnno{68} &\cr
& \supers &&  1 &&      54 &&              40552 &&  && \eqnno{50} &\cr
& \supers &&  3 &&      54 &&              64744 &&  && \eqnnol{223}{280} &\cr
& \supers &&  1 &&      55 &&             116555 &&  && \eqnno{118} &\cr
& \supers && 24 &&      56 &&               3552 &&  && \eqnnoltmp{248}{311} &\cr
& \supers &&  5 &&      59 &&              22503 &&  && \eqnnol{224}{281} &\cr
& \supers && 48 &&      60 &&                840 &&  && \eqnnoltmp{376}{436} &\cr
& \supers &&  3 &&      60 &&               1684 && \eqnnol{(A)$*$(g)}{137} && \eqnno{137} &\cr
& \supers &&  3 &&      60 &&              28820 &&  && \eqnnoltmp{255}{318} &\cr
& \supers &&  1 &&      60 &&             134292 &&  && {\small \eqnno{5}, \eqnno{90}, \eqnno{91}, \eqnno{93}, \eqnno{157}}&\cr
& \supers &&  1 &&      60 &&             307860 && \eqnnol{(B)$*$(i)}{194} &&  &\cr
\noalign{\hrule}
}}\hss}

\newpage

\hbox to\hsize{\hss\vbox{\offinterlineskip
\halign to127.5mm{\strut\tabskip=100pt minus 100pt
\strut\vrule\vphantom{\vrule height9.8pt}#&\hbox to4.8mm{\hss$#$\kern-1.5pt}&%
\vrule#&\hbox to4.2mm{\hss$#$\kern-1.2pt\hss}&%
\vrule#&\hbox to12.2mm{\hss$#$\kern-1.2pt\hss}&%
\vrule#&\hbox to32.2mm{\hss$#$\kern-1.2pt\hss}&%
\vrule#&\hbox to29.2mm{\hss#\kern-1.2pt\hss}&%
\vrule#&\hbox to25.2mm{\hss#\kern-1.2pt\hss}&%
\vrule#\tabskip=0pt\cr\noalign{\hrule}
& && N_0 && |N_1| && |N_3| && {\small Hadamard product}&%
&\kern-3pt{\small\# in Table~A}\kern-5pt&\cr
\noalign{\hrule\vskip1pt\hrule}
& \supers &&  1 &&      63 &&              96866 && \eqnnol{(f)$*$(h)}{172} && \eqnno{172} &\cr
& \supers &&  6 &&      64 &&              13504 &&  && \eqnnoltmp{377}{437} &\cr
& \supers &&  1 &&      64 &&              23360 &&  && \eqnno{116} &\cr
& \supers &&  2 &&      64 &&              32576 &&  && \eqnnoltmp{328}{388} &\cr
& \supers &&  1 &&      64 &&             131904 &&  && \eqnnoltmp{383}{443} &\cr
& \supers &&  1 &&      64 &&             246848 && \eqnnol{(C)$*$(c)}{69} && \eqnno{69} &\cr
& \supers &&  7 &&      66 &&               8716 &&  && \eqnnoltmp{379}{439} &\cr
& \supers &&  2 &&      66 &&              59386 && (B)$*(\delta)$ && \eqnno{151} &\cr
& \supers &&  1 &&      66 &&              69048 &&  && \eqnnoltmp{389}{449} &\cr
& \supers &&  4 &&      68 &&              95246 && \eqnnol{(c)$*$(g)}{175} && \eqnno{175} &\cr
& \supers && 11 &&      69 &&               8883 &&  && \eqnnoltmp{307}{367} &\cr
& \supers &&  5 &&      70 &&                980 &&  && \eqnnoltmp{356}{416} &\cr
& \supers && 30 &&      72 &&               1360 &&  && \eqnnoltmp{359}{419} &\cr
& \supers && 24 &&      72 &&               3496 &&  && \eqnnoltmp{286}{348} &\cr
& \supers && 12 &&      72 &&               3900 && \eqnnol{(a)$*$(f)}{160} && \eqnno{160} &\cr
& \supers &&  4 &&      72 &&              20708 && \eqnnol{(B)$*$(f)}{134} && \eqnno{134} &\cr
& \supers &&  3 &&      75 &&              52356 &&  && \eqnnoltmp{336}{396} &\cr
& \supers &&  5 &&      76 &&              10500 &&  && \eqnnoltmp{270}{333} &\cr
& \supers &&  5 &&      76 &&              24836 &&  && \eqnnoltmp{261}{324} &\cr
& \supers &&  3 &&      76 &&             144196 &&  && \eqnno{55} &\cr
& \supers &&  1 &&      76 &&             415420 && \eqnnol{(a)$*$(i)}{233}, (C)$*(\delta)$ && \eqnno{152} &\cr
& \supers && 28 &&      80 &&               2912 &&  && \eqnnol{212}{264} &\cr
& \supers && 23 &&      80 &&               4655 &&  && \eqnno{19} &\cr
& \supers &&  1 &&      80 &&             104976 &&  && \eqnnol{233}{297} &\cr
& \supers &&  1 &&      80 &&             174096 &&  && \eqnno{83} &\cr
& \supers &&  1 &&      80 &&             249872 &&  && \eqnnol{236}{300} &\cr
& \supers && 11 &&      84 &&               9052 &&  && \eqnnol{198}{259} &\cr
& \supers &&  6 &&      84 &&              20848 &&  && \eqnno{29} &\cr
& \supers &&  2 &&      84 &&              83412 &&  && \eqnnol{243}{222} &\cr
& \supers &&  3 &&      84 &&             113304 &&  && \eqnnol{291}{226} &\cr
& \supers && 13 &&      87 &&              21589 &&  && \eqnnoltmp{341}{401} &\cr
& \supers &&  1 &&      90 &&             151648 &&  && \eqnno{73} &\cr
& \supers &&  1 &&      92 &&             585396 && \eqnnol{(C)$*$(b)}{51} && \eqnno{51} &\cr
& \supers &&  9 &&      93 &&              43174 &&  && \eqnnoltmp{394}{454} &\cr
& \supers &&  1 &&      96 &&              12064 && (C)$*(\beta)$ && \eqnnol{\hbox{$7^*$}}{192}, \eqnno{43} &\cr
& \supers &&  6 &&      96 &&              15136 &&  && \eqnnoltmp{378}{438} &\cr
& \supers &&  3 &&      96 &&              26208 && \eqnnol{(c)$*$(e)}{120}, (A)$*(\alpha)$ && \eqnno{39}, \eqnno{120} &\cr
& \supers && 52 &&     100 &&               3500 &&  && \eqnnoltmp{311}{371} &\cr
& \supers &&  1 &&     100 &&             126580 && \eqnnol{(b)$*$(i)}{229}, (C)$*(\eta)$ &&  &\cr
& \supers &&  2 &&     104 &&              89544 &&  && \eqnnoltmp{348}{408} &\cr
& \supers && 12 &&     108 &&                968 && \eqnnol{(c)$*$(f)}{162} && \eqnno{162} &\cr
& \supers &&  6 &&     108 &&               3136 &&  && \eqnnoltmp{242}{306} &\cr
& \supers && 99 &&     108 &&               3213 &&  && \eqnnoltmp{345}{405} &\cr
& \supers && 12 &&     108 &&               4916 && \eqnnol{(b)$*$(f)}{161} && \eqnno{161} &\cr
\noalign{\hrule}
}}\hss}

\newpage

\hbox to\hsize{\hss\vbox{\offinterlineskip
\halign to127.5mm{\strut\tabskip=100pt minus 100pt
\strut\vrule\vphantom{\vrule height9.8pt}#&\hbox to4.8mm{\hss$#$\kern-1.5pt}&%
\vrule#&\hbox to4.2mm{\hss$#$\kern-1.2pt\hss}&%
\vrule#&\hbox to12.2mm{\hss$#$\kern-1.2pt\hss}&%
\vrule#&\hbox to32.2mm{\hss$#$\kern-1.2pt\hss}&%
\vrule#&\hbox to29.2mm{\hss#\kern-1.2pt\hss}&%
\vrule#&\hbox to25.2mm{\hss#\kern-1.2pt\hss}&%
\vrule#\tabskip=0pt\cr\noalign{\hrule}
& && N_0 && |N_1| && |N_3| && {\small Hadamard product}&%
&\kern-3pt{\small\# in Table~A}\kern-5pt&\cr
\noalign{\hrule\vskip1pt\hrule}
& \supers &&  4 &&     108 &&              10472 && \eqnnol{(B)$*$(g)}{138} && \eqnno{138} &\cr
& \supers && 12 &&     108 &&              12580 &&  && \eqnnoltmp{251}{314} &\cr
& \supers &&  6 &&     108 &&              19598 && \eqnnol{(f)$*$(g)}{178} && \eqnno{178} &\cr
& \supers &&  3 &&     108 &&              62596 && \eqnnol{(e)$*$(f)}{164}, (A)$*(\zeta)$ && \eqnno{164} &\cr
& \supers &&  2 &&     108 &&              81104 && \eqnnol{(c)$*$(h)}{169} && \eqnno{169} &\cr
& \supers &&  3 &&     108 &&             206716 && \eqnnol{(C)$*$(f)}{135} && \eqnno{135} &\cr
& \supers &&  1 &&     108 &&           49457556 && \eqnnol{(D)$*$(g)}{140} && \eqnno{140} &\cr
& \supers &&  5 &&     109 &&              16777 &&  && \eqnnoltmp{302}{362} &\cr
& \supers && 11 &&     110 &&               3740 &&  && \eqnnoltmp{355}{415} &\cr
& \supers &&  6 &&     112 &&              35408 &&  && \eqnnoltmp{400}{460} &\cr
& \supers &&  1 &&     112 &&             186800 &&  && \eqnnoltmp{331}{391} &\cr
& \supers &&  1 &&     112 &&             378800 &&  && \eqnno{71} &\cr
& \supers && 19 &&     113 &&               8515 &&  && \eqnnol{202}{272} &\cr
& \supers &&  5 &&     116 &&             186172 &&  && \eqnnoltmp{275}{338} &\cr
& \supers &&  1 &&     117 &&             713814 &&  && \eqnno{4} &\cr
& \supers &&  1 &&     117 &&             844872 &&  && \eqnnoltmp{280}{342} &\cr
& \supers && 91 &&     118 &&               1876 &&  && \eqnnoltmp{375}{435} &\cr
& \supers && 17 &&     126 &&              11700 &&  && \eqnnol{194}{255} &\cr
& \supers &&  1 &&     128 &&             263808 &&  && \eqnnoltmp{393}{453} &\cr
& \supers &&  1 &&     128 &&             382592 &&  && \eqnnol{220}{291} &\cr
& \supers &&  1 &&     128 &&             800384 &&  && \eqnnoltmp{256}{319} &\cr
& \supers &&  7 &&     129 &&              41441 &&  && \eqnnol{193}{254} &\cr
& \supers && 12 &&     132 &&               9736 && \eqnnol{(a)$*$(g)}{173} && \eqnno{173} &\cr
& \supers &&  4 &&     132 &&              52204 &&  && \eqnno{32} &\cr
& \supers &&  4 &&     132 &&             118772 && \eqnnol{(a)$*$(h)}{167} && \eqnno{167} &\cr
& \supers && 97 &&     136 &&               1768 &&  && \eqnnoltmp{374}{434} &\cr
& \supers && 36 &&     140 &&              12008 &&  && \eqnnoltmp{321}{381} &\cr
& \supers && 14 &&     140 &&              24136 &&  && \eqnno{26} &\cr
& \supers &&  3 &&     140 &&             198276 &&  && \eqnnoltmp{338}{398} &\cr
& \supers && 12 &&     144 &&               7312 && \eqnnol{(c)$*$(d)}{123} && \eqnno{123} &\cr
& \supers &&  6 &&     144 &&              30896 && \eqnnol{(d)$*$(g)}{176} && \eqnno{176} &\cr
& \supers &&  6 &&     146 &&              66714 &&  && \eqnnoltmp{306}{366} &\cr
& \supers && 15 &&     147 &&               6032 &&  && \eqnnoltmp{274}{337} &\cr
& \supers && 11 &&     148 &&              44108 &&  && \eqnnoltmp{324}{384} &\cr
& \supers && 13 &&     151 &&              26293 &&  && \eqnnoltmp{303}{363} &\cr
& \supers &&  6 &&     156 &&              29884 &&  && \eqnnoltmp{319}{379} &\cr
& \supers &&  1 &&     160 &&               9310 &&  && \eqnno{19} &\cr
& \supers &&  2 &&     160 &&             539680 && {\small \eqnnol{(e)$*$(e)}{115}, \eqnnol{(c)$*$(i)}{235}, (C)$*(\theta)$}&& {\small\eqnnol{\hbox{$\wh3$}}{208}, \eqnnol{\hbox{$6^*$}}{220}, \eqnno{115}, \eqnnol{190}{286}, \eqnnol{204}{200}} &\cr
& \supers &&  1 &&     160 &&            1956896 &&  && \eqnno{6}, \eqnno{75}, \eqnno{96}, \eqnno{146} &\cr
& \supers &&  1 &&     160 &&            5870688 &&  && \eqnno{76} &\cr
& \supers &&  3 &&     162 &&             197216 &&  && \eqnnoltmp{299}{359} &\cr
& \supers && 61 &&     163 &&               4795 &&  && \eqnno{124} &\cr
& \supers &&  7 &&     178 &&             129516 &&  && \eqnnoltmp{401}{461} &\cr
& \supers &&  4 &&     180 &&              28320 && \eqnnol{(b)$*$(h)}{168} && \eqnno{168} &\cr
\noalign{\hrule}
}}\hss}

\newpage

\hbox to\hsize{\hss\vbox{\offinterlineskip
\halign to127.5mm{\strut\tabskip=100pt minus 100pt
\strut\vrule\vphantom{\vrule height9.8pt}#&\hbox to4.8mm{\hss$#$\kern-1.5pt}&%
\vrule#&\hbox to4.2mm{\hss$#$\kern-1.2pt\hss}&%
\vrule#&\hbox to12.2mm{\hss$#$\kern-1.2pt\hss}&%
\vrule#&\hbox to32.2mm{\hss$#$\kern-1.2pt\hss}&%
\vrule#&\hbox to29.2mm{\hss#\kern-1.2pt\hss}&%
\vrule#&\hbox to25.2mm{\hss#\kern-1.2pt\hss}&%
\vrule#\tabskip=0pt\cr\noalign{\hrule}
& && N_0 && |N_1| && |N_3| && {\small Hadamard product}&%
&\kern-3pt{\small\# in Table~A}\kern-5pt&\cr
\noalign{\hrule\vskip1pt\hrule}
& \supers &&  4 &&     180 &&             110940 && \eqnnol{(B)$*$(h)}{142} && \eqnno{142} &\cr
& \supers &&  9 &&     180 &&             119332 &&  && \eqnnoltmp{361}{421} &\cr
& \supers &&  1 &&     180 &&           21847076 && \eqnnol{(D)$*$(f)}{136} && \eqnno{136} &\cr
& \supers && 15 &&     186 &&              20300 &&  && \eqnnol{226}{283} &\cr
& \supers &&  5 &&     188 &&             450516 &&  && \eqnnoltmp{260}{323} &\cr
& \supers && 81 &&     189 &&               4843 &&  && \eqnnoltmp{334}{394} &\cr
& \supers && 47 &&     189 &&               9277 &&  && \eqnnol{196}{257} &\cr
& \supers && 12 &&     192 &&                156 && \eqnnol{(b)$*$(g)}{174} && \eqnno{174} &\cr
& \supers &&  1 &&     192 &&             616896 &&  && \eqnnoltmp{384}{444} &\cr
& \supers && 11 &&     193 &&              48570 &&  && \eqnnoltmp{301}{361} &\cr
& \supers &&  1 &&     196 &&            2993772 &&  && \eqnno{33} &\cr
& \supers &&  3 &&     204 &&              18628 &&  && \eqnnol{228}{285} &\cr
& \supers &&  3 &&     204 &&             125636 &&  && \eqnnoltmp{387}{447} &\cr
& \supers &&  9 &&     205 &&              97622 &&  && \eqnnoltmp{309}{369} &\cr
& \supers &&  1 &&     207 &&             621972 &&  && \eqnnoltmp{363}{423} &\cr
& \supers &&  1 &&     208 &&            1218192 &&  && \eqnnol{237}{301} &\cr
& \supers &&  1 &&     208 &&            1863312 && \eqnnol{(d)$*$(i)}{231}, (C)$*(\varepsilon)$ &&  &\cr
& \supers &&  9 &&     209 &&              97622 &&  && \eqnnoltmp{298}{358} &\cr
& \supers && 63 &&     216 &&               7371 &&  && \eqnnoltmp{349}{409} &\cr
& \supers && 54 &&     216 &&               9900 &&  && \eqnnoltmp{343}{403} &\cr
& \supers &&  3 &&     220 &&             267636 &&  && \eqnnol{215}{267} &\cr
& \supers &&  3 &&     228 &&             278988 && \eqnnol{(e)$*$(g)}{177}, (A)$*(\gamma)$ && \eqnno{44}, \eqnno{177} &\cr
& \supers &&  4 &&     229 &&             297111 &&  && \eqnnoltmp{314}{374} &\cr
& \supers && 13 &&     231 &&              38037 &&  && \eqnnoltmp{240}{304} &\cr
& \supers && 10 &&     232 &&              59256 &&  && \eqnnoltmp{396}{456} &\cr
& \supers &&  5 &&     232 &&             122168 &&  && \eqnnoltmp{252}{315} &\cr
& \supers &&  9 &&     234 &&             103520 &&  && \eqnnol{214}{266} &\cr
& \supers && 34 &&     236 &&              22848 &&  && \eqnnol{213}{265} &\cr
& \supers && 56 &&     240 &&               6944 &&  && \eqnno{59} &\cr
& \supers &&  8 &&     240 &&             117056 &&  && \eqnno{74} &\cr
& \supers &&  1 &&     240 &&           19105840 && \eqnnol{(D)$*$(d)}{65} && \eqnno{65} &\cr
& \supers && 38 &&     241 &&              17458 &&  && \eqnnoltmp{284}{346} &\cr
& \supers &&  1 &&     243 &&             513936 &&  && \eqnnoltmp{278}{340} &\cr
& \supers && 29 &&     248 &&              38708 &&  && \eqnnoltmp{308}{368} &\cr
& \supers && 20 &&     252 &&              56064 &&  && \eqnnoltmp{273}{336} &\cr
& \supers &&  4 &&     252 &&             387464 && (B)$*(\zeta)$ &&  &\cr
& \supers &&  1 &&     252 &&            1162036 &&  && \eqnno{154} &\cr
& \supers &&  3 &&     252 &&            1522388 &&  && \eqnno{35} &\cr
& \supers &&  2 &&     270 &&             835370 && \eqnnol{(g)$*$(h)}{179} && \eqnno{53}, \eqnno{179} &\cr
& \supers &&  7 &&     274 &&             281388 &&  && \eqnnol{208}{289} &\cr
& \supers && 30 &&     276 &&              33780 &&  && \eqnnoltmp{318}{378} &\cr
& \supers && 20 &&     276 &&             116324 &&  && \eqnnol{222}{279} &\cr
& \supers && 29 &&     285 &&              40626 &&  && \eqnnol{195}{256} &\cr
& \supers &&  1 &&     288 &&            2339616 && \eqnnol{(e)$*$(h)}{171} && \eqnnol{\hbox{$\wh5$}}{210}, \eqnno{98}, \eqnno{171} &\cr
\noalign{\hrule}
}}\hss}

\newpage

\hbox to\hsize{\hss\vbox{\offinterlineskip
\halign to127.5mm{\strut\tabskip=100pt minus 100pt
\strut\vrule\vphantom{\vrule height9.8pt}#&\hbox to4.8mm{\hss$#$\kern-1.5pt}&%
\vrule#&\hbox to4.2mm{\hss$#$\kern-1.2pt\hss}&%
\vrule#&\hbox to12.2mm{\hss$#$\kern-1.2pt\hss}&%
\vrule#&\hbox to32.2mm{\hss$#$\kern-1.2pt\hss}&%
\vrule#&\hbox to29.2mm{\hss#\kern-1.2pt\hss}&%
\vrule#&\hbox to25.2mm{\hss#\kern-1.2pt\hss}&%
\vrule#\tabskip=0pt\cr\noalign{\hrule}
& && N_0 && |N_1| && |N_3| && {\small Hadamard product}&%
&\kern-3pt{\small\# in Table~A}\kern-5pt&\cr
\noalign{\hrule\vskip1pt\hrule}
& \supers &&  1 &&     288 &&           96055968 && \eqnnol{(D)$*$(e)}{112}, (A)$*(\kappa)$ && \eqnno{112} &\cr
& \supers &&  1 &&     291 &&            7935104 &&  && \eqnnol{230}{294} &\cr
& \supers && 23 &&     308 &&              70799 &&  && \eqnnoltmp{250}{313} &\cr
& \supers &&  1 &&     320 &&           19748928 &&  && \eqnnoltmp{241}{305} &\cr
& \supers &&  2 &&     324 &&            1502052 &&  && \eqnnol{290}{221} &\cr
& \supers &&  1 &&     324 &&           10792428 &&  && \eqnno{11}, \eqnno{95} &\cr
& \supers &&  4 &&     336 &&             595280 && \eqnnol{(d)$*$(h)}{170}, (B)$*(\varepsilon)$ && \eqnno{170} &\cr
& \supers &&  1 &&     336 &&            4761360 &&  && \eqnnoltmp{358}{418} &\cr
& \supers &&  1 &&     352 &&            3284448 &&  && \eqnnoltmp{330}{390} &\cr
& \supers && 13 &&     359 &&             393749 &&  && \eqnnoltmp{398}{458} &\cr
& \supers && 21 &&     361 &&             120472 &&  && \eqnnoltmp{287}{349} &\cr
& \supers && 22 &&     362 &&              94342 &&  && \eqnnoltmp{310}{370} &\cr
& \supers &&  5 &&     364 &&            6324580 &&  && \eqnnoltmp{282}{344} &\cr
& \supers && 19 &&     370 &&             140636 &&  && \eqnnoltmp{325}{385} &\cr
& \supers &&  1 &&     372 &&           71562236 && \eqnnol{(D)$*$(a)}{62} && \eqnno{62} &\cr
& \supers &&  5 &&     379 &&            1364199 &&  && \eqnnol{232}{296} &\cr
& \supers &&  1 &&     384 &&             164736 && \eqnnol{(m)$*$(m)}{277} && \eqnnol{13$^*$}{188}, \eqnno{47} &\cr
& \supers &&  3 &&     384 &&            1546624 && \eqnnol{(c)$*$(i)}{235}, (C)$*(\alpha)$ && \eqnno{37} &\cr
& \supers && 12 &&     400 &&             292444 &&  && \eqnnoltmp{323}{383} &\cr
& \supers && 21 &&     414 &&             128592 &&  && \eqnnol{218}{270} &\cr
& \supers && 14 &&     420 &&             159040 &&  && \eqnnol{189}{246} &\cr
& \supers &&  4 &&     428 &&             485244 &&  && \eqnnol{187}{248} &\cr
& \supers &&  1 &&     432 &&           78259376 && \eqnnol{(D)$*$(c)}{64} && \eqnno{64} &\cr
& \supers &&  9 &&     441 &&             173876 &&  && \eqnnoltmp{350}{410} &\cr
& \supers &&  5 &&     444 &&            1501908 &&  && \eqnnol{210}{263} &\cr
& \supers &&  1 &&     444 &&           19050964 &&  && \eqnnol{238}{302} &\cr
& \supers &&  3 &&     460 &&             894404 &&  && \eqnnol{231}{295} &\cr
& \supers &&  3 &&     468 &&            3687996 && \eqnnol{(f)$*$(i)}{227}, (C)$*(\zeta)$ &&  &\cr
& \supers &&  5 &&     468 &&           59427420 &&  && \eqnnoltmp{272}{335} &\cr
& \supers && 17 &&     478 &&             285760 &&  && \eqnnol{209}{290} &\cr
& \supers &&  1 &&     480 &&            4215904 &&  && \eqnnoltmp{258}{321} &\cr
& \supers &&  1 &&     480 &&           16034720 && (D)$*(\beta)$ && \eqnnol{9$^*$}{191}, \eqnnol{9$^{**}$}{190}, \eqnno{67} &\cr
& \supers &&  3 &&     484 &&             819404 &&  && \eqnnoltmp{340}{400} &\cr
& \supers &&  5 &&     492 &&             872164 &&  && \eqnnol{221}{292} &\cr
& \supers &&  1 &&     492 &&          136094428 && \eqnnol{(a)$*$(j)}{234}, (D)$*(\delta)$ && \eqnno{153} &\cr
& \supers && 15 &&     498 &&             360988 &&  && \eqnnol{219}{274} &\cr
& \supers && 20 &&     500 &&             343500 &&  && \eqnnoltmp{354}{414} &\cr
& \supers &&  1 &&     522 &&            9879192 && \eqnnol{(h)$*$(h)}{145} && \eqnnol{\hbox{$\wh4$}}{209}, \eqnno{145}, \eqnno{181} &\cr
& \supers && 36 &&     540 &&             325680 &&  && \eqnnoltmp{347}{407} &\cr
& \supers &&  3 &&     564 &&            2422620 &&  && \eqnnoltmp{283}{345} &\cr
& \supers && 13 &&     567 &&             512341 &&  && \eqnnoltmp{399}{459} &\cr
& \supers && 15 &&     570 &&             392025 &&  && \eqnnoltmp{315}{375} &\cr
& \supers &&  1 &&     575 &&           63441275 &&  && \eqnno{1}, \eqnno{79}, \eqnno{87}, \eqnno{128} &\cr
& \supers &&  1 &&     608 &&           22293216 &&  && \eqnnoltmp{247}{310} &\cr
\noalign{\hrule}
}}\hss}

\newpage

\hbox to\hsize{\hss\vbox{\offinterlineskip
\halign to127.5mm{\strut\tabskip=100pt minus 100pt
\strut\vrule\vphantom{\vrule height9.8pt}#&\hbox to4.8mm{\hss$#$\kern-1.5pt}&%
\vrule#&\hbox to4.2mm{\hss$#$\kern-1.2pt\hss}&%
\vrule#&\hbox to12.2mm{\hss$#$\kern-1.2pt\hss}&%
\vrule#&\hbox to32.2mm{\hss$#$\kern-1.2pt\hss}&%
\vrule#&\hbox to29.2mm{\hss#\kern-1.2pt\hss}&%
\vrule#&\hbox to25.2mm{\hss#\kern-1.2pt\hss}&%
\vrule#\tabskip=0pt\cr\noalign{\hrule}
& && N_0 && |N_1| && |N_3| && {\small Hadamard product}&%
&\kern-3pt{\small\# in Table~A}\kern-5pt&\cr
\noalign{\hrule\vskip1pt\hrule}
& \supers &&  1 &&     612 &&           51318900 && \eqnnol{(b)$*$(j)}{230}, (D)$*(\eta)$ &&  &\cr
& \supers &&  1 &&     624 &&           43406256 &&  && \eqnno{180} &\cr
& \supers && 13 &&     647 &&             942613 &&  && \eqnno{99} &\cr
& \supers && 72 &&     684 &&             398428 &&  && \eqnnol{200}{261} &\cr
& \supers &&  1 &&     684 &&          195638820 && \eqnnol{(D)$*$(b)}{63} && \eqnno{63} &\cr
& \supers &&  8 &&     713 &&            2286220 &&  && \eqnnoltmp{346}{406} &\cr
& \supers &&  1 &&     736 &&           26911072 && \eqnnol{(e)$*$(i)}{239} && \eqnnol{\hbox{$\wh6$}}{211}, \eqnno{77}, \eqnno{78}, \eqnno{97} &\cr
& \supers && 11 &&     741 &&            1526195 &&  && \eqnnoltmp{320}{380} &\cr
& \supers &&  3 &&     798 &&           11433160 &&  && \eqnnoltmp{388}{448} &\cr
& \supers && 48 &&     828 &&             344760 &&  && \eqnnoltmp{317}{377} &\cr
& \supers &&  5 &&     828 &&            4270932 &&  && \eqnnoltmp{268}{331} &\cr
& \supers && 11 &&     852 &&            1678156 &&  && \eqnnoltmp{316}{376} &\cr
& \supers &&  1 &&     864 &&          147560800 && \eqnnol{(c)$*$(j)}{236}, (D)$*(\alpha)$ && \eqnno{66} &\cr
& \supers &&  1 &&     900 &&            8364884 &&  && \eqnnol{227}{284} &\cr
& \supers &&  1 &&     928 &&          170809536 &&  && \eqnno{10}, \eqnno{54} &\cr
& \supers &&  3 &&     988 &&           14008436 &&  && \eqnnoltmp{367}{427} &\cr
& \supers &&  1 &&     992 &&           63721056 &&  && \eqnnoltmp{257}{320} &\cr
& \supers &&  1 &&    1008 &&          607849200 && \eqnnol{(D)$*$(h)}{143}, (B)$*(\kappa)$ && \eqnno{143} &\cr
& \supers &&  3 &&    1020 &&           15174100 && \eqnnol{(g)$*$(i)}{237}, (C)$*(\gamma)$ && \eqnno{52} &\cr
& \supers &&  3 &&    1056 &&           15001120 &&  && \eqnnol{229}{293} &\cr
& \supers &&  1 &&    1056 &&          138459552 &&  && \eqnnoltmp{265}{328} &\cr
& \supers &&  1 &&    1116 &&          349462868 && \eqnnol{(f)$*$(j)}{228}, (D)$*(\zeta)$ &&  &\cr
& \supers &&  7 &&    1162 &&              71127 &&  && \eqnnoltmp{392}{452} &\cr
& \supers &&  1 &&    1248 &&          683015008 &&  && {\small \eqnno{14}, \eqnno{85}, \eqnno{86}, \eqnno{156}}&\cr
& \supers &&  1 &&    1312 &&           58156708 &&  && \eqnnoltmp{263}{326} &\cr
& \supers &&  1 &&    1344 &&          109320512 && \eqnnol{(h)$*$(i)}{241} && \eqnnol{\hbox{$\wh{11}$}}{216}, \eqnno{94} &\cr
& \supers &&  7 &&    1434 &&           18676572 &&  && \eqnno{109} &\cr
& \supers &&  1 &&    1488 &&          517984144 && \eqnnol{(d)$*$(j)}{232}, (D)$*(\varepsilon)$ && &\cr
& \supers &&  1 &&    1584 &&          171534960 &&  && \eqnnoltmp{239}{303} &\cr
& \supers &&  1 &&    1616 &&          283183120 &&  && \eqnnoltmp{300}{360} &\cr
& \supers &&  1 &&    1818 &&          467810538 &&  && \eqnnoltmp{267}{330} &\cr
& \supers &&  5 &&    2043 &&           88982631 &&  && \eqnnoltmp{337}{397} &\cr
& \supers &&  4 &&    2300 &&          253765100 &&  && \eqnno{148} &\cr
& \supers &&  1 &&    2400 &&         2956977632 &&  && \eqnnol{211}{271} &\cr
& \supers &&  1 &&    2450 &&          623291900 &&  && \eqnnol{\hbox{$\wh1$}}{206}, \eqnno{80}, \eqnno{81}, \eqnno{131} &\cr
& \supers &&  1 &&    2484 &&         1327731388 && \eqnnol{(g)$*$(j)}{238}, (D)$*(\gamma)$ && \eqnno{149} &\cr
& \supers &&  1 &&    2592 &&           81451104 && \eqnnol{(C)$*$(j)}{198} &&  &\cr
& \supers &&  1 &&    2628 &&         3966805740 &&  && {\small \eqnno{8}, \eqnno{92}, \eqnno{125}, \eqnno{158}}&\cr
& \supers &&  1 &&    2656 &&         2493879008 &&  && \eqnnoltmp{362}{422} &\cr
& \supers &&  3 &&    2892 &&           85888580 &&  && \eqnnoltmp{391}{451} &\cr
& \supers && 44 &&    3180 &&            6378752 &&  && \eqnnoltmp{285}{347} &\cr
& \supers &&  1 &&    3488 &&         1142687008 && \eqnnol{(i)$*$(i)}{155} && \eqnnol{\hbox{$\wh{10}$}}{215}, \eqnno{155} &\cr
& \supers &&  1 &&    3616 &&          264403872 &&  && \eqnnoltmp{288}{350} &\cr
& \supers &&  1 &&    3936 &&        10892932064 && \eqnnol{(D)$*$(i)}{199} &&  &\cr
\noalign{\hrule}
}}\hss}

\newpage

\hbox to\hsize{\hss\vbox{\offinterlineskip
\halign to127.5mm{\strut\tabskip=100pt minus 100pt
\strut\vrule\vphantom{\vrule height9.8pt}#&\hbox to4.8mm{\hss$#$\kern-1.5pt}&%
\vrule#&\hbox to4.2mm{\hss$#$\kern-1.2pt\hss}&%
\vrule#&\hbox to12.2mm{\hss$#$\kern-1.2pt\hss}&%
\vrule#&\hbox to32.2mm{\hss$#$\kern-1.2pt\hss}&%
\vrule#&\hbox to29.2mm{\hss#\kern-1.2pt\hss}&%
\vrule#&\hbox to25.2mm{\hss#\kern-1.2pt\hss}&%
\vrule#\tabskip=0pt\cr\noalign{\hrule}
& && N_0 && |N_1| && |N_3| && {\small Hadamard product}&%
&\kern-3pt{\small\# in Table~A}\kern-5pt&\cr
\noalign{\hrule\vskip1pt\hrule}
& \supers &&  1 &&    4192 &&         2124587232 &&  && \eqnnol{277}{223} &\cr
& \supers &&  3 &&    4300 &&         1701817028 &&  && \eqnnoltmp{292}{352} &\cr
& \supers &&  1 &&    5408 &&         4296119968 &&  && \eqnnoltmp{295}{355} &\cr
& \supers &&  1 &&    5408 &&        22147077792 &&  && \eqnnoltmp{254}{317} &\cr
& \supers &&  1 &&    5472 &&         6444589536 && \eqnnol{(e)$*$(j)}{240} && \eqnnol{\hbox{$\wh{14}$}}{219}, \eqnno{88}, \eqnno{89} &\cr
& \supers &&  1 &&    7776 &&        66942277344 &&  && \eqnno{12} &\cr
& \supers &&  1 &&    8096 &&         9215266592 &&  && \eqnnoltmp{368}{428} &\cr
& \supers &&  1 &&    8224 &&        15542388128 &&  && \eqnnoltmp{289}{351} &\cr
& \supers &&  1 &&   10080 &&        24400330080 && \eqnnol{(h)$*$(j)}{242} && {\footnotesize \eqnnol{\hbox{$\wh8$}}{213}, \eqnno{82}, \eqnno{126}, \eqnno{127}, \eqnno{129}}&\cr
& \supers &&  1 &&   10912 &&        71557619232 &&  && \eqnnoltmp{271}{334} &\cr
& \supers && 68 &&   12676 &&           65175340 &&  && \eqnnoltmp{360}{420} &\cr
& \supers &&  1 &&   14752 &&       711860273440 &&  && \eqnno{7}, \eqnno{147} &\cr
& \supers &&  1 &&   26400 &&       230398034080 && \eqnnol{(i)$*$(j)}{243} && \eqnnol{\hbox{$\wh{12}$}}{217} &\cr
& \supers &&  1 &&   37216 &&       464865119712 &&  && \eqnnoltmp{264}{327} &\cr
& \supers &&  1 &&   41184 &&      5124430612320 && \eqnnol{(D)$*$(j)}{61}, (D)$*(\kappa)$ && \eqnno{61} &\cr
& \supers &&  1 &&   57760 &&      3869123234080 &&  && \eqnnol{\hbox{$\wh7$}}{212} &\cr
& \supers &&  1 &&   67104 &&     28583248229280 &&  && \eqnno{13}, \eqnno{57}, \eqnno{108} &\cr
& \supers &&  1 &&   70944 &&      3707752060576 &&  && \eqnnol{207}{225} &\cr
& \supers &&  1 &&   80416 &&     15561562691488 &&  && \eqnnoltmp{294}{354} &\cr
& \supers &&  1 &&   82450 &&     22323908689400 &&  && \eqnnoltmp{259}{322} &\cr
& \supers &&  1 &&   93984 &&     25265152551072 &&  && \eqnnol{225}{282} &\cr
& \supers &&  1 &&  177184 &&     45194569320864 &&  && \eqnnoltmp{351}{411} &\cr
& \supers &&  1 &&  188832 &&    101990911789344 &&  && \eqnnoltmp{276}{339} &\cr
& \supers &&  1 &&  201888 &&     40177844666400 && \eqnnol{(j)$*$(j)}{166} && \eqnnol{\hbox{$\wh{13}$}}{218}, \eqnno{166} &\cr
& \supers &&  1 &&  231200 &&   1700894366474400 &&  && \eqnno{2}, \eqnno{159} &\cr
& \supers &&  1 &&  549216 &&   5134247872650720 &&  && \eqnnoltmp{269}{332} &\cr
& \supers &&  1 &&  678816 &&  69080128815414048 &&  && \eqnno{9} &\cr
& \supers &&  1 &&  791200 &&   4288711075194400 &&  && \eqnnol{\hbox{$\wh2$}}{207} &\cr
& \supers &&  1 && 1565472 &&  28381748186959008 &&  && \eqnnoltmp{305}{365} &\cr
& \supers &&  1 && 2710944 && 302270555492914464 &&  && \eqnnol{\hbox{$\wh9$}}{214} &\cr
\noalign{\hrule}
}}\hss}

\begin{comments}
We list $|N_1|$ and $|N_3|$. The reason for not
using $N_2$ is that it is not invariant under the transformation
$q\mapsto-q$. There are many differential equations in the table that
are just transformations (by \cite{AZ}, Proposition~8)
$$
y_0(z)\mapsto\frac 1{1-pz}y_0\biggl(\biggl(\frac z{1-pz}\biggr)^r\biggr)
$$
which transforms the Yukawa coupling
$K(q)\mapsto K(q^r)$
These differential equations should not be in the table but we did not know
this transformation when we found them. In the Superseeker table we identify
$K(q)$ and $K(q^r)$ as the instanton numbers are just thinned out.
\end{comments}

\vskip5mm
\bgroup
\footnotesize\scshape

\hbox to80mm{\kern2mm\vbox{\hsize=80mm%
\leftline{Matematikcentrum}
\leftline{Lunds Universitet}
\leftline{Matematik MNF, Box 118}
\leftline{SE-221\,00 Lund, SWEDEN}
\leftline{\textit{E-mail address\/}: \rm\mailto{gert@maths.lth.se}}
}\hss}

\bigskip
\hbox to90mm{\kern2mm\vbox{\hsize=90mm%
\leftline{Fachbereich Mathematik 17}
\leftline{AG Algebraische Geometrie}
\leftline{Johannes Gutenberg-Universit\"at}
\leftline{D-55099 Mainz, GERMANY}
\leftline{\textit{E-mail address\/}: \rm\mailto{straten@mathematik.uni-mainz.de}}
}}

\bigskip
\hbox to80mm{\kern2mm\vbox{\hsize=80mm%
\leftline{School of Mathematical and Physical Sciences}
\leftline{The University of Newcastle}
\leftline{Callaghan 2308, NSW, AUSTRALIA}
\leftline{\textit{E-mail address\/}: \rm\mailto{wadim.zudilin@newcastle.edu.au}}
}\hss}

\egroup

\end{document}